\pdfoutput=1

\RequirePackage{fix-cm}

\documentclass[leqno, fleqn, centertags, 12pt]{article}

\usepackage{latexsym}
\usepackage{amsmath}
\usepackage{amssymb}
\usepackage{verbatim}

\usepackage[T1]{fontenc}

\usepackage[upright]{fourier}

\usepackage{bm}

\usepackage{fullpage}

\usepackage{tikz-cd}

\tikzcdset{row sep/special/.initial=12ex}
\tikzcdset{row sep/spe/.initial=10.5ex}
\tikzcdset{row sep/speh/.initial=13ex}

\tikzcdset{column sep/spec/.initial=5em}
\tikzcdset{column sep/spech/.initial=6em}
\tikzcdset{column sep/speci/.initial=6em}
\tikzcdset{column sep/specc/.initial=8em}
\tikzcdset{column sep/speccc/.initial=10em}

\tikzcdset{row sep/red/.initial=3.79em}

\tikzcdset{row sep/tri/.initial=5em}
\tikzcdset{column sep/tri/.initial=2.4em}

\tikzcdset{row sep/boom/.initial=5em}
\tikzcdset{column sep/boom/.initial=5.1em}

\tikzcdset{row sep/booom/.initial=5.5em}
\tikzcdset{column sep/booom/.initial=5.712em}

\tikzcdset{row sep/boomm/.initial=5.6em}
\tikzcdset{column sep/boomm/.initial=4.25em}

\tikzcdset{row sep/bboom/.initial=8em}
\tikzcdset{column sep/bboom/.initial=8em}

\tikzcdset{row sep/booms/.initial=4.25em}
\tikzcdset{column sep/booms/.initial=4.335em}

\tikzcdset{row sep/boomsss/.initial=4.45em}
\tikzcdset{column sep/boomsss/.initial=4.54em}

\tikzcdset{row sep/oom/.initial=4.5em}
\tikzcdset{column sep/oom/.initial=4.59em}

\tikzcdset{row sep/boomss/.initial=4.75em}
\tikzcdset{column sep/boomss/.initial=4.845em}

\tikzcdset{row sep/free/.initial=12ex}
\tikzcdset{column sep/freee/.initial=4em}

\makeatletter
\renewcommand{\@makefntext}[1]{\vspace*{0.5ex}\parindent=0em
\hspace*{-0.4em}
\hbox to 0.4em{\hss\@makefnmark}\hspace*{0.4em}{#1}
}
\makeatother

\newcounter{mysectionnumber}
\newcommand{\mysection}[2]{\setcounter{footnote}{0}
\setcounter{equation}{0}
\setcounter{myparnum}{0}
\refstepcounter{mysectionnumber}
\vspace{27pt}{\Large {\themysectionnumber.} {#1}}\label{#2}\vspace*{15pt}}

\numberwithin{equation}{section}

\newcommand{\myuppar}[1]{\vspace{\medskipamount}\textbf{#1}\hspace*{0.5em}}

\newcommand{\myit}[1]{\textbf{\textit{#1}}\hspace{0.0em}}

\newcounter{myparnum}[mysectionnumber]
\renewcommand{\themyparnum}{\themysectionnumber.\arabic{myparnum}}
\newcommand{\mypar}[2]{\refstepcounter{myparnum}{\vspace{\medskipamount}\textbf{{\themyparnum. #1}\label{#2}}\hspace{0.5em}}}

\newcounter{mylemmanum}[myparnum]
\newcounter{myaparnum}
\newcommand{\myappend}[2]{\setcounter{footnote}{0}
\setcounter{myaparnum}{0}
\vspace{27pt}{\Large A\dff.\oss {#1}}\label{#2}\vspace*{15pt}}

\newcommand{\myapar}[2]{\refstepcounter{myaparnum}{\vspace{\medskipamount}\textbf{{\themyaparnum. #1}\label{#2}}\hspace{0.5em}}}
\renewcommand{\themyaparnum}{A\halfff\fff.\fff\arabic{myaparnum}}

\newcommand{\proof}{\vspace{\medskipamount}{\textbf{{\emph{Proof}.}}\hspace*{1em}}}

\newcommand{\eproof}{ $\blacksquare$}

\newcommand{\dis}{\displaystyle}

\def\sss{\hspace{0.05em}\ }
\def\dss{\hspace{0.1em}\ }
\def\trs{\hspace{0.15em}\ }
\def\qss{\hspace{0.2em}\ }
\def\pss{\hspace{0.3em}\ }
\def\oss{\hspace{0.4em}\ }

\def\halfff{\hspace*{0.025em}}
\def\fff{\hspace*{0.05em}}
\def\dff{\hspace*{0.1em}}
\def\trf{\hspace*{0.15em}}
\def\qff{\hspace*{0.2em}}
\def\pff{\hspace*{0.3em}}
\def\off{\hspace*{0.4em}}

\def\ttff{{\hspace*{-0.05em}--\hspace*{0.15em}}}

\newcommand{\hnsp}{\hspace*{-0.05em}}
\newcommand{\nsp}{\hspace*{-0.1em}}
\newcommand{\nnsp}{\hspace*{-0.15em}}

\newcommand{\dnsp}{\hspace*{-0.2em}}

\renewcommand{\leq}{\leqslant}
\renewcommand{\geq}{\geqslant}

\newcommand{\zzz}{\mathbf{Z}}

\newcommand{\rrr}{\mathbf{R}}

\newcommand{\image}{\operatorname{Im}\trf}
\newcommand{\kernel}{\operatorname{Ker}}

\newcommand{\inte}{\operatorname{int}\trf}
\newcommand{\cat}{\operatorname{cat}\trf}
\newcommand{\sub}{\operatorname{\textit{sub}}\qff}
\newcommand{\sta}{\operatorname{St}\trf}

\newcommand{\num}[1]{|\qff #1 \qff|}

\newcommand{\norm}[1]{\|\qff #1 \qff\|}

\newcommand{\ttoo}{\hspace*{0.2em}\longrightarrow\hspace*{0.2em}}

\newcommand{\particular}{par\-tic\-u\-lar}

\begin{document}

\setlength{\baselineskip}{12pt plus 0pt minus 0pt}
\setlength{\parskip}{12pt plus 0pt minus 0pt}
\setlength{\abovedisplayskip}{12pt plus 0pt minus 0pt}
\setlength{\belowdisplayskip}{12pt plus 0pt minus 0pt}

\newskip\smallskipamount \smallskipamount=3pt plus 0pt minus 0pt
\newskip\medskipamount   \medskipamount  =6pt plus 0pt minus 0pt
\newskip\bigskipamount   \bigskipamount =12pt plus 0pt minus 0pt

\author{Nikolai\qss V.\qss Ivanov}
\title{Leray\qss theorems\qss in\qss bounded\qss cohomology\qss theory}
\date{}

\footnotetext{\hspace*{-0.65em}\copyright\oss 
Nikolai\qss V.\qss Ivanov,\oss 2020.\trs 
Neither\dss the work reported\sss in\dss the present\dss paper\halfff,\qss
nor its preparation were supported\dss by any\dss corporate entity.}

\maketitle

\renewcommand{\baselinestretch}{1}
\selectfont

\vspace*{12ex}

\myit{\hspace*{0em}\large Contents}\vspace*{1ex} \vspace*{\bigskipamount}\\ 
\hbox to 0.8\textwidth{\myit{\phantom{A.}1.}\hspace*{0.5em} Introduction\hfil  2}\hspace*{0.5em} \vspace*{0.25ex}\\
\hbox to 0.8\textwidth{\myit{\phantom{A.}2.}\hspace*{0.5em} Cohomological\trs Leray\trs theorems\hfil 5}\hspace*{0.5em} \vspace*{0.25ex}\\
\hbox to 0.8\textwidth{\myit{\phantom{A.}3.}\hspace*{0.5em} Homological\trs Leray\trs theorems\hfil 15}\hspace*{0.5em} \vspace*{0.25ex}\\
\hbox to 0.8\textwidth{\myit{\phantom{A.}4.}\hspace*{0.5em} Extensions of\dss coverings and\dss bounded cohomology\hfil 19}\hspace*{0.5em} \vspace*{0.25ex}\\
\hbox to 0.8\textwidth{\myit{\phantom{A.}5.}\hspace*{0.5em} A\trs Leray\trs theorem\dss for $l_{\dff 1}$\dnsp-homology\hfil 23}\hspace*{0.5em} \vspace*{0.25ex}\\
\hbox to 0.8\textwidth{\myit{\phantom{A.}6.}\hspace*{0.5em} Uniqueness of\qss Leray\trs homomorphisms\hfil 26}\hspace*{0.5em} \vspace*{0.25ex}\\
\hbox to 0.8\textwidth{\myit{\phantom{A.}7.}\hspace*{0.5em} Nerves of\trs families and\dss paracompact\sss spaces\hfil 30}\hspace*{0.5em} \vspace*{0.25ex}\\
\hbox to 0.8\textwidth{\myit{\phantom{A.}8.}\hspace*{0.5em} Closed subspaces and\dss fundamental\dss groups\hfil 38}\hspace*{0.5em} \vspace*{0.25ex}\\
\hbox to 0.8\textwidth{\myit{\phantom{A.}9.}\hspace*{0.5em} Closed subspaces and\dss homology\dss groups\hfil 42}\hspace*{0.5em} \vspace*{1ex}\\
\hbox to 0.8\textwidth{\myit{Appendix.}\hspace*{0.5em} Double complexes\hfil 49}\hspace*{0.5em}
\vspace*{1ex}\\
\hbox to 0.8\textwidth{\myit{References}\hspace*{0.5em}\hfil 57}\hspace*{0.5em}  \vspace*{0.25ex}

\newpage
\mysection{Introduction}{introduction}

\myuppar{Leray\trs theory.}
Let\sss $\mathcal{U}$\sss be a covering of\dss a\sss topological\sss space\sss $X$\nnsp,\oss
and\dss let\sss $\mathcal{U}^{\dff \cap}$\sss be\sss the collection of\dss all\sss
non-empty\sss finite intersection of\dss elements of\dss $\mathcal{U}$\dnsp.\oss
A classical\dss theorem of\qss Leray\trs relates\sss the cohomology of\dss $X$\sss
with\sss cohomology\sss of\trs the sets\sss $U\qff \in\qff \mathcal{U}^{\dff \cap}$\sss
and\dss the combinatorial\sss structure of\trs the covering\sss $\mathcal{U}$\dnsp.\oss
The\sss latter\dss is\dss encoded\sss in a simplicial\sss complex\sss $N_{\trf \mathcal{U}}$\nsp,\oss
the nerve of\dss $\mathcal{U}$\dss 
(see\dss Section\qss \ref{cohomological-leray}\qss
for\sss the definition).\oss
If\dss every\sss element\sss of\dss $\mathcal{U}^{\dff \cap}$\sss
is\qss \emph{acyclic},\oss
i.e.\qss has\sss the same cohomology\sss as a point,\oss
Leray\qss theorem\dss implies\sss that\dss the cohomology of\dss $X$\sss are equal\dss to\sss
that\sss of\dss $N_{\trf \mathcal{U}}$\nsp.\oss
Here\qss ``cohomology''\qss are understood\sss very\sss broadly.\oss
Leray\trs theorem applies\sss to\sss the cohomology\sss of\dss sheaves on\sss $X$\nnsp,\oss
as also\sss to\sss the singular cohomology\dss theory\qss (under moderate assumptions).\oss
Morally,\oss Leray\trs theorem applies\sss to every\sss cohomology\dss theory\sss
which can\sss be\sss locally\sss defined.\oss
For example,\oss singular cohomology can be defined\sss in\sss terms of\dss
arbitrarily\sss small\sss simplices.\oss

\myuppar{Bounded cohomology.}
Suppose\sss that\sss elements of\dss $\mathcal{U}^{\dff \cap}$\sss are\dss
in an appropriate sense\qss ``acyclic''\qss with respect\dss to\sss the bounded cohomology.\oss
Gromov's\qss \emph{Vanishing\dss theorem}\qss asserts\sss that\sss if\dss $\mathcal{U}$\sss is\dss
also open,\oss then\dss the image of\trs the canonical\dss homomorphism\sss
$\widehat{H}^{\dff *}\dff(\trf X\trf)
\qff \ttoo\qff
H^{\dff *}\dff(\trf X\trf)$\sss
vanishes in dimensions bigger\sss than\sss the dimension of\dss $N$\nnsp.\oss
See\qss \cite{gro},\pss Section\qss 3.1.\oss

As\dss is\dss well\dss known,\oss the bounded cohomology\sss theory\sss 
cannot\dss be\sss locally\sss defined.\oss
In an attempt\dss to find a proof\dss and a conceptual\dss framework\sss for\sss the\sss 
Vanishing\dss theorem,\oss
the author\qss \cite{i1},\oss \cite{i2}\qss discovered\dss that\sss a part\sss of\qss Leray\trs theory\sss
survives in\sss the non-local\sss setting of\trs the bounded cohomology\sss theory\sss
and\dss leads\sss to stronger\sss results.\oss
Namely,\oss under moderate\sss assumptions 
about\sss $X$\sss and\sss $\mathcal{U}$\nnsp,\oss
the canonical\dss map\sss
$\widehat{H}^{\dff *}\dff(\trf X\trf)
\qff \ttoo\qff
H^{\dff *}\dff(\trf X\trf)$\sss
factors\sss through\dss the natural\dss homomorphism\vspace{1.5pt}\vspace{0.625pt}
\[
\quad
l_{\trf \mathcal{U}}\dff \colon\dff
H^{\dff *}\fff(\trf N_{\trf \mathcal{U}} \dff)
\qff \ttoo\qff 
H^{\dff *}\fff(\trf X \trf)
\]

\vspace{-10.5pt}\vspace{0.625pt}
from\trs Leray's\dss theory.\oss
We call\sss such a result\sss a\qss \emph{Leray\trs theorem},\oss
and call $l_{\trf \mathcal{U}}${\nsp} a\qss \emph{Leray\dss homomorphism}.\oss 
The fact\dss that\dss $H^{\dff *}\fff(\trf X \trf)$\sss can be\sss locally\sss defined\sss
plays a crucial\dss role.\oss
The proofs in\qss \cite{i1},\qss \cite{i2}\qss were inspired\sss by\dss the sheaf\trs theory\sss
and\sss phrased\sss in\sss its\sss language.\oss

The goal\sss of\trs the present\sss paper\dss is\dss to generalize\sss these results\sss 
from\qss \cite{i1},\oss \cite{i2}\qss and\sss at\dss the same\sss time\sss provide 
elementary\sss and\dss transparent\dss proofs.\oss
The sheaf\trs theory\dss is\dss invoked only\dss to deal\sss with fairly\sss
bad spaces and coverings.\oss
In\sss particular\halfff,\oss no sheaf\trs theory\dss is\dss needed\sss to deal\sss
with open coverings of\dss arbitrary spaces.\oss
Our main\sss tools are\sss the double complexes of\dss coverings.\oss
See\trs Section\qss \ref{cohomological-leray}.\oss
In order\sss to stress\sss the elementary\sss nature of\trs these\sss tools,\oss
direct\sss and elementary\sss proofs of\dss required\dss properties
are included\sss in\trs Sections\qss \ref{cohomological-leray},\pss \ref{homological-leray},\oss 
and\sss an\dss Appendix.

At\dss the same\sss time\sss our proofs 
do not\sss depend on any\sss machinery\sss specific\sss to\sss the bounded
cohomology\sss theory.\oss
We will\sss use\sss the\sss fact\dss that\dss
the bounded cohomology of\dss a path-connected space depend only\sss on\sss its
fundamental\sss group,\oss but\sss not\sss any\sss ideas involved\sss in\sss its proofs.\oss
Also,\oss in order\sss to relate our main\sss theorems\sss to\sss the\sss 
Vanishing\dss theorem,\oss we will\sss use\sss the vanishing of\trs the
bounded cohomology of\sss path-connected spaces with amenable fundamental\sss group.\oss

\myuppar{Bounded\sss acyclicity\sss and\sss amenability.}
Let\dss us\sss call\sss a\sss topological\sss space\qss
\emph{boundedly\dss acyclic}\pss if\dss its bounded cohomology\sss are isomorphic\sss
to\sss the bounded cohomology\sss of\dss a point.\oss
This\dss is\dss the most\sss natural\dss notion of\pss ``acyclicity''\qss
with respect\dss to\sss the bounded cohomology.\oss
Let\dss us\sss say\dss that\dss the covering\sss $\mathcal{U}$\sss is\qss
\emph{boundedly\dss acyclic}\pss if\dss every\sss element\sss of\dss
$\mathcal{U}^{\dff \cap}$\sss is\dss boundedly\sss acyclic.\oss
Since\sss the bounded cohomology\sss of\dss a path-connected space are equal\dss to\sss the
bounded cohomology\sss of\dss its fundamental\sss group,\oss
this\sss means\sss that\dss for\sss $U\qff \in\qff \mathcal{U}^{\dff \cap}$\sss 
the fundamental\dss group\sss $\pi_{\dff 1}\dff(\trf V\trf)$\sss
is\qss \emph{boundedly\sss acyclic}\pss in\sss the obvious sense.\oss
Our\sss main\sss results are concerned\sss with boundedly\sss acyclic coverings.\pss
Moreover\halfff,\pss
using an argument\sss from\trs \cite{i1},\qss \cite{i2},\pss
the condition of\dss being\sss boundedly\sss acyclic can be relaxed.\oss
Namely,\oss it\dss is\dss sufficient\dss to assume\sss that\sss every\sss 
$U\qff \in\qff \mathcal{U}^{\dff \cap}$\sss
is\qss \emph{weakly\dss boundedly\sss acyclic}\pss in\dss the sense\sss that image of\trs
$\pi_{\dff 1}\dff(\trf U\trf)$\sss in\sss $\pi_{\dff 1}\dff(\trf X\trf)$\sss
are boundedly\sss acyclic.\oss
See\dss Section\qss \ref{extensions-coverings}.\oss\vspace{-1.75pt}

Starting\sss with\trs Gromov\qss \cite{gro},\oss
the bounded acyclicity\dss is\dss almost\sss always replaced\dss by\dss the stronger property\sss
of\trs being\dss boundedly\sss acyclic\qss ``by\dss the good\dss reason'',\oss namely,\oss
being\sss a path-connected space with amenable fundamental\sss group.\oss
In\dss particular,\oss Gromov's\trs proof\dss
of\trs the\dss Vanishing\dss theorem\dss depends on averaging over\sss amenable groups.\pss
Amenable groups are boundedly\sss acyclic,\oss but\dss the converse\dss is\dss not\dss true.\oss
Sh.\dss Matsumoto\dss and\qss Sh.\dss Morita\qss \cite{mm}\qss provided
examples of\dss groups which are boundedly\sss acyclic by\dss reasons completely\sss
different\dss from\sss being amenable.\oss
So,\oss our\sss results are stronger\dss than\trs Gromov's\dss ones\sss in\dss
this respect\sss also.\oss\vspace{-1.75pt} 

Gromov\dss observed\dss that\sss in\sss the context\sss of\trs the\dss Vanishing\dss theorem\sss
it\dss is\dss sufficient\dss to assume\sss that\dss the images of\trs the inclusion\sss
homomorphisms\sss
$\pi_{\dff 1}\dff(\trf U\trf)
\qff \ttoo\qff
\pi_{\dff 1}\dff(\trf X\trf)$\sss
are amenable for every\sss set\sss $U\qff \in\qff \mathcal{U}$\dnsp.\oss
Cf.\qss Theorem\qss \ref{vanishing}.\oss
Since a subgroup of\dss an amenable group\dss is\dss amenable,\oss
under\dss this assumption\sss the image of\dss
$\pi_{\dff 1}\dff(\trf U\trf)
\qff \ttoo\qff
\pi_{\dff 1}\dff(\trf X\trf)$\sss
is\dss amenable also for every\sss
$U\qff \in\qff \mathcal{U}^{\dff \cap}$\dnsp.\oss
The notion of\dss weakly\dss boundedly\sss acyclic subsets\sss was suggested\dss
by\dss this idea of\qss Gromov.\oss\vspace{-1.75pt}

\myuppar{Open and closed coverings.}
In order\sss to deal\dss with\qss Leray\dss maps\sss
$l_{\trf \mathcal{U}}\dff \colon\dff
H^{\dff *}\fff(\trf N \trf)
\qff \ttoo\qff 
H^{\dff *}\fff(\trf X \trf)$\sss
one needs\sss to assume\sss that\dss the covering\sss $\mathcal{U}$\sss
behaves sufficiently\sss nicely\sss with\sss respect\dss to\sss the singular cohomology\dss theory.\oss
A classical\dss theorem of\qss Eilenberg\qss \cite{e}\qss ensures\sss that\sss
all\sss open coverings are sufficiently\sss nice.\oss
See\trs Theorem\qss \ref{eilenberg}.\oss
Theorem\qss \ref{covering-theorem-open}\qss is\dss 
our\trs Leray\trs theorem\sss for open coverings.\oss
Gromov's\qss Vanishing\dss theorem\dss can\sss be proved\sss in\sss the same way.\oss
See\dss Theorem\qss \ref{vanishing}.\oss\vspace{-1.75pt}

While\dss Gromov's\qss Vanishing\dss theorem\dss is\dss concerned only\dss with open coverings,\oss
Leray\trs theory\sss suggests\sss to consider also closed\dss locally\dss finite\sss coverings.\oss
Suppose\sss that\sss $\mathcal{U}$\sss is\dss a closed\dss locally\dss finite\sss covering,\oss
and\dss that\sss $\mathcal{U}$\sss is\dss weakly\dss boundedly\sss acyclic.\oss
We prove\sss a\trs Leray\trs theorem\dss for\sss such coverings\sss 
$\mathcal{U}$\sss in\sss two different\sss situations.\oss
In\dss both situations we need\dss to assume\sss that\sss $X$\sss is\trs
Hausdorff\trs and\dss paracompact,\oss 
but\dss further assumptions differ.\oss\vspace{-1.75pt}

In\dss Section\qss \ref{nerves-groups}\qss we assume\sss that\sss $\mathcal{U}$\sss
behaves nicely\sss with respect\dss to\sss fundamental\dss groups and covering spaces.\oss
More precisely,\oss we assume\sss that\sss $X$\sss is\dss 
path connected,\oss locally\dss path connected,\oss
and semilocally\sss simply\sss connected,\oss
and\dss that\sss subsets\sss
$U\qff \in\qff \mathcal{U}$\sss 
are path connected and\dss locally\dss path connected.\oss
It\dss turns out\dss that\sss in\dss this case\sss $\mathcal{U}$\sss
can\sss be replaced\dss by\sss an open covering\sss with\sss the same nerve.\oss 
See\dss Theorem\qss \ref{pi-one-coverings-extension}.\oss
This\sss theorem depends on subtle properties of\dss paracompact\sss spaces.\oss
It\dss is\dss hard\dss to extract\dss the full\dss proofs of\trs the needed\sss 
results\sss from\dss the\sss literature,\oss
so we pre\-sent\-ed\dss them\sss  
in\dss Section\qss \ref{paracompact-spaces}.\oss
Theorem\qss \ref{covering-theorem-closed}\qss is\dss 
our\trs Leray\trs theorem\dss in\dss this situation.\oss

In\dss Section\qss \ref{closed-homology}\qss we assume\sss that\sss $\mathcal{U}$\sss
behaves nicely\sss with respect\dss to\sss the singular\sss homology\dss theory.\oss
More precisely,\oss we assume\sss that\dss the space\sss $X$\sss
and all\sss elements of\dss $\mathcal{U}^{\dff \cap}$\sss are\qss
\emph{ho\-mo\-log\-i\-cal\-ly\dss locally\dss connected}.\oss
See\dss Section\qss \ref{closed-homology}\qss for\sss the definition.\oss
Only\sss in\sss this section\sss we resort\dss to\sss the sheaf\trs theory.\oss
Probably,\oss this can\sss be avoided,\oss
but\dss at\dss the cost\sss of\dss obscuring\sss the underlying\sss ideas.\oss
Theorem\qss \ref{acyclic-closed-homology}\qss is\dss our\trs Leray\trs theorem\sss
in\sss this situation.\oss

\myuppar{Abstract\dss Leray\dss theorems.}
As we already\dss pointed out,\oss the proofs do not\dss rely\sss on any\dss tools
from\dss the bounded cohomology\dss theory.\oss
In\sss fact,\oss the basic results hold\sss for any\sss cohomology\dss theory\sss
arising\sss from cochain\sss complexes\sss $A^{\bullet}\dff(\trf Z\trf)$\sss
functorially\sss depending on subspaces\dss $Z\qff \subset\qff X$\dss
and equipped\sss with a natural\dss transformation\dss
$A^{\bullet}\dff(\trf Z\trf)
\qff \ttoo\qff
C^{\bullet}\dff(\trf Z\trf)$\sss
to\sss the complexes of\dss singular cochains.\oss
Theorems\qss \ref{A-acyclic-open}\qss and\qss \ref{A-acyclic-closed}\qss
are such\qss \emph{abstract\trs Leray\trs theorems}\pss
for\sss open and closed coverings respectively.\oss
These results are stated and\sss proved\sss in such an abstract\dss form\dss
not\dss for\sss the sake of\dss generality,\oss 
but\dss in\sss order\dss to make\sss their\sss nature more\sss transparent.\oss

\myuppar{\hspace*{-0.2em}$l_{\dff 1}$\dnsp-homology.}
These results and\sss proofs admit\sss a straightforward dualization,\oss
leading\dss to\trs Leray\trs theorems for $l_{\dff 1}$\dnsp-homology.\oss
Namely,\oss under appropriate assumptions\sss the natural\dss homomorphism\sss
$H_{\dff *}\fff(\trf X \trf)
\qff \ttoo\qff
H^{\dff l_{\dff 1}}_{\dff *}\dff(\trf X\trf)$
can be\dss factored\dss through\dss the map
$H_{\dff *}\fff(\trf X \trf)
\qff \ttoo\qff 
H_{\dff *}\fff(\trf N_{\trf \mathcal{U}} \dff)$
from\trs Leray\trs theory.\oss
We\sss limited ourselves by\dss the case of\dss open coverings.\oss
See\trs Theorem\qss \ref{homology-covering-theorem-open}.\oss
This\sss theorem\dss is\dss deduced\sss from an abstract\dss
homological\qss Leray\trs theorem,\oss Theorem\qss \ref{e-acyclic-open}.\oss
Similar\sss results for closed coverings also can\sss be proved\sss
by\sss dualization of\dss cohomological\dss proofs.\oss

Theorem\qss \ref{homology-covering-theorem-open}\qss is\dss a strengthening of\dss
a recent\sss result\sss of\qss R.\sss Frigerio\qss \cite{f},\oss 
who proved an analogue for $l_{\dff 1}$\dnsp-homology\sss of\qss Gromov's\qss
Vanishing\dss theorem.\oss
Frigerio's\dss proofs are based on\trs Gromov's\qss theories\sss of\dss
multicomplexes and of\trs the diffusion of\dss chains,\oss
recently\dss reconstructed\dss by\trs R.\dss Frigerio\dss and\trs M.\dss Moraschini\qss \cite{fmo},\oss
and are far\sss from\sss being elementary.\oss
A\trs Leray\trs theorem\dss for $l_{\dff 1}$\dnsp-homology\sss was recently\sss
proved\dss by\trs Cl.\dss L\"{o}h\dss and\trs R.\dss Sauer\qss \cite{ls}.\oss
The\sss methods of\qss \cite{ls}\qss are based on\sss author's
homological\sss approach\sss to\sss bounded cohomology\qss \cite{i1},\oss \cite{i2}\qss
and assume amenability.\oss 
It\sss seems\sss that\dss the bounded acyclicity\dss is\dss not\sss sufficient\sss
for\sss methods of\trs \cite{f}\qss and\qss \cite{ls}.\oss

\myuppar{Uniqueness of\qss Leray\trs maps.}
In\dss this paper\dss the\dss Leray\trs maps\sss
$l_{\trf \mathcal{U}}\dff \colon\dff
H^{\dff *}\fff(\trf N_{\trf \mathcal{U}} \dff)
\qff \ttoo\qff 
H^{\dff *}\fff(\trf X \trf)$\sss
for open coverings $\mathcal{U}$\sss 
are defined\sss in\sss terms of\dss the double complex of\dss $\mathcal{U}$\dnsp.\oss
For closed\dss locally\sss finite coverings\sss the definition\dss is\dss
similar\sss in spirit,\oss but\dss is\dss more complicated.\oss

If\dss $\mathcal{U}$\sss is\dss open and\dss the space\sss $X$\sss is\dss paracompact,\oss
one can also use partitions of\dss unity\sss in order\sss to define a natural\dss map\sss
$H^{\dff *}\fff(\trf N_{\trf \mathcal{U}} \dff)
\qff \ttoo\qff 
H^{\dff *}\fff(\trf X \trf)$\nnsp.\oss
R.\dss Frigerio\dss and\dss A.\dss Maffei\qss \cite{fma}\qss
recently\dss proved\dss that\dss for open coverings of\dss
paracompact\sss spaces\sss the\sss two definitions agree.\oss
In\dss Section\qss \ref{uniqueness}\qss we generalize\sss their\sss
result\dss by\sss proving\dss that\sss every\dss two reasonable definitions
of\qss Leray\dss maps agree for open coverings of\dss paracompact\sss spaces.\oss
In\sss more details,\oss following\qss M.\dss Barr\qss \cite{ba},\oss
we consider\dss a category\dss  
having as objects pairs\dss
$(\trf X\fff,\pff \mathcal{U}\dff)$\nnsp,\oss 
where\dss $\mathcal{U}$\dss is\dss a covering\sss of\dss a\sss topological\sss space\dss $X$\nnsp,\oss
and define a\qss \emph{Leray\trs transformation}\qss as a natural\dss transformation\sss
of\dss functors on\sss this category.\oss
See\dss Section\qss \ref{uniqueness}.\oss
It\dss turns out\dss that\sss on\sss the subcategory\sss of\dss
open coverings of\dss paracompact\sss spaces all\qss Leray\trs transformations agree.\oss
See\trs Theorem\qss \ref{coverings-of-paracompact}.\oss

\newpage
\mysection{Cohomological\pss Leray\qss theorems}{cohomological-leray}

\myuppar{The nerve of\dss a covering.}
Let\dss $\mathcal{U}$\trs be a covering of\dss a\sss topological\sss space\sss $X$\sss 
by subsets of\dss $X$\nnsp.\oss
The\qss \emph{nerve}\dss 
$N\off =\off N_{\trf \mathcal{U}}$\dss 
of\dss $\mathcal{U}$\dss 
is\dss a simplicial\sss complex in the sense,\oss for example,\oss 
of\qss \cite{sp},\oss Section\qss 3.1,\oss such\dss that\dss its set\sss of\dss vertices\dss
is\dss in a one-to-one correspondence with\dss  
$\mathcal{U}$\dss and\dss its\sss simplices
are finite non-empty\sss sets of\dss vertices 
such\dss that\dss the intersection of\trs the corresponding
elements of\dss $\mathcal{U}$\dss is\dss non-empty.\oss
We will\sss assume\sss that\sss a\sss linear order\dss $<$\dss on\dss 
the set\sss of\dss vertices of\dss $N$\sss is\dss fixed.\oss 
For each simplex $\sigma$ of\dss $N$ we denote by\sss $\num{\sigma}$\sss 
the intersection of\trs the elements of\dss $\mathcal{U}$ 
corresponding\dss to\sss the vertices of\dss $\sigma$\nnsp.\oss
We\sss denote by\sss $\mathcal{U}^{\dff \cap}$\sss the collection of\dss
all\sss sets $\num{\sigma}$\nnsp.\oss

\myuppar{Classical\qss Leray\trs theorems.}
Let\dss $\mathcal{U}$\trs be a covering\sss of\trs $X$\nnsp.\oss
It\dss is\dss well\dss known\dss that\dss 
under moderate\qss ``niceness''\qss assumptions about\trs $X$\dss
and\dss $\mathcal{U}$\dss
there\dss is\dss a canonical\dss homomorphism\vspace{0pt}\vspace{0.375pt}
\[
\quad
l_{\trf \mathcal{U}}\dff \colon\dff
H^{\dff *}\fff(\trf N \trf)
\qff \ttoo\qff 
H^{\dff *}\fff(\trf X \trf)
\qff,
\]

\vspace{-12pt}\vspace{0.375pt}
where\dss $H^{\dff *}\fff(\trf N \trf)$\dss is\dss the simplicial\sss cohomology\dss
of\trs $N$\dss and\dss $H^{\dff *}\fff(\trf X \trf)$\dss is\dss some appropriate cohomology\sss 
of\trs $X$\nnsp.\oss
For example,\oss if\trs $\mathcal{U}$\dss is\dss open and\dss 
$H^{\dff *}\fff(\trf X \trf)$\dss is\dss the\qss \v{C}ech\sss cohomology,\oss
then such a\sss homomorphism exists\sss by\dss the definition of\trs the\sss latter\halfff.\oss
By\sss a classical\dss theorem of\pss Leray\qss $l_{\trf \mathcal{U}}$\dss is\dss
an\dss isomorphism\dss if\trs $\mathcal{U}$\dss is\qss \emph{acyclic},\oss
i.e.\qss if\trs every\dss non-empty\dss element\sss of\dss $\mathcal{U}^{\dff \cap}$\dss 
has\sss the same cohomology\sss as a\sss point\halfff.\oss
If\trs $\mathcal{U}$\dss is\dss not\sss acyclic,\oss then\dss $H^{\dff *}\fff(\trf N \trf)$\dss
and\dss $H^{\dff *}\fff(\trf X \trf)$\dss are related\dss by\sss a\trs Leray\dss
spectral\sss sequence and\dss $l_{\trf \mathcal{U}}$\dss is\dss one of\trs its\dss two edge 
homomorphisms.\oss 
There are similar results\sss in\sss the homology\dss theory.\oss
Cf.\dss Section\qss \ref{homological-leray}.\oss

Leray\trs theory\dss applies\sss to cohomology\dss theory\sss
$X\off \longmapsto\off
H^{\dff *}\fff(\trf X \trf)$\sss
which can\dss be\qss ``locally\dss defined''.\oss
Orig\-i\-nal\-ly\halfff,\oss the cohomology\dss $H^{\dff *}\fff(\trf X \trf)$\dss
was\sss the sheaf\dss cohomology\sss of\dss $X$\dss 
with coefficients\sss in\sss a sheaf\dss on\dss $X$\nnsp,\oss
and\dss the acyclicity\dss was understood\sss in\dss terms of\dss
cohomology\sss with coefficients in\dss the same sheaf\halfff.\oss
While\sss the definition of\dss singular cohomology\sss of\dss $X$\dss is\dss not\dss local,\oss 
singular cohomology\sss can\sss be\qss ``localized''\qss
in a sense made precise in\trs Theorem\qss \ref{eilenberg}\qss below,\oss
and\qss Leray\trs theory\sss applies\sss to\sss singular cohomology 
at\dss least\dss when\dss the covering\dss $\mathcal{U}$\dss is\dss open.\oss

\myuppar{Leray\trs theorems\sss for\dss non-local\sss cohomology.}
The goal\sss of\trs this section\dss is\dss to show\dss that\sss a part\sss
of\qss Leray\trs theory\sss survives even\dss when\dss the cohomology\dss groups
involved cannot\dss be\qss ``localized''\qss at\sss all.\oss
The motivating\dss example\dss is\dss the bounded cohomology\dss theory\halfff,\oss
but\dss the results are more\sss transparent\dss if\dss stated\dss and\dss proved\dss 
in\dss a\sss general\dss form.\oss

Suppose\sss that\dss
$Z\qff \longmapsto\qff \widetilde{H}^{\dff *}\fff(\trf Z \trf)$\dss
is\trs a\qss ``cohomology\dss theory''\qss defined\dss in\dss terms of\dss
natural\sss cochain complexes,\oss
which are equipped\dss with a\sss natural\dss transformation\dss
to\sss the singular cochain complexes\qss
(all\dss this will\dss be made precise in a moment\halfff).\pss
This natural\dss transformation\dss leads\sss to canonical\dss maps\dss
$\widetilde{H}^{\dff *}\fff(\trf Z \trf)
\qff \ttoo\qff
H^{\dff *}\fff(\trf Z \trf)$\dss
from\dss $\widetilde{H}^{\dff *}\fff(\trf Z \trf)$\dss to\sss
the singular cohomology\dss groups\dss $H^{\dff *}\fff(\trf Z \trf)$\nnsp.\oss

Let\dss us\dss say\dss that\dss a covering\dss $\mathcal{U}$\dss of\trs $X$\dss
is\dss \emph{$\widetilde{H}$\dnsp-acyclic}\pss
if\qss for\dss every\dss non-empty\dss
$U\qff \in\qff \mathcal{U}^{\dff \cap}$\dss
the cohomology\dss
$\widetilde{H}^{\dff *}\fff(\trf U \trf)$\sss
is\trs naturally\dss isomorphic\sss to\sss the singular\sss cohomology\sss 
of\dss a\sss point\halfff.\oss

It\dss turns out\dss that\qss
if\trs $\mathcal{U}$\dss is\qss ``nice''\qss and\dss $\widetilde{H}$\dnsp-acyclic,\oss
then\dss the diagram of the solid arrows\vspace{0pt}\vspace{-0.625pt}
\begin{equation*}
\quad
\begin{tikzcd}[column sep=normal, row sep=huge]\dis
\widetilde{H}^{\fff *}\fff(\dff X \dff) \arrow[rr]
\arrow[rd, dashed]
&
& 
{H}^{\fff *}\fff(\dff X \dff) 
\\
&
{H}^{\fff *}\fff(\trf N \trf) \arrow[ru, "\dis l_{\trf \mathcal{U}}"']
& 
\end{tikzcd}
\end{equation*}

\vspace*{-12pt}
can\dss be completed\dss to a commutative\sss triangle\sss by a dashed arrow.\oss 
In other\dss words,\oss 
the canonical\dss map\dss
$\widetilde{H}^{\dff *}\fff(\trf X \trf)
\qff \ttoo\qff
H^{\dff *}\fff(\trf X \trf)$\dss
factors\sss through\dss the canonical\dss homomorphism\vspace{1.5pt}
\[
\quad 
l_{\trf \mathcal{U}}\dff \colon\dff
H^{\dff *}\fff(\trf N \trf)
\qff \ttoo\qff 
H^{\dff *}\fff(\trf X \trf)
\pff.
\]

\vspace{-10.5pt}
The existence of\dss such\sss factorization\dss could\dss be called\dss a\qss
\emph{Leray\qss theorem}\pss for\dss the cohomology\dss theory\dss
$Z\qff \longmapsto\qff \widetilde{H}^{\dff *}\fff(\trf Z \trf)$\nnsp.\oss
In\dss particular\halfff,\oss there\dss is\dss a\trs Leray\trs theorem\trs
for\dss the bounded cohomology\dss theory\dss
$\widetilde{H}^{\dff *}\fff(\trf Z \trf)
\off =\off
\widehat{H}^{\dff *}\fff(\trf Z \trf)$\nnsp,\oss
and\dss this\sss theorem\dss implies\qss
\emph{Vanishing\dss theorem}\pss
of\qss Gromov\qss \cite{gro}\qss 
(see\trs Theorem\qss \ref{vanishing}\qss below)\qss
and\dss is\dss a\sss vast\dss generalization of\trs the\sss latter\halfff.\oss

Naively\halfff,\oss one could\dss expect\dss that\dss
$\widetilde{H}^{\dff *}\fff(\trf X \trf)$\dss
is\dss isomorphic\dss to\dss
$H^{\dff *}\fff(\trf N \trf)$\dss
if\trs the covering\dss $\mathcal{U}$\dss is\qss ``nice''\qss and\dss $\widetilde{H}$\dnsp-acyclic.\oss
In\dss fact\halfff,\oss this\dss is\dss very\sss unlikely\dss unless
the cohomology\dss theory\dss 
$Z\qff \longmapsto\qff \widetilde{H}^{\dff *}\fff(\trf Z \trf)$\dss
is\dss locally\dss defined.\oss
In\dss particular\halfff,\oss this\dss is\dss not\dss true for\dss the
bounded cohomology\dss theory\halfff.\oss
A\sss natural\dss transformation\dss from\dss $\widetilde{H}^{\dff *}\fff(\trf Z \trf)$\dss
to a\sss locally\dss defined\dss theory\dss
appears\sss to be an\dss inevitable feature.\oss

\myuppar{Generalized\sss cochains.}
In\dss fact\halfff,\oss 
for our purposes\sss it\dss is\dss sufficient\dss to have\dss cohomology\dss
$\widetilde{H}^{\dff *}\fff(\trf Z \trf)$\dss
and\dss the corresponding\sss cochain complexes\sss defined only\dss
for subspaces\dss $Z\qff \subset\pff X$\nnsp.\oss
Let\dss us\dss consider\dss the category\dss $\sub X$\dss having\dss subspaces of\trs $X$\dss
as objects and\dss the inclusions\dss
$Y\qff \subset\pff Z$\dss as morphisms.\oss
Let\trs $A^{\bullet}$\sss be\dss a contravariant\dss functor\dss from\dss
$\sub X$\dss to
augmented cochain complexes of\dss modules over some ring\dss $R$\nnsp.\oss
In our applications\sss $R\off =\off \rrr$\nnsp,\oss
but\sss assuming\dss this does\sss not\sss simplifies anything.\oss
The\sss functor\dss $A^{\bullet}$\dss assigns\sss 
to a subspace\dss $Z\qff \subset\pff X$\dss a\sss complex\vspace{-1.5pt}
\begin{equation*}
\quad
\begin{tikzcd}[column sep=large, row sep=normal]\dis
0 \arrow[r]
& 
R \arrow[r, "\dis d_{\trf -\dff 1}\qff"]
& 
A^{\fff 0}\dff(\trf Z \trf) \arrow[r, "\dis d_{\trf 0}\off"]
& 
A^{\fff 1}\dff(\trf Z \trf) \arrow[r, "\dis d_{\trf 1}\off"]
&   
A^{\fff 2}\dff(\trf Z \trf) \arrow[r, "\dis d_{\trf 2}\off"]
&
\off \ldots \off\off
\end{tikzcd}
\end{equation*}

\vspace{-9pt}
of\trs $R$\dnsp-modules,\oss
and\dss to every\dss inclusion\dss $Y\qff \subset\pff Z$\dss a\qss
\emph{restriction\dss morphism}\qss
$A^{\bullet}\dff(\trf Z \trf)
\dff \ttoo\dff
A^{\bullet}\dff(\trf Y \trf)$\nnsp.\oss
Elements\sss of\dss $A^{\fff q}\dff(\trf Z \trf)$
are\sss thought\sss as\qss \emph{generalized $q$\dnsp-cochains}\qss of\dss $Z$\nnsp.\oss
The action of\trs restriction\dss morphisms\dss will\dss be denoted\dss by\vspace{1.5pt}
\[
\quad
c\off \longmapsto\off c\qff\bigl|_{\qff Z} 
\pff.
\]

\vspace{-10.5pt}
We will\dss need\dss mostly\dss the complexes\dss $A^{\bullet}\dff(\trf Z \trf)$\dss
for\dss $Z\off =\off X$\sss or\dss 
$Z\qff \in\qff \mathcal{U}^{\dff \cap}$\dnsp,\oss
but\dss it\dss is\dss hard\dss to imagine such a\sss functor defined only\dss for\dss
such\dss $Z$\nnsp.\oss
The examples are\sss the functor\dss
$Z\off \longmapsto\off C^{\dff \bullet}\dff(\trf Z \trf)$\dss
of\dss singular cochains with coefficients in $R$ and\dss the functor\dss
$Z\off \longmapsto\off B^{\dff \bullet}\dff(\trf Z \trf)$\dss
of\dss bounded\dss real-valued singular cochains\sss in\dss the case of\trs
$R\off =\off \rrr$\nnsp.\oss
We will\sss say\dss that\dss the covering\dss $\mathcal{U}$\dss is\qss 
\emph{$A^{\dff \bullet}$\dnsp\dnsp-acyclic}\oss if\pss
$A^{\bullet}\dff(\trf Z\trf)$\dss is\dss exact\dss for every\dss
$Z\qff \in\qff \mathcal{U}^{\dff \cap}$\dnsp,\oss
and\dss will\sss say\dss that\dss $\mathcal{U}$\dss is\qss
\emph{boundedly\dss acyclic}\pss if\qss it\trs is\dss
$B^{\dff \bullet}$\dnsp\dnsp-acyclic.\oss

\myuppar{The double complex of\dss a covering.}
Let\dss $\mathcal{U}$\trs be a covering\sss of\trs $X$\dss and\dss let\trs
$N$\trs be\sss its\sss nerve.\oss
For\qss $p\qff \geq\qff 0$\qss let\dss $N_{\dff p}$\dss 
be\sss the set\sss of\dss $p$\dnsp-dimensional\sss simplices of\trs $N$\nnsp.\oss
For\qss $p\fff,\pff q\qff \geq\qff 0$\pss let\vspace{3.75pt}
\[
\quad
C^{\dff p}\dff(\trf N\fff,\pff A^{q} \trf)
\off\off =\off\off
\prod\nolimits_{\qff \sigma\qff \in\qff N_{\dff p}}\qff 
A^{q}\fff\left(\trf \num{\sigma} \trf\right)
\pff.
\]

\vspace{-8.25pt}
An element\dss 
$c\qff \in\qff
C^{\dff p}\dff(\trf N\fff,\pff A^{q} \trf)$\dss
is\dss a\sss map assigning\dss to each\sss
$\sigma\qff \in\pff N_{\dff p}$\dss 
a generalized $q$\dnsp-cochain\vspace{1.5pt}
\[
\quad
c_{\dff \sigma}
\qff \in\pff
A^{q}\fff\left(\trf \num{\sigma} \trf\right)
\pff.
\]

\vspace{-10.5pt}
The functor\dss $A^{\bullet}$ assigns\dss to each\dss 
$\sigma\qff \in\qff N_{\dff p}$\dss the\sss complex\dss $A^{\bullet}\dff(\trf \num{\sigma} \trf)$,\oss
i.e.\qss the complex\vspace{0pt}
\begin{equation}
\label{b-sigma}
\quad
\begin{tikzcd}[column sep=large, row sep=huge]\dis
0 \arrow[r]
& 
R \arrow[r, "\dis d_{\trf -\dff 1}\qff"]
& 
A^{\fff 0}\dff(\trf \num{\sigma} \trf) \arrow[r, "\dis d_{\trf 0}"] 
&   
A^{\fff 1}\dff(\trf \num{\sigma} \trf) \arrow[r, "\dis d_{\trf 1}"] 
&
\off \ldots \off\off.
\end{tikzcd}
\end{equation}

\vspace{-9pt}
For each\qss $p\qff \geq\qff 0$\qss
the\sss term-wise direct\dss product\dss of\trs complexes\qss (\ref{b-sigma})\qss over\qss
$\sigma\qff \in\qff N_{\dff p}$\qss has\sss the form\vspace{0pt}
\begin{equation}
\label{c-nerve}
\quad
\begin{tikzcd}[column sep=large, row sep=huge]\dis
0 \arrow[r]
& 
C^{\dff p}\dff(\trf N \trf) \arrow[r, "\dis d_{\trf -\dff 1}\qff"]
& 
C^{\dff p}\dff(\trf N\fff,\pff A^{\fff 0} \trf) \arrow[r, "\dis d_{\trf 0}"] 
&   
C^{\dff p}\dff(\trf N\fff,\pff A^{\fff 1} \trf) \arrow[r, "\dis d_{\trf 1}"] 
&
\off \ldots \off,
\end{tikzcd}
\end{equation}

\vspace{-9pt}
where\dss $C^{\dff p}\dff(\trf N \trf)$\dss is\dss the space of\dss simplicial
$p$\dnsp-cochains of\dss $N$\dss with coefficients in\sss $R$\nnsp.\oss
As every\sss simplicial\sss complex,\oss 
the nerve\dss $N$\dss has one empty\sss simplex.\oss
By\dss the definition,\oss
the dimension of\trs the simplex\dss $\varnothing$\dss
is\dss $-\qff 1$\dss and\dss
$\num{\varnothing}\off =\off X$\nnsp.\oss
So,\pss the complex\qss (\ref{c-nerve})\qss 
is\dss defined also for\dss $p\off =\off -\qff 1$\dss 
and\dss is\dss equal\dss to\dss 
$A^{\bullet}\dff(\trf X\trf)$\dss in\dss this case.\oss 

One can also set\trs
$A^{\fff -\dff 1}\dff(\trf Z\trf)\off =\off R$\dss for every\sss subspace\dss $Z$\dnsp.\oss
Then\sss one can\dss replace\sss the\sss term\sss $R$\sss in\qss (\ref{b-sigma})\qss
by\dss $A^{\fff -\dff 1}\dff(\trf \num{\sigma} \trf)$\dss
and\dss interpret\dss the\sss term\dss $C^{\dff p}\dff(\dff N \dff)$\dss
in\qss (\ref{c-nerve})\qss as\dss
$C^{\dff p}\dff(\trf N\fff,\pff A^{\fff -\dff 1} \trf)$\nnsp.\oss
We will\sss denote\sss the complex\qss (\ref{c-nerve})\qss
by\dss $C^{\dff p}\dff(\dff N\fff,\pff A^{\bullet} \dff)$\nnsp.\oss
For every\dss $p\qff \geq\qff -\qff 1$\dss 
there\dss is\dss a canonical\dss morphism\vspace{3pt}
\[
\quad
\delta_{\fff p}\dff \colon\dff
C^{\dff p}\dff(\dff N\fff,\pff A^{\bullet} \dff)
\qff \ttoo\qff
C^{\dff p\dff +\dff 1}\dff(\dff N\fff,\pff A^{\bullet} \dff)
\pff
\]

\vspace{-9pt}
defined as follows.\oss
Suppose\sss that\dss $\sigma\qff \in\qff N_{\dff p\dff +\dff 1}$\dss and\dss let\trs
$v_{\dff 0}\off <\off v_{\dff 1}\off <\off \ldots\off <\off v_{\dff p\dff +\dff 1}$\dss
be\sss the vertices of\dss $\sigma$\dss listed\dss in\dss the increasing order\halfff.\oss 
For each\qss
$i\qff =\qff 0\fff,\pff 1\fff,\pff \ldots\fff,\pff p\qff +\qff 1$\qss
let\oss 
$\partial_{\dff i}\dff\sigma
\off =\off 
\sigma\qff \smallsetminus\qff \{\qff v_{\fff i} \qff\}$\dss
be\sss the $i${\dnsp}th\sss face of\dss $\sigma$\nnsp.\oss
Clearly\halfff,\pss
$\num{\sigma}
\pff \subset\off
\num{\partial_{\dff i}\dff\sigma}$\dss
and\dss there\dss is\dss a\sss restriction\dss
morphism\vspace{3pt}
\[
\quad
A^{\bullet}\dff(\trf \num{\partial_{\dff i}\dff\sigma} \trf)
\qff \ttoo\qff
A^{\bullet}\dff(\trf \num{\sigma} \trf)
\pff.
\]

\vspace{-9pt}
If\qss 
$c\qff \in\qff
C^{\dff p}\dff(\trf N\fff,\pff A^{q} \trf)$\nnsp,\oss
then\dss $\delta_{\fff p}\dff(\dff c\trf)$\dss assigns\sss to\sss $\sigma$\dss
the\sss generalized $q$\dnsp-cochain\vspace{6pt}
\[
\quad
\left(\trf \delta_{\fff p}\dff(\dff c\trf)\trf\right)_{\dff \sigma}
\off\qff =\off\off
\sum\nolimits_{\dff i\qff =\qff 0}^{\dff p\qff +\qff 1}\off 
(\dff -\qff 1 \dff)^{\dff i}\off
c_{\pff \partial_{\dff i}\dff\sigma}\qff\bigl|_{\qff \num{\sigma}}
\off\pff \in\off\qff
A^{q}\fff\left(\trf \num{\sigma} \trf\right)
\pff.
\]

\vspace{-6pt}
The fact\dss that\sss each\dss  
$A^{\bullet}\dff(\trf \num{\partial_{\dff i}\dff\sigma} \trf)
\qff \ttoo\qff
A^{\bullet}\dff(\trf \num{\sigma} \trf)$\dss
is\dss
a\sss morphism of\dss complexes implies\sss that\dss\vspace{4.5pt}
\[
\quad
\delta_{\fff p}\dff \colon\dff
C^{\dff p}\dff(\dff N\fff,\pff A^{\bullet} \dff)
\qff \ttoo\qff
C^{\dff p\dff +\dff 1}\dff(\dff N\fff,\pff A^{\bullet} \dff)
\pff
\]

\vspace{-7.5pt}
is\dss a\sss morphism of\dss complexes also.\oss
Equivalently\halfff,\pss 
$\delta_{\fff p}\dff \circ\qff d_{\trf q}
\off =\off
d_{\trf q}\dff \circ\qff \delta_{\fff p}$\dss
on\dss
$C^{\dff p}\dff(\trf N\fff,\pff A^{q} \trf)$\dss
for every\trs
$p\fff,\pff q\off \geq\off -\qff 1$\nnsp.\oss
Let\dss us\dss collect\dss 
$C^{\dff p}\dff(\trf N\fff,\pff A^{q} \trf)$\dss
and\dss $d_{\trf q}\dff,\off \delta_{\fff p}$\dss
into a single diagram\vspace{3pt}
\begin{equation}
\label{big-diagram}
\qquad
\begin{tikzcd}[column sep=large, row sep=huge]\dis
R \arrow[r]
\arrow[d]
& 
A^{\fff 0}\dff(\trf X \trf) \arrow[r] 
\arrow[d]
&   
A^{\fff 1}\dff(\trf X \trf) \arrow[r] 
\arrow[d]
&
\off \ldots \off
\\
C^{\trf 0}\dff(\trf N \trf) \arrow[r]
\arrow[d]
& 
C^{\trf 0}\dff(\trf N\fff,\pff A^{\fff 0} \trf) \arrow[r]
\arrow[d]
&   
C^{\trf 0}\dff(\trf N\fff,\pff A^{\fff 1} \trf) \arrow[r]
\arrow[d]
&
\off \ldots \off
\\
C^{\dff 1}\dff(\trf N \dff) \arrow[r]
\arrow[d]
& 
C^{\dff 1}\dff(\trf N\fff,\pff A^{\fff 0} \trf) \arrow[r]
\arrow[d]
&   
C^{\dff 1}\dff(\trf N\fff,\pff A^{\fff 1} \trf) \arrow[r]
\arrow[d]
&
\off \ldots \off
\\
C^{\trf 2}\dff(\trf N \dff) \arrow[r]
\arrow[d]
& 
C^{\trf 2}\dff(\trf N\fff,\pff A^{\fff 0} \trf) \arrow[r]
\arrow[d]
&   
C^{\trf 2}\fff(\trf N\fff,\pff A^{\fff 1} \trf) \arrow[r]
\arrow[d]
&
\off \ldots \off
\\
\ldots \vphantom{C^{\trf 2}\fff(\trf N\fff,\pff A^{\trf 1} \trf)} &
\ldots \vphantom{C^{\trf 2}\fff(\trf N\fff,\pff A^{\trf 1} \trf)} &
\ldots \vphantom{C^{\trf 2}\fff(\trf N\fff,\pff A^{\trf 1} \trf)} &
\ldots \vphantom{C^{\trf 2}\fff(\trf N\fff,\pff A^{\trf 1} \trf)} &
\quad \vphantom{C^{\trf 2}\fff(\trf N\fff,\pff A^{\trf 1} \trf)} .
\end{tikzcd}
\end{equation}

\vspace{-9pt}\vspace{-1.125pt}
Since\dss 
$\delta_{\fff p}\dff \circ\qff d_{\trf q}
\off =\off
d_{\trf q}\dff \circ\qff \delta_{\fff p}$\dss
on\dss
$C^{\dff p}\dff(\trf N\fff,\pff A^{q} \trf)$\nnsp,\oss
this diagram\dss is\dss commutative.\oss
By\dss the construction,\pss
$d_{\trf q\dff +\dff 1}\dff \circ\qff d_{\trf q}\off =\off 0$\dss
for every\dss $q$\nnsp,\oss
i.e.\qss the rows of\trs this diagram are complexes.\oss
A standard computation shows that\dss
$\delta_{\fff p\dff +\dff 1}\dff \circ\qff \delta_{\fff p}
\off =\off
0$\dss
for every\dss $p$\nnsp,\oss
i.e.\qss the columns of\trs this diagram are also complexes.\oss
It\dss turns out\dss that\dss the\sss top row and\dss the\sss left\sss column
of\trs this diagram should\dss be\sss treated differently\dss from\dss the rest\halfff.\oss
The part\vspace{3pt}
\begin{equation}
\label{double-complex-covering}
\qquad
\begin{tikzcd}[column sep=large, row sep=huge]\dis
C^{\trf 0}\dff(\trf N\fff,\pff A^{\fff 0} \trf) \arrow[r]
\arrow[d]
&   
C^{\trf 0}\dff(\trf N\fff,\pff A^{\fff 1} \trf) \arrow[r]
\arrow[d]
&
\off \ldots \off
\\
C^{\dff 1}\dff(\trf N\fff,\pff A^{\fff 0} \trf) \arrow[r]
\arrow[d]
&   
C^{\dff 1}\dff(\trf N\fff,\pff A^{\fff 1} \trf) \arrow[r]
\arrow[d]
&
\off \ldots \off
\\
C^{\trf 2}\dff(\trf N\fff,\pff A^{\fff 0} \trf) \arrow[r]
\arrow[d]
&   
C^{\trf 2}\fff(\trf N\fff,\pff A^{\fff 1} \trf) \arrow[r]
\arrow[d]
&
\off \ldots \off
\\
\ldots \vphantom{C^{\trf 2}\fff(\trf N\fff,\pff A^{\trf 1} \trf)} &
\ldots \vphantom{C^{\trf 2}\fff(\trf N\fff,\pff A^{\trf 1} \trf)} &
\quad \vphantom{C^{\trf 2}\fff(\trf N\fff,\pff A^{\trf 1} \trf)} .
\end{tikzcd}
\end{equation}

\vspace{-9pt}\vspace{-1.125pt}
of\trs the diagram\qss (\ref{big-diagram})\qss
is\dss the\qss \emph{double complex of\qss the covering}\qss $\mathcal{U}$\nnsp.\oss
See\sss the\qss Appendix\qss for\dss a discussion of\dss
double complexes and\dss related\dss notions.\oss
We will\sss denote\sss this double complex\dss by\dss
$C^{\trf \bullet}\fff(\trf N\fff,\pff A^{\bullet} \trf)$\sss
and assume from\dss now\sss on\dss that\dss $\bullet$\dss 
\emph{stands\dss for\sss non-negative\dss integers}.\oss
Let\trs 
$T^{\trf \bullet}\dff(\trf N\fff,\pff A \trf)$\dss
be\sss the\sss total\sss complex of\trs the double complex\dss
$C^{\trf \bullet}\fff(\trf N\fff,\pff A^{\bullet} \trf)$\nnsp.\oss
By\dss the definition,\oss\vspace{0pt}
\[
\quad
T^{\dff n}\dff(\trf N\fff,\pff A \dff)
\off\off =\off\off
\bigoplus_{p\qff =\qff 0}^n\qff
C^{\dff p}\dff(\trf N\fff,\pff A^{\fff n\dff -\dff p} \trf)
\]

\vspace{-12pt}
and\dss the differential\qss
$D\dff \colon\dff T^{\dff n}\dff(\trf N\fff,\pff A \dff)
\qff \ttoo\qff 
T^{\dff n\dff +\dff 1}\dff(\trf N\fff,\pff A \dff)$\qss
is\dss given\dss by\dss the formula\oss\vspace{3pt}
\[
\quad
D\qff \bigl|\qff C^{\dff p}\dff(\trf N\fff,\pff A^{q} \trf)
\off =\off\dff 
d_{\dff q}\qff +\qff (\dff -\qff 1 \fff)^{\dff p}\qff \delta_{\dff p} 
\pff.\oss
\]

\vspace{-9pt}
The homomorphisms\dss $d_{\dff -\dff 1}$\dss from\qss (\ref{c-nerve})\qss
lead\dss to a morphism\dss\vspace{3pt}
\[
\quad
\lambda_{\trf A}\qff \colon\qff
C^{\dff \bullet}\dff(\trf N \trf) 
\qff \ttoo\qff
T^{\dff \bullet}\dff(\trf N\fff,\pff A \trf)
\]

\vspace{-9pt}
Similarly\halfff,\oss the morphism\qss
$\delta_{\dff -\dff 1}\dff \colon\dff
A^{\bullet}\dff(\trf X\trf)
\off =\off
C^{\dff -\dff 1}\dff(\trf N\fff,\pff A^{\bullet} \trf)
\qff \ttoo\qff
C^{\dff 0}\dff(\trf N\fff,\pff A^{\bullet} \trf)$\qss
leads\sss to a morphism\dss\vspace{3pt}
\[
\quad
\tau_{\dff A}\qff \colon\qff
A^{\bullet}\dff(\trf X\trf)
\qff \ttoo\qff
T^{\dff \bullet}\dff(\trf N\fff,\pff A \trf)
\pff,
\]

\vspace{-9pt}
where\dss it\dss is\dss understood\dss that\dss the augmentation\dss term
$A^{\fff -\dff 1}\dff(\trf X\trf)$ 
is\dss  removed\dss from\dss 
$A^{\bullet}\dff(\trf X\trf)$\nnsp.\oss

\mypar{Lemma.}{acyclic-coverings}
\emph{If\pss $\mathcal{U}$\dss is\pss 
$A^{\bullet}$\dnsp\dnsp-acyclic,\oss 
then\dss the morphism\qss
$\lambda_{\trf A}\dff \colon\dff
C^{\dff \bullet}\dff(\trf N \trf) 
\qff \ttoo\qff
T^{\dff \bullet}\dff(\trf N\fff,\pff A \trf)$\qss
induces an\dss isomorphism of\qss cohomology\dss groups.\oss}

\proof
If\qss $\mathcal{U}$\dss is\pss 
$A^{\bullet}$\dnsp\dnsp-acyclic,\oss
then\dss for every\dss  
simplex\dss $\sigma\off \neq\off \varnothing$\dss the complex\qss
(\ref{b-sigma})\qss is\dss exact\halfff.\oss
This implies\sss that\dss for every\dss $p\qff \geq\qff 0$\dss the complex\qss
(\ref{c-nerve})\qss is\dss exact\halfff.\oss
It\dss follows\dss that\dss the rows of\trs the double complex\trs
$C^{\dff \bullet}\dff(\trf N\fff,\pff A^{\bullet} \trf)$\dss
are exact\sss and\dss $d_{\trf -\dff 1}$\dss 
induces an\dss isomorphism of\trs the complex\dss
$C^{\dff \bullet}\dff(\trf N \trf)$\dss 
with\dss the kernel\dss of\trs the morphism of\dss complexes\vspace{3pt}
\[
\quad
\hspace{-1.4em}
\begin{tikzcd}[column sep=large, row sep=huge]\dis
d_{\trf 0}\qff \colon\qff\hspace{-6em}
&
C^{\dff \bullet}\dff(\trf N\fff,\pff A^{\fff 0} \trf) \arrow[r]
& 
C^{\dff \bullet}\dff(\trf N\fff,\pff A^{\fff 1} \trf)\pff.
\end{tikzcd}
\]

\vspace{-9pt}
It\dss remains\sss to apply\qss Theorem\qss \ref{double-complex}\qss
from\dss the\dss Appendix.\oss  \eproof

\myuppar{Singular cochains.}
Let\dss $Y$\dss be a\sss topological\sss space.\oss
Following\dss Godement\qss \cite{go},\oss
let\dss us\dss fix a\sss topological\sss space\sss $\Delta$\sss
and call\sss a continuous map\dss
$s\dff \colon\dff \Delta\qff \ttoo\qff Y$\dss a\qss \emph{singular\dss simplex}\pss
in\dss $Y$\nnsp.\oss
Cf\halfff.\qss \cite{go},\oss Chapter\qss II,\oss Example\qss 3.9.1.\oss
The main examples are\sss $\Delta\off =\off \Delta^{\fff q}$\trs for\dss
$q\qff \geq\qff 0$\nnsp.\oss
Let\dss $S\dff(\trf Y\trf)$\dss be\sss the set\sss of\dss singular\sss simplices in\dss $Y$\dnsp.\oss
A\qss \emph{singular\sss cochain}\qss with coefficients\sss in\sss $R$\sss
is\dss defined as 
a\sss map\dss $S\dff(\trf Y\trf)\qff \ttoo\qff R$\nnsp.\oss
Let\dss $C\dff(\trf Y\trf)$\dss be\sss the $R$\dnsp-module of\dss singular cochains\dss
in\dss $Y$\dnsp.\oss
For\qss $p\qff \geq\qff 0$\pss let\vspace{6pt}
\[
\quad
C^{\dff p}\dff(\trf N\fff,\pff C \trf)
\off\off =\off\off
\prod_{\sigma\qff \in\qff N_{\dff p}\vphantom{N^N}}\qff 
C\fff\left(\trf \num{\sigma} \trf\right)
\pff.
\]

\vspace{-6pt}
So,\oss an element\dss 
$c\qff \in\qff
C^{\dff p}\dff(\trf N\fff,\pff C \trf)$\dss
is a map assigning\dss to a simplex\sss
$\sigma\qff \in\pff N_{\dff p}$\dss 
a cochain\vspace{3pt}
\[
\quad
c_{\dff \sigma}
\qff \colon\pff
S\fff\left(\trf \num{\sigma} \trf\right)
\qff \ttoo\qff
R
\pff.
\]

\vspace{-9pt}
The maps\qss 
$\delta_{\fff p}\dff \colon\dff
C^{\dff p}\dff(\trf N\fff,\pff C \trf)
\qff \ttoo\qff
C^{\dff p\dff +\dff 1}\dff(\trf N\fff,\pff C \trf)$\qss
are defined as before.\oss
In\dss more details,\oss\vspace{1.5pt}
\begin{equation}
\label{classical-boundary}
\quad
\left(\trf \delta_{\fff p}\dff c\trf\right)_{\dff \sigma}
\qff(\dff s\trf)
\off\qff =\off\qff
\sum_{\dff i\qff =\qff 0}^{\dff p\qff +\qff 1}\off 
(\dff -\qff 1 \dff)^{\dff i}\off
c_{\pff \partial_{\dff i}\dff\sigma}
\qff(\dff s\trf)
\end{equation}

\vspace{-12pt}
if\trs $p\qff \geq\qff 0$\nnsp,\qss
$\sigma$\dss
is\dss a $p$\dnsp-simplex\sss of\trs $N$\nnsp,\pss
$c\qff \in\qff
C^{\dff p}\dff(\trf N\fff,\pff C \trf)$\nnsp,\oss
and\dss
$s\dff \colon\dff
\Delta\qff \ttoo\qff \num{\sigma}$\dss singular simplex.\oss

A singular simplex\sss
$s\dff \colon\dff
\Delta\qff \ttoo\qff X$\dss
is\dss said\dss to be\qss
\emph{small}\oss if\trs $s\trf(\dff \Delta\dff)$\dss
is\dss contained\dss in some\dss $U\qff \in\pff \mathcal{U}$\dnsp.\oss
Let\dss $S\dff(\qff X\fff,\pff \mathcal{U} \trf)$\dss
be\sss the set\sss of\dss small\sss singular simplices in\dss $X$\nnsp.\oss
A\qss \emph{small\sss singular\sss cochain}\pss is\dss defined as 
a\sss map\dss
$S\dff(\qff X\fff,\pff \mathcal{U} \trf)
\qff \ttoo\qff
R$\nnsp.\oss
Let\dss $C\dff(\qff X\fff,\pff \mathcal{U} \trf)$\dss 
be\sss the $R$\dnsp-module of\dss small\sss singular cochains\dss
in\dss $X$\nnsp.\oss
The restriction of\dss singular cochains\sss to\dss
$S\dff(\qff X\fff,\pff \mathcal{U} \trf)$\dss
leads\sss to a\sss map\dss
$C\dff(\trf X \trf)
\qff \ttoo\qff
C\dff(\qff X\fff,\pff \mathcal{U} \trf)$\nnsp.\oss
Sim\-i\-lar\-ly\halfff,\oss 
restrictions of\dss small\sss singular cochains\dss
to\dss $S\dff(\dff U\dff)$\nnsp,\pss $U\qff \in\qff \mathcal{U}$\dnsp,\oss
lead\dss to a\sss map\dss\vspace{3pt}
\[
\quad
\overline{\delta}_{\trf -\dff 1}\dff \colon\dff
C\dff(\qff X\fff,\pff \mathcal{U} \trf)
\qff \ttoo\qff
C^{\dff 0}\dff(\trf N\fff,\pff C \trf)
\pff.
\]

\vspace{-9pt}
As before,\pss
$\delta_{\fff p\dff +\dff 1}\dff \circ\qff \delta_{\fff p}
\off =\off
0$\dss
for every\dss $p\qff \geq\qff 0$\dss
and also\dss
$\delta_{\trf 0}\dff \circ\qff \overline{\delta}_{\trf -\dff 1}
\off =\off
0$\nnsp.\oss

\mypar{Lemma.}{columns-are-exact}
\emph{The\qss following\dss sequence\dss is\dss exact\dff:}\vspace{0pt}
\[
\quad
\begin{tikzcd}[column sep=large, row sep=huge]\dis
0 \arrow[r]
&
C\dff(\qff X\fff,\pff \mathcal{U} \trf) \arrow[r, "\dis \delta_{\trf -\dff 1}\qff"]
& 
C^{\trf 0}\dff(\trf N\fff,\pff C \trf) \arrow[r, "\dis \delta_{\trf 0}\qff"]
& 
C^{\trf 1}\dff(\trf N\fff,\pff C \trf) \arrow[r, "\dis \delta_{\trf 1}"] 
&   
\off \ldots\off \off.
\end{tikzcd}
\]

\vspace{-9pt}
\proof
Obviously\halfff,\pss  
$\overline{\delta}_{\trf -\dff 1}$\dss 
is\dss injective.\oss
Let\dss us\dss prove\sss the exactness\sss in\dss the\sss term\dss
$C^{\trf 0}\dff(\trf N\fff,\pff C \trf)$\nnsp.\oss
An element\dss
$c\qff \in\qff
C^{\trf 0}\dff(\trf N\fff,\pff C \trf)$\dss
is\dss a\sss family\sss of\dss cochains\dss
$c_{\trf U}\qff \in\pff C\dff(\trf U\trf)$\nnsp,\pss
$U\qff \in\pff \mathcal{U}$\dnsp.\oss
It\dss belongs\sss to\sss the\sss kernel\sss of\dss $\delta_{\trf 0}$\dss
if\trs and\dss only\trs if\trs for every\dss 
$U\fff,\pff V\qff \in\pff \mathcal{U}$\dss
the cochains\dss $c_{\trf U}$\dss and\dss $c_{\qff V}$\dss
are equal\sss on singular simplices\sss with\dss the image
in\dss $U\qff \cap\qff V$\dnsp.\oss
If\trs this\dss is\dss the case,\oss
then\dss there exists a unique\sss map\dss
$\gamma\dff \colon\dff
S\dff(\qff X\fff,\pff \mathcal{U} \trf)
\qff \ttoo\qff
R$\dss
such\dss that\dss the restriction of\dss $\gamma$\dss
to\dss $S\dff(\trf U \trf)$\dss is\dss equal\dss to\dss $c_{\trf U}$\dss
for every\dss $U\qff \in\pff \mathcal{U}$\dnsp.\oss
Clearly\halfff,\pss
$\overline{\delta}_{\trf -\dff 1}\dff(\qff \gamma\qff)
\off =\off
c$\nnsp.\oss
It\dss follows\dss that\dss $\overline{\delta}_{\trf -\dff 1}$\dss
is\dss an\sss isomorphism onto\sss the kernel of\dss $\delta_{\trf 0}$\nsp,\oss
i.e.\qss the sequence\dss is\dss exact\dss at\dss the\sss terms $0$ and\dss
$C^{\trf 0}\dff(\trf N\fff,\pff C \trf)$\nnsp.\oss
In order\dss to prove\sss the exactness at\dss the\sss terms\dss
$C^{\trf p}\dff(\trf N\fff,\pff C \trf)$\dss with\dss $p\qff \geq\qff 1$\nnsp,\oss
it\dss is\dss sufficient\dss to\sss construct\sss a 
contracting\sss chain\dss homotopy\vspace{0.5pt}
\[
\quad
\begin{tikzcd}[column sep=large, row sep=huge]\dis
C^{\trf p\dff -\dff 1}\dff(\trf N\fff,\pff C \trf) 
& 
C^{\trf p}\dff(\trf N\fff,\pff C \trf)\pff,\oss 
p\qff \geq\qff 1\pff. 
\arrow[l, "\dis\pff k_{\dff p}"'] 
\end{tikzcd}
\]

\vspace{-9pt}
Let\dss us choose for every\sss small\sss singular simplex\dss 
$s\dff \colon\dff
\Delta\qff \ttoo\qff X$\dss
an element\dss $U_{\dff s}\qff \in\qff \mathcal{U}$\dss
such\dss that\dss 
$s\trf(\dff \Delta\dff) 
\qff \subset\pff 
U_{\dff s}$\dss
and\dss denote by\dss $u_{\fff s}$\dss the corresponding\dss vertex of\trs $N$\nnsp.\oss
Let\dss
$c\qff \in\qff
C^{\dff p}\dff(\trf N\fff,\pff C \trf)$\nnsp.\oss
We need\dss to define a cochain\qss
$k_{\trf p}\dff(\dff c\trf)_{\dff \tau}\qff \in\qff
C^{\trf p\dff -\dff 1}\fff\left(\trf \num{\tau} \trf\right)$\qss
for every\dss
$\tau\dff \in\pff N_{\dff p\dff -\dff 1}$\nsp.\oss
Suppose\sss that\vspace{1.5pt}
\[
\quad
r\dff \colon\dff
\Delta\qff \ttoo\qff \num{\tau}
\pff
\]

\vspace{-10.5pt}
is\dss a singular simplex\sss
and\dss let\trs
$\sigma
\off =\off
\tau\dff \cup\dff \{\trf u_{\dff r}\trf\}$\nnsp.\oss
Clearly\halfff,\pss
$r\dff(\dff \Delta\trf)
\qff \subset\qff
\num{\tau}\qff \cap\qff U_{\dff r}
\off =\off
\num{\sigma}$\nnsp.\oss
In\dss \particular\halfff,\pss $\sigma$\dss is\dss a simplex.\oss
If\trs $u_{\dff r}\qff \in\pff \tau$\nnsp,\oss
then\dss $\sigma$\dss is\dss a $(\dff p\qff -\qff 1\dff)$\dnsp-simplex
and\dss we set\dss\vspace{3pt}
\[
\quad
k_{\trf p}\dff(\dff c\trf)_{\dff \tau}\trf(\dff r\trf)
\off =\off
0 
\pff.
\]

\vspace{-9pt}
Otherwise,\pss $\sigma$\dss is\dss a $p$\dnsp-simplex,\pss
$\tau\off =\off \partial_{\fff a}\trf \sigma$\dss for some\dss $a$\nnsp,\oss 
and\dss we set\vspace{3pt}
\[
\quad
k_{\dff p}\dff(\dff c\trf)_{\dff \tau}\trf(\dff r\trf)
\off =\off
(\qff -\qff 1\trf)^{\dff a}\trf c_{\trf \sigma}\dff(\trf r\trf)
\pff.
\]

\vspace{-9pt}
Let\dss 
$c\qff \in\qff
C^{\dff p}\dff(\trf N\fff,\pff C \trf)$\nnsp,\oss
$\sigma\qff \in\pff N_{\dff p}$\nsp,\oss 
and\dss
$s\dff \colon\dff
\Delta\qff \ttoo\qff \num{\sigma}$\dss
is\dss a singular simplex.\oss
In order\dss to prove\sss that\dss $k_{\dff \bullet}$\dss is\dss
a contracting chain\dss homotopy\halfff,\oss
we need\dss to prove\sss that\vspace{6pt}
\begin{equation}
\label{homotopy-check}
\quad
\delta_{\dff p\dff -\dff 1}\trf \bigl(\trf k_{\dff p}\dff(\dff c\trf)
\qff\bigr)_{\dff \sigma}
\qff(\dff s\trf)
\off +\off
k_{\dff p\dff +\dff 1}\trf \bigl(\qff \delta_{\dff p}\dff(\dff c\trf)
\qff\bigr)_{\dff \sigma}
\qff(\dff s\trf)
\off =\off
c_{\trf \sigma}\dff(\dff s\trf)
\end{equation}

\vspace{-6pt}
for every\sss such\dss $c\fff,\pff \sigma$\nnsp,\oss
and\dss $s$\nnsp.\oss
The formula\qss (\ref{classical-boundary})\qss implies\sss that\vspace{6pt}
\[
\quad
\delta_{\dff p\dff -\dff 1}\trf \bigl(\trf k_{\dff p}\dff(\dff c\trf)
\qff\bigr)_{\dff \sigma}
\qff(\dff s\trf)
\off\qff =\off\qff
\sum\nolimits_{\dff i\qff =\qff 0}^{\dff p}\off (\dff -\qff 1 \dff)^{\dff i}\qff
k_{\dff p}\dff(\dff c\trf)_{\qff \partial_{\dff i}\dff\sigma}
\qff(\dff s\trf)
\pff.
\]

\vspace{-6pt}
Suppose first\dss that\dss
$u_{\dff s}\qff \in\qff \sigma$\dss
and\dss let\dss
$\tau\off =\off \sigma\qff \smallsetminus\qff \{\trf u_{\dff s}\trf\}$\nnsp.\oss
Then\dss
$\tau\off =\off \partial_{\fff a}\trf \sigma$\dss for some\dss $a$\nnsp.\oss
Clearly\halfff,\pss
$u_{\dff s}\qff \in\qff \partial_{\dff i}\dff \sigma$\qss
if\qss $i\off \neq\off a$\dss
and\dss hence\dss
$k_{\dff p}\dff(\dff c\trf)_{\qff \partial_{\dff i}\dff\sigma}
\qff(\dff s\trf)
\off =\off
0$\qss
if\qss $i\off \neq\off a$\nnsp.\qff\oss
It\dss follows\dss that\vspace{4.5pt}
\[
\quad
\delta_{\dff p\dff -\dff 1}\trf \bigl(\trf k_{\dff p}\dff(\dff c\trf)
\qff\bigr)_{\dff \sigma}
\qff(\dff s\trf)
\off =\off
(\dff -\qff 1 \dff)^{\dff a}\qff
k_{\dff p}\dff(\dff c\trf)_{\qff \partial_{\dff a}\dff\sigma}
\qff(\dff s\trf)
\off =\off
c_{\trf \sigma}\dff(\dff s\trf)
\pff.
\]

\vspace{-7.5pt}
Since also\qss 
$k_{\dff p\dff +\dff 1}\trf (\trf d\qff)_{\dff \sigma}
\qff(\dff s\trf)
\off =\off
0$\qss
for every\sss $d$\nnsp,\oss this implies\qss (\ref{homotopy-check})\qss 
in\dss the case when\dss
$u_{\dff s}\qff \in\qff \sigma$\nnsp.\oss

Suppose now\dss that\trs
$u_{\dff s}\qff \not\in\qff \sigma$\nnsp.\oss
Let\trs
$v_{\dff 0}\off <\off v_{\dff 1}\off <\off \ldots\off <\off v_{\dff p}$\dss
be\sss the vertices of\dss $\sigma$\dss listed\dss in\dss the increasing order\halfff,\oss
and\dss let\dss $a$\dss be such\dss that\trs
$v_{\dff a\dff -\dff 1}\off <\off u_{\dff s}\off <\off v_{\dff a}$\qss
({\fff}for\dss $a\off =\off 0$\dss or\dss $p$\dss one of\trs these inequalities\dss
is\dss vacuous).\oss
Let\trs
$\rho
\off =\off
\sigma\dff \cup\dff \{\trf u_{\dff s}\trf\}$\dss
and\dss
$\rho_{\dff i}
\off =\off
\partial_{\dff i}\dff \sigma\dff \cup\dff \{\trf u_{\dff s}\trf\}$\nnsp.\oss
Then\vspace{6pt}
\[
\quad
\rho_{\dff i}\off =\off \partial_{\dff i}\trf \rho
\quad
\mbox{for}\quad 
i\qff <\qff a
\quad 
\mbox{and}\quad
\rho_{\dff i}\off =\off \partial_{\dff i\dff +\dff 1}\trf \rho 
\quad 
\mbox{for}\quad  
i\qff \geq\qff a\dff.\quad
\]

\vspace{-7.5pt}
Also\dss $\sigma\off =\off \partial_{\fff a}\trf \rho$\nnsp,\oss\vspace{4.5pt}
\[
\quad
\partial_{\dff i}\trf \sigma\off =\off \partial_{\fff a\dff -\dff 1}\qff \rho_{\dff i}
\quad
\mbox{for}\quad 
i\qff <\qff a\dff,
\quad 
\mbox{and}\quad
\partial_{\dff i}\trf \sigma\off =\off \partial_{\fff a}\qff \rho_{\dff i}
\quad
\mbox{for}\quad 
i\qff \geq\qff a\dff.\quad
\]

\vspace{-6pt}
Similarly\dss to\sss the above,\oss
(\ref{classical-boundary})\qss implies\sss that\vspace{6pt}
\[
\quad
\delta_{\dff p\dff -\dff 1}\trf \bigl(\trf k_{\dff p}\dff(\dff c\trf)
\qff\bigr)_{\dff \sigma}
\qff(\dff s\trf)
\off\qff =\off\qff
\sum\nolimits_{\dff i\qff =\qff 0}^{\dff p}\off (\dff -\qff 1 \dff)^{\dff i}\qff
k_{\dff p}\dff(\dff c\trf)_{\qff \partial_{\dff i}\dff\sigma}
\qff(\dff s\trf)\pff
\]

\vspace{-27pt}
\[
\quad
\phantom{\delta_{\dff p\dff -\dff 1}\trf \bigl(\trf k_{\dff p}\dff(\dff c\trf)
\qff\bigr)_{\dff \sigma}
\qff(\dff s\trf)
\off\qff }
=\off\qff
\sum\nolimits_{\dff i\qff =\qff 0}^{\dff a\dff -\dff 1}\off 
(\dff -\qff 1 \dff)^{\dff i\dff +\dff a\dff -\dff 1}\qff
c_{\qff \rho_{\dff i}}
\qff(\dff s\trf)
\off +\off
\sum\nolimits_{\dff i\qff =\qff a}^{\dff p}\off 
(\dff -\qff 1 \dff)^{\dff i\dff +\dff a}\qff
c_{\qff \rho_{\dff i}}
\qff(\dff s\trf)
\pff.
\]

\vspace{-12pt}\vspace{1.75pt}
At\dss the same\sss time\vspace{3pt}
\[
\quad
k_{\dff p\dff +\dff 1}\trf \bigl(\qff \delta_{\dff p}\dff(\dff c\trf)
\qff\bigr)_{\dff \sigma}
\qff(\dff s\trf)
\off\qff =\off\qff
(\qff -\qff 1\trf)^{\dff a}\pff 
\delta_{\dff p}\dff(\dff c\trf)_{\trf \rho}\dff(\dff s\trf)
\off\qff =\off\qff
(\qff -\qff 1\trf)^{\dff a}\pff 
\sum\nolimits_{\dff i\qff =\qff 0}^{\dff p\qff +\qff 1}\off 
(\dff -\qff 1 \dff)^{\dff i}\qff
c_{\qff \partial_{\dff i}\dff\rho}
\dff(\dff s\trf)
\pff
\]

\vspace{-27pt}
\[
\quad
=\off\qff
\sum\nolimits_{\dff i\qff =\qff 0}^{\dff a\qff -\qff 1}\off 
(\dff -\qff 1 \dff)^{\dff i\dff +\dff a}\qff
c_{\qff \rho_{\dff i}}
\dff(\dff s\trf)
\off\qff +\off\qff
c_{\qff \sigma}
\dff(\dff s\trf)
\off\qff +\off\qff
\sum\nolimits_{\dff i\qff =\qff a\dff +\dff 1}^{\dff p\dff +\dff 1}\off 
(\dff -\qff 1 \dff)^{\dff i\dff +\dff a}\qff
c_{\qff \rho_{\dff i\dff -\dff 1}}
\dff(\dff s\trf)
\pff
\]

\vspace{-27pt}\vspace{0.25pt}
\[
\quad
=\off\qff
\sum\nolimits_{\dff i\qff =\qff 0}^{\dff a\qff -\qff 1}\off 
(\dff -\qff 1 \dff)^{\dff i\dff +\dff a}\qff
c_{\qff \rho_{\dff i}}
\dff(\dff s\trf)
\off\qff +\off\qff
c_{\qff \sigma}
\dff(\dff s\trf)
\off\qff +\off\qff
\sum\nolimits_{\dff i\qff =\qff a}^{\dff p}\off 
(\dff -\qff 1 \dff)^{\dff i\dff +\dff a\dff +\dff 1}\qff
c_{\qff \rho_{\dff i\dff -\dff 1}}
\dff(\dff s\trf)
\pff.
\]

\vspace{-6pt}
By\dss adding\dss the results of\trs these calculations we see\sss
that\qss (\ref{homotopy-check})\qss holds when\dss
$u_{\dff s}\qff \not\in\qff \sigma$\dss also.\oss  \eproof

\myuppar{The complex of\dss small\sss singular\sss cochains.}
Taking\dss $\Delta\off =\off \Delta^{q}$\dss
leads\sss to\sss the notions of\dss small\sss singular $q$\dnsp-simplices
and\sss small\sss singular 
$q$\dnsp-cochains.\oss
Let\dss $S_{\fff q}\dff(\qff X\fff,\pff \mathcal{U} \trf)$\dss
be\sss the set\sss of\dss small\sss singular $q$\dnsp-simplices
and\dss
$C^{\fff q}\dff(\qff X\fff,\pff \mathcal{U} \trf)$\dss
be\sss the $R$\dnsp-module of\dss small\sss singular $q$\dnsp-cochains.\oss
Clearly\halfff,\oss every\dss face of\dss a small\sss singular $q$\dnsp-simplex\dss
is\dss small.\oss
Therefore,\oss one can define\sss the\qss \emph{coboundary\dss maps}\vspace{1.5pt}
\[
\quad
d\dff \colon\dff
C^{\fff q}\dff(\qff X\fff,\pff \mathcal{U} \trf)
\qff \ttoo\qff
C^{\fff q\dff +\dff 1}\dff(\qff X\fff,\pff \mathcal{U} \trf) 
\]

\vspace{-10.5pt}
in\dss the usual\sss way\halfff.\oss
As usual,\oss
$d\dff \circ\dff d\off =\off 0$\dss 
and\dss 
$C^{\fff q}\dff(\qff X\fff,\pff \mathcal{U} \trf)$\dss
together\dss with\sss $d$\trs form a complex,\oss
which we\sss denote by\trs
$C^{\dff \bullet}\dff(\qff X\fff,\pff \mathcal{U} \trf)$\nnsp.\oss
The restriction of\dss singular $q$\dnsp-cochains\sss to\dss
$S_{\fff q}\dff(\qff X\fff,\pff \mathcal{U} \trf)$\dss
leads\sss to a morphism\dss
$C^{\dff \bullet}\dff(\trf X \trf)
\qff \ttoo\qff
C^{\dff \bullet}\dff(\qff X\fff,\pff \mathcal{U} \trf)$\nnsp.\oss 
The following\dss fundamental\dss theorem\dss is\dss well\dss known.\oss
It\dss is\dss due\sss to\qss Eilenberg\qss \cite{e},\oss
but\dss this\dss is\dss hardly\dss ever\sss mentioned\dss nowadays.\oss

\mypar{Theorem.}{eilenberg}
\emph{The morphism\pss
$C^{\dff \bullet}\dff(\trf X \trf)
\qff \ttoo\qff
C^{\dff \bullet}\dff(\qff X\fff,\pff \mathcal{U} \trf)$\qss
induces an\dss isomorphism of\dss cohomology\dss groups\dss
if\trs the interiors\qss $\inte U$\dss of\qss sets\qss $U\qff \in\pff \mathcal{U}$\dss
cover\pss $X$\nnsp.\oss}  \eproof

\mypar{Lemma.}{usual-cochains-isomorphism}
\emph{The homomorphisms\qss
$\overline{\delta}_{\trf -\dff 1}\dff \colon\dff
C^{\fff q}\dff(\qff X\fff,\pff \mathcal{U} \trf)
\qff \ttoo\qff 
C^{\trf 0}\dff(\trf N\fff,\pff C^{\fff q} \trf)$\qss
define a morphism}\qss\vspace{1.5pt}
\[
\quad
\tau_{\qff C\dff,\dff \mathcal{U}}\qff \colon\qff
C^{\dff \bullet}\dff(\qff X\fff,\pff \mathcal{U} \trf)
\qff \ttoo\qff
T^{\dff \bullet}\dff(\trf N\fff,\pff C \trf)
\pff
\]

\vspace{-10.5pt}
\emph{inducing\sss an\dss isomorphism\dss of\trs cohomology\dss groups.\oss
The morphism}\vspace{1.5pt}
\[
\quad
\tau_{\qff C}\qff \colon\qff
C^{\dff \bullet}\dff(\trf X \trf)
\qff \ttoo\qff
T^{\dff \bullet}\dff(\trf N\fff,\pff C \trf)
\pff
\]

\vspace{-10.5pt}
\emph{induces an\dss isomorphism\dss of\trs cohomology\dss groups\dss
if\pss the interiors of\qss the elements of\pss $\mathcal{U}$\dss cover\trs $X$\nnsp.\oss}

\proof
Clearly\halfff,\pss $\overline{\delta}_{\trf -\dff 1}$\dss commutes with\dss
the coboundary\dss maps.\oss
Lemma\qss \ref{columns-are-exact}\qss implies\sss that\dss 
$\overline{\delta}_{\trf -\dff 1}$\dss induces an\dss isomorphism of\trs the complex\dss 
$C^{\dff \bullet}\dff(\qff X\fff,\pff \mathcal{U} \trf)$\dss
with\dss the kernel\sss of\trs the morphism\vspace{1.5pt}
\[
\quad
\delta_{\trf 0}\dff \colon\dff
C^{\trf 0}\fff(\trf N\fff,\pff C^{\trf \bullet} \trf)
\qff \ttoo\qff
C^{\trf 1}\fff(\trf N\fff,\pff C^{\trf \bullet} \trf)
\pff.
\]

\vspace{-10.5pt}
At\dss the same\sss time,\oss Lemma\qss \ref{columns-are-exact}\qss
implies\sss that\dss the columns of\trs the double complex\dss
$C^{\trf \bullet}\fff(\trf N\fff,\pff C^{\trf \bullet} \trf)$\dss
are exact\halfff.\oss
Therefore\trs Theorem\qss \ref{double-complex}\qss
with\dss the rows and columns interchanged\dss
implies\dss the first\sss statement\sss of\trs the\sss lemma.\oss  
The first\sss statement\dss implies\sss the second\dss in\sss view\sss of\trs
Lemma\qss \ref{eilenberg}.\oss  \eproof

\myuppar{The homomorphism\dss
$H^{\dff *}\fff(\trf N \trf)
\qff \ttoo\qff 
H^{\dff *}\fff(\trf X \trf)$\dss for open coverings.}
The morphisms\vspace{-1.5pt}
\[
\quad
\begin{tikzcd}[column sep=large, row sep=huge]\dis
C^{\dff \bullet}\dff(\trf N \trf) 
\arrow[r, "\dis \lambda_{\qff C}"]
& 
T^{\dff \bullet}\dff(\trf N\fff,\pff C \dff)
&
C^{\dff \bullet}\dff(\trf X \trf)\pff
\arrow[l, "\dis \off \tau_{\trf C}"']
\end{tikzcd}
\]

\vspace{-9pt}
lead\dss to homomorphisms of\dss cohomology\dss groups,\oss\vspace{-1.5pt}
\[
\quad
\begin{tikzcd}[column sep=large, row sep=huge]\dis
H^{\dff *}\dff(\trf N \trf) 
\arrow[r, "\dis \lambda_{\qff C\dff *}"]
& 
H^{\dff *}\dff(\trf N\fff,\pff C \dff)
&
H^{\dff *}\dff(\trf X \trf)\pff,
\arrow[l, "\dis \off \tau_{\trf C\dff *}"']
\end{tikzcd}
\]

\vspace{-9pt}
where\dss 
$H^{\dff *}\dff(\trf N\fff,\pff C \dff)$\dss
denotes\sss the cohomology\sss of\trs
$T^{\dff \bullet}\dff(\trf N\fff,\pff C \dff)$\nnsp.\oss
If\trs the\sss interiors\sss of\qss the\sss elements of\trs $\mathcal{U}$\dss cover\trs $X$\nnsp,\oss
then\dss $\tau_{\trf C\dff *}$\dss is\dss an\dss isomorphism\dss by\trs
Lemma\qss \ref{usual-cochains-isomorphism}\qss
and\dss we\sss take\sss the composition\vspace{4.5pt}
\[
\quad
\tau_{\trf C\dff *}^{\dff -\dff 1}\qff \circ\qff \lambda_{\qff C\dff *}\qff \colon\qff
H^{\dff *}\fff(\trf N \trf)
\off \ttoo\off 
H^{\dff *}\fff(\trf X \trf)
\pff,
\]

\vspace{-7.5pt}
as\sss the\qss \emph{canonical\qss homomorphism}\pss
$l_{\trf \mathcal{U}}\dff \colon\dff
H^{\dff *}\fff(\trf N \trf)
\qff \ttoo\qff 
H^{\dff *}\fff(\trf X \trf)$\nnsp.\oss

\myuppar{The\dss $A^{\bullet}$\dnsp\dnsp-cohomology\sss and\dss the singular cohomology\halfff.}
The\dss \emph{$A^{\bullet}$\dnsp\dnsp-cohomology\dss groups}\qss\vspace{3pt}\vspace{-0.25pt}
\[
\quad
\widetilde{H}^{\dff *}\dff(\trf X\trf)
\off =\off\dff
H_{\dff A}^{\dff *}\dff(\trf X\trf)
\]

\vspace{-9pt}\vspace{-0.25pt} 
of\trs $X$\dss
are simply\dss the cohomology\dss groups of\trs the complex\dss $A^{\bullet}\dff(\trf X\trf)$\dss
with\dss the\sss term\dss $R$\dss omitted.\oss 
Suppose\sss that\dss the functor\dss $A^{\bullet}$\dss is\dss equipped\dss with a natural\dss
transformation\dss 
$\varphi^{\dff \bullet}\dff \colon\dff
A^{\bullet}\qff \ttoo\qff C^{\dff \bullet}$\dnsp.\oss
This natural\dss transformation\dss leads\sss to a\sss homomorphism\dss\vspace{3pt}\vspace{-0.25pt}
\[
\quad
\widetilde{H}^{\dff *}\dff(\trf X\trf)
\off =\off\dff
H_{\dff A}^{\dff *}\dff(\trf X\trf)
\off \ttoo\off
H^{\dff *}\dff(\trf X\trf)
\pff.
\]

\vspace{-9pt}\vspace{-0.25pt} 
\mypar{Theorem.}{A-acyclic-open}
\emph{If\pss $\mathcal{U}$\dss is\qss $A^{\bullet}$\dnsp\dnsp-acyclic 
and\trs the interiors of\qss the elements of\qss $\mathcal{U}$\dss cover\trs $X$\nnsp,\oss
then\dss the\sss homomorphism\qss
$\widetilde{H}^{\dff *}\dff(\trf X\trf)
\qff \ttoo\qff
H^{\fff *}\fff(\trf X \trf)$\qss
can be factored\dss through\dss 
the\sss canonical\dss 
homomorphism\qss
$l_{\trf \mathcal{U}}\dff \colon\dff
H^{\dff *}\fff(\trf N \trf)
\qff \ttoo\qff 
H^{\dff *}\fff(\trf X \trf)$\dnsp.\oss}

\proof
The homomorphisms\qss\vspace{2pt}
\[
\quad
\varphi^{\dff q}\dff\left(\trf \num{\sigma} \trf\right)\dff \colon\dff
A^{q}\fff\left(\trf \num{\sigma} \trf\right)
\qff \ttoo\qff
C^{\dff q}\fff\left(\trf \num{\sigma} \trf\right)
\pff
\]

\vspace{-10pt}
lead\dss to a morphism\vspace{2pt}
\[
\quad
\varphi^{\dff \bullet\dff \bullet}\dff \colon\dff
C^{\dff \bullet}\dff(\trf N\fff,\pff A^{\bullet} \trf)
\off \ttoo\off
C^{\dff \bullet}\dff(\trf N\fff,\pff C^{\dff \bullet} \trf)
\pff
\]

\vspace{-10.5pt}
of\dss double complexes.\oss 

In\dss turn,\pss $\varphi^{\dff \bullet\dff \bullet}$\dss leads\sss to a morphism\vspace{3pt}
\[
\quad
\Phi^{\dff \bullet}\dff \colon\dff
T^{\dff \bullet}\dff(\trf N\fff,\pff A \trf)
\off \ttoo\off
C^{\dff \bullet}\dff(\trf N\fff,\pff C \trf)
\pff
\]

\vspace{-9pt}
Clearly\halfff,\oss the diagram\vspace{4.5pt}
\[
\quad
\begin{tikzcd}[column sep=large, row sep=boom]\dis
C^{\dff \bullet}\dff(\trf N \trf) 
\arrow[r]
\arrow[d, "\dis \qff ="]
& 
T^{\dff \bullet}\dff(\trf N\fff,\pff A \trf)
\arrow[d, "\dis \qff \Phi^{\dff \bullet}"]
&
A^{\bullet}\dff(\trf X \trf)
\arrow[l]
\arrow[d, "\dis \qff \varphi^{\dff \bullet}"]
\\
C^{\dff \bullet}\dff(\trf N \trf) \arrow[r]
& 
T^{\dff \bullet}\dff(\trf N\fff,\pff C \trf)
&
C^{\dff \bullet}\dff(\trf X \trf)\pff.
\arrow[l]
\end{tikzcd}
\]

\vspace{-7.5pt}
is\dss commutative and\dss leads\sss
to\sss the following\sss commutative diagram of\dss cohomology\dss
groups\vspace{4.5pt}
\[
\quad
\begin{tikzcd}[column sep=large, row sep=boom]\dis
H^{\dff *}\dff(\trf N \trf) 
\arrow[r, red, line width=0.8pt]
\arrow[d, "\dis \qff ="]
& 
H^{\dff *}\dff(\trf N\fff,\pff A \trf)
\arrow[d]
&
\widetilde{H}^{\dff *}\fff(\trf X \trf)
\arrow[l]
\arrow[d]
\\
H^{\dff *}\dff(\trf N \trf) \arrow[r]
& 
H^{\dff *}\dff(\trf N\fff,\pff C \trf)
&
H^{\dff *}\dff(\trf X \trf)\pff,
\arrow[l, red, line width=0.8pt]
\end{tikzcd}
\]

\vspace{-7.5pt}
where\dss
$H^{\dff *}\dff(\trf N\fff,\pff A \trf)$\dss
denotes\sss the cohomology\sss of\trs the\sss total\sss complex\dss
$T^{\dff \bullet}\dff(\trf N\fff,\pff A \trf)$\nnsp.\oss

The red arrows are isomorphisms.\oss
Indeed,\oss
since\dss the covering\dss 
$\mathcal{U}$\dss is\dss $A^{\bullet}$\dnsp\dnsp-acyclic,\oss
the arrow\dss
$H^{\dff *}\dff(\trf N \trf)
\qff \ttoo\qff
H^{\dff *}\dff(\trf N\fff,\pff A \trf)$\dss
is\dss an\dss isomorphism\dss by\qss Lemma\qss \ref{acyclic-coverings}.\oss
Since\sss the interiors of\qss the elements of\qss $\mathcal{U}$\dss cover\trs $X$\nnsp,\oss
the arrow\dss
$H^{\dff *}\dff(\trf X \trf)
\qff \ttoo\qff
H^{\dff *}\dff(\trf N\fff,\pff C \trf)$\dss
is\dss an\dss isomorphism\dss by\qss Lemma\qss \ref{usual-cochains-isomorphism}.\oss
By\dss inverting\dss these\sss two arrows we get\dss the commutative diagram\vspace{4.5pt}
\[
\quad
\begin{tikzcd}[column sep=large, row sep=boom]\dis
H^{\dff *}\dff(\trf N \trf) 
\arrow[d, "\dis \qff ="]
& 
H^{\dff *}\dff(\trf N\fff,\pff A \trf)
\arrow[l, red, line width=0.8pt]
\arrow[d]
&
\widetilde{H}^{\dff *}\fff(\trf X \trf)
\arrow[l]
\arrow[d]
\\
H^{\dff *}\dff(\trf N \trf) \arrow[r]
& 
H^{\dff *}\dff(\trf N\fff,\pff C \trf)
\arrow[r, red, line width=0.8pt]
&
H^{\dff *}\dff(\trf X \trf)\pff.
\end{tikzcd}
\]

\vspace{-7.5pt}
It\dss follows\dss that\dss 
$\widetilde{H}^{\dff *}\fff(\trf X \trf)
\qff \ttoo\qff
H^{\dff *}\dff(\trf X \trf)$\dss
factors\sss through\dss the
canonical\dss homomorphism\vspace{4.5pt}\vspace{-0.625pt}
\[
\quad
\begin{tikzcd}[column sep=large, row sep=boom]\dis
H^{\dff *}\dff(\trf N \trf) \arrow[r]
& 
H^{\dff *}\dff(\trf N\fff,\pff C \trf)
\arrow[r]
&
H^{\dff *}\dff(\trf X \trf)\pff.
\end{tikzcd}
\]

\vspace{-7.5pt}\vspace{-0.625pt}
This completes\sss the proof\halfff.\oss  \eproof

\newpage
\mysection{Homological\pss Leray\qss theorems}{homological-leray}

\myuppar{Generalized\sss chains.}
Let\dss $e_{\dff \bullet}$\dss be\dss a covariant\dss functor\dss from\dss
$\sub X$\dss to
augmented chain complexes of\dss modules over some ring\dss $R$\nnsp.\oss
The\sss functor\dss $e_{\dff \bullet}$\dss assigns\sss to a subspace $Z\qff \subset\pff X$ 
a\sss complex\vspace{0pt}
\begin{equation*}
\quad
\begin{tikzcd}[column sep=large, row sep=normal]\dis
0 \arrow[r]
& 
R  \arrow[l]
& 
e_{\trf 0}\dff(\trf Z \trf) 
\arrow[l, "\dis \off d_{\trf 0}\qff"']
& 
e_{\trf 1}\dff(\trf Z \trf)
\arrow[l, "\dis \off d_{\trf 1}\off"']
&   
e_{\trf 2}\dff(\trf Z \trf)
\arrow[l, "\dis \off d_{\trf 2}\off"']
&
\off \ldots \off\off,
\arrow[l, "\dis \off d_{\trf 3}\off"']
\end{tikzcd}
\end{equation*}

\vspace{-9pt}
For every\dss $Y\qff \subset\qff Z$\dss there\dss is\dss a\qss
\emph{inclusion\dss morphism}\qss
$e_{\trf \bullet}\dff(\trf Z \trf)
\dff \ttoo\dff
e_{\trf \bullet}\dff(\trf Y \trf)$\dnsp.\oss
For\qss $p\fff,\pff q\qff \geq\qff 0$\pss let\vspace{4pt}
\[
\quad
c_{\dff p}\dff(\trf N\fff,\pff e_{\dff q} \trf)
\off\off =\off\off
\bigoplus\nolimits_{\qff \sigma\qff \in\qff N_{\dff p}}\qff 
e_{\dff q}\fff\left(\trf \num{\sigma} \trf\right)
\pff.
\]

\vspace{-8pt}
So,\oss an element 
$c\qff \in\qff
c_{\dff p}\dff(\trf N\fff,\pff e_{\dff q} \trf)$
is\dss a\sss direct\sss sum of\qss \emph{generalized
$q$\dnsp-chains}\vspace{4pt}
\begin{equation*}
\quad
c
\off =\off\
\bigoplus\nolimits_{\qff \sigma\qff \in\qff N_{\dff p}}\qff c_{\dff \sigma}
\off\off \in\off\off\dff
\bigoplus\nolimits_{\qff \sigma\qff \in\qff N_{\dff p}}\qff
e_{\dff q}\fff\left(\trf \num{\sigma} \trf\right)
\pff.
\end{equation*}

\vspace{-8pt}
For every\dss $p\qff \geq\qff 0$\dss 
there\dss is\dss a canonical\dss morphism\qss 
$\delta_{\fff p}\dff \colon\dff
c_{\dff p}\dff(\dff N\fff,\pff e_{\dff \bullet} \dff)
\qff \ttoo\qff
c_{\dff p\dff -\dff 1}\dff(\dff N\fff,\pff e_{\dff \bullet} \dff)$\nnsp,\oss
defined as follows.\oss
Let\dss $\sigma\qff \in\qff N_{\dff p}$\nsp.\oss
For each\dss face\dss $\partial_{\dff i}\dff\sigma$\dss
there\dss is\dss an\dss inclusion\dss morphism\dss\vspace{4pt}
\[
\quad
\Delta_{\qff \sigma\fff,\dff i}
\pff \colon\pff
e_{\dff \bullet}\dff(\trf \num{\sigma} \trf)
\qff \ttoo\qff
e_{\dff \bullet}\dff(\trf \num{\partial_{\dff i}\dff\sigma} \trf)
\pff. 
\]

\vspace{-8pt}
For\dss 
$c_{\dff \sigma}
\qff \in\qff 
e_{\dff q}\fff\left(\trf \num{\sigma} \trf\right)$\dss
we set\vspace{4pt}
\[
\quad
\delta_{\fff p}\dff \left(\trf c_{\dff \sigma}\trf\right)
\off =\off\dff
\bigoplus\nolimits_{\qff i\qff =\qff 0}^{\qff p}\off 
(\dff -\qff 1 \dff)^{\dff i}\qff
\Delta_{\qff \sigma\fff,\dff i}\dff\left(\trf
c_{\dff \sigma}
\trf\right)
\off\pff \in\off\pff\dff
\bigoplus\nolimits_{\qff i\qff =\qff 0}^{\qff p}\off
e_{\dff q}\fff\left(\trf \num{\partial_{\dff i}\dff\sigma} \trf\right)
\pff
\]

\vspace{-8pt}
and extend\dss $\delta_{\fff p}$\dss 
to\sss the direct\sss sum\dss
$c_{\dff p}\dff(\trf N\fff,\pff e_{\dff q} \trf)$\dss
by\dss linearity\halfff.\oss
Since\dss
$\Delta_{\qff \sigma\fff,\dff i}$\dss 
are\dss
morphisms of\dss complexes,\pss 
$\delta_{\fff p}$\dss is\dss a\sss morphism also.\oss
The\qss \emph{double complex}\pss
$c_{\trf \bullet}\fff(\trf N\fff,\pff e_{\trf \bullet} \trf)$\qss 
\emph{of\qss the covering}\qss $\mathcal{U}$\qss
is\dss\vspace{3.75pt}
\begin{equation}
\label{double-complex-covering-homology}
\qquad
\begin{tikzcd}[column sep=large, row sep=huge]\dis
c_{\trf 0}\dff(\trf N\fff,\pff e_{\trf 0} \trf) 
&   
c_{\trf 0}\dff(\trf N\fff,\pff e_{\trf 1} \trf) \arrow[l]
&
\off \ldots \off \arrow[l]
\\
c_{\dff 1}\dff(\trf N\fff,\pff e_{\trf 0} \trf) 
\arrow[u]
&   
c_{\dff 1}\dff(\trf N\fff,\pff e_{\trf 1} \trf) \arrow[l]
\arrow[u]
&
\off \ldots \off \arrow[l]
\\
c_{\trf 2}\dff(\trf N\fff,\pff e_{\trf 0} \trf) 
\arrow[u]
&   
c_{\trf 2}\fff(\trf N\fff,\pff e_{\trf 1} \trf) \arrow[l]
\arrow[u]
&
\off \ldots \off \arrow[l]
\\
\ldots \vphantom{C_{\trf 2}\fff(\trf N\fff,\pff e_{\trf 1} \trf)} \arrow[u] &
\ldots \vphantom{C_{\trf 2}\fff(\trf N\fff,\pff e_{\trf 1} \trf)} \arrow[u]&
\quad \vphantom{C_{\trf 2}\fff(\trf N\fff,\pff e_{\trf 1} \trf)} ,
\end{tikzcd}
\end{equation}

\vspace{-12pt}
where\sss the horizontal\sss arrows are\sss the direct\sss sums of\trs the maps\dss $d_{\dff i}$\dss
and\dss the vertical\sss arrows are\sss $\delta_{\dff i}$\nsp.\oss
Let\trs 
$t_{\trf \bullet}\dff(\trf N\fff,\pff e \trf)$\dss
be\sss the\sss total\sss complex of\dss
$c_{\trf \bullet}\fff(\trf N\fff,\pff e_{\trf \bullet} \trf)$\dss
and\dss let\dss $C_{\trf \bullet}\fff(\trf N\trf)$\dss
be\sss the complex of\dss simplicial\sss chains of\trs $N$\dss with coefficients\sss
in\dss $R$\nnsp.\oss
The boundary\dss maps\dss $d_{\trf 0}$\dss and\dss $\delta_{\trf 0}$\dss 
lead\dss to morphisms\dss\vspace{2pt}
\[
\quad
\lambda_{\dff e}\qff \colon\qff
t_{\dff \bullet}\dff(\trf N\fff,\pff e \trf) 
\qff \ttoo\qff
C_{\dff \bullet}\dff(\trf N \trf)
\quad\
\mbox{and}\quad\
\tau_{\dff e}\qff \colon\qff
t_{\dff \bullet}\dff(\trf N\fff,\pff e \trf) 
\qff \ttoo\qff
e_{\dff \bullet}\dff(\trf X \trf)
\pff,
\]

\vspace{-10pt}
where\dss it\dss is\dss understood\dss that\dss the augmentation\dss term\dss 
is\dss  removed\dss from\dss 
$e_{\dff \bullet}\dff(\trf X\trf)$\nnsp.\oss
The covering\dss $\mathcal{U}$\dss is\dss said\dss to be\dss 
\emph{$e_{\fff \bullet}$\nsp\dnsp-acyclic}\oss if\pss
$e_{\dff \bullet}\dff(\trf Z\trf)$\dss is\dss exact\dss for every\dss
$Z\qff \in\qff \mathcal{U}^{\dff \cap}$\dnsp.\oss
If\qss the covering\dss $\mathcal{U}$\dss is\dss 
$e_{\dff \bullet}$\nsp\dnsp-acyclic,\oss 
then\sss $\lambda_{\dff e}$ induces an\dss isomorphism
of\dss homology\dss groups.\oss
The proof\dss is\dss similar\dss to\sss the proof\dss of\qss
Lemma\qss \ref{acyclic-coverings},\oss
using\qss Theorem\qss \ref{homology-double-complex}\qss
instead\dss of\qss Theorem\qss \ref{double-complex}.\oss

\myuppar{Singular\sss chains.}
Suppose\sss that\dss a\sss space $\Delta$ is\dss fixed and\dss
maps\dss $s\dff \colon\dff \Delta\qff \ttoo\qff Y$\dss
are\sss treated as singular simplices.\oss 
A\qss \emph{singular\dss chain}\pss is\dss a finite formal\sss sum
of\dss singular\sss simplices with coefficients in\dss $R$\nnsp.\oss
Let\dss $c\trf(\trf Y\trf)$\dss be\sss the $R$\dnsp-module of\dss singular chains\dss
in\dss $Y$\dss
and\dss let\vspace{1.5pt}
\[
\quad
C_{\dff p}\dff(\trf N\fff,\pff c \trf)
\off\off =\off\off
\bigoplus\nolimits_{\qff \sigma\qff \in\qff N_{\dff p}}\qff 
c\dff\left(\trf \num{\sigma} \trf\right)
\pff,
\]

\vspace{-10.5pt}
where\dss $p\qff \geq\qff -\qff 1$\nnsp.\oss
The maps\qss 
$\delta_{\fff p}\dff \colon\dff
C_{\dff p}\dff(\trf N\fff,\pff c \trf)
\qff \ttoo\qff
C_{\dff p\dff -\dff 1}\dff(\trf N\fff,\pff c \trf)$\nnsp,\pss
$p\qff \geq\qff 0$\nnsp,\oss
are defined as before.\oss
Similarly\dss to\sss the cohomological\sss situation,\pss
$\delta_{\fff p\dff -\dff 1}\dff \circ\qff \delta_{\fff p}
\off =\off
0$\dss
for every\dss $p\qff \geq\qff 1$\nnsp.\oss
A\qss \emph{small\sss singular\sss chain}\pss is\dss defined as 
a\sss formal\sss sum of\dss small\sss singular\sss simplices with coefficients in\dss $R$\nnsp.\oss
Let\dss $c\trf(\qff X\fff,\pff \mathcal{U} \trf)$\dss 
be\sss the $R$\dnsp-module of\dss small\sss singular chains\dss
in\dss $X$\nnsp.\oss
There\dss is\dss an obvious map\dss
$c\trf(\qff X\fff,\pff \mathcal{U} \trf)
\qff \ttoo\qff
c\trf(\trf X \trf)$\nnsp.\oss
The inclusion\dss maps\dss
$c\trf(\trf U\trf)\qff \ttoo\qff c\trf(\trf X\trf)$\nnsp,\pss
$U\qff \in\qff \mathcal{U}$\dnsp,\oss
lead\dss to a\sss map\dss\vspace{1.5pt}
\[
\quad
\overline{\delta}_{\trf 0}\dff \colon\dff
C_{\trf 0}\dff(\trf N\fff,\pff c \trf)
\qff \ttoo\qff
c\trf(\qff X\fff,\pff \mathcal{U} \trf)
\pff.
\]

\vspace{-10.5pt}
Clearly\halfff,\pss
$\delta_{\dff 0}$\dss is\dss equal\dss to\sss the composition of\dss
$\overline{\delta}_{\trf 0}$\dss 
and\dss the inclusion\dss
$c\trf(\qff X\fff,\pff \mathcal{U} \trf)
\qff \ttoo\qff
c\trf(\trf X \trf)$\nnsp.\oss

\mypar{Lemma.}{homology-columns-are-exact}
\emph{The\qss following\dss sequence\dss is\dss exact\dff:}\vspace{-3.5pt}\vspace{-1pt}
\[
\quad
\begin{tikzcd}[column sep=large, row sep=normal]\dis
0
&
c\trf(\trf X\fff,\pff \mathcal{U}\trf) \arrow[l]
&
C_{\trf 0}\dff(\trf N\fff,\pff c \trf) \arrow[l, "\dis \off \overline{\delta}_{\trf 0}\qff"']
& 
C_{\trf 1}\dff(\trf N\fff,\pff c \trf) \arrow[l, "\dis \off \delta_{\trf 1}\qff"'] 
&   
\off \ldots\off . \arrow[l, "\dis \off \delta_{\trf 2}"'] 
\end{tikzcd}
\]

\vspace{-12pt}
\proof
It\dss is\dss sufficient\dss to\sss construct\sss 
a contracting\sss chain\dss homotopy\vspace{3pt}
\[
\quad
k_{\qff -\dff 1}\dff \colon\dff
c\trf(\trf X\fff,\pff \mathcal{U}\trf) 
\qff \ttoo\qff
C_{\trf 0}\dff(\trf N\fff,\pff c \trf)\dff,\quad\ 
k_{\trf q}\dff \colon\dff
C_{\trf q}\dff(\trf N\fff,\pff c \trf) 
\qff \ttoo\qff
C_{\trf q\dff +\dff 1}\dff(\trf N\fff,\pff c \trf) 
\dff, 
\]

\vspace{-9pt}
where\dss $q\qff \geq\qff 0$\nnsp.\oss
For every\sss small\sss singular simplex\dss 
$s\dff \colon\dff
\Delta\qff \ttoo\qff X$\dss
let\dss us\sss choose a subset\dss $U_{\dff s}\qff \in\qff \mathcal{U}$\dss
such\dss that\dss 
$s\trf(\dff \Delta\dff) 
\qff \subset\pff 
U_{\dff s}$\dss
and\sss denote by\sss $u_{\fff s}$\sss be\sss the corresponding\dss vertex of\sss $N$\nnsp.\oss
If\trs 
$s\dff(\dff \Delta\trf)
\qff \subset\qff
\num{\sigma}$\dss
for some\dss $\sigma\qff \in\qff N_{\dff p}$\nsp,\oss
we will\sss denote by\dss $s\dff *\dff \sigma$\dss 
the singular simplex\dss $s$\dss
considered as an element\dss of\trs $c\trf(\trf \num{\sigma}\trf)$\nnsp.\oss
In order\dss to define\sss 
$k_{\trf p}$\nsp,\oss 
it\dss is\dss sufficient\dss to define\sss $k_{\trf p}$\sss on\dss the chains of\trs the form\dss
$s\dff *\dff \sigma$\dss and\dss then extend\dss $k_{\trf p}$\dss 
first\dss to\dss
$c\trf(\trf \num{\sigma} \trf)$\dss
and\dss then\dss to\dss
$C_{\trf p}\trf(\trf N\fff,\pff c \trf)$\dss
by\dss linearity\halfff.\oss
Suppose\sss that\dss
$\sigma\qff \in\qff N_{\dff p}$
and\sss $s$\dss be a singular $q$\dnsp-simplex \dss in\dss $X$\dss
such\dss that\dss
$s\dff(\dff \Delta^{\fff q}\trf)
\qff \subset\qff
\num{\sigma}$\nnsp.\oss
Let\dss
$\rho
\off =\off
\sigma\dff \cup\dff \{\trf u_{\dff s}\trf\}$\nnsp.\oss
Then\dss
$s\dff(\dff \Delta\trf)
\qff \subset\qff
\num{\sigma}\qff \cap\qff U_{\dff s}
\off =\off
\num{\rho}$\dss
and,\oss in\dss \particular\halfff,\pss $\rho$\dss is\dss a simplex.\oss
If\trs $u_{\dff s}\qff \in\pff \sigma$\nnsp,\oss
then\dss $\rho$\dss is\dss a $p$\dnsp-simplex.\oss
Otherwise,\pss $\rho$\dss is\dss a $(\dff p\qff +\qff 1\dff)$\dnsp-simplex and\dss
$\sigma\off =\off \partial_{\fff a}\trf \rho$\dss for some\dss $a$\nnsp.\oss 
Let\vspace{3pt}
\[
\quad
k_{\trf p}\dff(\dff s\dff *\dff \sigma\trf)
\off =\off
0
\quad\
\mbox{if}\quad\
u_{\dff s}\qff \in\pff \sigma
\pff,
\]

\vspace{-36pt}
\[
\quad
k_{\trf p}\dff(\dff s\dff *\dff \sigma\trf)
\off =\off
(\qff -\qff 1\trf)^{\dff a}\qff 
s\dff *\dff \rho\off \in\off c\trf(\trf \num{\rho}\trf)
\quad\
\mbox{if}\quad\
u_{\dff s}\qff \not\in\pff \sigma
\pff.
\]

\vspace{-9pt}
In order\dss to prove\sss that\dss $k_{\dff \bullet}$\dss is\dss
a contracting chain\dss homotopy\halfff,\oss
it\dss is\dss sufficient\dss to prove\sss that\vspace{4.5pt}
\begin{equation}
\label{homotopy-check-homology}
\quad
\delta_{\dff p\dff +\dff 1}\trf \bigl(\trf k_{\dff p}\dff(\trf s\dff *\dff \sigma\trf)
\qff\bigr)
\off +\off
k_{\dff p\dff -\dff 1}\trf \bigl(\qff \delta_{\dff p}\dff(\trf s\dff *\dff \sigma\trf)
\qff\bigr)
\off =\off
s\dff *\dff \sigma
\pff.
\end{equation}

\vspace{-7.5pt}
By\dss rewriting\dss the definition of\dss $\delta_{\dff p}$\dss in\dss terms of\trs
the notations\dss
$s\dff *\dff \sigma$\dss we see\sss that\vspace{4.5pt}
\[
\quad
\delta_{\fff p}\dff(\trf s\dff *\dff \sigma\trf)
\off =\off\dff
\bigoplus\nolimits_{\qff i\qff =\qff 0}^{\qff p}\off 
(\dff -\qff 1 \dff)^{\dff i}\qff
s\dff *\dff \partial_{\dff i}\dff \sigma
\pff.
\]

\vspace{-7.5pt}
Suppose\sss first\dss that\dss
$u_{\dff s}\qff \in\qff \sigma$\nnsp.\oss
Let\dss
$\tau\off =\off \sigma\qff \smallsetminus\qff \{\trf u_{\dff s}\trf\}$\nnsp.\oss
Then\dss
$\sigma\off =\off \tau\dff \cup\dff \{\trf u_{\dff s}\trf\}$\dss
and\dss
$\tau\off =\off \partial_{\fff a}\trf \sigma$\dss for some\dss $a$\nnsp.\oss
Clearly\halfff,\pss
$u_{\dff s}\qff \in\qff \partial_{\dff i}\dff \sigma$\qss
if\qss $i\off \neq\off a$\dss
and\dss hence\dss 
$k_{\trf p\dff -\dff 1}\dff \left(\qff
s\dff *\dff \partial_{\dff i}\dff \sigma
\qff\right)
\off =\off
0$\dss
if\qss $i\off \neq\off a$\nnsp.\oss
Therefore\vspace{3pt}
\[
\quad
k_{\dff p\dff -\dff 1}\trf 
\bigl(\qff \delta_{\dff p}\dff(\trf s\dff *\dff \sigma\trf)\qff\bigr)
\off =\off
k_{\dff p\dff -\dff 1}\trf 
\left(\qff (\dff -\qff 1 \dff)^{\dff a}\qff
s\dff *\dff \partial_{\dff a}\dff \sigma
\qff\right)
\]

\vspace{-36pt}
\[
\quad
\phantom{k_{\dff p\dff -\dff 1}\trf 
\bigl(\qff \delta_{\dff p}\dff(\trf s\dff *\dff \sigma\trf)\qff\bigr)
\off }
=\off
(\dff -\qff 1 \dff)^{\dff a}\qff
k_{\dff p\dff -\dff 1}\trf 
\left(\qff 
s\dff *\dff \tau
\qff\right)
\off =\off
(\dff -\qff 1 \dff)^{\dff a}\qff (\dff -\qff 1 \dff)^{\dff a}\qff s\dff *\dff \sigma
\off =\off
s\dff *\dff \sigma
\pff.
\]

\vspace{-9pt}
Since also\qss 
$k_{\dff p}\trf (\trf s\dff *\dff \sigma\trf)
\off =\off
0$\nnsp,\oss 
this implies\qss (\ref{homotopy-check-homology})\qss 
in\dss the case when\dss
$u_{\dff s}\qff \in\qff \sigma$\nnsp.\oss
Suppose now\dss that\dss
$u_{\dff s}\qff \not\in\qff \sigma$\nnsp.\oss
Let\trs
$\rho
\off =\off
\sigma\dff \cup\dff \{\trf u_{\dff s}\trf\}$\dss
and\dss
let\sss $a$\sss be such\dss that\dss
$\sigma\off =\off \partial_{\fff a}\trf \rho$\nnsp.\oss
Let\dss
$\rho_{\dff i}
\off =\off
\partial_{\dff i}\dff \sigma\dff \cup\dff \{\trf u_{\dff s}\trf\}$\nnsp.\oss
Then\vspace{3pt}
\[
\quad
\rho_{\dff i}\off =\off \partial_{\dff i}\trf \rho
\quad
\mbox{for}\quad 
i\qff <\qff a
\quad 
\mbox{and}\quad
\rho_{\dff i}\off =\off \partial_{\dff i\dff +\dff 1}\trf \rho 
\quad 
\mbox{for}\quad  
i\qff \geq\qff a\pff;\quad
\]

\vspace{-36pt}
\[
\quad
\partial_{\dff i}\trf \sigma\off =\off \partial_{\fff a\dff -\dff 1}\qff \rho_{\dff i}
\quad
\mbox{for}\quad 
i\qff <\qff a
\quad 
\mbox{and}\quad
\partial_{\dff i}\trf \sigma\off =\off \partial_{\fff a}\qff \rho_{\dff i}
\quad
\mbox{for}\quad 
i\qff \geq\qff a\dff.\quad
\]

\vspace{-9.5pt}
Therefore\qss\vspace{-0.5pt}
\[
\quad
\delta_{\dff p\dff +\dff 1}\trf \bigl(\trf k_{\dff p}\dff(\trf s\dff *\dff \sigma\trf)
\qff\bigr)
\off =\off
\delta_{\dff p\dff +\dff 1}\dff\left(\trf 
(\dff -\qff 1 \dff)^{\dff a}\pff s\dff *\dff \rho\trf\right)
\off =\off
\bigoplus\nolimits_{\qff i\qff =\qff 0}^{\qff p\dff +\dff 1}\off 
(\dff -\qff 1 \dff)^{\dff i\dff +\dff a}\off
s\dff *\dff \partial_{\dff i}\trf \rho
\]

\vspace{-27pt}
\[
\quad
\off =\off
\bigoplus\nolimits_{\qff i\qff =\qff 0}^{\qff a\dff -\dff 1}\off 
(\dff -\qff 1 \dff)^{\dff i\dff +\dff a}\qff
s\dff *\dff \rho_{\dff i}
\off\off \oplus\off\off 
s\dff *\dff \sigma
\off\off \oplus\off\off 
\bigoplus\nolimits_{\qff i\qff =\qff a}^{\qff p\dff +\dff 1}\off 
(\dff -\qff 1 \dff)^{\dff i\dff +\dff a\dff -\dff 1}\qff
s\dff *\dff \rho_{\dff i}
\pff.
\]

\vspace{-9pt}
At\dss the same\sss time\vspace{1.5pt}
\[
\quad
k_{\dff p\dff -\dff 1}\trf \bigl(\qff \delta_{\dff p}\dff(\trf s\dff *\dff \sigma\trf)
\qff\bigr)
\off =\off
\bigoplus\nolimits_{\qff i\qff =\qff 0}^{\qff p}\off 
k_{\dff p\dff -\dff 1}\trf \bigl(\qff 
(\dff -\qff 1 \dff)^{\dff i}\qff
s\dff *\dff \partial_{\dff i}\trf \sigma
\qff\bigr)
\]

\vspace{-27pt}
\[
\quad
\off =\off
\bigoplus\nolimits_{\qff i\qff =\qff 0}^{\qff a\dff -\dff 1}\off 
(\dff -\qff 1 \dff)^{\dff i\dff +\dff a}\qff
s\dff *\dff \rho_{\dff i}
\off\off \oplus\off\off
\bigoplus\nolimits_{\qff i\qff =\qff a}^{\qff p\dff +\dff 1}\off 
(\dff -\qff 1 \dff)^{\dff i\dff +\dff a\dff -\dff 1}\qff
s\dff *\dff \rho_{\dff i}
\pff.
\]

\vspace{-7.5pt}
By\sss combining\dss the\sss last\dss two
calculations we see\sss that\qss (\ref{homotopy-check-homology})\qss
holds for\dss
$u_{\dff s}\qff \in\qff \sigma$\dss also.\oss  \eproof

\myuppar{The homomorphism
$H_{\dff *}\fff(\trf X \trf)
\qff \ttoo\qff 
H_{\dff *}\fff(\trf N \trf)$ for open coverings.}
Let\dss $C_{\dff q}\trf(\qff X\fff,\pff \mathcal{U} \trf)$\dss 
be\sss the $R$\dnsp-module of\dss small\sss singular $q$\dnsp-chains\dss
in\dss $X$\nnsp.\oss
Since\sss the boundary\sss of\dss a small\sss singular $q$\dnsp-simplex\dss
is\dss obviously\sss a small\sss chain,\oss
the modules\sss $C_{\dff q}\trf(\qff X\fff,\pff \mathcal{U} \trf)$\sss
together\dss with\dss the restrictions of\trs the usual\sss boundary\dss maps
form a chain complex,\oss
which\dss we\sss denote\sss by\dss 
$C_{\dff \bullet}\trf(\trf X\fff,\pff \mathcal{U}\trf)$\nnsp.\oss
Recall\dss that\dss $C_{\dff \bullet}\trf(\trf X\trf)$\dss
is\dss the usual\sss chain complex of\trs $X$\nnsp.\oss
The morphisms\vspace{-3pt}
\[
\quad
\begin{tikzcd}[column sep=large, row sep=huge]\dis
C_{\dff \bullet}\dff(\trf N \trf) 
& 
t_{\dff \bullet}\dff(\trf N\fff,\pff C \dff)
\arrow[l, "\dis \off\qff \lambda_{\qff C}"']
\arrow[r, "\dis \tau_{\trf C}\off"]
&
C_{\dff \bullet}\dff(\trf X \trf)\pff
\end{tikzcd}
\]

\vspace{-10pt}
lead\dss to homomorphisms of\dss cohomology\dss groups,\oss\vspace{-3pt}
\[
\quad
\begin{tikzcd}[column sep=large, row sep=huge]\dis
H_{\dff *}\dff(\trf N \trf) 
& 
H_{\dff *}\dff(\trf N\fff,\pff C \dff)
\arrow[l, "\dis \off\qff \lambda_{\qff C\dff *}"']
\arrow[r, "\dis \tau_{\trf C\dff *}\off"]
&
H_{\dff *}\dff(\trf X \trf)\pff
\end{tikzcd}
\]

\vspace{-10pt}
where\dss 
$H_{\dff *}\dff(\trf N\fff,\pff C \dff)$\dss
is\sss the homology\sss of\trs
$t_{\dff \bullet}\dff(\trf N\fff,\pff C \dff)$\nnsp.\oss
If\trs $\tau_{\trf C\dff *}$\dss is\dss an\sss isomorphism,\oss
we can\dss take\vspace{3pt}
\[
\quad
\lambda_{\qff C\dff *}\dff \circ\dff \tau_{\trf C\dff *}^{\dff -\dff 1}
\pff \colon\dff
H_{\dff *}\dff(\trf X \trf)
\off \ttoo\off
H_{\dff *}\dff(\trf N \trf)
\]

\vspace{-9pt}
as\sss the canonical\dss homomorphism\dss
$H_{\dff *}\fff(\trf X \trf)
\qff \ttoo\qff 
H_{\dff *}\fff(\trf N \trf)$\nnsp.\oss
This\dss is\dss the case,\oss for example,\oss when\dss $\mathcal{U}$\dss
is\dss open,\pss or\halfff,\oss at\dss least\halfff,\oss
he interiors\qss $\inte U$\dss of\qss sets\qss $U\qff \in\pff \mathcal{U}$\dss
cover\pss $X$\nnsp.\oss
Indeed,\oss
Lemma\qss \ref{homology-columns-are-exact}\qss 
with\dss $\Delta\off =\off \Delta^{q}$\dss implies\sss that\dss
the columns of\qss 
(\ref{double-complex-covering-homology})\qss
are exact\sss and\dss that\dss $\overline{\delta}_{\trf 0}$\dss
induces an\dss isomorphism of\trs the complex\dss
$C_{\dff \bullet}\trf(\trf X\fff,\pff \mathcal{U}\trf)$\dss
with\dss the cokernel\sss of\trs the morphism\vspace{1.5pt}
\[
\quad
\delta_{\dff 1}\dff \colon\dff
C_{\trf 1}\fff(\trf N\fff,\pff C_{\dff \bullet} \trf)
\qff \ttoo\qff
C_{\trf 0}\fff(\trf N\fff,\pff C_{\dff \bullet} \trf)
\pff.
\]

\vspace{-10.5pt}
Theorem\qss \ref{homology-double-complex}\qss
with\dss the rows and columns interchanged\dss
implies\sss that\dss the resulting\dss morphism\vspace{1.5pt}
\[
\quad
t_{\dff \bullet}\dff(\trf N\fff,\pff C \trf) 
\qff \ttoo\qff
C_{\dff \bullet}\trf(\trf X\fff,\pff \mathcal{U}\trf)
\pff
\]

\vspace{-10.5pt}
induces an\dss isomorphism\dss in\dss homology\halfff.\oss
If\trs the interiors\qss $\inte U$\dss of\qss sets\qss $U\qff \in\pff \mathcal{U}$\dss
cover\pss $X$\nnsp,\oss
then\dss the homological\dss version of\qss Theorem\qss \ref{eilenberg}\qss
implies\sss that\dss the inclusion\dss 
$C_{\dff \bullet}\trf(\trf X\fff,\pff \mathcal{U}\trf)
\qff \ttoo\qff
C_{\dff \bullet}\trf(\trf X\trf)$\dss
induces an\dss isomorphism\dss in\dss homology\sss
and\dss hence\dss
$\tau_{\dff C}\qff \colon\qff
t_{\dff \bullet}\dff(\trf N\fff,\pff C \trf) 
\qff \ttoo\qff
C_{\dff \bullet}\dff(\trf X \trf)$\dss
induces an\dss isomorphism\dss in\dss homology\halfff.\oss

\myuppar{The\dss $e_{\dff \bullet}$\dnsp-homology\sss and\dss the singular homology\halfff.}
The\dss \emph{$e_{\dff \bullet}$\dnsp-homology\dss groups}\qss
$\widetilde{H}_{\dff *}\dff(\trf X\trf)
\off =\off
H^{\trf e}_{\dff *}\dff(\trf X\trf)$\dss of\trs $X$\dss
are simply\dss the cohomology\dss groups of\trs the complex\dss $e_{\dff \bullet}\dff(\trf X\trf)$\dss
with\dss the\sss term\dss $R$\dss omitted.\oss 
Suppose\sss that\dss the functor\dss $e_{\dff \bullet}$\dss is\dss equipped\dss with a natural\dss
transformation\dss 
$\varphi\dff \colon\dff
C_{\dff \bullet}\qff \ttoo\qff e_{\dff \bullet}$\dnsp.\oss
This natural\dss transformation\dss leads\sss to a homomorphism\dss
$H_{\dff *}\dff(\trf X\trf)
\qff \ttoo\qff
\widetilde{H}_{\dff *}\dff(\trf X\trf)$\nnsp.\oss

\mypar{Theorem.}{e-acyclic-open}
\emph{If\qss $\mathcal{U}$\dss is\qss $e_{\dff \bullet}$\dnsp\dnsp-acyclic 
and\trs the interiors of\qss the elements of\qss $\mathcal{U}$\dss cover\trs $X$\nnsp,\oss
then\dss the\sss homomorphism\qss
$H_{\dff *}\dff(\trf X\trf)
\qff \ttoo\qff
\widetilde{H}_{\dff *}\dff(\trf X\trf)$\qss
can be factored\dss through\dss the\sss canonical\dss homomorphism\qss
$H_{\dff *}\dff(\trf X\trf)
\qff \ttoo\qff
H_{\dff *}\dff(\trf N\trf)$\dnsp.\oss}

\proof
It\dss differs\sss from\dss the proof\dss of\qss Theorem\qss \ref{A-acyclic-open}\qss
only\dss by\dss the directions of\dss arrows.\oss  \eproof

\newpage
\mysection{Extensions\qss of\pss coverings\qss and\qss bounded\qss cohomology}{extensions-coverings}

\myuppar{Weakly\sss boundedly\sss acyclic coverings.}
By\sss applying\qss Theorem\qss \ref{A-acyclic-open}\qss 
to\sss the functor\dss $B^{\dff \bullet}$\dss in\dss the role of\dss $A^{\bullet}$,\oss
we see\sss that\dss
if\trs $\mathcal{U}$\dss is\dss boundedly\sss acyclic\qss
and\dss the interiors of\trs the elements of\dss $\mathcal{U}$ cover\dss $X$\nnsp,\oss
then\dss the\sss map\dss
$\widehat{H}^{\dff *}\fff(\trf X \trf)
\qff \ttoo\qff
H^{\fff *}\fff(\trf X \trf)$\qss
can\sss be factored\dss through\dss 
$l_{\trf \mathcal{U}}\dff \colon\dff
H^{\dff *}\fff(\trf N \trf)
\qff \ttoo\qff 
H^{\dff *}\fff(\trf X \trf)$\dnsp.\oss
The fact\dss that\dss the bounded cohomology\sss
depend only\sss on\dss the fundamental\dss groups 
allows\sss to prove\sss that\dss
the same conclusion\sss holds under\sss a\sss weaker\sss assumption.\oss
Namely\halfff,\oss it\dss holds for\qss \emph{weakly\dss boundedly\sss acyclic}\qss coverings,\oss
to be defined\sss in\sss a\sss moment\halfff.\oss
See\trs Theorem\qss \ref{pi-one-extensions}\qss below.\oss

A group $\pi$ is\dss said\sss to be\qss \emph{boundedly\sss acyclic}\qss
if\trs the bounded cohomology\sss of\qss Eilenberg-MacLane space\dss 
$K\dff(\trf \pi\fff,\qff 1\dff)$\sss 
are\sss the same as\sss the bounded cohomology\sss of\dss a point.\oss 
A\sss path connected subset\sss $Z$\sss of\trs  $X$\sss 
is\dss said\dss to be\qss \emph{weakly\dss boundedly\sss acyclic}\pss if\trs 
the image of\trs the inclusion homomorphism\dss
$\pi_{\dff 1}\dff(\trf Z\fff,\qff z \trf)
\trf \ttoo\trf 
\pi_{\dff 1}\dff(\trf X\fff,\qff z \trf)$\dss
is\dss boundedly\sss acyclic,\pss
and\qss \emph{amenable}\qss if\trs this image\dss is\sss amenable.\oss
Since amenable\sss groups are boundedly acyclic,\oss
an amenable subset\dss is\dss weakly\dss boundedly\sss acyclic.\oss
A\sss covering\sss $\mathcal{U}$\sss  
is\dss said\dss to be\qss \emph{weakly\dss boundedly\sss acyclic}\pss
if\dss every\sss element\sss of\dss $\mathcal{U}^{\dff \cap}$\trs
is\dss weakly\dss boundedly\sss acyclic,\oss
and\qss \emph{amenable}\pss if\dss
every\dss element\sss of\trs $\mathcal{U}$\trs
is\dss amenable and\sss
every\dss element\sss of\trs $\mathcal{U}^{\dff \cap}$\sss
is\dss path-con\-nec\-ted.\oss
Since every\sss subgroup of\dss an amenable\sss group\dss is\dss amenable,\oss
in\dss this case every\sss element\sss of\trs $\mathcal{U}^{\dff \cap}$\sss
is\dss amenable and\dss hence
an amenable covering\dss is\dss
weakly\dss boundedly\sss acyclic.\oss

We would\dss like\sss to be able\sss to\sss turn 
a weakly\dss boundedly\sss acyclic covering\dss into
a boundedly\sss acyclic one\sss without\sss affecting\dss its\sss nerve.\oss
This can\dss be done after\dss replacing\dss $X$\dss
by\sss a\sss larger space.\oss

\myuppar{Extensions of\dss coverings.}
Let\dss $\mathcal{U}$\dss be a covering of\trs $X$\nnsp.\oss
Suppose\sss that\dss $X$\dss is\dss a subspace of\dss some other space\dss $X'$\dnsp.\oss
An\qss \emph{extension}\pss of\dss $\mathcal{U}$\dss to\dss $X'$\dss is\dss defined
as a map\dss $U\qff \longmapsto\qff U'$\dss from\dss $\mathcal{U}$\dss to\sss the
set\sss of\dss subsets of\trs $X'$\dss
such\dss that\dss $U'\qff \cap\qff X\off =\off U$\dss for\sss every\dss
$U\qff \in\qff \mathcal{U}$\dss and\dss the collection\vspace{3pt}
\[
\quad
\mathcal{U}\fff'
\off =\off
\bigl\{\pff
U'
\pff |\off
U\qff \in\qff \mathcal{U}
\pff\bigr\}
\]

\vspace{-9pt}
is\dss a covering\sss of\trs $X'$\dnsp.\oss
Since\dss $U'\qff \cap\qff X\off =\off U$ for\dss
$U\qff \in\qff \mathcal{U}$\nnsp,\oss 
the covering\dss $\mathcal{U}\fff'$\dss determines\sss the extension\dss
$U\qff \longmapsto\qff U'$\dss
and\dss we may\dss identify\dss this extension\dss 
with\dss $\mathcal{U}\fff'$\dnsp.\oss
We will\sss say\dss that\sss an extension\dss $\mathcal{U}\fff'$\sss
of\pss $\mathcal{U}$\dss is\qss \emph{nerve-preserving}\pss if\trs the conditions\vspace{3pt}
\[
\quad
\bigcap_{U\qff \in\qff \mathcal{V}} U
\off \neq\off
\varnothing
\quad\
\mbox{and}\quad\
\bigcap_{U\qff \in\qff \mathcal{V}} U'
\off \neq\off
\varnothing
\]

\vspace{-9pt}
are equivalent\dss for every\pss \emph{finite}\pss subset\dss
$\mathcal{V}\qff \subset\qff \mathcal{U}$\nnsp.\oss
Clearly\halfff,\oss if\qss $\mathcal{U}\fff'$\dss is\dss a\sss nerve-preserving\sss
extension of\qss $\mathcal{U}$\nnsp,\oss
then\dss the nerves\sss of\qss $\mathcal{U}\fff'$\dss and\qss $\mathcal{U}$\dss
are\sss the same,\oss
or\halfff,\oss rather\halfff,\oss are canonically\dss isomorphic.\oss

Suppose\sss that\dss $\mathcal{U}\fff'$\sss is\dss a\sss nerve-preserving\dss extension
of\qss $\mathcal{U}$\dss
and\dss let\trs $N$\dss be\sss the common\dss nerve of\qss
$\mathcal{U}$\dss and\trs $\mathcal{U}\fff'$\dnsp.\oss
Recall\dss that\dss
for a simplex\dss $\sigma$\dss of\dss $N$\dss we denote by\dss $\num{\sigma}$\dss 
the intersection of\dss the elements of\trs $\mathcal{U}$\dss 
corresponding\dss to\sss the vertices of\dss $\sigma$\nnsp,\oss
and\dss let\dss us\sss denote by\dss $\num{\sigma}\fff'$\dss 
the intersection of\dss the elements of\trs $\mathcal{U}\fff'$\dss 
corresponding\dss to\sss the vertices of\dss $\sigma$\nnsp.\oss
Also,\oss let\dss $\num{\varnothing}\fff'\off =\off X'$\dnsp.\oss

\mypar{Theorem.}{pi-one-extensions}
\emph{Let\trs $\mathcal{U}$\sss be a covering\dss of\pss $X$\sss
such\dss that\sss every\dss element\sss of\pss $\mathcal{U}^{\dff \cap}$\sss
is\dss path-con\-nec\-ted.\oss  
Let\qss $N$\dss be\sss the nerve of\pss $\mathcal{U}$\nnsp.\oss
Then\dss there exists a space\qss $X'\pff \supset\pff X$\dss 
and\sss a nerve-preserving\dss extension\dss
$\mathcal{U}\fff'$\dss of\oss $\mathcal{U}$\dss to\pss $X'$\qss
such\dss that\dss every\dss element\dss of\oss $\mathcal{U}\fff'^{\qff \cap}$\dss
is\dss path-con\-nec\-ted,\oss
the inclusion\dss map}\vspace{3pt}\vspace{-0.4pt}
\[
\quad
\pi_{\trf 1}\dff(\trf X\fff,\qff x\trf)
\off \ttoo\off
\pi_{\trf 1}\dff(\trf X'\fff,\qff x\trf)
\pff,
\]

\vspace{-9pt}\vspace{-0.4pt}
\emph{where\qss $x\qff \in\qff X$\nnsp,\oss
is\dss an\dss isomorphism,\oss
and\qss for\dss every\dss
simplex\dss $\sigma$\dss of\oss $N$\dss the inclusion\dss maps}\qss\vspace{3pt}\vspace{-0.4pt}
\[
\quad
\pi_{\trf 1}\dff(\trf \num{\sigma},\qff z\trf)
\off \ttoo\off
\pi_{\trf 1}\dff(\trf \num{\sigma}\fff'\fff,\qff z\trf)
\quad\
\mbox{\emph{and}}\dff\quad\
\pi_{\trf 1}\dff(\trf \num{\sigma}\fff'\fff,\qff z\trf)
\off \ttoo\off
\pi_{\trf 1}\dff(\trf X'\fff,\qff z\trf)
\pff,
\]

\vspace{-9pt}\vspace{-0.4pt}
\emph{where\qss $z\qff \in\qff \num{\sigma}$\nnsp,\qff\oss
are,\qss respectively\halfff,\qss surjective\dss and\dss injective.\oss
If\qss $\mathcal{U}$\dss is\dss an open\dss 
(respectively\halfff,\qss closed\halfff)\dss
covering,\oss
then\dss $\mathcal{U}\fff'$\dss can\dss be assumed\dss to be open\dss 
(respectively\halfff,\oss closed\halfff).\oss}

\proof
Let\dss us\sss choose for every\sss simplex\dss $\sigma$\dss of\trs $N$\dss 
a\sss collection of\qss loops\sss in\dss $\num{\sigma}$\dss 
such\dss that\dss the homotopy\sss
classes of\trs these\sss loops normally\dss generate\sss the kernel\sss of\qss
$\pi_{\dff 1}\dff(\trf \num{\sigma}\fff,\qff z \trf)
\qff \ttoo\qff 
\pi_{\dff 1}\dff(\trf X\fff,\qff z \trf)$\nnsp,\oss
where\dss
$z\qff \in\qff \num{\sigma}$\nnsp.\oss
For every\dss $U\qff \in\qff \mathcal{U}$\dss
let\dss $U'$\dss be\sss the result\sss of\dss attaching\sss a\sss two-dimensional\sss 
disc\sss to\dss $U$\dss along each\dss loop\dss from\dss this collection
contained\dss in\dss $U$\nnsp.\oss
A\sss loop\dss from\dss this collection\dss may\qss (and\dss usually\dss will\fff)\qss 
be contained\dss in several\sss sets\dss $U\qff \in\qff \mathcal{U}$\nnsp.\oss
In\dss this case we attach\dss the same disc\sss to all\dss sets\dss $U\qff \in\qff \mathcal{U}$\dss
containing\dss this\sss loop.\oss
Let\trs $X'$\dss be\sss the result\sss of\dss attaching\sss all\dss these
discs\dss to\dss $X$\nnsp.\oss

Clearly\halfff,\oss the sets\dss $U'$\dss form\sss a\sss covering\dss
$\mathcal{U}\fff'$\dss of\trs $X'$\nnsp.\oss
By\dss the construction,\pss
$\mathcal{U}\fff'$\dss is\dss an\sss extension of\qss $\mathcal{U}$\dss
to\dss $X'$\dnsp.\oss
Moreover\halfff,\oss several\sss elements of\dss $\mathcal{U}\fff'$\dss
intersect\dss if\trs and\dss only\trs if\trs the corresponding\sss
elements of\dss $\mathcal{U}$\dss intersect\halfff,\oss
i.e.\dss $\mathcal{U}\fff'$\dss
is\dss a\sss nerve-preserving\sss extension\sss of\dss $\mathcal{U}$\nnsp.\oss
Since\sss the discs are attached\dss to\dss $X$\dss along\dss loops
contractible in\dss $X$\nnsp,\oss
the inclusion\dss map\dss
$\pi_{\dff 1}\dff(\trf X\fff,\qff x \trf)
\qff \ttoo\qff 
\pi_{\dff 1}\dff(\trf X'\fff,\qff x \trf)$\dss
is\dss an\sss isomorphism.\oss 
We\sss will\dss use\dss it\dss to\sss identify\dss 
the fundamental\dss groups\dss
$\pi_{\dff 1}\dff(\trf X\fff,\qff x \trf)$\dss
and\dss
$\pi_{\dff 1}\dff(\trf X'\fff,\qff x \trf)$\nnsp.\oss

Let\sss $\sigma$\sss be a simplex\sss of\dss $N$\nnsp.\oss
Clearly\halfff,\oss the inclusion\dss map\dss
$\pi_{\trf 1}\dff(\trf \num{\sigma}\fff,\qff z\trf)
\off \ttoo\off
\pi_{\trf 1}\dff(\trf  \num{\sigma}\fff'\fff,\qff z\trf)$\dss
is\dss surjective.\oss
Among\dss the\sss loops used\dss to attach discs\sss there are\sss
loops\sss in\dss $\num{\sigma}$\dss such\dss that\dss their\dss homotopy\sss
classes generate\sss the kernel\sss of\qss
$\pi_{\dff 1}\dff(\trf \num{\sigma}\fff,\qff z \trf)
\qff \ttoo\qff 
\pi_{\dff 1}\dff(\trf X\fff,\qff z \trf)$\nnsp.\oss
The discs attached\dss to\sss these loops are contained\dss in\dss $U'$\dss
for each\dss $U\qff \in\qff \mathcal{U}$\dss corresponding\dss
to a vertex of\dss $\sigma$\nnsp.\oss
Therefore\sss these discs are contained\dss in\dss $\num{\sigma}\fff'$\dnsp.\oss
It\dss follows\dss that\dss the\sss kernel\sss 
of\trs the inclusion\dss map\dss\vspace{3pt}\vspace{-0.4pt}
\[
\quad
\pi_{\trf 1}\dff(\trf \num{\sigma}\fff,\qff z\trf)
\off \ttoo\off
\pi_{\trf 1}\dff(\trf X\fff,\qff z\trf)
\off\qff =\off\qff
\pi_{\trf 1}\dff(\trf X'\fff,\qff z\trf)
\pff
\]

\vspace{-9pt}\vspace{-0.4pt}
is\dss contained\dss in\dss the\sss kernel\sss of\trs the inclusion\dss map\dss 
$\pi_{\trf 1}\dff(\trf \num{\sigma}\fff,\qff z\trf)
\off \ttoo\off
\pi_{\trf 1}\dff(\trf  \num{\sigma}\fff'\fff,\qff z\trf)$\nnsp.\oss
Since\sss the\sss latter\dss is\dss surjective,\oss
this implies\sss the injectivity\sss of\vspace{3pt}\vspace{-0.4pt}
\[
\quad
\pi_{\trf 1}\dff(\trf \num{\sigma}\fff'\fff,\qff z\trf)
\off \ttoo\off
\pi_{\trf 1}\dff(\trf  X'\fff,\qff z\trf)
\pff.
\]

\vspace{-9pt}\vspace{-0.4pt}
It\dss remains\sss to prove\sss the\sss last\sss statement\sss of\trs the\sss theorem.\oss
Clearly\halfff,\oss if\qss $\mathcal{U}$\dss is\dss a\sss closed\sss covering,\oss
then\dss $\mathcal{U}\fff'$\dss is\dss also closed.\oss
In\dss general,\pss $\mathcal{U}\fff'$\dss is\dss not\sss an\sss open covering
even\dss if\trs $\mathcal{U}$\dss is\dss open.\oss
But\dss if\trs $\mathcal{U}$\dss is\dss open,\oss
one can\dss turn\dss $\mathcal{U}'$\dss into an\sss open covering\dss
without\dss affecting\dss $X'$\dss and\dss the fundamental\dss groups of\dss intersections.\oss
In order\dss to do\sss this,\oss let\dss us\dss for each attached disc 
remove from\dss $X'$\dss
a closed disc with\dss the same center and smaller radius.\oss
Let\dss $X''$\dss be\sss the result\halfff.\oss
Clearly\halfff,\pss $X''$\dss is\dss open\dss in\dss  
$X'$\dss
and\dss there\dss is\dss a\sss natural\qss ``radial''\qss retraction\dss
$r\dff \colon\dff
X''\qff \ttoo\qff X$\nnsp.\oss
For every\dss $U\qff \in\qff \mathcal{U}$\dss
let\dss\vspace{1.5pt}
\[
\quad
U''
\off =\off\dff
U'\qff \cup\pff r^{\dff -\dff 1}\dff(\trf U\trf)
\pff.
\]

\vspace{-10.5pt} 
Then\dss $U''$\dss is\dss open\sss and\dss $U'$\dss is\dss a deformation\dss
retract\sss of\dss $U''$\dnsp.\oss
Clearly\halfff,\oss
every\dss finite intersection of\dss sets\dss $U''$\dss
has\sss the corresponding\dss finite intersection of\dss sets\dss $U'$\dss
as a deformation\dss retract\halfff.\oss
Hence\vspace{1.5pt}
\[
\quad
\mathcal{U}\fff''
\off =\off
\bigl\{\pff
U''
\pff |\off
U\qff \in\qff \mathcal{U}
\pff\bigr\}
\]

\vspace{-10.5pt}
is\dss an open covering of\trs $X'$\dss and\dss has all\sss already\sss
established\dss properties of\trs $\mathcal{U}\fff'$\dnsp.\oss  \eproof

\mypar{Corollary.}{amenable-extensions}
\emph{If\trs $\mathcal{U}$ be\dss a\sss weakly\sss boundedly\sss acyclic covering\dss of\pss $X$\nnsp,\pss
then\dss there exists\qss $X'\pff \supset\pff X$\dss 
and\sss a\sss nerve-preserving\dss extension\qss
$\mathcal{U}\fff'$\dss of\oss $\mathcal{U}$\dss to\pss $X'$\qss
such\dss that\qss $\mathcal{U}\fff'$\dss is\dss boundedly\sss acyclic\dss and\dss}\vspace{3pt}
\[
\quad
\pi_{\trf 1}\dff(\trf X\fff,\qff x\trf)
\off \ttoo\off
\pi_{\trf 1}\dff(\trf X'\fff,\qff x\trf)
\pff,
\]

\vspace{-9pt}
\emph{where\qss $x\qff \in\qff X$\nnsp,\oss
is\dss an\dss isomorphism.\oss
If\qss $\mathcal{U}$\dss is\dss an open\dss 
(respectively\halfff,\qss closed\halfff)\dss
covering,\oss
then\dss $\mathcal{U}\fff'$\dss can\dss be assumed\dss to be open\dss 
(respectively\halfff,\oss closed\halfff).\oss}

\proof
Let\dss $X'$\dss and\dss $\mathcal{U}\fff'$\dss be\sss the space and\dss
the covering\dss provided\dss by\qss Theorem\qss \ref{pi-one-extensions}.\oss
Then\dss the inclusion\dss map\dss
$\pi_{\dff 1}\dff(\trf X\fff,\qff x \trf)
\qff \ttoo\qff 
\pi_{\dff 1}\dff(\trf X'\fff,\qff x \trf)$\dss
is\dss an\sss isomorphism and we can use\sss it\dss to identify\dss
$\pi_{\dff 1}\dff(\trf X\fff,\qff x \trf)$\dss
and\dss
$\pi_{\dff 1}\dff(\trf X'\fff,\qff x \trf)$\nnsp.\oss
Let\sss $\sigma$\sss be\dss an arbitrary\sss 
simplex of\trs the nerve of\pss $\mathcal{U}$\nnsp.\oss
Since\vspace{3pt}
\[
\quad
\pi_{\trf 1}\dff(\trf \num{\sigma},\qff z\trf)
\off \ttoo\off
\pi_{\trf 1}\dff(\trf \num{\sigma}\fff'\fff,\qff z\trf)
\]

\vspace{-9pt}
is\dss surjective,\oss
the images of\trs
$\pi_{\trf 1}\dff(\trf \num{\sigma},\qff z\trf)$\dss
and\dss
$\pi_{\trf 1}\dff(\trf \num{\sigma}\fff'\fff,\qff z\trf)$\dss
in\dss
$\pi_{\trf 1}\dff(\trf X\fff'\fff,\qff z\trf)$\dss
are equal.\oss
Since\dss $\mathcal{U}$\dss is\dss weakly\sss boundedly\sss acyclic,\oss
the image of\trs $\pi_{\trf 1}\dff(\trf \num{\sigma},\qff z\trf)$\dss
is\dss boundedly\sss acyclic.\oss
It\dss follows\dss that\dss the image\sss $\pi_{\dff \sigma}$\sss of\trs
$\pi_{\trf 1}\dff(\trf \num{\sigma}\fff'\fff,\qff z\trf)$\dss
is\dss boundedly\sss acyclic.\oss
But\sss since\vspace{3pt}
\[
\quad
\pi_{\trf 1}\dff(\trf \num{\sigma}\fff',\qff z\trf)
\off \ttoo\off
\pi_{\trf 1}\dff(\trf X'\fff,\qff z\trf)
\]

\vspace{-9pt}
is\dss injective,\qss
$\pi_{\dff \sigma}$
is\dss isomorphic\sss to
$\pi_{\trf 1}\dff(\trf \num{\sigma}\fff',\qff z\trf)$\nnsp.\pss
Since\sss $\widehat{H}^{\dff *}\fff(\trf \num{\sigma}\fff' \trf)$\sss
is\dss determined\sss by 
$\pi_{\trf 1}\dff(\trf \num{\sigma}\fff',\qff z\trf)$\nnsp,\oss
the subset\dss $\num{\sigma}\fff'$\dss
is\dss a\sss boundedly\sss acyclic.\oss 
It\dss follows\sss that\dss $\mathcal{U}\fff'$\dss 
is\dss boundedly\sss acyclic.\oss  \eproof

\mypar{Theorem.}{covering-theorem-open}
\emph{If\qss $\mathcal{U}$\dss is\dss an\sss 
open\sss weakly\sss boundedly\sss acyclic\sss covering\halfff,\oss
then\dss the\sss canonical\dss homomorphism\qss
$\widehat{H}^{\dff *}\fff(\trf X \trf)
\qff \ttoo\qff
H^{\fff *}\fff(\trf X \trf)$\qss
can be factored\dss through\qss
$H^{\dff *}\fff(\trf N \trf)
\qff \ttoo\qff 
H^{\dff *}\fff(\trf X \trf)$\dnsp.\oss
In\dss particular\halfff,\oss this\dss is\dss true\dss if\qss
$\mathcal{U}$\dss is\dss an\sss open amenable covering.\oss }

\proof
Theorem\qss \ref{A-acyclic-open}\qss
this\sss implies\sss that\dss the conclusion of\trs the\sss theorem\dss
holds\dss if\dss $\mathcal{U}$\dss is\dss boundedly\sss acyclic.\oss
If\dss $\mathcal{U}$\dss is\dss only\sss weakly\sss boundedly\sss acyclic,\oss
Corollary\qss \ref{amenable-extensions}\qss provides\sss us\dss
with a space\dss $X'\qff \supset\qff X$\dss and an open\sss boundedly\sss acyclic
covering\trs $\mathcal{U}\fff'$\sss of\trs $X'$\dss with\dss the same nerve\sss $N$\nnsp.\oss
Moreover\halfff,\oss the inclusion\dss 
$i\dff \colon\dff
X\qff \ttoo\qff X'$\dss induces an\dss isomorphism of\trs
fundamental\dss groups and\dss hence an\dss isomorphism of\trs
bounded cohomology\dss groups.\oss
The already\dss proved case of\trs the\sss theorem\dss 
implies\sss that\dss the canonical\dss map\dss 
${\widehat{H}^{\dff *}\fff(\trf X' \trf)
\qff \ttoo\qff
H^{\dff *}\dff(\trf X'\trf)}$\dss 
is\dss equal\dss to\sss the composition\vspace{1.5pt}
\[
\quad
\begin{tikzcd}[column sep=large, row sep=huge]\dis
\widehat{H}^{\dff *}\fff(\trf X' \trf) 
\arrow[r]
&
H^{\dff *}\dff(\trf N \trf)
\arrow[r]
&
H^{\dff *}\dff(\trf X '\trf)\pff.
\end{tikzcd}
\]

\vspace{-10.5pt}
Next\halfff,\oss
the construction of\trs the horizontal\sss arrows\sss 
in\dss the square\vspace{1.25pt}\vspace{-0.4pt}
\[
\quad
\begin{tikzcd}[column sep=large, row sep=boom]\dis
H^{\dff *}\dff(\trf N \trf)
\arrow[r]
\arrow[d, "\dis \qff = "]
&
H^{\dff *}\dff(\trf X '\trf)
\arrow[d]
\\
H^{\dff *}\dff(\trf N \trf) 
\arrow[r]
&
H^{\dff *}\dff(\trf X \trf)
\end{tikzcd}
\]

\vspace{-10.75pt}\vspace{-0.4pt}
shows\sss that\dss this square\dss is\dss commutative.\oss
This\sss leads\sss to\sss the diagram\vspace{1.25pt}\vspace{-0.4pt}
\[
\quad
\begin{tikzcd}[column sep=large, row sep=boom]\dis
\widehat{H}^{\dff *}\fff(\trf X' \trf) 
\arrow[r]
\arrow[d, "\dis \qff \widehat{i}^{\qff *}"]
&
H^{\dff *}\dff(\trf N \trf)
\arrow[r]
\arrow[d, "\dis \qff ="]
&
H^{\dff *}\dff(\trf X '\trf)
\arrow[d, "\dis \qff i^{\fff *}"]
\\
\widehat{H}^{\dff *}\fff(\trf X \trf)
\arrow[r, dashed]
& 
H^{\dff *}\dff(\trf N \trf) 
\arrow[r]
&
H^{\dff *}\dff(\trf X \trf)
\end{tikzcd}
\]

\vspace{-10.75pt}\vspace{-2.75pt}
of\dss solid arrows\sss with commutative right\sss square.\oss
Since\sss $\widehat{i}^{\qff *}$\dss is\dss an\dss isomorphism,\oss
this diagram can\sss be completed\dss by\sss a dashed arrow\dss to a commutative diagram.\oss
The\sss theorem\dss follows.\oss  \eproof

\myuppar{Coverings amenable in\dss the sense of\qss Gromov.}
Gromov's\pss \cite{gro}\qss notion of\dss an amenable subset\dss is\dss
different\halfff.\oss
Namely\halfff,\pss he calls\dss $Z$\dss amenable if\trs every\dss path-connected
component\sss of\trs $Z$\dss is\dss amenable\sss in\dss
our sense.\oss
Let\dss us\sss call\sss such subsets\qss \emph{amenable\sss in\dss the sense of\qss Gromov},\oss
and\dss call\dss $\mathcal{U}$\qss 
\emph{amenable\sss in\dss the sense of\qss Gromov}\pss
if\dss elements\sss of\dss $\mathcal{U}$\sss
are amenable in\dss the sense of\qss Gromov.\oss

\mypar{Theorem\qss ({\fff}Vanishing\dss theorem).}{vanishing}
\emph{If\pss $\mathcal{U}$\dss is\dss an open covering amenable in\dss the sense of\pss Gromov,\oss
then\dss the canonical\dss homomorphism\oss
$\widehat{H}^{\dff p}\fff(\dff X \dff)
\qff \ttoo\qff
H^{\dff p}\fff(\dff X \dff)$\qss
vanishes\qss for\qss $p\qff >\qff \dim\dff N$\nnsp.\oss}

\proof
Now we know only\dss that\dss path components\sss of\dss sets\sss $\num{\sigma}$\sss
are boundedly acyclic.\oss 
Let\dss us\sss replace in\sss the complex\qss (\ref{b-sigma})\qss with\sss
$A^{\dff \bullet}\off =\off B^{\dff \bullet}$\sss the\sss term $R\off =\off \rrr$\sss
by\sss the product\sss of\dss copies of\dss $\rrr$\sss 
corresponding\sss to path components of\dss $\num{\sigma}$\nnsp.\oss
This forces us\sss to replace\sss the spaces\sss $C^{\dff p}\dff(\trf N\trf)$\sss
in\qss (\ref{c-nerve})\qss and\qss (\ref{big-diagram})\qss by\sss 
some other spaces which are still\sss equal\dss to $0$\sss for\sss $p\qff >\qff \dim\dff N$\nnsp.\oss
Now the cohomology\sss spaces\sss $H^{\dff p}$\sss 
of\trs the\sss left\sss column of\qss (\ref{big-diagram})\qss
may\sss be not\sss equal\dss to\sss $H^{\dff p}\dff(\trf N \trf)$\nnsp,\oss but,\oss
obviously,\pss 
$H^{\dff p}\off =\off 0$\sss if\dss $p\qff >\qff \dim\dff N$\nnsp.\oss 
The rest\sss of\dss our arguments apply\sss and show\sss that\dss
$\widehat{H}^{\dff p}\fff(\dff X \dff)
\qff \ttoo\qff
H^{\dff p}\fff(\dff X \dff)$\sss
factors\sss through\sss $H^{\dff p}$\dnsp.\oss
The\sss theorem\sss follows.\oss  \eproof

\newpage
\mysection{A\qss Leray\qss theorem\pss for\qss $l_{\dff 1}$\dnsp-homology}{l-one-homology}

\myuppar{The\dss $l_{\dff 1}$\dnsp-norm\sss of\dss chains.}
The\dss $l_{\fff 1}$\dnsp-norm\qss $\norm{c\fff}_{\dff 1}$\qss
of\dss a singular chain\qss\vspace{3pt} 
\[
\quad 
c
\off\off =\off\off
\sum_{\qff \sigma} \dff c_{\dff \sigma}\dff \sigma
\off \in\off C_{\dff m}\fff(\dff X \dff)\dff,
\]

\vspace{-12pt} 
where\qss $c_{\dff \sigma}\qff \in\qff \rrr$\qss 
and\dss the sum\dss is\dss taken over all\sss singular $m$\dnsp-simplices\sss $\sigma$\sss in\dss 
$X$\nnsp,\oss
is defined as\oss\vspace{3pt} 
\[
\quad 
\norm{c\fff}_{\dff 1}
\off\off =\off\off
\sum_{\qff \sigma} \qff \num{c_{\dff \sigma}}
\]

\vspace{-12pt}
By\dss the definition of\dss singular chains,\oss 
the sums above involve only\sss a finite number
of\dss non-zero coefficients\dss $c_{\dff \sigma}$\nnsp,\oss
and\dss hence\dss $\norm{c\fff}_{\dff 1}$\dss is\dss well\sss defined.\oss 
The normed space\dss $C_{\dff m}\dff(\trf X\trf)$\dss is\dss not\sss complete.\oss
Let\dss $L_{\dff m}\dff(\trf X\trf)$\dss be its completion\dss with respect\dss to\sss the
$l_{\dff 1}$\dnsp-norm.\oss
The boundary\sss operators\vspace{3pt} 
\[
\quad 
\partial_{\fff m}\dff \colon\dff
C_{\dff m}\dff(\trf X\trf)
\qff \ttoo\qff
C_{\dff m\dff -\dff 1}\dff(\trf X\trf)
\]

\vspace{-9pt}
are obviously\dss bounded and\dss hence\sss extend\dss by\dss continuity\dss to
bounded\sss operators\dss\vspace{3pt} 
\[
\quad 
d_{\dff m}\dff \colon\dff
L_{\dff m}\dff(\trf X\trf)
\qff \ttoo\qff
L_{\dff m\dff -\dff 1}\dff(\trf X\trf)
\pff.
\]

\vspace{-9pt}
Since\dss
$C_{\dff m}\dff(\trf X\trf)$\dss is\dss dense\sss in\dss $L_{\dff m}\dff(\trf X\trf)$\nnsp,\oss
the identity\dss
$\partial_{\fff m\dff -\dff 1}\fff \circ\qff \partial_{\fff m}
\off =\off
0$\qss
implies\dss
the identity\dss\vspace{3pt} 
\[
\quad 
d_{\dff m\dff -\dff 1}\fff \circ\qff d_{\dff m}
\off =\off
0
\pff.
\]

\vspace{-9pt}
Therefore\sss the sequence\vspace{0pt}
\[
\quad
\begin{tikzcd}[column sep=large, row sep=normal]\dis
0 
& 
L_{\trf 0}\dff(\trf X \trf) \arrow[l]
& 
L_{\trf 1}\dff(\trf X \trf) \arrow[l, "\dis \off \partial_{\dff 1}"']
&   
L_{\trf 2}\dff(\trf X \trf) \arrow[l, "\dis \off \partial_{\dff 2}"']
&
\off \ldots \off\off \arrow[l, "\dis \off \partial_{\dff 3}"']
\end{tikzcd}
\]

\vspace{-9pt}
is\dss a chain complex.\oss
The\dss \emph{$l_{\dff 1}$\dnsp-homology\dss spaces}\pss of\dss $X$\dss are defined as 
the homology spaces of\dss this complex
and are denoted\dss by\dss
$H^{\dff l_{\dff 1}}_{\dff *}\dff(\dff X \dff)$\dnsp.\oss 
In more details,\vspace{3pt}\vspace{0.5pt}
\[
\quad
H^{\dff l_{\dff 1}}_{\dff m}\dff(\trf X \trf) 
\off =\off
\kernel\qff d_{\dff m}
\qff\bigl/\qff
\image\qff
d_{\dff m\dff +\dff 1}
\pff.
\]

\vspace{-9pt}\vspace{0.5pt}
The $l_{\dff 1}$\dnsp-homology\dss spaces\dss $H^{\dff l_{\dff 1}}_{\dff m}\dff(\trf X \trf)$\dss
are real vector spaces carrying a\qss \emph{canonical\dss semi-norm}\qss
induced\dss by\dss the $l_{\fff 1}$\dnsp-norm on\sss $L_{\trf m}\dff(\trf X \trf)$\nnsp.\oss
If\trs the image of\trs the boundary\sss operator\vspace{3pt}\vspace{0.5pt} 
\[
\quad 
d_{\dff m\dff +\dff 1}\dff \colon\dff
L_{\dff m\dff +\dff 1}\dff(\trf X\trf)
\qff \ttoo\qff
L_{\dff m}\dff(\trf X\trf)
\pff
\]

\vspace{-9pt}\vspace{0.5pt}
is\dss not\sss closed,\oss
then\dss 
$H^{\dff l_{\dff 1}}_{\dff m}\dff(\trf X\trf)$\dss
contains non-zero homology\sss classes with\dss the norm\dss $0$\dss and\dss
the canonical\sss semi-norm on\dss
$H^{\dff l_{\dff 1}}_{\dff m}\dff(\trf X\trf)$\dss
is\dss not\sss a norm.\oss

\myuppar{The duality\dss between\dss
$L_{\fff \bullet}\dff(\trf X\trf)$\dss
and\dss
$B^{\fff \bullet}\dff(\trf X\trf)$\nnsp.}
The space\dss
$L_{\dff m}\dff(\trf X\trf)$\dss
is\dss nothing else\sss but\dss the space of\dss real-valued
$l_{\dff 1}$\hnsp\dnsp-functions on\dss the set $S_{\dff m}\dff(\trf X\trf)$
of\dss singular $m$\dnsp-simplices in\dss $X$\nnsp.\oss
Similarly\halfff,\pss
$B^{\dff m}\dff(\trf X\trf)$
is\dss the space of\dss real-valued
$l_{\dff \infty}$\hnsp\dnsp-functions on\sss $S_{\dff m}\dff(\trf X\trf)$\nnsp.\oss
Let\sss $L_{\dff m}\dff(\trf X\trf)^{\fff *}$\dss
be\sss the\dss Banach space dual\dss to\dss $L_{\dff m}\dff(\trf X\trf)$\nnsp.\oss
The standard\dss pairing\dss between $m$\dnsp-chains and $m$\dnsp-cochains extends\sss
to a pairing\dss\vspace{1.5pt} 
\[
\quad 
\langle\qff \bullet\fff,\pff \bullet \qff\rangle
\qff \colon\qff
L_{\dff m}\dff(\trf X\trf)
\pff \times\off
B^{\dff m}\dff(\trf X\trf)
\qff \ttoo\qff
\rrr
\pff
\] 

\vspace{-10.5pt}
and\dss leads\sss to a map\dss 
$L_{\dff m}\dff(\trf X\trf)^{\fff *}
\qff \ttoo\qff
B^{\fff m}\dff(\trf X\trf)$\nnsp.\oss
By\dss the well\dss known duality\dss between\dss $l_{\dff 1}$\dss
and\dss $l_{\dff \infty}$\dss spaces\sss this map\dss
is\dss an\dss isometry\sss of\trs Banach spaces.\oss
Let\dss us\sss use\sss this map\sss to\sss identify\dss $B^{\fff m}\dff(\trf X\trf)$\dss
with\dss $L_{\dff m}\dff(\trf X\trf)^{\fff *}$\dnsp.\oss
By\dss the very\sss definition,\oss this\sss turns\sss the coboundary\sss operator\vspace{1.5pt} 
\[
\quad 
\partial^{\dff m\dff -\dff 1}\dff \colon\dff
B^{\fff m\dff -\dff 1}\dff(\trf X\trf)
\qff \ttoo\qff
B^{\fff m}\dff(\trf X\trf)
\pff
\]

\vspace{-10.5pt}
into\sss the\dss Banach\dss space adjoint\sss
of\trs the boundary\sss operator\dss
$\partial_{\dff m}\dff \colon\dff
L_{\fff m}\dff(\trf X\trf)
\qff \ttoo\qff
L_{\fff m\dff -\dff 1}\dff(\trf X\trf)$\nnsp.\oss
Therefore\sss the complex\dss $B^{\fff \bullet}\dff(\trf X\trf)$\dss
is\dss the\dss Banach\sss dual\dss of\trs the complex\dss
$L_{\fff \bullet}\dff(\trf X\trf)$\nnsp.\oss

This duality\dss between\dss 
$L_{\fff \bullet}\dff(\trf X\trf)$\dss
and\dss
$B^{\fff \bullet}\dff(\trf X\trf)$\dss
leads\sss to an\dss imperfect\sss duality between\dss the
$l_{\dff 1}$\dnsp-ho\-mol\-o\-gy\sss and\dss the bounded cohomology\halfff.\oss
In\dss fact\halfff,\oss already\dss the duality\dss between\dss 
$l_{\dff 1}$\dss
and\dss $l_{\dff \infty}$\dss spaces\dss is\dss not\dss perfect\halfff.\oss
While\dss $l_{\dff \infty}$\dss spaces are\sss the duals of\trs the corresponding\dss
$l_{\dff 1}$\dss spaces,\pss
$l_{\dff 1}$\dss spaces are\sss not\sss the duals\dss of\dss
$l_{\dff \infty}$\dss spaces.\oss
When one passes\sss to\sss the cohomology\sss spaces,\oss
an additional\sss difficulty\sss arises.\oss
The canonical\dss pairing\dss
$\langle\qff \bullet\fff,\pff \bullet \qff\rangle$\dss
between\dss $L_{\dff m}\dff(\trf X\trf)$\dss and\dss $B^{\fff m}\dff(\trf X\trf)$\dss
leads\sss to a pairing\vspace{1.5pt} 
\[
\quad 
H^{\dff l_{\dff 1}}_{\dff m}\dff(\trf X\trf)
\pff \times\off
\widehat{H}^{\dff m}\dff(\trf X\trf)
\qff \ttoo\qff
\rrr
\pff,
\]

\vspace{-10.5pt}
but\halfff,\oss in\dss general,\oss
both\dss $H^{\dff l_{\dff 1}}_{\dff m}\dff(\trf X\trf)$\dss
and\dss
$\widehat{H}^{\dff m}\dff(\trf X\trf)$\dss
contain\dss non-zero elements with\dss the canonical\sss semi-norm equal\dss to $0$\nnsp,\oss
and\dss the pairing of\dss such elements with every\sss element\sss of\trs the
other space\dss is\dss $0$\nnsp.\oss
Nevertheless,\oss sometimes\dss
$H^{\dff l_{\dff 1}}_{\dff m}\dff(\trf X\trf)$\dss
and\dss
$\widehat{H}^{\dff m}\dff(\trf X\trf)$\dss
behave\sss like a\sss true duality\dss between\dss them exists.\oss
C.\qss L\"{o}h\qss \cite{l}\qss provided a systematic exploration of\trs this duality\halfff.\oss
We will\dss limit\sss ourselves by\dss the following observation of\pss
Sh.\dss Matsumoto\dss and\qss Sh.\dss Morita\qss \cite{mm}.\oss

\mypar{Theorem.}{zeros-are-dual}
\emph{Let\qss $m\qff \geq\qff 1$\nnsp.\oss
If\oss $\widehat{H}^{\dff m}\dff(\trf X\trf)
\off =\off
\widehat{H}^{\dff m\dff +\dff 1}\dff(\trf X\trf)
\off =\off 
0$\nnsp,\oss
then\qss
$H^{\dff l_{\dff 1}}_{\dff m}\dff(\trf X\trf)\off =\off 0$\nnsp.\oss
Also,\oss if\oss $H^{\dff l_{\dff 1}}_{\dff m}\dff(\trf X\trf)
\off =\off 
H^{\dff l_{\dff 1}}_{\dff m\dff +\dff 1}\dff(\trf X\trf)
\off =\off
0$\nnsp,\oss
then\qss
$\widehat{H}^{\dff m\dff +\dff 1}\dff(\trf X\trf)\off =\off 0$\nnsp.\oss}

\proof
Given\dss vector subspaces\qss
$U\qff \subset\pff L_{\fff m}\dff(\trf X\trf)$\qss and\qss 
$V\qff \subset\pff B^{\fff m}\dff(\trf X\trf)$\nnsp,\qff\oss
let\vspace{3pt}\vspace{-0.125pt} 
\[
\quad
U^{\dff \perp}
\off =\off
\bigl\{\off
u\qff \in\qff B^{\fff m}\dff(\trf X\trf)
\off \bigl|\off
\langle\qff c\fff,\pff u \qff\rangle
\off =\off
0
\off\pff
\mbox{for\dss every}\off\qff
c\qff \in\pff U
\off\bigr\}
\pff,
\]

\vspace{-36pt}
\[
\quad
{}^{\perp}\dff V
\off\trf =\off
\bigl\{\off
c\pff\fff \in\qff L_{\dff m}\dff(\trf X\trf)
\off \bigl|\off
\langle\qff c\fff,\pff u \qff\rangle
\off =\off
0
\off\pff
\mbox{for\dss every}\off\qff
u\qff \in\pff V
\off\bigr\}
\pff
\]

\vspace{-9pt}\vspace{-0.125pt}
be\sss their orthogonal\sss complements with respect\dss to\sss 
$\langle\qff \bullet\fff,\pff \bullet \qff\rangle$\nnsp.\oss
Recall\dss that\dss the complements\dss $U^{\dff \perp}$\dss and\dss ${}^{\perp}\dff V$\dss
are always closed.\oss
Also,\oss
the double complement\dss
${}^{\perp}\dff (\qff U^{\dff \perp}\trf)$\dss is\dss equal\dss to\sss the closure of\dss $U$\dss
with\dss respect\dss to\sss the $l_{\dff 1}$\dnsp-norm,\oss
but\dss $(\trf {}^{\perp}\trf V\qff)^{\dff \perp}$\dss is\dss
the closure of\dss $V$\dss only\sss in\sss a\sss weaker sense.\oss
Namely\halfff,\pss $(\trf {}^{\perp}\trf V\qff)^{\dff \perp}$\dss
is\dss the\sss weak\dnsp$^*$\dnsp-closure of\dss $V$\nnsp.\oss
See\qss \cite{r},\oss Theorem\qss 4.7.\oss

Since\sss boundary\sss and coboundary\sss operators are bounded\sss operators,\oss
their\dss kernels are closed.\oss
Moreover\halfff,\oss
since\dss $\partial^{\dff m}$\dss is\dss 
the adjoint\sss operator of\dss $\partial_{\dff m\dff +\dff 1}$\nsp,\oss\vspace{3pt}
\[
\quad
\kernel\qff \partial^{\dff m}
\off =\off
\kernel\qff \partial_{\dff m\dff +\dff 1}^{\dff *}
\off =\off
\bigl(\qff \image \partial_{\dff m\dff +\dff 1}\trf\bigr)^{\dff \perp}
\quad\
\mbox{and}
\]

\vspace{-36pt}
\[
\quad
\kernel\qff \partial_{\dff m\dff +\dff 1}
\off =\off
{}^{\perp}\dff \left(\qff \image \partial^{\dff m}\trf\right)
\pff.
\]

\vspace{-9pt}
See\qss \cite{r},\oss Theorem\qss 4.12.\oss
Also,\pss $\image \partial^{\dff m}$\dss is\dss closed\dss 
with\dss respect\dss to\sss the $l_{\dff \infty}$\dnsp-norm
if\trs and\dss only\trs if\trs $\image \partial_{\dff m\dff +\dff 1}$\dss
is\dss closed\dss with\dss respect\dss to\sss the $l_{\dff 1}$\dnsp-norm.\oss
In\dss this case\dss $\image \partial^{\dff m}$\dss is\dss also 
weak\dnsp$^*$\dnsp-closed.\oss
This\dss is\dss a special\sss case of\trs the closed\dss range\sss theorem.\oss
See\qss \cite{r},\oss Theorem\qss 4.14.\oss

Suppose\sss that\qss 
$H^{\dff l_{\dff 1}}_{\dff m}\dff(\trf X\trf)\off =\off 0$\nnsp.\oss
Then\dss
$\image \partial_{\dff m\dff +\dff 1}
\off =\off
\kernel\trf \partial_{\dff m}$\dss and\dss hence\dss
$\image \partial_{\dff m\dff +\dff 1}$\dss is\dss closed.\oss
By\dss the previous paragraph\dss this\sss implies\sss that\dss
$\image \partial^{\dff m}$\dss is\dss also closed\sss
and\dss weak\dnsp$^*$\dnsp-closed.\oss
Suppose\sss that\sss also\dss 
$H^{\dff l_{\dff 1}}_{\dff m\dff +\dff 1}\dff(\trf X\trf)\off =\off 0$\nnsp.\oss
Then\dss
$\image \partial_{\dff m\dff +\dff 2}
\off =\off
\kernel\trf \partial_{\dff m\dff +\dff 1}$\dss
and\dss hence\vspace{4.5pt}\vspace{1.25pt}
\[
\quad
\kernel\trf \partial^{\dff m\dff +\dff 1}
\off =\off
\bigl(\qff \image \partial_{\dff m\dff +\dff 2}\trf\bigr)^{\dff \perp}
\off =\off
\bigl(\qff \kernel\trf \partial_{\dff m\dff +\dff 1}\trf\bigr)^{\dff \perp}
\off =\off
\bigl(\qff {}^{\perp}\dff \left(\qff \image \partial^{\dff m}\trf\right)\qff\bigr)^{\dff \perp}
\pff.
\]

\vspace{-7.5pt}\vspace{1.25pt}
The\sss last\sss space\dss is\dss equal\dss to\sss 
the\sss weak\dnsp$^*$\dnsp-closure of\dss
$\image \partial^{\dff m}$\dnsp.\oss
But\dss $\image \partial^{\dff m}$\dss is\dss weak\dnsp$^*$\dnsp-closed\sss
and\dss hence\qss
$\kernel\trf \partial^{\dff m\dff +\dff 1}
\off =\off
\image \partial^{\dff m}$\dnsp.\oss
It\dss follows\dss that\trs
$\widehat{H}^{\dff m\dff +\dff 1}\dff(\trf X\trf)\off =\off 0$\nnsp.\oss

Suppose\sss now\dss that\qss 
$\widehat{H}^{\dff m\dff +\dff 1}\dff(\trf X\trf)\off =\off 0$\nnsp.\oss
Then\dss 
$\image \partial^{\dff m}
\off =\off
\kernel\trf \partial^{\dff m\dff +\dff 1}$\dss
and\dss hence\dss
$\image \partial^{\dff m}$\dss 
is\dss closed.\oss
As\sss explained\dss above,\pss then\dss 
$\image \partial_{\dff m\dff +\dff 1}$ is\dss closed.\oss
If\trs also\dss 
$\widehat{H}^{\dff m}\dff(\trf X\trf)\off =\off 0$\nnsp,\pss
then
$\image \partial^{\dff m\dff -\dff 1}
\off =\off
\kernel\trf \partial^{\dff m}$\dss
and\dss hence\dss $\image \partial^{\dff m\dff -\dff 1}$\dss is\dss closed.\oss
It\dss follows\dss that\vspace{4.5pt}\vspace{1.25pt}
\[
\quad
\kernel\qff \partial_{\dff m}
\off =\off
{}^{\perp}\dff \left(\qff \image \partial^{\dff m\dff -\dff 1}\trf\right)
\off =\off
{}^{\perp}\dff \left(\qff \kernel\trf \partial^{\dff m}\trf\right)
\off =\off
{}^{\perp}\dff \left(\qff \bigl(\qff \image \partial_{\dff m\dff +\dff 1}\trf\bigr)^{\dff \perp}\trf\right)
\pff.
\]

\vspace{-7.5pt}\vspace{1.25pt}
The\sss last\sss space\dss is\dss equal\dss to\sss 
the\sss closure of\dss
$\image \partial_{\dff m\dff +\dff 1}$\dnsp.\oss
Since\dss $\image \partial_{\dff m\dff +\dff 1}$\dss is\dss closed,\oss
it\dss follows\dss that\qss
$\kernel\qff \partial_{\dff m}
\off =\off
\image \partial_{\dff m\dff +\dff 1}$\dss
and\dss hence\trs
$H^{\dff l_{\dff 1}}_{\dff m}\dff(\trf X\trf)\off =\off 0$\nnsp.\oss
This completes\sss the proof\halfff.\oss  \eproof

\mypar{Corollary.}{homology-amenable}
\emph{If\pss $X$\dss is\dss path-connected and\dss the fundamental\dss group\qss 
$\pi_{\dff 1}\dff(\trf X\fff,\qff x\trf)$\dss
is\dss boundedly\sss acyclic,\oss
then\pss
$H^{\dff l_{\dff 1}}_{\dff m}\dff(\trf X\trf)\off =\off 0$\qss
for every\dss $m\qff \geq\qff 1$\nnsp.\oss}

\proof
Since\qss $\widehat{H}^{\dff m}\dff(\trf X\trf)\off =\off 0$\qss
for every\dss $m\qff \geq\qff 1$\nnsp,\oss
this\sss follows\sss from\qss Theorem\qss \ref{zeros-are-dual}.\oss  \eproof

\mypar{Theorem.}{homology-covering-theorem-open}
\emph{If\qss $\mathcal{U}$\dss be\dss an\sss open\sss 
weakly\dss boundedly\sss acyclic covering of\pss $X$\nnsp,\oss
then\dss the\sss canonical\dss homomorphism\dss
$H_{\dff *}\fff(\trf X \trf)
\qff \ttoo\qff
H^{\dff l_{\dff 1}}_{\dff *}\dff(\trf X\trf)$\dss
can be\dss factored\dss through\qss
$H_{\dff *}\fff(\trf X \trf)
\qff \ttoo\qff 
H_{\dff *}\fff(\trf N \trf)$\dnsp.}

\proof
Clearly\halfff,\pss
$H^{\dff l_{\dff 1}}_{\trf 0}\dff(\trf Z\trf)$\dss 
is\dss canonically\dss isomorphic\sss to\dss $\rrr$\sss for every\dss
path-connected space\dss $Z$\nnsp.\oss
Together\dss with\trs
Corollary\qss \ref{homology-amenable}\qss
this implies\sss that\dss path-connected boundedly\sss acyclic subsets
are\sss $L_{\dff \bullet}$\nsp\dnsp-acyclic.\oss
The rest\sss of\trs the proof\trs differs\sss from\dss 
the proof\dss of\qss
Theorem\qss \ref{covering-theorem-open}\qss
by\dss using\trs Theorem\qss \ref{e-acyclic-open}\qss
instead of\qss Theorem\qss \ref{A-acyclic-open}\qss 
and\dss inverting\dss the directions of\dss arrows.\oss  \eproof

\newpage
\mysection{Uniqueness\qss of\pss Leray\qss homomorphisms}{uniqueness}

\myuppar{Leray\qss homomorphisms.}
For every\dss topological\sss space\dss $X$\dss and an open covering\dss
$\mathcal{U}$\dss of\trs $X$\dss we constructed\dss
in\trs Section\qss \ref{cohomological-leray}\qss a canonical\dss
homomorphism\vspace{1.5pt}
\[
\quad
l_{\trf \mathcal{U}}\dff \colon\dff
H^{\dff *}\fff(\trf N_{\trf \mathcal{U}} \trf)
\qff \ttoo\qff 
H^{\dff *}\fff(\trf X \trf)
\pff,
\]

\vspace{-10.5pt}
where\trs $N_{\trf \mathcal{U}}$\dss is\dss the nerve of\dss $\mathcal{U}$\nnsp,\pss
$H^{\dff *}\fff(\trf N_{\trf \mathcal{U}} \trf)$\dss 
is\dss the simplicial\sss cohomology\sss of\trs $N_{\trf \mathcal{U}}$\nsp,\oss
and\dss $H^{\dff *}\fff(\trf X \trf)$\dss
is\dss the singular cohomology\sss of\trs $X$\nnsp.\oss
The homomorphisms $l_{\trf \mathcal{U}}$ are called\qss
\emph{Leray\qss homomorphisms}.\oss
The goal\sss of\trs this section\dss is\dss to provide an axiomatic characterization
of\trs them.\oss

\myuppar{The category\sss of\dss coverings.}
Following\trs M.\dss Barr\qss \cite{ba},\oss
let\dss us consider\dss the category\dss $\mathfrak{Co}$\dss  
having as objects pairs\dss
$(\trf X\fff,\pff \mathcal{U}\dff)$\nnsp,\oss 
where\dss $\mathcal{U}$\dss is\dss a covering\sss of\dss a\sss topological\sss space\dss $X$\nnsp,\oss
and as morphisms\dss
$(\trf X\fff,\pff \mathcal{U}\dff)
\qff \ttoo\qff
(\trf Y\fff,\pff \mathcal{V}\dff)$\dss
continuous maps\dss
$f\dff \colon\dff
X\qff \ttoo\qff Y$\dss
such\dss that\sss every\sss $U\qff \in\qff \mathcal{U}$\dss
is\dss contained\dss in\dss the preimage\dss $f^{\dff -\dff 1}\dff(\trf V\trf)$\dss
for some\dss  $V\qff \in\qff \mathcal{V}$\dnsp,\oss
i.e.\qss the covering\dss $\mathcal{U}$\dss is\dss a\dss refinement\sss
of\trs the covering\dss\vspace{1.5pt}
\[
\quad
f^{\dff -\dff 1}\dff(\trf \mathcal{V}\trf)
\off =\off
\left\{\pff f^{\dff -\dff 1}\dff(\trf V\trf)
\pff \mid\off
V\qff \in\qff \mathcal{V}
\pff\right\}
\pff.
\]

\vspace{-10.5pt}
As\sss before,\pss 
$N_{\trf \mathcal{U}}$
is\dss the nerve of\trs the covering\dss $\mathcal{U}$\nnsp.\oss
A morphism\dss
$f\dff \colon\dff
(\trf X\fff,\pff \mathcal{U}\dff)
\qff \ttoo\qff
(\trf Y\fff,\pff \mathcal{V}\dff)$\dss
induces a homomorphism\dss
$f^{\dff *}\dff \colon\dff
H^{\dff *}\fff(\trf Y \trf)
\qff \ttoo\qff
H^{\dff *}\fff(\trf X \trf)$\nnsp,\oss
turning\dss 
$(\trf X\fff,\pff \mathcal{U}\dff)\off \longmapsto\off
H^{\dff *}\fff(\trf X \trf)$\dss
into a functor\dss from\dss $\mathfrak{Co}$\dss to graded abelian\dss groups.\oss
Such a\sss morphism $f$ also
induces a homomorphism\vspace{1.5pt}
\[
\quad
f^{\dff *}\dff \colon\dff
H^{\dff *}\fff(\trf N_{\trf \mathcal{V}} \trf)
\qff \ttoo\qff
H^{\dff *}\fff(\trf N_{\trf \mathcal{U}} \trf)
\pff.
\]

\vspace{-10.5pt}
Indeed,\oss since\sss the covering\dss $\mathcal{U}$\dss refines\dss
$f^{\dff -\dff 1}\dff(\trf \mathcal{V}\trf)$\nnsp,\oss
one can choose for each\dss $U\qff \in\qff \mathcal{U}$\dss
some\dss
$\varphi\dff(\trf U\trf)\qff \in\qff \mathcal{V}$\dss
such\dss that\dss
$f\dff(\trf U\trf)\qff \subset\qff \varphi\dff(\trf U\trf)$\nnsp.\oss
Clearly\halfff,\oss if\trs the intersection of\dss several\sss sets\dss
$U_{\dff i}\qff \in\qff \mathcal{U}$\dss is\dss non-empty\halfff,\oss
then\dss the intersection of\trs the images\dss 
$\varphi\dff(\trf U_{\dff i}\trf)$\dss is\dss also non-empty\halfff.\oss
Therefore\dss $\varphi$\dss is\dss a\sss simplicial\dss map\dss
$N_{\trf \mathcal{U}}\qff \ttoo\qff N_{\trf \mathcal{V}}$\dss
and\dss hence defines an\dss induced\dss homomorphism\vspace{1.5pt}
\[
\quad
\varphi^{\dff *}\dff \colon\dff
H^{\dff *}\fff(\trf N_{\trf \mathcal{V}} \trf)
\qff \ttoo\qff
H^{\dff *}\fff(\trf N_{\trf \mathcal{U}} \trf)
\pff.
\]

\vspace{-10.5pt}
Suppose\sss that\dss
$\varphi'\dff \colon\dff
\mathcal{U}\qff \ttoo\qff \mathcal{V}$\dss
is\dss another\dss map such\dss that\dss
$f\dff(\trf U\trf)\qff \subset\qff \varphi'\dff(\trf U\trf)$\dss
for every\dss $U\qff \in\qff \mathcal{U}$\nnsp.\oss
If\trs the intersection of\dss several\sss sets\dss
$U_{\dff i}\qff \in\qff \mathcal{U}$\dss is\dss non-empty\halfff,\oss
then\dss the set\vspace{3pt}\vspace{-0.625pt}
\[
\quad
\bigcap\nolimits_{\dff i}\dff \varphi\trf(\trf U_{\dff i}\trf)
\off \cap\off\qff
\bigcap\nolimits_{\dff i}\dff \varphi'\trf(\trf U_{\dff i}\trf)
\off\off \supset\off\off\qff
f\dff\left(\pff \bigcap\nolimits_{\dff i}\dff U_{\dff i}\trf\right)
\pff
\]

\vspace{-9pt}\vspace{-0.625pt}
is\dss also non-empty\halfff.\oss
It\dss follows\dss that\dss if\trs $\sigma$\dss is\dss a simplex of\trs
$N_{\qff \mathcal{U}}$\nsp,\oss
then\dss 
$\varphi\dff(\trf \sigma\trf)\qff \cup\qff \varphi'\dff(\trf \sigma\trf)$\dss
is\dss a simplex of\trs $N_{\qff \mathcal{V}}$\nsp.\oss
Therefore\dss $\varphi$\dss and\dss $\varphi'$\dss are connected\dss
by\sss an elementary\sss simplicial\dss homotopy\sss
and\dss hence\dss
$\varphi^{\dff *}\off =\off \varphi'^{\dff *}$\nnsp.\oss
We see\sss that\dss the induced\dss homomorphism\dss $\varphi^{\dff *}$\dss
does not\sss depend on\dss the choice of\dss $\varphi$\dss
and\dss hence we can\dss take\sss it\sss as\sss the induced\dss map\dss
$f^{\dff *}\dff \colon\dff
H^{\dff *}\fff(\trf N_{\trf \mathcal{V}} \trf)
\qff \ttoo\qff
H^{\dff *}\fff(\trf N_{\trf \mathcal{U}} \trf)$\nnsp.\oss
A\sss trivial\dss verification shows\sss that\dss these induced\dss maps\sss
turn\dss
$(\trf X\fff,\pff \mathcal{U}\dff)\off \longmapsto\off
H^{\dff *}\fff(\trf N_{\trf \mathcal{U}} \trf)$\dss 
into a functor\dss from\dss $\mathfrak{Co}$\dss to graded abelian\dss groups.\oss

\myuppar{Leray\qss transformations.}
For a\sss space\dss $X$\trs let\dss $\mathcal{U}\dff(\trf X\trf)$\dss
be\sss the covering of\trs $X$\dss by\dss the single set\sss $X$\nnsp.\oss
Clearly\halfff,\pss $N_{\trf \mathcal{U}\dff(\trf X\trf)}$\dss
consists of\dss only\sss one vertex.\oss 
Therefore\dss
$H^{\trf m}\dff(\trf N_{\trf \mathcal{U}\dff(\trf X\trf)} \trf)
\off =\off 
0$\dss
for\dss $m\qff \geq\qff 1$\dss and\dss
$H^{\trf 0}\dff(\qff N_{\trf \mathcal{U}\dff(\trf X\trf)} \trf)$\dss
is\dss equal\dss to\sss the group of\dss coefficients.\oss
If\trs $X$\dss is\dss path-connected,\oss then\dss
$H^{\trf 0}\dff(\trf X \trf)$\dss
is\dss also equal\dss to\sss the group of\dss coefficients
and\dss there\dss is\dss a canonical\dss isomorphism\vspace{2.5pt}
\[
\quad
i_{\qff X}
\qff \colon\qff
H^{\trf 0}\dff(\trf N_{\trf \mathcal{U}\dff(\trf X\trf)} \trf)
\off \ttoo\off
H^{\trf 0}\dff(\trf X \trf)
\pff.
\]

\vspace{-9.5pt}
Let $\mathfrak{C}$ be a\sss full\sss subcategory\dss of\sss $\mathfrak{Co}$
containing\dss all\dss coverings 
$(\trf X\fff,\pff \mathcal{U}\dff(\trf X\trf)\trf)$
with\dss path-connected $X$\nnsp.\oss
A natural\dss transformation\dss from\dss the functor\dss
$(\trf X\fff,\pff \mathcal{U}\dff)
\off \longmapsto\off
H^{\dff *}\fff(\trf N_{\trf \mathcal{U}} \trf)$\dss
on\dss the category\dss $\mathfrak{C}$\dss
to\sss the functor\dss
$(\trf X\fff,\pff \mathcal{U}\dff)
\off \longmapsto\off
H^{\dff *}\fff(\trf X \trf)$\dss
assigns\sss to each\dss object\dss
$(\trf X\fff,\pff \mathcal{U}\dff)$\dss
of\trs $\mathfrak{C}$\dss a\sss homomorphism\vspace{2.5pt}
\[
\quad
l_{\trf \mathcal{U}}\dff \colon\dff
H^{\dff *}\fff(\trf N_{\trf \mathcal{U}} \trf)
\qff \ttoo\qff 
H^{\dff *}\fff(\trf X \trf)
\pff
\]

\vspace{-9.5pt}
in such a\sss way\dss that\dss for every\dss morphism\dss
$f\dff \colon\dff
(\trf X\fff,\pff \mathcal{U}\dff)
\qff \ttoo\qff
(\trf Y\fff,\pff \mathcal{V}\dff)$\dss
the square\vspace{0pt}\vspace{-0.75pt}
\begin{equation}
\label{leray-transformation}
\quad
\begin{tikzcd}[column sep=boom, row sep=boom]\dis
H^{\dff *}\dff(\trf N_{\trf \mathcal{V}} \dff)
\arrow[r, "\dis l_{\trf \mathcal{V}}"]
\arrow[d, "\dis \qff f^{\dff *}"]
&
H^{\dff *}\dff(\trf Y\trf)
\arrow[d, "\dis \qff f^{\dff *}"]
\\
H^{\dff *}\dff(\trf N_{\trf \mathcal{U}} \dff) 
\arrow[r, "\dis l_{\trf \mathcal{U}}"]
&
H^{\dff *}\dff(\trf X \trf)
\end{tikzcd}
\end{equation}

\vspace{-9pt}\vspace{-0.75pt}
is\dss commutative.\oss
Let\dss us\dss call\sss such a\sss natural\dss transformation\sss a\pss
\emph{Leray\dss transformation}\pss on\dss $\mathfrak{C}$\qss
if\qss for every\sss path-connected space\dss $X$\dss the homomorphism\sss
$l_{\trf \mathcal{U}\dff(\trf X\trf)}$\sss
is\dss equal\sss in dimension $0$\sss to\sss $i_{\qff X}$\nsp.\oss

\mypar{Theorem.}{double-complex-functoriality}
\emph{The canonical\dss homomorphisms\qss
$l_{\trf \mathcal{U}}\dff \colon\dff
H^{\dff *}\fff(\trf N_{\trf \mathcal{U}} \trf)
\qff \ttoo\qff 
H^{\dff *}\fff(\trf X \trf)$\qss
from\dss Section\qss \ref{cohomological-leray}\qss
form a\qss Leray\qss transformation on\dss the category\sss of\dss
open coverings.\oss}

\proof
The only\dss not\dss quite functorial\sss element\sss 
constructions of\trs the construction of\dss 
$l_{\trf \mathcal{U}}$
is\dss the choice of\dss linear orders on\sss sets of\dss vertices.\oss
As\dss is\dss well\dss known,\oss
the cohomology\dss
$H^{\dff *}\fff(\trf N_{\trf \mathcal{U}} \trf)$\sss
is\dss independent\sss on\sss this choice up\sss to canonical\dss isomorphisms,\oss
as also\sss the induced\dss maps.\oss
Given\sss $f$\sss as above,\oss one can choose\sss 
the orders of\dss vertices of\trs $N_{\trf \mathcal{U}}$
and\sss $N_{\trf \mathcal{V}}$ in such a way\dss that\sss $f$\sss is\qss
(non-strictly)\qss order-preserving.\oss
A routine check\sss shows\sss that\dss then\qss (\ref{leray-transformation})\qss
is\dss commutative.\oss 
Another\sss routine check shows\sss that\sss
$l_{\trf \mathcal{U}\dff(\trf X\trf)}$\sss
is\dss equal\sss in dimension $0$\sss to\sss $i_{\qff X}$\nsp.\oss  \eproof

\myuppar{Geometric realizations of\dss simplicial\dss complexes.}
Let\trs $N$\dss is\dss a simplicial\sss complex,\pss 
$\num{N}$\dss be its geometric realization,\oss
and\dss 
$\mathcal{U}\dff(\trf N\trf)$\dss be\sss the covering of\trs $\num{N}$\dss 
by\dss the open stars of\dss vertices.\oss
By\dss a\sss well\dss known\dss theorem of\qss Eilenberg\qss \cite{e},\pss 
$N_{\qff \mathcal{U}\dff(\trf N\trf)}\off =\off N$\dss
and\dss there\dss is\dss a\sss canonical\dss isomorphism\vspace{1.5pt}
\[
\quad
i_{\qff N}
\qff \colon\qff
H^{\dff *}\dff(\trf N \trf)
\off \ttoo\off
H^{\dff *}\dff(\trf \num{N} \trf)
\pff.
\]

\vspace{-10.5pt}
For any\dss reasonable\trs Leray\dss transformation
one should\dss have\dss
$l_{\trf \mathcal{U}\dff(\trf N\trf)}\off =\off i_{\qff N}$\nsp.\oss

\mypar{Theorem.}{Leray-realizations}
\emph{$l_{\trf \mathcal{U}\dff(\trf N\trf)}\off =\off i_{\qff N}$\pss
for every\sss simplicial\sss complex\dss $N$\nnsp.\oss}

\proof
The idea\dss is\dss to prove\sss that\sss
$N\off \longmapsto\off l_{\trf \mathcal{U}\dff(\trf N\trf)}$\sss
is\dss a natural\dss transformation of\dss cohomology\dss theories
and apply\dss the\qss Eilenberg--Steenrod\dss uniqueness\sss theorem.\oss
First\sss of\dss all,\oss we need\sss 
to extend\dss the construction of\dss
$l_{\trf \mathcal{U}\dff(\trf N\trf)}$\sss
to pairs of\dss simplicial\sss complexes.\oss
Let\vspace{3pt}
\[
\quad
\mathcal{T}^{\dff \bullet}\dff(\trf N\trf)
\off =\off
T^{\dff \bullet}\dff(\trf N_{\qff \mathcal{U}\dff(\trf N\trf)}\fff,\pff C\trf)
\]

\vspace{-9pt}
be\sss the\sss total\sss complex associated\dss with\dss
the covering\dss $\mathcal{U}_{\qff N}$\dss of\trs $\num{N}$\dss by\sss open stars.\oss
For a simplicial\sss subcomplex\sss $L$\sss of\trs $N$\dss let\dss
$\mathcal{T}^{\dff \bullet}\dff(\trf N,\pff L\trf)$\sss
be\sss the kernel\sss of\dss morphism\sss
$\mathcal{T}^{\dff \bullet}\dff(\trf N\trf)
\qff \ttoo\qff
\mathcal{T}^{\dff \bullet}\dff(\trf L\trf)$\sss
induced\dss by\dss the inclusion\sss
$L\qff \ttoo\qff N$\nnsp.\oss
As we know,\oss the morphisms\vspace{3pt}
\[
\quad
\tau_{\dff C}\dff \colon\dff
C^{\dff \bullet}\dff(\trf \num{N}\trf)
\off \ttoo\off
\mathcal{T}^{\dff \bullet}\dff(\trf N\trf)
\quad\
\mbox{and}\quad\
C^{\dff \bullet}\dff(\trf \num{L}\trf)
\off \ttoo\off
\mathcal{T}^{\dff \bullet}\dff(\trf L\trf)
\]

\vspace{-9pt}
induce isomorphisms\sss in cohomology\dss groups.\oss
It\dss follows\dss that\dss the natural\dss morphism\vspace{3pt}
\[
\quad
C^{\dff \bullet}\dff(\trf \num{N}\fff,\pff \num{L}\trf)
\off \ttoo\off
\mathcal{T}^{\dff \bullet}\dff(\trf N\fff,\pff L\trf)
\]

\vspace{-9pt}
also induces isomorphism\sss in cohomology\dss groups.\oss
This allows\sss to define\sss the map\vspace{3pt}
\[
\quad
l_{\qff \mathcal{U}\dff(\trf N\trf)\fff,\trf \mathcal{U}\dff(\trf L\trf)}
\qff \colon\qff
H^{\dff *}\fff(\trf N\fff,\pff L \trf)
\qff \ttoo\qff 
H^{\dff *}\fff(\trf \num{N}\fff,\pff \num{L} \trf)
\pff
\]

\vspace{-9pt}
exactly\sss as in\dss the absolute case.\oss

Next,\oss we need\sss to check\dss the functoriality\sss of\trs the maps\sss
$l_{\trf \mathcal{U}\dff(\trf N\trf)}$\sss
and\sss
$l_{\qff \mathcal{U}\dff(\trf N\trf)\fff,\trf \mathcal{U}\dff(\trf L\trf)}$\nnsp.\oss
We will\dss limit\sss ourselves by\dss the absolute case,\oss
the case of\dss pairs being completely\sss similar\halfff.\oss
For a vertex $v$ of\dss a simplicial\sss complex we will\sss denote by\sss $U_{\dff v}$\dss
the open star of\dss $v$\nnsp.\oss
Let\dss
$f\dff \colon\dff
M\qff \ttoo\qff N$\dss
be a simplicial\dss map and\dss 
$\num{f}\qff \colon\qff
\num{M}\qff \ttoo\qff \num{N}$\dss
be its geometric realization.\oss
Then\dss
$\num{f}\qff(\trf U_{\dff v}\trf)
\off \subset\off 
U_{\dff f\dff(\dff v\trf)}$
for every\sss vertex $v$ of\dss $M$\nnsp.\oss
It\dss follows\dss that\dss $\num{f}$\dss is\dss 
a morphism of\dss coverings\dss\vspace{3pt}
\[
\quad
(\trf \num{N}\fff,\pff \mathcal{U}\dff(\trf N\trf)\trf)
\off \ttoo\off
(\trf \num{M}\fff,\pff \mathcal{U}\dff(\trf M\trf)\trf)
\pff
\]

\vspace{-9pt}
and\dss hence induces a morphism\dss
$\num{f}^{\dff *}\dff \colon\dff
\mathcal{T}^{\dff \bullet}\dff(\trf N\trf)
\qff \ttoo\qff
\mathcal{T}^{\dff \bullet}\dff(\trf M\trf)$\nnsp.\oss
The diagram\vspace{0pt}
\[
\quad
\begin{tikzcd}[column sep=boom, row sep=boom]\dis
C^{\dff \bullet}\dff(\trf N \trf) 
\arrow[r, "\dis \lambda_{\qff C}"]
\arrow[d, "\dis \qff f^{\dff *}"]
& 
\mathcal{T}^{\dff \bullet}\dff(\trf N\trf)
\arrow[d, "\dis \qff \num{f}^{\dff *}"]
&
C^{\dff \bullet}\dff(\trf \num{N} \trf)\pff
\arrow[l, "\dis \off \tau_{\trf C}"']
\arrow[d, "\dis \qff \num{f}^{\dff *}"]
\\
C^{\dff \bullet}\dff(\trf M \trf) 
\arrow[r, "\dis \lambda_{\qff C}"]
& 
\mathcal{T}^{\dff \bullet}\dff(\trf M\trf)
&
C^{\dff \bullet}\dff(\trf \num{M} \trf)\pff
\arrow[l, "\dis \off \tau_{\trf C}"']
\end{tikzcd}
\]

\vspace{-7.5pt}
is\dss commutative,\oss
as a routine check shows.\oss
This\sss implies\sss the commutativity\sss of\trs
the corresponding diagram of\trs the (co)homology\dss groups.\oss
In\dss turn,\oss this implies\sss that\dss the diagram\vspace{1.4pt}
\[
\quad
\begin{tikzcd}[column sep=boom, row sep=boom]\dis
H^{\dff *}\dff(\trf N \trf) 
\arrow[r, "\dis l_{\trf \mathcal{U}\dff(\trf N\trf)}"]
\arrow[d, "\dis \qff f^{\dff *}"]
&
H^{\dff *}\dff(\trf \num{N} \trf)
\arrow[d, "\dis \qff \num{f}^{\dff *}"]
\\
H^{\dff *}\dff(\trf M \trf) 
\arrow[r, "\dis l_{\trf \mathcal{U}\dff(\trf M\trf)}"]
& 
H^{\dff *}\dff(\trf \num{M} \trf)
\end{tikzcd}
\]

\vspace{-7.5pt}
is\dss commutative.\oss
This\dss is\dss the required\dss functoriality\sss of\dss 
$l_{\trf \mathcal{U}\dff(\trf N\trf)}$\nnsp.\oss
We\sss leave\sss to\sss the reader\sss to check\dss the functoriality\sss of\dss 
$l_{\qff \mathcal{U}\dff(\trf N\trf)\fff,\trf \mathcal{U}\dff(\trf L\trf)}$\nnsp,\oss
as also\sss to check\dss that\dss the maps\sss
$l_{\trf \mathcal{U}\dff(\trf N\trf)}$\sss
and\sss
$l_{\qff \mathcal{U}\dff(\trf N\trf)\fff,\trf \mathcal{U}\dff(\trf L\trf)}$\sss
commute with\dss the connecting\sss
homomorphisms\sss $\partial$\sss of\trs the cohomological\dss 
long exact\sss sequence of\trs the pair\sss
$(\trf N\fff,\pff L\trf)$\nnsp.\oss
Clearly,\oss in dimension $0$\sss these maps are equal\dss to\sss the
isomorphism\sss $i_{\qff N}$\sss and\dss its analogue for\sss the pair\sss
$(\trf N\fff,\pff L\trf)$\sss respectively.\oss
Now\dss the cohomological\sss version of\qss Eilenberg--Steenrod\dss uniqueness\sss theorem\dss
implies\sss that\dss
$l_{\trf \mathcal{U}\dff(\trf N\trf)}\off =\off\dff i_{\qff N}$\nsp.\oss
See\qss \cite{es},\oss Theorem\qss VI.8.1.\oss  \eproof

\mypar{Theorem.}{coverings-of-paracompact}
\emph{Let\qss $\mathfrak{P}$ be\sss the\sss full\sss subcategory\dss of\qss $\mathfrak{Co}$
having\sss open coverings of\qss paracompact\sss spaces as\dss its\sss objects.\oss
There\sss exists\sss exactly\sss one\trs Leray\trs transformation on\dss $\mathfrak{P}$\nnsp.\oss}

\proof
We already\dss constructed a\dss Leray\trs transformation on a bigger category.\oss
It\dss remains\sss to prove\sss the uniqueness.\oss
Suppose\sss that\sss $\mathcal{U}$\sss is\dss an open covering\sss o
f\dss a paracompact\sss space\sss $X$\nnsp,\oss
and\dss let\dss $N\off =\off N_{\qff \mathcal{U}}$\nsp.\oss
Since\sss $X$\sss is\dss paracompact,\oss 
there exists a partition of\dss unity\dss $t_{\trf U}$\nsp,\pss
$U\qff \in\pff \mathcal{U}$\dss subordinated\dss to\dss $\mathcal{U}$\nnsp.\oss
It\trs leads\sss to a continuous map\dss
$t\dff \colon\dff
X\qff \ttoo\qff \num{N}$\nnsp.\oss
The open star of\dss a vertex\dss $u$\dss corresponding\dss to\dss
$U\qff \in\pff \mathcal{U}$\dss
consists of\dss points with\dss the barycentric coordinate\dss 
$t_{\trf u}\qff >\qff 0$\nnsp.\oss
It\dss follows\dss that\dss the image\dss $t\trf(\trf U\trf)$\dss is\dss contained\dss
in\dss the open star of\dss $u$\nnsp.\oss
Therefore\sss $t$\sss is\dss a\sss morphism\vspace{3pt}
\[
\quad
(\trf X\fff,\pff \mathcal{U}\dff)
\qff \ttoo\qff
(\trf \num{N}\fff,\pff \mathcal{U}\dff(\trf N\trf)\trf)
\pff.
\]

\vspace{-9pt}
Clearly,\pss
$N_{\qff \mathcal{U}\dff(\trf N\trf)}
\off =\off 
N_{\qff \mathcal{U}}$\sss
and one can\sss take\sss the identity\dss map as\sss the map\sss
$\varphi\dff \colon\dff
N_{\qff \mathcal{U}}
\qff \ttoo\qff
N_{\qff \mathcal{U}\dff(\trf N\trf)}$\sss
used\dss to construct\sss
$t^{\dff *}\dff \colon\dff
H^{\dff *}\fff(\trf N_{\trf \mathcal{U}\dff(\trf N\trf)} \trf)
\qff \ttoo\qff
H^{\dff *}\fff(\trf N_{\trf \mathcal{U}} \trf)$\nnsp.\oss
Hence\sss the\sss latter\sss map\dss is\dss equal\dss to identity.\oss
If\dss 
$\mathcal{U}\off \longmapsto\off l_{\trf \mathcal{U}}$\sss
is\dss a\dss Leray\dss transformation\sss on\sss $\mathfrak{P}$\nnsp,\oss
the\sss the diagram\vspace{1.4pt}
\[
\quad
\begin{tikzcd}[column sep=boom, row sep=boom]\dis
H^{\dff *}\dff(\trf N \trf)
\arrow[r, "\dis l_{\trf \mathcal{U}\dff(\trf N\trf)}"]
\arrow[d, "\dis \qff ="]
&
H^{\dff *}\dff(\trf \num{N}\trf)
\arrow[d, "\dis \qff t^{\dff *}"]
\\
H^{\dff *}\dff(\trf N \trf) 
\arrow[r, "\dis l_{\trf \mathcal{U}}"]
&
H^{\dff *}\dff(\trf X \trf)
\end{tikzcd}
\]

\vspace{-7.5pt}
is\dss commutative.\oss
Together\sss with\qss Theorem\qss \ref{Leray-realizations}\qss
this implies\sss that\sss
$l_{\trf \mathcal{U}}
\off =\off
t^{\dff *}\dff \circ\qff i_{\qff N}$\nsp.\oss
The uniqueness\sss follows,\oss
as also\sss the fact\dss that\sss $t^{\dff *}$\sss
does not\sss depend on\sss the partition of\dss unity.\oss  \eproof

\newpage
\mysection{Nerves\qss of\pss families\qss and\qss paracompact\qss spaces}{paracompact-spaces}

\myuppar{Families of\dss subsets.}
Let\dss $S$\dss be a\sss set\halfff.\oss
A\qss \emph{family\sss of\dss subsets}\qss 
$\mathcal{F}
\off =\off
\{\qff
F_{\dff i}
\qff\}_{\dff i\dff \in\qff I}$\qss
of\trs $S$\dss is\dss a map\dss\vspace{1.5pt}
\[
\quad 
i\off \longmapsto\off F_{\dff i}\pff \subset\off S
\]

\vspace{-10.5pt}
from a set\sss $I$\sss to\sss the set\sss of\dss all\sss subsets of\trs $S$\nnsp.\oss
Usually\dss the nature of\trs the set\dss $I$\dss is\dss of\dss no importance 
and\dss $\mathcal{F}$\dss is\dss treated almost\sss as a collection of\dss
subsets of\trs $S$\nnsp.\oss
But\dss our\sss main\dss interest\dss is\dss in\dss the\sss families indexed\dss by\dss the same
set\dss $I$\nnsp.\oss
If\qss 
$\mathcal{G}
\off =\off
\{\qff
G_{\dff i}
\qff\}_{\dff i\dff \in\qff I}$\dss
is\dss another\dss family\dss indexed\dss by\trs 
$I$\dss
and\dss $F_{\dff i}\qff \subset\qff G_{\dff i}$\dss for every\dss
$i\qff \in\pff I$\nnsp,\oss
we say\dss that\dss $\mathcal{F}$\sss is\qss \emph{combinatorially\dss refining}\pss
$\mathcal{G}$\nnsp.\oss

The\qss \emph{nerve}\pss of\dss a\sss family\trs 
$\mathcal{F}
\off =\off
\{\qff
F_{\dff i}
\qff\}_{\dff i\dff \in\qff I}$\dss 
is\dss an abstract\sss
simplicial\sss complex\dss having\dss the set\dss $I$\dss
as\sss the set\sss of\dss vertices.\oss
Its simplices are\dss finite non-empty\sss subsets\dss
$\sigma\qff \subset\pff I$\dss such\dss that\vspace{3pt}
\[
\quad
\bigcap\nolimits_{\qff i\qff \in\qff \sigma}\qff F_{\dff i}
\off \neq\off
\varnothing
\pff.
\]

\vspace{-9pt}
When\dss two families\qss
$\mathcal{F}
\off =\off
\{\qff
F_{\dff i}
\qff\}_{\dff i\dff \in\qff I}$\qss
and\qss
$\mathcal{G}
\off =\off
\{\qff
G_{\dff i}
\qff\}_{\dff i\dff \in\qff I}$\qss
are indexed\dss by\dss the same set\dss $I$\nnsp,\oss
it\dss make sense\sss to say\dss that\dss their\dss nerves are\qss
\emph{equal}.\oss

\myuppar{Families of\dss subspaces.}
Suppose now\dss that\dss $S$\dss is\dss a\sss topological\sss space.\oss
A\dss family\dss 
$\mathcal{F}
\off =\off
\{\qff
F_{\dff i}
\qff\}_{\dff i\dff \in\qff I}$\dss 
of\dss subsets,\oss or\halfff,\oss
what\dss is\dss the same,\oss
of\dss subspaces of\trs $S$\dss is\dss said\dss to be\qss \emph{open}\qss
if\trs the sets\dss $F_{\dff i}$\dss are open,\oss
and\qss \emph{closed}\pss if\trs the sets\dss $F_{\dff i}$\dss are closed.\oss
The\sss family\dss 
$\mathcal{F}
\off =\off
\{\qff
F_{\dff i}
\qff\}_{\dff i\dff \in\qff I}$\dss 
is\dss said\dss to be\qss 
\emph{locally\dss finite}\pss if\dss for every\dss $x\qff \in\qff S$\dss
there exists an open set\dss $U$\dss such\dss that\dss $x\qff \in\qff U$\dss
and\dss $U\qff \cap\qff F_{\dff i}\off \neq\off \varnothing$\dss
for only\sss a\sss finite number of\dss indices\dss $i\qff \in\pff I$\nnsp.\oss

Given a closed\dss family\dss 
$\mathcal{F}
\off =\off
\{\qff
F_{\dff i}
\qff\}_{\dff i\dff \in\qff I}$\nsp,\oss
we would\dss like\sss to know\dss when\dss there exists an open\dss family\dss
$\mathcal{U}$\dss such\dss that\dss $\mathcal{F}$\dss is\dss 
combinatorially\dss refining\dss $\mathcal{G}$\dss and\dss
the nerves of\trs $\mathcal{F}$\dss and\dss $\mathcal{G}$\dss 
are equal.\oss
Since\sss the nerves of\dss families are involved,\oss
it\dss is\dss only\dss natural\dss to assume\sss that\dss $\mathcal{F}$\dss is\trs
locally\dss finite.\oss
If\trs $I$\dss is\dss countable,\oss
then\dss such a\sss family\dss $\mathcal{U}$\dss exists\dss if\trs
$\mathcal{F}$\dss is\trs locally\dss finite and\dss $S$\dss is\dss a\sss normal\sss space.\oss
See\qss Proposition\qss \ref{countable-extensions}.\oss
But\dss the analogue of\trs this result\dss for\dss general\trs $I$\dss involves
an additional\sss assumption about\dss $\mathcal{F}$\dss
satisfied\dss when\dss $S$\dss is\dss not\sss only\dss normal,\oss
but\dss is,\oss moreover\halfff,\oss \emph{paracompact}\qss (see below).\oss
Now\dss we\sss turn\dss to\sss the\sss key\dss property\sss 
of\trs locally\dss finite families.\oss

\mypar{Lemma.}{closure}
\emph{Let\dss $\mathcal{F}$\dss be a\sss locally\dss finite\sss family\dss of\dss
subspaces of\trs a\sss topological\dss space\dss $S$\nnsp.\oss
Then}\vspace{3pt}\vspace{0.375pt}
\[
\quad
\overline{\dff\bigcup\nolimits_{\qff i\dff \in\qff I}\qff F_{\dff i}}
\off\qff =\off\pff
\bigcup\nolimits_{\qff i\dff \in\qff I}\qff \overline{F}_{\dff i}
\pff,
\]

\vspace{-9pt}\vspace{0.375pt}
\emph{where\qss $\overline{A}$\qss denotes\sss the closure\sss of\qss
a\sss subset\dss $A\qff \subset\qff S$\nnsp,\oss
and\dss the family\qss
$\bigl\{\qff
\overline{F}_{\dff i}
\qff\bigr\}_{\dff i\dff \in\qff I}$\dss is\dss locally\qss finite.}\oss

\proof
If\dss $U$\dss is\dss an open set\dss intersecting\dss $F_{\dff i}$\dss 
for only\sss a\sss finite number of\trs $i\qff \in\pff I$\nnsp,\oss
then\dss $U$\dss intersects\sss the closures\dss $\overline{F}_{\dff i}$\dss
for only\sss a\sss finite number of\trs $i\qff \in\pff I$\nnsp.\oss
The\dss lemma\sss follows.\oss  \eproof

\mypar{Proposition.}{countable-extensions}
\emph{Let\dss $Z$\dss be a\sss normal\sss space.\oss
Let\pss 
$\mathcal{F}
\off =\off
\{\qff
F_{\dff i}
\qff\}_{\dff i\dff \in\qff I}$\dss
be a closed\dss family\dss 
combinatorially\dss refining\sss an\sss open\dss family\qss
$\mathcal{U}
\off =\off
\{\qff
U_{\dff i}
\qff\}_{\dff i\dff \in\qff I}$\qss
of\qss subsets of\pss $Z$\nnsp.\oss
If\qss the family\qss 
$\mathcal{F}$\dss is\dss locally\dss finite and\qss $I$\dss
is\dss countable,\oss
then\dss there\dss exists\dss an\dss open\dss family\qss
$\mathcal{G}
\off =\off
\{\qff
G_{\dff i}
\qff\}_{\dff i\dff \in\qff I}$\dss
such\dss that\dss}\vspace{1.5pt}
\[
\quad
F_{\dff i}
\off \subset\off 
G_{\dff i}
\off \subset\off
\overline{G}_{\dff i}
\off \subset\off 
U_{\dff i}
\]

\vspace{-10.5pt}
\emph{for\dss every\dss $i\qff \in\pff I$\dss
and\dss the\sss nerves\sss of\pss $\mathcal{F}$\dss
and\qss $\mathcal{G}$\dss are equal.\oss}

\proof
We may\sss assume\sss that\qss
$I\off =\off \{\qff 0\fff,\pff 1\fff,\pff 2\fff,\pff \ldots\off \qff\}$\dss
is\dss the set\sss of\dss non-negative integers.\oss
Let\dss us\dss first\dss ignore\sss the condition\dss
$\overline{G}_{\dff i}\off \subset\off U_{\dff i}$\nsp.\oss
Equivalently\halfff,\oss let\dss us\sss assume\sss that\dss
$U_{\dff i}\off =\off Z$\qss for every\dss $i$\nnsp.\oss
Let\dss us\sss consider\dss the family\dss of\dss 
intersections\dss\vspace{1.5pt}
\[
\quad
\bigcap_{\qff i\qff \in\qff K}\qff F_{\dff i}
\quad\
\mbox{such\dss that}\quad\
F_{\trf 0}
\off \cap\pff
\bigcap_{\qff i\qff \in\qff K}\qff F_{\dff i}
\off =\off
\varnothing
\]

\vspace{-10.5pt}
with\sss finite\dss $K\qff \subset\pff I$\nnsp.\oss
Since\dss $\mathcal{F}$\dss is\dss locally\dss finite,\oss
this\sss family\dss is\dss also\sss locally\dss finite.\oss
By\qss Lemma\qss \ref{closure}\qss the union\dss $C_{\trf 0}$\dss
of\trs this family\dss is\dss closed.\oss
Clearly\halfff,\pss $F_{\trf 0}\qff \cap\pff C_{\trf 0}\off =\off \varnothing$\nnsp.\oss
Since\dss $Z$\dss is\dss normal,\oss
there exists an open set\dss $G_{\trf 0}$\dss such\dss that\qss 
$F_{\trf 0}
\off \subset\off 
G_{\trf 0}$\qss
and\qss
$\overline{G}_{\trf 0}
\off \cap\off 
C_{\trf 0}
\off =\off 
\varnothing$\nnsp.\oss
Let\dss us\sss consider\dss the\sss family\vspace{1.5pt}
\[
\quad
\mathcal{E}_{\trf 1}
\off =\off
\bigl\{\off \overline{G}_{\trf 0}\qff,\off
F_{\trf 1}\dff,\off
F_{\trf 2}\dff,\off
\ldots\off
\qff\bigr\}
\]

\vspace{-10.5pt}
Clearly\halfff,\oss the\sss nerve of\dss $\mathcal{E}_{\trf 1}$\dss
is\dss equal\dss to\sss the nerve of\trs
$\mathcal{F}$\nnsp.\oss
If\dss an open set\dss $U$\dss intersects an\dss infinite number of\dss sets\sss from\dss
$\mathcal{E}_{\trf 1}$\nsp,\oss
then\dss $U$\dss intersects an\dss infinite number of\trs the sets\dss
$F_{\dff i}$\dss with\dss $i\qff \geq\qff 1$\nnsp,\oss
contrary\dss to\sss the assumption\dss that\dss $\mathcal{F}$\dss is\dss
locally\dss finite.\oss
Hence\dss $\mathcal{E}_{\trf 1}$\dss is\dss locally\dss finite.\oss
Arguing\dss by\dss induction,\oss
we may\sss assume\sss that\dss the open sets\dss
$G_{\trf 0}\dff,\off
G_{\trf 1}\dff,\off
\ldots\dff,\off
G_{\dff k\dff -\dff 1}$\dss
are already\sss defined,\oss
that\trs
$F_{\trf i}
\off \subset\off 
G_{\trf i}$\qss
for\qss $i\off =\off 0\fff,\pff 1\fff,\pff 2\fff,\pff \ldots\fff,\pff k\qff -\qff 1$\nnsp,\oss 
the\sss family\qss\vspace{1.5pt}
\[
\quad
\mathcal{E}_{\dff k}
\off =\off
\bigl\{\off
\overline{G}_{\trf 0}\qff,\off
\ldots\qff,\off
\overline{G}_{\dff k\dff -\dff 1}\qff,\dff\off
F_{\dff k}\dff,\off
F_{\dff k\dff +\dff 1}\dff,\off
\ldots\off
\qff\bigr\}
\]

\vspace{-10.5pt}
is\dss locally\dss finite,\oss
and\trs the nerves of\trs $\mathcal{E}_{\dff k}$\dss
and\dss $\mathcal{F}$\dss are equal.\oss
By\dss applying\dss the above arguments\sss to\dss
$\mathcal{E}_{\dff k}$\dss and\dss $F_{\dff k}$\dss
in\dss the roles\sss of\trs 
$\mathcal{F}$\sss
and\dss 
$F_{\trf 0}$\dss respectively\halfff,\oss
we will\dss get\sss the next\dss family\vspace{1.5pt}
\[
\quad
\mathcal{E}_{\dff k\dff +\dff 1}
\off =\off
\bigl\{\off
\overline{G}_{\trf 0}\qff,\off
\ldots\qff,\off
\overline{G}_{\dff k\dff -\dff 1}\qff,\dff\off
\overline{G}_{\dff k}\qff,\dff\off
F_{\dff k\dff +\dff 1}\dff,\off
\ldots\off
\qff\bigr\}
\]

\vspace{-10.5pt}
with similar\dss properties.\oss
By\dss continuing\dss this\dss process indefinitely\dss
we define\dss an open subset\dss $G_{\dff i}\qff \supset\qff F_{\dff i}$\dss
for every\dss $i\off =\off 0\fff,\pff 1\fff,\pff 2\fff,\pff \ldots \off$\nnsp.\oss
Let\dss us\dss consider\dss the family\vspace{1.5pt}
\[
\quad
\mathcal{G}
\off =\off
\bigl\{\qff
G_{\trf 0}\dff,\off
G_{\trf 1}\dff,\off
G_{\trf 2}\dff,\off
\ldots\off
\qff\bigr\}
\pff.
\]

\vspace{-10.5pt}
Suppose\sss that\qss
$\bigcap\nolimits_{\qff i\qff \in\qff K}\qff G_{\dff i}
\off \neq\off
\varnothing$\qss 
for a finite subset\qss $K\qff \subset\pff I$\nnsp.\oss
Then\qss
$\bigcap\nolimits_{\qff i\qff \in\qff K}\qff \overline{G}_{\dff i}
\off \neq\off
\varnothing$\qss 
and\dss if\sss $k$\sss is\dss the maximal\sss element\sss of\trs $K$\nnsp,\oss
then\dss the sets\dss $\overline{G}_{\dff i}$\dss
with\dss $i\qff \in\pff K$\dss occur\sss in\sss the family\dss
$\mathcal{E}_{\dff k\dff +\dff 1}$\nsp.\oss
Since\sss the nerves of\trs the families\sss
$\mathcal{E}_{\dff k\dff +\dff 1}$
and\dss
$\mathcal{F}$\dss are equal,\oss
it\dss follows\dss that\dss
$\bigcap\nolimits_{\qff i\qff \in\qff K}\qff F_{\dff i}
\off \neq\off
\varnothing$\nnsp.\oss
We see\sss that\qss
$\bigcap\nolimits_{\qff i\qff \in\qff K}\qff G_{\dff i}
\off \neq\off
\varnothing$\qss
implies\qss
$\bigcap\nolimits_{\qff i\qff \in\qff K}\qff F_{\dff i}
\off \neq\off
\varnothing$\nnsp.\oss
The converse implication\dss holds\sss because\dss
$F_{\dff i}\qff \subset\qff G_{\dff i}$\dss
for every\dss $i$\nnsp.\oss
Therefore\sss the nerves of\trs
$\mathcal{F}$\dss and\dss $\mathcal{G}$\dss are equal.\oss

This proves\sss the\sss theorem\dss when\dss
$U_{\dff i}\off =\off Z$\qss for\dss every\dss $i$\nnsp.\oss
The extension\dss to\sss the general\sss case\dss is\dss easy\halfff.\oss
By\dss the already\dss proved\dss part\sss of\trs the\sss theorem\dss
there exists a\sss family\trs
$\mathcal{G}'
\off =\off
\{\qff
G_{\dff i}'
\qff\}_{\dff i\dff \in\qff I}$\dss of\dss open sets such\dss that\trs
$F_{\dff i}\qff \subset\qff G_{\dff i}'$\dss
for every\dss $i$\dss and\dss the nerves of\trs
$\mathcal{F}$\dss and\dss $\mathcal{G}'$\dss are equal.\oss 
For every\dss $i$\qss let\qss 
$V_{\fff i}\off =\off G_{\dff i}'\qff \cap\qff U_{\dff i}$\nsp.\oss
Clearly\halfff,\pss
$F_{\dff i}\qff \subset\qff V_{\fff i}$\nsp.\oss 
Since\dss $Z$\dss is\dss normal,\oss
for every\dss $i$\dss there exists an open set\dss $G_{\dff i}$\dss
such\dss that\trs
$F_{\dff i}\off \subset\off G_{\dff i}$\dss
and\dss
$\overline{G}_{\dff i}
\off \subset\off 
V_{\fff i}
\off \subset\off 
U_{\fff i}$\nsp.\oss
Also,\pss
$G_{\dff i}
\off \subset\off 
G_{\dff i}'$\dss 
for every\dss $i$\nnsp.\oss
Since\sss the nerves of\trs
$\mathcal{F}$\dss and\dss $\mathcal{G}'$\dss are equal,\oss
this implies\sss that\dss the nerves of\trs
$\mathcal{F}$\dss and\dss $\mathcal{G}$\dss are also equal.\oss  \eproof

\mypar{Theorem.}{general-extensions}
\emph{Let\dss $Z$\dss be a\sss normal\sss space.\oss
Let\pss 
$\mathcal{F}
\off =\off
\{\qff
F_{\dff i}
\qff\}_{\dff i\dff \in\qff I}$\dss
be a closed\dss family\dss  
combinatorially\dss refining\sss an\sss open\dss family\qss
$\mathcal{U}
\off =\off
\{\qff
U_{\dff i}
\qff\}_{\dff i\dff \in\qff I}$\qss
of\qss subsets of\pss $Z$\nnsp.\oss
If\qss both\trs families\qss $\mathcal{F}$\dss and\qss $\mathcal{U}$\dss
are\dss locally\qss finite,\oss
then\dss there\dss exists\dss an\dss open\dss locally\dss finite\dss family\qss
$\mathcal{G}
\off =\off
\{\qff
G_{\dff i}
\qff\}_{\dff i\dff \in\qff I}$\qss
such\dss that}\vspace{2.5pt}
\[
\quad
F_{\dff i}
\off \subset\off 
G_{\dff i}
\off \subset\off
\overline{G}_{\dff i}
\off \subset\off
U_{\dff i}
\]

\vspace{-9.5pt}
\emph{for\dss every\dss $i\qff \in\pff I$\dss
and\dss the\sss nerves\sss of\pss $\mathcal{F}$\dss
and\qss $\mathcal{G}$\dss are equal.\oss}

\proof
When\dss the set\sss $I$\sss is\dss uncountable,\oss one needs\sss to replace\sss
the usual\dss induction\dss by\trs the\sss transfinite\sss induction
or\dss some equivalent\dss tool.\oss
We\sss will\dss use\qss Zorn\qss lemma.\oss
More importantly\halfff,\oss in\dss this case\sss
there\dss is\dss no easy\dss way\dss to\ ensure\dss that\dss the families\dss
$\mathcal{E}_{\dff k}$\dss 
from\dss the proof\dss of\qss
Proposition\qss \ref{countable-extensions}\qss are\sss
locally\dss finite.\oss
By\dss this\sss reason one needs\sss to assume\dss that\dss the family\dss $\mathcal{U}$\dss
is\dss locally\dss finite.\oss

Let\dss us\dss consider\dss the set\dss $\mathfrak{A}$\dss of\dss all\dss families\dss
$\mathcal{A}
\off =\off
\{\qff
A_{\dff i}
\qff\}_{\dff i\dff \in\qff I}$\dss
of\dss subspaces of\trs $Z$\dss
such\dss that\dss for every\dss $i\qff \in\pff I$\dss either\dss
$A_{\dff i}\off =\off F_{\dff i}$\nsp,\oss
or\dss $A_{\dff i}$\dss is\dss open and\dss\vspace{2pt}
\[
\quad
F_{\dff i}
\off \subset\off 
A_{\dff i} 
\off \subset\off
\overline{A}_{\dff i}
\off \subset\off
U_{\dff i} 
\pff.
\]

\vspace{-10pt}
In order\dss to apply\qss Zorn\qss lemma,\oss
let\dss us\dss define a partial\sss order\dss $\leq$\dss on\dss $\mathfrak{A}$\dss
as follows.\oss
Suppose\sss that\dss
$\mathcal{A}
\off =\off
\{\qff
A_{\dff i}
\qff\}_{\dff i\dff \in\qff I}$\dss
and\dss
$\mathcal{B}
\off =\off
\{\qff
B_{\dff i}
\qff\}_{\dff i\dff \in\qff I}$\dss 
are\dss two elements of\trs $\mathfrak{A}$\nnsp.\oss
Then\dss $\mathcal{A}\qff \leq\pff \mathcal{B}$\qss
if\trs the nerves of\dss $\mathcal{A}$\dss and\dss $\mathcal{B}$\dss
are equal\sss and\dss
for every\dss $i\qff \in\pff I$\dss either\dss
$A_{\dff i}\off =\off B_{\dff i}$\nsp,\oss
or\dss $A_{\dff i}\off =\off F_{\dff i}$\nsp,\oss or\dss both.\oss
Clearly\halfff,\pss $\leq$\dss is\dss a\sss partial\sss order on\dss $\mathfrak{A}$\nnsp.\oss
We claim\dss that\dss $\leq$\dss satisfies\sss the assumptions of\pss Zorn\qss lemma.\oss

Suppose\sss that\dss 
$\mathfrak{B}\off \subset\off \mathfrak{A}$\dss is\dss a\sss linearly\dss ordered\dss by\dss $\leq$\dss
subset\halfff.\oss
Let\dss $i\qff \in\pff I$\nnsp.\oss
If\qss $A_{\dff i}\off \neq\off F_{\dff i}$\qss for some family\dss
$\mathcal{A}
\off =\off
\{\qff
A_{\dff i}
\qff\}_{\dff i\dff \in\qff I}$\nnsp,\oss
then\dss
$B_{\dff i}\off =\off A_{\dff i}$\dss
for every\sss family\dss
$\mathcal{B}
\off =\off
\{\qff
B_{\dff i}
\qff\}_{\dff i\dff \in\qff I}$\dss
such\dss that\dss
$\mathcal{A}\qff \leq\pff \mathcal{B}$\nnsp.\oss 
Let\qss $E_{\dff i}\off =\off F_{\dff i}$\qss
if\qss $A_{\dff i}\off =\off F_{\dff i}$\qss
for every\dss $\mathcal{A}\qff \in\qff \mathfrak{B}$\nnsp,\oss
and\qss $E_{\dff i}\off =\off A_{\dff i}$\qss
if\qss $A_{\dff i}\off \neq\off F_{\dff i}$\qss
for some\qss $\mathcal{A}\qff \in\qff \mathfrak{B}$\nnsp.\oss
Since\dss $\mathfrak{B}$\dss is\dss linearly\sss ordered,\oss
this definition\dss is\dss correct\halfff.\oss
Let\qss 
$\mathcal{E}
\off =\off
\{\qff
E_{\dff i}
\qff\}_{\dff i\dff \in\qff I}$\nsp.\oss
If\trs $\mathcal{E}\qff \in\qff \mathfrak{A}$\nnsp,\oss
then,\oss obviously\halfff,\pss
$\mathcal{A}\qff \leq\qff \mathcal{B}$\dss for every\dss
$\mathcal{A}\qff \in\qff \mathfrak{B}$\nnsp,\oss
and\dss hence\dss $\mathfrak{B}$\dss admits an upper\sss bound.\oss

Let\dss us\dss prove\sss that\trs
$\mathcal{E}\qff \in\qff \mathfrak{A}$\nnsp.\oss
Clearly\halfff,\oss
for every\dss $i\qff \in\pff I$\dss either\qss
$E_{\dff i}\off =\off F_{\dff i}$\nsp,\oss
or\qss $E_{\dff i}$\dss is\dss open and\dss\vspace{2pt}
\[
\quad
F_{\dff i}
\off \subset\off 
E_{\dff i} 
\off \subset\off
\overline{E}_{\dff i}
\off \subset\off
U_{\dff i} 
\pff.
\]

\vspace{-10pt}
Suppose\sss that\qss
$\bigcap\nolimits_{\qff i\qff \in\qff K}\qff E_{\dff i}
\off \neq\off
\varnothing$\qss 
for\sss a\sss finite subset\qss $K\qff \subset\pff I$\nnsp.\oss
Each set\trs $E_{\dff i}$\dss with\dss $i\qff \in\pff K$\dss occurs\dss in\sss some\sss family\dss 
$\mathcal{B}_{\fff i}\qff \in\qff \mathfrak{B}$\nnsp.\oss
Since\dss $K$\dss is\dss finite and\dss $\mathfrak{B}$\dss is\dss linearly\sss ordered,\oss
the set\sss of\dss families\dss $\mathcal{B}_{\fff i}$\nnsp,\oss 
where\dss $i\qff \in\pff K$\nnsp,\oss
has a maximal\sss element\halfff.\oss
If\trs
$\mathcal{B}
\off =\off
\{\qff
B_{\dff i}
\qff\}_{\dff i\dff \in\qff I}$\dss
is\dss this maximal\dss element\halfff,\oss
then\dss $E_{\dff i}\off =\off B_{\dff i}$\dss for every\dss $i\qff \in\pff K$\nnsp.\oss
Since\sss the nerves\sss of\dss $\mathcal{B}$\dss and\dss $\mathcal{F}$\dss are equal,\oss
this implies\sss that\trs
$\bigcap\nolimits_{\qff i\qff \in\qff K}\qff F_{\dff i}
\off \neq\off
\varnothing$\nnsp.\oss
It\dss follows\dss that\dss the nerves\sss of\dss
$\mathcal{E}$\dss and\dss $\mathcal{F}$\dss are equal.\oss
Therefore\dss $\mathcal{E}\qff \in\qff \mathfrak{A}$\nnsp.\oss

We see\sss that\sss every\dss linearly\sss ordered subset\sss of\trs $\mathfrak{A}$\dss
admits an upper bound.\oss
Therefore\qss Zorn\qss lem\-ma applies
and\dss there exists a\sss family\dss
$\mathcal{G}
\off =\off
\{\qff
G_{\dff i}
\qff\}_{\dff i\dff \in\qff I}$\dss
in\dss $\mathfrak{A}$\dss maximal\dss 
with\sss respect\dss to\dss $\leq$\nnsp.\oss
As above,\oss
for every\dss $i\qff \in\pff I$\dss either\dss
$G_{\dff i}\off =\off F_{\dff i}$\nsp,\oss
or\dss $G_{\dff i}$\dss is\dss open and\dss\vspace{1.5pt}
\[
\quad
F_{\dff i}
\off \subset\off 
G_{\dff i} 
\off \subset\off
\overline{G}_{\dff i}
\off \subset\off
U_{\dff i} 
\pff.
\]

\vspace{-10.5pt}
The family\dss $\mathcal{G}$\dss is\dss locally\dss finite because\dss $\mathcal{U}$\dss is.\oss
The arguments\sss of\trs the previous paragraph show\dss that\dss
the nerves of\dss $\mathcal{B}$\dss and\dss $\mathcal{F}$\dss are equal.\oss 
It\dss remains\sss to show\dss that\dss
$G_{\dff i}$\dss is\dss open\dss
for every\dss $i\qff \in\pff I$\nnsp.\oss
Suppose\sss that\dss $G_{\dff k}$\dss is\dss not\sss open\dss
for some\dss $k\qff \in\pff I$\nnsp.\oss
Then\dss
$G_{\dff k}\off =\off F_{\dff k}$\nsp.\oss
In\dss particular\halfff,\pss
$G_{\dff k}$\dss is\dss closed.\oss
Since\sss the family\dss $\mathcal{U}$\dss is\dss locally\dss finite
and\qss
$\overline{G}_{\dff i}
\off \subset\off
U_{\dff i}$\qss
for every\dss $i\qff \in\pff I$\nnsp,\oss 
the family\sss of\dss closures\dss\vspace{1.5pt}
\[
\quad
\bigl\{\pff 
\overline{G}_{\dff i}
\pff\bigr\}_{\dff i\dff \in\qff I}
\]

 \vspace{-10.5pt}
is\dss locally\dss finite.\oss
It\dss follows\dss that\dss the family\dss of\dss 
intersections\dss\vspace{3pt}
\[
\quad
\bigcap_{\qff i\qff \in\qff K}\qff \overline{G}_{\dff i}
\quad\
\mbox{such\dss that}\quad\
G_{\dff k}
\off \cap\pff
\bigcap_{\qff i\qff \in\qff K}\qff \overline{G}_{\dff i}
\off =\off
\varnothing
\]

\vspace{-9pt}
with\sss finite\dss $K\qff \subset\pff I$\dss
is\dss also\sss locally\dss finite.\oss
By\qss Lemma\qss \ref{closure}\qss the union\dss $C_{\trf k}$\dss
of\trs this family\dss is\dss closed.\oss
Clearly\halfff,\qss 
$G_{\dff k}\pff \cap\pff C_{\trf k}\off =\off \varnothing$\nnsp.\oss
Since\sss the space\dss $Z$\dss is\dss normal,\oss
there exists an open set\dss $H_{\dff k}$\dss such\dss that\vspace{1.5pt}
\[
\quad
 F_{\dff k}
\off \subset\off 
H_{\dff k}
\qff,\qff\quad
\overline{H}_{\dff k}
\off \subset\off 
U_{\dff k}
\qff,\quad\
\mbox{and}\qff\quad\
\overline{H}_{\dff k}
\off \cap\off 
C_{\dff k}
\off =\off 
\varnothing
\pff.
\]

\vspace{-10.5pt}
Let\dss us\sss define a\sss family\dss 
$\mathcal{E}
\off =\off
\{\qff
E_{\dff i}
\qff\}_{\dff i\dff \in\qff I}$\qss
as\dss follows\fff:\vspace{1.5pt}
\[
\quad
E_{\dff i}\off =\off G_{\dff i}
\quad
\mbox{if}\quad
i\off \neq\off k
\quad\
\mbox{and}\quad\
E_{\dff k}\off =\off H_{\dff k}
\pff.
\] 

\vspace{-10.5pt}
Then\dss the nerve of\trs $\mathcal{E}$\dss is\dss equal\dss to\sss the nerve of\trs
$\mathcal{G}$\dss and\dss hence\sss
to\sss the nerve of\trs
$\mathcal{F}$\dss and\dss
for every\dss $i\qff \in\pff I$\dss either\dss
$E_{\dff i}\off =\off F_{\dff i}$\nsp,\oss
or\dss $E_{\dff i}$\dss is\dss open and\dss\vspace{1.5pt}
\[
\quad
F_{\dff i}
\off \subset\off 
E_{\dff i} 
\off \subset\off
\overline{E}_{\dff i}
\off \subset\off
U_{\dff i} 
\pff.
\]

\vspace{-10.5pt}
Hence\dss $\mathcal{E}\qff \in\qff \mathfrak{A}$\nnsp.\oss
But\sss since\dss $G_{\dff k}\off =\off F_{\dff k}$\dss
is\dss assumed\dss to\sss be not\sss open and\dss 
$E_{\dff k}\off =\off H_{\dff k}$\dss is\dss open,\pss
$E_{\dff k}\off \neq\off F_{\dff k}$\dss
and\dss
$\mathcal{G}\off <\off \mathcal{E}$\nnsp,\oss
contrary\dss to\sss the maximality\sss of\trs $\mathcal{G}$\nnsp.\oss
The contradiction shows\sss that\trs
$G_{\dff i}$\dss is\dss open\dss
for every\dss $i\qff \in\pff I$\nnsp.\oss  \eproof

\myuppar{Remarks.}
Theorem\qss \ref{general-extensions}\qss was\dss proved\dss by\qss
K.\dss Morita\qss \cite{mor}.\oss 
See\qss \cite{mor},\oss 
Theorem\qss 1.3.\oss
The above proof\trs is\dss a version of\qss Morita's\trs one,\oss but\halfff,\oss
following\pss J.\dss van\sss Mill\qss \cite{vm},\oss
we\sss included\sss a\sss proof\dss of\qss its\sss simplified\dss version,\oss
Proposition\qss  \ref{countable-extensions}.\oss
Cf.\qss \cite{vm},\oss Proposition\qss 3.2.1\qss and\dss Corollary\qss 3.2.2.\oss
The main\dss idea\dss is\dss more\sss transparent\dss in\dss the countable case,\oss
and\dss the separation of\trs the countable and\dss the general\sss cases
shows why\dss one needs\sss to assume\sss that\dss $\mathcal{U}$\dss is\dss
locally\dss finite.\oss

\myuppar{Paracompactness.}
Suppose\sss that\qss
$\mathcal{F}
\off =\off
\{\qff
F_{\dff i}
\qff\}_{\dff i\dff \in\qff I}$\qss
and\qss 
$\mathcal{U}
\off =\off
\{\qff
U_{\dff k}
\qff\}_{\dff k\dff \in\qff K}$\qss
are\sss two families of\dss subsets of\dss a\sss set\dss $S$\nnsp.\oss
If\trs for every\dss $i\qff \in\pff I$\dss there exits\dss
$k\qff \in\pff K$\dss such\dss that\trs
$F_{\dff i}\qff \subset\qff U_{\dff k}$\nsp,\oss
we say\dss that\dss $\mathcal{F}$\sss is\qss \emph{refining}\qss
$\mathcal{U}$\nnsp.\oss
Obviously\halfff,\oss if\qss
$\mathcal{F}$\dss is\dss combinatorially\dss refining\dss $\mathcal{U}$\nnsp,\oss
then\dss $\mathcal{F}$\dss is\dss refining\dss $\mathcal{U}$\nnsp.\oss
The family\dss $\mathcal{F}$\dss  
is\dss said\dss to be a\qss \emph{covering}\qss
of\dss $S$\dss if\trs the union of\dss all\sss sets\dss $F_{\dff i}$\dss is\dss
equal\dss to\dss $S$\nnsp.\oss

A\sss topological\dss space $X$\dss is\dss called\qss \emph{paracompact}\pss
if\trs for every open covering\dss
$\mathcal{F}$\trs of\trs $X$\dss there exists\sss a\dss locally\dss finite\sss open\dss
covering\dss refining\dss
$\mathcal{F}$\nnsp.\oss
For\dss the rest\sss of\trs this section\dss we will\sss assume\sss
that\qss \emph{$X$\dss is\dss a\dss paracompact\qss Hausdorff\qss space.\oss}

\mypar{Lemma.}{normal}
\emph{The space\dss $X$\dss is\dss normal.}\oss

\proof
Let\dss $A\qff \subset\qff X$\dss be a closed set\sss and\dss let\dss
$b\qff \in\qff X\qff \smallsetminus\qff A$\nnsp.\oss
Since\dss $X$\dss is\dss a\dss Hausdorff\dss space,\oss for every\dss $a\qff \in\qff A$\dss
there are\sss open sets\dss $U_{\dff a}\qff \ni\qff a$\dss and\dss $V_{\dff a}\qff \ni\qff b$\dss 
such\dss that\dss
$U_{\dff a}\qff \cap\qff V_{\dff a}
\off =\off
\varnothing$\nnsp.\oss
The family\qss\vspace{0pt}
\[
\quad
\bigl\{\qff
U_{\dff a}
\qff\bigr\}_{\dff a\dff \in\qff A} 
\]

\vspace{-12pt}
together\dss with\dss $X\qff \smallsetminus\qff B$\dss
forms\sss an open covering of\trs $X$\nnsp.\oss
Let\dss
$\{\qff
F_{\dff i}
\qff\}_{\dff i\dff \in\qff I}$\dss be a\sss locally\dss finite open\dss
covering\dss refining\dss
this cover\halfff,\oss
and\dss let\dss
$I_{\dff A}\off =\off
\{\qff 
i\qff \in\qff I
\off |\off 
F_{\dff i}\qff \cap\qff A\off \neq\off \varnothing
\qff\}$\nnsp.\oss
Since\dss
$U_{\dff a}\qff \cap\qff V_{\dff a}
\off =\off
\varnothing$\nnsp,\oss
the closure\dss
$\overline{\dff F_{\dff i}}$\dss
does not\dss contain\dss $b$\dss if\trs $i\qff \in\pff I_{\dff A}$\nsp.\oss
Together\dss with\trs Lemma\qss \ref{closure}\qss this implies\sss that\vspace{3pt}
\[
\quad
\bigcup\nolimits_{\dff i\dff \in\qff I_{\dff A}}\dff F_{\dff i}
\quad\
\mbox{and}\quad\
X\off \smallsetminus\off
\overline{\dff\bigcup\nolimits_{\dff i\dff \in\qff I_{\dff A}}\dff F_{\dff i}}
\]

\vspace{-9pt}
are disjoint\sss open neighborhoods of\dss $A$\dss and\dss $b$\dss respectively\halfff.\oss 
Therefore\sss the space\dss $X$\dss is\dss regular\halfff.\oss
Now\sss a similar\sss argument\sss allows\sss to prove\sss that\dss $X$\dss is\dss normal.\oss
Let\dss $A\fff,\pff B\qff \subset\qff X$\dss be closed sets\sss such\dss that\dss
$A\qff \cap\qff B\off =\off \varnothing$\nnsp.\oss
By\dss the previous paragraph\dss $X$\dss is\dss regular\halfff.\oss
Hence\sss for every\dss $a\qff \in\qff A$\dss
there are\sss open sets\dss $U_{\dff a}$\dss and\dss $V_{\dff a}$\dss 
such\dss that\dss
$a\qff \in\qff U_{\dff a}$,\pss $B\qff \subset\qff V_{\dff a}$\nsp,\oss
and\dss
$U_{\dff a}\qff \cap\qff V_{\dff a}
\off =\off
\varnothing$\nnsp.\oss
The families\vspace{3pt}
\[
\quad
\bigl\{\qff
U_{\dff a}
\qff\bigr\}_{\dff a\dff \in\qff A}
\quad\
\mbox{and}\quad\
\bigl\{\qff
V_{\dff a}
\qff\bigr\}_{\dff a\dff \in\qff A}
\]

\vspace{-9pt}
together\dss with\dss $X\qff \smallsetminus\qff B$\dss
forms\sss an open covering of\trs $X$\nnsp.\oss
Let\dss
$\{\qff
F_{\dff i}
\qff\}_{\dff i\dff \in\qff I}$\dss be a\sss locally\dss finite open\dss covering\dss refining\dss
this covering\halfff,\oss
and\dss let\dss
$I_{\dff A}\off =\off
\{\qff 
i\qff \in\qff I
\off |\off 
F_{\dff i}\qff \cap\qff A\off \neq\off \varnothing
\qff\}$\nnsp.\oss
Since\dss
$U_{\dff a}\qff \cap\qff V_{\dff a}
\off =\off
\varnothing$\nnsp,\oss
the closure\dss
$\overline{\dff F_{\dff i}}$\dss
is\dss disjoint\dss from\dss $B$\dss if\trs $i\qff \in\pff I_{\dff A}$\nsp.\oss
Together\dss with\trs Lemma\qss \ref{closure}\qss this implies\sss that\vspace{3pt}
\[
\quad
\bigcup\nolimits_{\dff i\dff \in\qff I_{\dff A}}\dff F_{\dff i}
\quad\
\mbox{and}\quad\
X\off \smallsetminus\off
\overline{\dff\bigcup\nolimits_{\dff i\dff \in\qff I_{\dff A}}\dff F_{\dff i}}
\]

\vspace{-9pt}
are disjoint\sss open neighborhoods of\dss $A$\dss and\dss $B$\dss respectively\halfff.\oss  \eproof

\mypar{Theorem.}{closed-refinement}
\emph{If\pss 
$\mathcal{U}
\off =\off
\{\qff
U_{\dff i}
\qff\}_{\dff i\dff \in\qff I}$\dss
is\dss an open covering of\pss $X$\nnsp,\oss
then\dss there\dss exists\dss a\dss locally\dss finite closed\sss covering\qss
$\mathcal{F}
\off =\off
\{\qff
F_{\dff i}
\qff\}_{\dff i\dff \in\qff I}$\dss
of\pss $X$\dss such\dss that\pss
$\mathcal{F}$\dss is\dss combinatorially\dss refining\dss
$\mathcal{U}$\dnsp.\oss}

\proof
Lemma\qss \ref{normal}\qss implies\sss that\dss the space\dss $X$\dss is\dss regular\halfff.\oss
Hence\dss if\trs $x\qff \in\qff U_{\dff i}$\nsp,\oss
then\dss $x\qff \in\qff G$\dss for an open set\dss $G$\dss
such\dss that\dss $\overline{\dff G\dff}\qff \subset\qff U_{\dff i}$\nsp.\oss
Such open sets\dss $G$\dss form\sss an open covering\oss
$\mathcal{G}
\off =\off
\{\qff
G_{\dff k}
\qff\}_{\dff k\qff \in\qff K}$\dss 
of\trs $X$\dss such\dss that\dss the\sss family\sss of\trs the closures\vspace{1.5pt}
\[
\quad
\bigl\{\pff
\overline{G}_{\dff k}
\qff\bigr\}_{\dff k\qff \in\qff K}
\]

\vspace{-10.5pt}
is\dss refining\trs
$\{\qff
U_{\dff i}
\qff\}_{\dff i\dff \in\qff I}$\nsp.\oss
Let\dss
$\{\qff
A_{\dff s}
\qff\}_{\dff s\qff \in\qff S}$\dss
be a\sss locally\dss finite covering\dss refining\dss $\mathcal{G}$\dss
(at\dss this point\dss we don't\dss even\dss need\dss it\trs to\sss be open).\oss
Then\dss for every\dss $s\qff \in\qff S$\dss there exists\dss
$i\dff(\dff s\trf)\qff \in\pff I$\dss such\dss that\vspace{3pt}
\[
\quad
\overline{\dff A_{\dff s}}
\off\qff \subset\off\qff
U_{\dff i\dff(\dff s\trf)}
\pff.
\]

\vspace{-10pt}
For every\dss $i\qff \in\pff I$\pss
let\qss\vspace{2pt}
\[
\quad
F_{\dff i}
\off =\off
\bigcup\nolimits_{\qff i\dff(\dff s\trf)\qff =\qff i}\qff
\overline{\dff A_{\dff s}}
\pff.
\]

\vspace{-9pt}
By\qss Lemma\qss \ref{closure}\qss the family\dss
$\{\qff
F_{\dff i}
\qff\}_{\dff i\dff \in\qff I}$\dss
is\dss locally\dss finite.\oss
Since\dss
$\{\qff
A_{\dff s}
\qff\}_{\dff s\qff \in\qff S}$\dss
is\dss a covering,\pss
$\{\qff
F_{\dff i}
\qff\}_{\dff i\dff \in\qff I}$\dss
is\dss a covering also.\oss
Clearly\halfff,\pss
$F_{\dff i}$\dss is\dss closed\sss and\dss
$F_{\dff i}\qff \subset\qff U_{\dff i}$\pss
for every\dss $i\qff \in\pff I$\nnsp.\oss  \eproof

\myuppar{Remark.}
In\dss the situation of\pss Theorem\qss \ref{closed-refinement}\qss
it\dss may\dss happen\dss that\dss $F_{\dff i}\off =\off \varnothing$\dss
for some\dss $i\qff \in\qff I$\nnsp.\oss
Note\sss that\dss adding\sss empty\sss subsets\sss to a family\sss
does not\sss affect\dss the property\sss of\trs being\dss locally\dss finite.\oss

\myuppar{Stars.}
Let\dss
$\mathcal{U}
\off =\off
\{\qff
U_{\dff i}
\qff\}_{\dff i\dff \in\qff I}$\dss
be a covering of\dss a set\trs $Z$\nnsp.\oss
For\dss $A\qff \subset\qff Z$\qss let\vspace{3pt}
\[
\quad
\sta(\trf A\fff,\pff \mathcal{U}\trf)
\off =\off
\bigcup\qff U_{\dff i}
\pff,
\]

\vspace{-9pt}
where\sss the union\dss is\dss taken over\dss $i\qff \in\qff I$\dss
such\dss that\dss
$A\qff \cap\qff U_{\dff i}
\off \neq\off
\varnothing$\nnsp.\oss
For\dss $z\qff \in\qff Z$\dss we set\vspace{3pt}
\[
\quad
\sta(\trf z\fff,\pff \mathcal{U}\trf)
\off =\off
\sta(\trf \{\trf z\trf\}\fff,\pff \mathcal{U}\trf)
\qff.
\]

\vspace{-9pt}\vspace{-0.875pt}
Clearly\halfff,\oss for every\dss $A\qff \in\qff Z$\vspace{3pt}\vspace{-0.875pt}
\[
\quad
\sta(\trf A\fff,\pff \mathcal{U}\trf)
\off =\off
\bigcup\nolimits_{\dff a\qff \in\qff A}\qff \sta(\trf a\fff,\pff \mathcal{U}\trf)
\pff.
\]

\vspace{-9pt}
A covering\dss
$\mathcal{A}
\off =\off
\{\qff
A_{\dff s}
\qff\}_{\dff s\dff \in\qff S}$\dss
is\dss said\dss to be a\qss \emph{barycentric\dss refinement}\pss of\trs $\mathcal{U}$\dss
if\dss for every\dss $z\qff \in\qff Z$\dss
there exists\dss $i\qff \in\qff I$\dss such\dss that\dss
$\sta(\trf z\fff,\pff \mathcal{A}\trf)
\qff \subset\qff
U_{\dff i}$\nsp,\oss
and\dss is\dss to be a\qss \emph{star\dss refinement}\pss of\trs $\mathcal{U}$\dss
if\dss for every\dss $s\qff \in\qff S$\dss
there exists\dss $i\qff \in\qff I$\dss such\dss that\dss
$\sta(\trf A_{\dff s}\fff,\pff \mathcal{A}\trf)
\qff \subset\qff
U_{\dff i}$\nsp.\oss

\mypar{Lemma.}{barycentric-refinement}
\emph{Let\qss $\mathcal{U}$\dss
be\sss an open covering of\pss $X$\nnsp.\oss
Then\dss there\sss exists\sss an\dss open covering of\pss $X$\dss
barycentrically\dss refining\dss $\mathcal{U}$\nnsp.\oss}

\proof
Let\trs
$\mathcal{U}
\off =\off
\{\qff
U_{\dff i}
\qff\}_{\dff i\qff \in\qff I}$\nsp.\oss
Then\dss by\qss Theorem\qss \ref{closed-refinement}\qss
there\dss exists\dss a\dss locally\dss finite closed\sss covering\qss
$\mathcal{F}
\off =\off
\{\qff
F_{\dff i}
\qff\}_{\dff i\dff \in\qff I}$\dss
of\pss $X$\dss such\dss that\qss
$F_{\dff i}\qff \subset\qff U_{\dff i}$\qss
for every\dss $i\qff \in\pff I$\nnsp.\oss
For every\dss $x\qff \in\qff X$\trs let\vspace{3pt}
\[
\quad
I\trf(\dff x\trf)
\off =\off
\bigl\{\qff
i\qff \in\pff I
\off |\off
x\qff \in\qff F_{\dff i}
\qff\bigr\}
\pff.
\]

\vspace{-9pt}
Since\dss $\mathcal{F}$\dss is\dss locally\dss finite,\pss
$I\trf(\dff x\trf)$\dss
is\dss finite for every\dss $x\qff \in\qff X$\nnsp.\oss 
By\qss Lemma\qss \ref{closure}\qss the union\vspace{0pt}
\[
\quad
\bigcup\nolimits_{\qff k\qff \in\qff K}
\qff F_{\dff i}
\]

\vspace{-12pt}
is\dss closed\dss for every\dss $K\qff \subset\pff I$\nnsp.\oss
It\dss follows\dss that\dss the set\vspace{1.5pt}
\[
\quad
V_{\dff x}
\off =\off
\left(\qff
\bigcap_{\qff i\qff \in\qff I\trf(\dff x\trf)}
\qff U_{\dff i}
\qff\right)
\off \cap\off
\left(\qff
X\qff \smallsetminus\qff
\bigcup_{\qff i\qff \not\in\qff I\trf(\dff x\trf)}
\qff F_{\dff i}
\qff\right)
\]

\vspace{-10.5pt}
is\dss open for every\dss $x\qff \in\qff X$\nnsp.\oss
Therefore\dss
$\mathcal{V}
\off =\off
\{\qff
V_{\dff z}
\qff\}_{\dff z\qff \in\qff X}$\dss
is\dss an open covering of\trs $X$\nnsp.\oss

It\dss is\dss sufficient\dss to prove\sss that\dss
$\mathcal{V}$\dss is\dss a\sss barycentric\sss refinement\sss of\dss $\mathcal{U}$\nnsp.\oss
Suppose\sss that\dss $z\qff \in\qff X$\dss 
and\dss let\dss us\sss choose some\dss $j\qff \in\pff I\trf(\dff z\trf)$\nnsp,\oss
so\sss that\dss $z\qff \in\qff F_{\fff j}$\nsp.\oss 
By\dss the definition of\trs the sets\dss $V_{\dff x}$\nsp,\oss
if\dss $z\qff \in\qff V_{\dff x}$\nnsp,\oss
then\dss $z\qff \not\in\qff F_{\dff i}$\dss
for every\dss $i\qff \not\in\qff I\trf(\dff x\trf)$\nnsp.\oss
Therefore in\dss this case\dss $j\qff \in\pff I\trf(\dff x\trf)$\dss
and\dss hence\dss $V_{\dff x}\qff \subset\qff U_{\fff j}$\nnsp.\oss
It\dss follows\dss that\dss
$\sta(\trf z\fff,\pff \mathcal{V}\trf)
\qff \subset\qff
U_{\fff j}$\nsp.\oss
Since\dss $z\qff \in\qff X$\dss was arbitrary\halfff,\pss
$\mathcal{V}$\dss is\dss a\sss barycentric refinement\sss of\dss 
$\mathcal{U}$\dnsp.\oss  \eproof

\mypar{Lemma.}{double-refinement}
\emph{Let\dss $\mathcal{A}\dff,\off \mathcal{B}\dff,\off \mathcal{C}$\dss
be coverings of\trs a\sss set\halfff.\oss
If\trs $\mathcal{A}$\sss is\dss a\sss barycentric\sss refinement\sss of\trs $\mathcal{B}$\dss and\dss
$\mathcal{B}$\dss is\dss a\sss barycentric\sss refinement\sss of\qss $\mathcal{C}$\nnsp,\oss
then\dss $\mathcal{A}$\dss is\dss a\sss star\sss refinement\sss of\qss $\mathcal{C}$\nnsp.\oss}

\proof
Let\trs
$\mathcal{A}
\off =\off
\{\qff
A_{\dff s}
\qff\}_{\dff s\qff \in\qff S}$\qss
and\qss
$\mathcal{B}
\off =\off
\{\qff
B_{\dff t}
\qff\}_{\dff t\qff \in\qff T}$\nsp.\oss
Let\dss us\dss fix\dss $s\qff \in\qff S$\dss
and\dss for each\dss $a\qff \in\qff A_{\dff s}$\dss choose\dss
$t\trf(\dff a\trf)\qff \in\qff T$\dss
such\dss that\trs
$\sta(\dff a\fff,\pff \mathcal{A}\trf)
\qff \subset\qff
B_{\dff t\trf(\dff a\trf)}$\nsp.\oss
Then\vspace{3pt}
\[
\quad
\sta(\trf A_{\dff s}\fff,\pff \mathcal{U}\trf)
\off =\off\dff
\bigcup\nolimits_{\dff a\qff \in\qff A_{\dff s}}\qff \sta(\dff a\fff,\pff \mathcal{U}\trf)
\off \subset\off\dff
\bigcup\nolimits_{\dff a\qff \in\qff A_{\dff s}}\qff B_{\dff t\trf(\dff a\trf)}
\pff.
\]

\vspace{-9pt}
Let\dss us\dss fix\dss now\sss some\dss $z\qff \in\qff A_{\dff s}$\nsp.\oss
If\dss also $a\qff \in\qff A_{\dff s}$\nsp,\oss
then\dss\vspace{3pt}
\[
\quad
z\pff \in\pff 
\sta(\dff a\fff,\pff \mathcal{A}\trf)
\pff \subset\off
B_{\dff t\trf(\dff a\trf)}
\]

\vspace{-9pt}
and\dss hence\dss
$z\qff \in\qff B_{\dff t\trf(\dff a\trf)}$\nsp.\oss
It\dss follows\dss that\vspace{3pt}
\[
\quad
\bigcup\nolimits_{\qff a\qff \in\qff A_{\dff s}}\qff B_{\dff t\trf(\dff a\trf)}
\off \subset\off
\sta(\dff z\fff,\pff \mathcal{B}\trf)
\]

\vspace{-9pt}
and\dss hence\dss
$\sta(\trf A_{\dff s}\fff,\pff \mathcal{U}\trf)
\pff \subset\off
\sta(\dff z\fff,\pff \mathcal{B}\trf)$\nnsp.\oss
Since\dss $s\qff \in\qff S$\dss was arbitrary\sss
and\dss $\mathcal{B}$\dss is\dss a\sss barycentric refinement\sss of\trs
$\mathcal{C}$\nnsp,\oss
this\sss implies\dss that\dss $\mathcal{A}$\dss
is\dss a star\dss refinement\sss of\trs $\mathcal{C}$\nnsp.\oss  \eproof

\mypar{Theorem.}{star-refinement}
\emph{Let\qss $\mathcal{U}$\dss
be\sss an open covering of\pss $X$\nnsp.\oss
Then\dss there\sss exists\sss an\sss open\dss covering of\pss $X$\dss
star\dss refining\trs $\mathcal{U}$\nnsp.\oss}

\proof
It\dss is\dss sufficient\dss to apply\dss twice\trs
Lemma\qss \ref{barycentric-refinement}\qss
and\dss then\sss apply\qss
Lemma\qss \ref{double-refinement}.\oss  \eproof

\mypar{Theorem.}{closed-to-open-extension}
\emph{Let\pss 
$\mathcal{F}$\dss
be\sss a\sss locally\dss finite closed\dss family\dss of\qss subsets of\pss $X$\nnsp.\oss
Then\dss there\sss exists\dss a\dss locally\trs finite open\dss family\qss
$\mathcal{U}$\dss
of\qss subsets\sss of\pss $X$\dss
such\dss that\pss
$\mathcal{F}$\dss is\dss combinatorially\dss refining\dss
$\mathcal{U}$\dnsp.\oss}

\proof
Let\dss
$\mathcal{F}
\off =\off
\{\qff
F_{\dff i}
\qff\}_{\dff i\dff \in\qff I}$\nsp.\oss
For each\dss $x\qff \in\qff X$\dss there\dss is\dss an open set\dss
$V_{\fff x}$\dss such\dss that\dss $x\qff \in\qff V_{\fff x}$\dss and\dss
$V_{\fff x}\qff \cap\qff F_{\dff i}\off \neq\off \varnothing$\dss
for only\sss a\sss finite number of\dss $i\qff \in\pff I$\nnsp.\oss
Then\dss
$\mathcal{V}
\off =\off
\{\qff
V_{\fff x}
\qff\}_{\dff x\qff \in\qff X}$\dss
is\dss an open covering of\trs $X$\nnsp.\oss
By\qss The\-o\-rem\qss \ref{star-refinement}\qss
there\sss 
is\dss an open covering\dss
$\mathcal{W}
\qff =\pff
\{\qff
W_{\fff s}
\qff\}_{\dff s\qff \in\qff S}$\dss
star\dss refining\qss $\mathcal{V}$\nnsp.\oss
For 
$i\qff \in\pff I$\trs 
let\vspace{2.25pt}
\[
\quad
U_{\dff i}
\off =\off
\sta(\qff F_{\dff i}\fff,\pff \mathcal{W}\trf)
\pff.
\]

\vspace{-9.75pt}
Clearly\halfff,\pss
$F_{\dff i}\qff \subset\qff U_{\dff i}$\dss
and\dss $U_{\dff i}$\dss is\dss open\dss for every\dss $i$\nnsp.\oss 
Therefore,\oss it\dss is\dss sufficient\dss to prove\sss that\dss the family\trs
$\mathcal{U}
\off =\off
\{\qff
U_{\dff i}
\qff\}_{\dff i\dff \in\qff I}$\dss
is\dss locally\dss finite.\oss
Since\dss $\mathcal{W}$\dss is\dss an open covering,\oss
it\dss is\dss sufficient\dss to prove\sss that\dss
for every\dss $s\qff \in\qff S$\dss
the intersection\dss $W_{\fff s}\qff \cap\qff U_{\dff i}$\dss
is\dss non-empty\dss 
for no more\sss than a\sss finite number of\dss
$i\qff \in\pff I$\nnsp.\oss
Let\dss $s\qff \in\qff S$\nnsp.\oss
Since\dss $\mathcal{W}$\dss is\dss star\dss refining\dss $\mathcal{V}$\nnsp,\oss
there exists\dss $x\qff \in\qff X$\dss
such\dss that\qss\vspace{2.25pt}
\[
\quad
\sta(\qff W_{\fff s}\fff,\pff \mathcal{W}\trf)
\off \subset\off\dff
V_{\fff x}
\pff.
\]

\vspace{-9.25pt}
Suppose\sss that\qss $W_{\dff s}\qff \cap\qff U_{\dff i}\off \neq\off \varnothing$\nnsp.\oss
Then\dss by\dss the definition of\dss $U_{\dff i}$\dss there exists\dss $t\qff \in\qff S$\dss 
such\dss that\dss
$W_{\dff s}\qff \cap\qff W_{\dff t}
\off \neq\off 
\varnothing$\qss
and\qss
$W_{\dff t}\qff \cap\qff F_{\dff i}
\off \neq\off 
\varnothing$\nnsp.\oss
By\dss the\sss definition,\pss
$W_{\dff t}
\qff \subset\qff 
\sta(\qff W_{\fff s}\fff,\pff \mathcal{W}\trf)$\nnsp.\oss
It\dss follows\dss that\dss\vspace{3pt}\vspace{-0.375pt}
\[
\quad
\sta(\qff W_{\fff s}\fff,\pff \mathcal{W}\trf)
\qff \cap\qff 
F_{\dff i}
\off \neq\off \varnothing
\pff.
\]

\vspace{-9pt}\vspace{-0.375pt}
But\qss
$\sta(\qff W_{\fff s}\fff,\pff \mathcal{W}\trf)
\off \subset\off\dff
V_{\fff x}$\qss
and\dss hence\qss
$V_{\fff x}
\qff \cap\pff 
F_{\dff i}
\off \neq\off \varnothing$\nnsp.\oss
By\dss the choice of\trs $V_{\fff x}$\dss this\sss may\dss happen
only\dss for a finite number of\dss $i\qff \in\pff I$\nnsp.\oss
This completes\sss the proof\halfff.\oss  \eproof

\myuppar{Remark.}
This\sss lemma\dss is\dss due\sss to\dss C.H.\dss Dowker\qss \cite{do}.\oss
See\qss \cite{do},\oss Lemma\qss 1.\oss

\mypar{Theorem.}{paracompact-extension}
\emph{Let\dss $X$\dss be a\sss paracompact\sss space.\oss
Let\pss 
$\mathcal{F}
\off =\off
\{\qff
F_{\dff i}
\qff\}_{\dff i\dff \in\qff I}$\dss
be a closed\dss family\dss  
combinatorially\dss refining\sss an\sss open\qss family\qss
$\mathcal{U}
\off =\off
\{\qff
U_{\dff i}
\qff\}_{\dff i\dff \in\qff I}$\qss
of\qss subsets of\pss $X$\nnsp.\oss
If\pss $\mathcal{F}$\dss
is\dss locally\qss finite,\oss
then\dss there\dss exists\dss an\dss open\dss locally\dss finite\dss family\qss
$\mathcal{G}
\off =\off
\{\qff
G_{\dff i}
\qff\}_{\dff i\dff \in\qff I}$\qss
such\dss that}\vspace{2.5pt}
\[
\quad
F_{\dff i}
\off \subset\off 
G_{\dff i}
\off \subset\off
\overline{G}_{\dff i}
\off \subset\off
U_{\dff i}
\]

\vspace{-9.5pt}
\emph{for\dss every\dss $i\qff \in\pff I$\dss
and\dss the\sss nerves\sss of\pss $\mathcal{F}$\dss
and\qss $\mathcal{G}$\dss are equal.\oss}

\proof
Theorem\qss \ref{closed-to-open-extension}\qss implies\sss that\dss
$\mathcal{F}$\dss is\dss combinatorially\dss refining\dss
a\sss locally\dss finite open\dss family\dss 
$\mathcal{V}
\off =\off
\{\qff
V_{\fff i}
\qff\}_{\dff i\qff \in\qff I}$\nsp.\oss
Let\qss $W_{\fff i}\off =\off U_{\dff i}\qff \cap\qff V_{\fff i}$\nsp.\oss
Then\qss
$\mathcal{W}
\off =\off
\{\qff
W_{\fff i}
\qff\}_{\dff i\qff \in\qff I}$\dss
is\dss a\dss locally\dss finite open\dss family\sss and\dss 
$\mathcal{F}$\trs is\dss combinatorially\dss refining\dss $\mathcal{W}$\dnsp.\oss
Since\dss $X$\dss is\dss normal\dss by\qss Lemma\qss \ref{normal},\oss
we can apply\qss Theorem\qss \ref{general-extensions}\qss
to\dss $Z\off =\off X$\nnsp,\oss the family\dss $\mathcal{F}$\nnsp,\oss
and\dss the family\dss $\mathcal{W}$\dss is\dss the role of\dss $\mathcal{U}$\nnsp.\oss
Since\dss $W_{\fff i}\qff \subset\qff U_{\dff i}$\dss
for every\dss $i\qff \in\pff I$\nnsp,\oss
the family\qss $\mathcal{G}$\qss from\qss Theorem\qss \ref{general-extensions}\qss
has\sss the required\dss properties.\oss  \eproof

\myuppar{Remarks.}
The above presentation of\dss a\sss fragment\sss of\trs the\sss theory\sss
of\dss paracompact\sss spaces\dss is\dss large\-ly\dss based on\trs
Chapter\qss 5\qss of\trs the book\trs \cite{en}\trs of\pss R.\dss Engelking.\oss
Engelking's\dss priorities are quite different\dss from ours.\oss
For example,\oss Theorem\qss \ref{closed-refinement}\qss is\dss not\sss stated,\oss
but\sss only\sss mentioned\dss in\dss Remark\qss 5.1.7,\oss
and\qss Theorem\qss \ref{closed-to-open-extension}\qss
appears only\sss as\dss Problem\qss 5.5.17(a).\oss
Although\qss Theorem\qss \ref{paracompact-extension}\qss easily\dss follows
from\qss Theorem\qss \ref{closed-to-open-extension},\oss it\dss may\dss be new.\oss

\newpage
\mysection{Closed\qss subspaces\qss and\qss fundamental\qss groups}{nerves-groups}

\myuppar{Introduction.}
Let\dss $X$\dss be a\sss topological\sss space.\oss
The main\dss goal\sss of\trs this section\dss is\dss to prove an analogue of\trs
Theorem\qss \ref{covering-theorem-open}\qss 
for closed\dss locally\dss finite coverings of\trs $X$\nnsp.\oss
In\dss this case we need\dss to assume\sss that\dss 
$X$\dss is\dss Hausdorff\dss and\dss paracompact\halfff.\oss
In addition,\oss we need\dss to assume\sss that\dss
the main\dss results of\trs the\sss theory\sss of\dss covering spaces 
apply\dss to\dss $X$\nnsp,\oss
although\dss we will\dss not\dss use\sss this\sss theory\halfff,\oss
at\dss least\dss not\sss explicitly\halfff.\oss
Since\sss the bounded cohomology\dss theory\dss is\dss intimately\dss
related\dss to fundamental\dss groups and coverings,\oss
such an assumption\dss seems\sss to be\sss only\dss natural.\oss
Somewhat\sss surprisingly\halfff,\oss
a\sss weaker assumption about\sss elements of\dss coverings\dss 
turns out\dss to be sufficient\halfff.\oss\vspace{1pt}

Before\sss stating\sss our assumptions explicitly\halfff,\pss
let\dss us\dss review\sss some standard\sss definitions.\oss
By\sss a\sss neighborhood\dss 
we will\dss always understand an open set\halfff.\oss
So,\oss a subset\dss $U\qff \subset\qff Z$\dss of\dss a\sss topological\sss space\dss $Z$\dss is\dss 
a\qss \emph{neighborhood}\qss of\dss $z\qff \in\qff Z$\qss
if\trs $z\qff \in\qff U$\dss and\dss $U$\dss is\dss open.\oss
A\sss topological\sss space\dss $Z$\dss is\dss called\qss \emph{locally\dss path-connected}\pss 
if\trs every\dss point\dss $z\qff \in\qff Z$\qss has a path-connected\dss neighborhood.\oss\vspace{1pt}

It\dss is\dss well\dss known\dss that\dss the relation of\dss 
being connected\dss by\sss a\sss path\dss is\dss 
an equivalence relation on\dss $Z$\nnsp.\oss
The\qss \emph{path components}\qss of\trs $Z$\dss are\sss
the equivalence classes with respect\dss to\sss this equivalence relation.\oss
It\dss is\dss well\dss known\dss that\trs
$Z$\dss is\dss locally\dss path-connected\trs
if\trs and\dss only\trs if\trs for every\sss open set\dss $U\qff \subset\qff Z$\dss
each path component\dss of\pss $U$\sss is\dss open\dss in\dss $Z$\nnsp.\oss
See\qss \cite{mu},\oss Theorem\dss 25.4.\oss\vspace{1pt}

A\sss space\dss $Z$\dss is\dss called\qss \emph{semilocally\dss simply-connected}\pss 
if\dss every\dss point\dss
$z\qff \in\qff Z$\dss has\dss a\sss path connected\dss 
neighborhood\dss $U$\dss 
such\dss that\dss the inclusion\dss homomorphism\dss
$\pi_{\dff 1}\dff(\trf U\fff,\qff z\trf)
\qff \ttoo\qff
\pi_{\dff 1}\dff(\trf Z\fff,\qff z\trf)$\dss
is\dss trivial.\oss
A\sss space\dss $Z$\dss has a universal\sss covering\sss space\dss
if\trs and\dss only\trs if\qss $Z$\dss is\dss path connected,\oss
locally\dss path connected,\oss
and semilocally\sss simply\sss connected.\oss
See\qss \cite{mu},\oss Corollary\qss 82.2.\oss

\myuppar{The assumptions.}
For\dss the rest\sss of\trs this section\dss $X$\dss will\dss be
a\sss paracompact\trs Hausdorff\trs space\dss which\dss is\dss 
path connected,\oss
locally\dss path connected,\oss
and semilocally\sss simply\sss connected.\oss

\myuppar{Simple subsets.}
A subset\dss $U\qff \subset\qff X$\trs
is\dss said\dss to be\qss \emph{simple}\pss if\trs $U$\dss is\dss open and\dss for every\dss
path component\dss $V$\dss of\trs $U$\dss 
the inclusion\dss homomorphism\dss
$\pi_{\dff 1}\dff(\trf V\fff,\qff v\trf)
\qff \ttoo\qff
\pi_{\dff 1}\dff(\trf X\fff,\qff v\trf)$\dss
is\dss trivial,\oss
where\dss $v$\dss is\dss an arbitrary\dss point\sss of\trs $V$\dnsp.\oss
Equivalently\halfff,\pss $U$\dss is\qss \emph{simple}\pss if\trs $U$\dss is\dss open
and every\dss loop\sss in\dss $U$\dss is\dss contractible\sss in\dss $X$\nnsp.\oss
Clearly\halfff,\oss an open subset\sss of\dss a simple subset\dss is\dss simple.\oss

\mypar{Lemma.}{subspace-paracompact}
\emph{A closed subspace of\dss a\sss paracompact\dss space\dss is\dss paracompact\halfff.\oss}

\proof
Let\dss $Z$\dss be a\sss paracompact\sss space and\dss let\trs
$F\qff \subset\qff Z$\dss a closed subset\halfff.\oss
Let\dss
$\mathcal{V}
\off =\off
\{\qff
V_{\fff i}
\qff\}_{\dff i\dff \in\qff I}$\dss
be a covering of\trs $F$\dss by\dss subsets\dss $V_{\fff i}\qff \subset\pff F$\dss open in\dss $F$\dnsp.\oss
For every\dss $i\qff \in\pff I$\trs there exists a subset\dss $U_{\dff i}\qff \subset\qff Z$\dss
such\dss that\trs $V_{\fff i}\off =\off F\qff \cap\qff U_{\dff i}$\dss
and\dss $U_{\dff i}$\dss is\dss open\dss in\dss $Z$\nnsp.\oss
The family\trs
$\{\qff
U_{\dff i}
\qff\}_{\dff i\dff \in\qff I}$\dss
together\dss with\dss $Z\qff \smallsetminus\qff F$\dss
is\dss an open covering\dss of\trs $Z$\nnsp.\oss
Since\dss $Z$\dss is\dss paracompact\halfff,\oss
there exists a\sss locally\dss finite open covering\dss
$\mathcal{W}
\off =\off
\{\qff
W_{\dff k}
\qff\}_{\dff k\dff \in\qff K}$\dss
refining\dss this covering.\oss
Clearly\halfff,\oss
the family\dss
$\{\qff
W_{\dff k}\qff \cap\qff F
\qff\}_{\dff k\dff \in\qff K}$\dss
is\dss a\sss locally\dss finite covering of\trs $F$\dss by\sss sets open\dss in\dss $F$\dss
which\dss is\dss refining\dss $\mathcal{V}$\dnsp.\oss

\mypar{Theorem.}{pi-one-extension}
\emph{Let\dss $A\qff \subset\qff X$\dss be a\sss path connected and\dss locally\dss path connected
subset\sss and\dss $a\qff \in\qff A$\nnsp.\oss
Then\dss there exists a\sss path connected\sss open set\qss $U\qff \subset\qff X$\dss 
such\dss that\qss
$A\qff \subset\qff U$\dss
and\dss 
the images of\pss the inclusion\dss homomorphisms\dss 
$\pi_{\dff 1}\dff(\trf A\fff,\qff a\trf)
\qff \ttoo\qff
\pi_{\dff 1}\dff(\trf X\fff,\qff a\trf)$\dss
and\qss
$\pi_{\dff 1}\dff(\trf U\fff,\qff a\trf)
\qff \ttoo\qff
\pi_{\dff 1}\dff(\trf X\fff,\qff a\trf)$\dss
are equal.}

\proof
The proof\dss deals with several\dss families simultaneously\halfff.\oss
In order\dss not\dss to overburden\dss the\sss text\dss with subscripts,\oss
let\dss us\sss pretend\dss that\dss families are collections of\dss subsets.\oss
More precisely\halfff,\oss given a\sss family\dss
$\mathcal{F}
\off =\off
\{\qff
F_{\dff i}
\qff\}_{\dff i\dff \in\qff I}$\nsp,\oss
let\dss us agree\sss that\trs
$F\qff \in\qff \mathcal{F}$\dss means\sss that\trs
$F\off =\off F_{\dff i}$\dss 
for some\dss $i\qff \in\pff I$\nnsp.\oss 

Since\dss $X$\dss is\dss semilocally\dss simply-connected,\oss
there\dss exists\sss a covering\dss 
$\mathcal{S}$\dss 
of\trs $X$\trs
by\dss path connected simple sets.\oss
Since\dss $X$\dss is\dss paracompact\halfff,\oss
by\qss Theorem\qss \ref{star-refinement}\qss there exists an open covering\dss
$\mathcal{T}$\dss 
of\trs $X$\dss
star\dss refining\dss $\mathcal{S}$\dnsp.\qff\oss
If\pss $T\fff,\pff T\fff'\qff \in\qff \mathcal{T}$\trs
and\qss
$T\qff \cap\qff T\fff'\off \neq\off \varnothing$\nnsp,\oss
then\qss
$T\fff'
\off \subset\off 
\sta(\trf T\fff,\pff \mathcal{T}\trf)$\dss
and\dss hence\qss
$T\qff \cup\qff T\fff'
\off \subset\off
S$\dss
for some\dss $S\qff \in\qff \mathcal{S}$\nnsp.\oss
It\dss follows\dss that\dss the subset\qss
$T\qff \cup\qff T\fff'$\qss is\dss simple.\oss

Let\dss $A\qff \cap\qff \mathcal{T}$\dss be\sss the family\sss
of\dss intersections\dss $A\qff \cap\qff T$\dss with\dss
$T\qff \in\qff \mathcal{T}$\dnsp.\oss
Clearly\halfff,\pss 
$A\qff \cap\qff \mathcal{T}$\dss
is\dss a covering of\dss $A$\dss be sets open\sss in\dss $A$\nnsp.\oss
Since\dss $A$\dss is\dss locally\dss path-connected,\oss
there exists a covering\dss 
$\mathcal{P}$\dss 
of\dss $A$\dss
by\dss path-connected subsets open\dss in\dss $A$\dss such\dss that\dss
$\mathcal{P}$\dss is\dss refining\dss $A\qff \cap\qff \mathcal{T}$\dnsp.\oss
Then\dss for every\dss $P\qff \in\qff \mathcal{P}$\dss
there exists\dss $T\qff \in\qff \mathcal{T}$\dss
such\dss that\trs
$P\qff \subset\pff T$\dnsp.\oss
Let\dss us\dss choose such a set\trs 
$T\off =\off T\dff(\trf P\trf)$\dss
for each\trs $P\qff \in\qff \mathcal{P}$\dnsp.\oss

Since\dss $X$\dss is\dss locally\dss path-connected,\oss
for every\trs $P\qff \in\qff \mathcal{P}$\dss
there exists an\sss open set\dss $V\trf(\trf P\trf)$\dss
such\dss that\qss
$P\qff \subset\qff
V\trf(\trf P\trf)
\qff \subset\qff
T\trf(\trf P\trf)$\qss
and\dss every\dss point\dss in\dss
$V\trf(\trf P\trf)$\dss
is\dss connected\dss by\sss a path with some point\sss in\dss $P$\nnsp.\oss
Since\dss $P$\dss is\dss path connected,\pss
$V\trf(\trf P\trf)$\dss is\dss also path connected.\oss

By\qss Theorem\qss \ref{closed-refinement}\qss there\dss is\dss
a closed\dss locally\dss finite covering\dss $\mathcal{F}$\dss
of\dss $A$\dss combinatorially\dss refining\dss $\mathcal{P}$\dnsp.\oss
For each\qss $F\qff \in\qff \mathcal{F}$\qss
let\trs $P\trf(\trf F\trf)\qff \in\qff \mathcal{P}$\qss
be\sss the subset\sss corresponding\dss to\dss $F$\dnsp,\oss
so\sss that\trs $F\qff \subset\qff P\trf(\trf F\trf)$\nnsp.\oss
Let\vspace{3pt}
\[
\quad
V\dff(\trf F\trf)\off =\off
V\trf\bigl(\trf P\dff(\trf F\trf)\trf\bigr)
\quad\
\mbox{and}\quad\
T\dff(\trf F\trf)\off =\off
T\trf\bigl(\trf P\dff(\trf F\trf)\trf\bigr)
\pff
\]

\vspace{-9pt}
for each\qss $F\qff \in\qff \mathcal{F}$\dnsp.\oss
Then\qss\vspace{3pt}
\[
\quad
F
\qff \subset\qff
V\dff(\trf F\trf)
\qff \subset\pff
T\dff(\trf F\trf)
\]

\vspace{-9pt}
for every\qss $F\qff \in\qff \mathcal{F}$\dnsp.\oss
By\qss Theorem\qss \ref{closed-to-open-extension}\qss
there exists\dss a\dss locally\trs finite open\dss family\qss
$\mathcal{W}$\dss
of\qss subsets\sss of\pss $X$\dss
such\dss that\dss
$\mathcal{F}$\dss is\dss combinatorially\dss refining\dss
$\mathcal{W}$\nsp\dnsp.\dff\oss
For\sss each\qss $F\qff \in\qff \mathcal{F}$\qss
let\trs $W\dff(\trf F\trf)\qff \in\qff \mathcal{W}$\qss
be\sss the subset\sss corresponding\dss to\dss $F$\nnsp,\oss
so\sss that\trs $F\qff \subset\qff W\dff(\trf F\trf)$\nnsp.\oss
Let\qss
$U\trf(\trf F\trf)
\off =\off
V\dff(\trf F\trf)\qff \cap\qff W\dff(\trf F\trf)$\nnsp.\oss
Then\dss\vspace{3pt}\vspace{0.875pt}
\[
\quad
 F
\qff \subset\qff 
U\trf(\trf F\trf)
\qff \subset\qff 
W\dff(\trf F\trf)
\]

\vspace{-9pt}\vspace{0.875pt}
and\dss hence\sss sets\qss $U\trf(\trf F\trf)$\nnsp,\qss
$F\qff \in\qff \mathcal{F}$\nnsp,\qff\oss 
form an open\dss locally\dss finite\sss family\dss $\mathcal{U}$\dnsp.\oss
By\qss Theorem\qss \ref{paracompact-extension}\qss
there exists a family\sss of\dss open\dss in\dss $X$\dss sets\dss
$G\trf(\trf F\trf)$\nnsp,\oss where\qss
$F\qff \in\qff \mathcal{F}$\nnsp,\qff\oss
such\dss that\vspace{3pt}
\[
\quad
F
\off \subset\off
G\trf(\trf F\trf)
\off \subset\off
\overline{G\trf(\trf F\trf)}
\off \subset\off
U\trf(\trf F\trf)
\]

\vspace{-9pt}
for each\qss $F\qff \in\qff \mathcal{F}$\qss
and\dss the family\trs $\mathcal{G}$\dss of\trs these sets\sss has\sss the nerve
equal\dss to\sss the\sss nerve\dss of\qss $\mathcal{F}$\nnsp.\oss

In\dss particular\halfff,\pss
$G\trf(\trf F\trf)\qff \cap\qff G\trf(\trf F\fff'\trf)
\off \neq\off
\varnothing$\qss
if\qss and\dss only\qss if\pss
$F\qff \cap\qff F\fff'
\off \neq\off
\varnothing$\nnsp.\oss
Since\vspace{3pt}
\[
\quad
G\trf(\trf F\trf)
\off \subset\off
U\trf(\trf F\trf)
\off \subset\off
V\trf(\trf F\trf)
\off \subset\off
T\trf(\trf F\trf)
\]

\vspace{-9pt}
for every\qss $F\qff \in\qff \mathcal{F}$\dnsp,\oss
it\dss follows\dss that\trs
$T\trf(\trf F\trf)\qff \cap\qff T\trf(\trf F\fff'\trf)
\off \neq\off
\varnothing$\trs
if\qss
$F\qff \cap\qff F\fff'
\off \neq\off
\varnothing$\nnsp.\oss
This,\oss in\dss turn,\oss implies\sss that\dss the union\qss
$T\trf(\trf F\trf)\qff \cup\qff T\trf(\trf F\fff'\trf)$\qss
is\dss simple\dss if\qss
$F\qff \cap\qff F\fff'
\off \neq\off
\varnothing$\nnsp.\oss

Since\dss $\mathcal{F}$\dss is\dss a covering\sss of\dss $A$\nnsp,\oss
i.e.\qss the union of\trs the sets\dss $F\qff \in\qff \mathcal{F}$\dss
is\dss equal\dss to\dss $A$\nnsp,\oss
the union of\trs the sets\dss
$G\trf(\trf F\trf)$\nnsp,\qss
$F\qff \in\qff \mathcal{F}$\nnsp,\qff\oss
contains\dss $A$\nnsp.\oss
Let\dss $U$\dss be\sss the path component\sss of\trs this union
containing\dss $A$\nnsp.\oss
It\dss is\dss sufficient\dss to show\dss that\dss every\dss loop\sss in\dss $U$\dss
based at\dss $a$\dss is\dss homotopic\sss in\dss $X$\trs relatively\dss to\dss $a$\dss 
to a\sss loop\sss in\dss $A$\nnsp.\oss
Let\dss
$p\dff \colon\dff
[\trf 0\fff,\qff 1\trf]
\qff \ttoo\qff
U$\trs
be\sss such a\sss loop.\oss
Let\dss us\dss partition\dss
$[\trf 0\fff,\qff 1\trf]$\dss
by\sss several\dss points\qss\vspace{3pt}\vspace{-1.125pt}
\[
\quad
0\off =\off x_{\trf 0}\qff <\qff
x_{\trf 1}\qff <\qff
\ldots\qff <\qff
x_{\dff n}\off =\off 1
\]

\vspace{-9pt}\vspace{-1.125pt}
into\sss subintervals\dss
$I_{\dff i}\off =\off [\trf x_{\dff i\dff -\dff 1}\dff,\qff x_{\dff i}\trf]$\nnsp.\oss
As\dss is\dss well\dss known,\oss one can choose\sss this\sss partition\dss
in\sss such a\sss way\dss that\dss
$p\trf(\qff I_{\dff i}\trf)\qff \subset\pff G\trf(\trf F_{\dff i}\trf)$\dss
for\dss some\dss $F_{\dff i}\qff \in\qff \mathcal{F}$\dss for every\dss $i$\nnsp.\oss
Let\dss 
$y_{\trf 1}\dff,\pff
y_{\trf 2}\dff,\pff
\ldots\dff,\pff
y_{\trf n}
\qff \in\pff
[\trf 0\fff,\qff 1\trf]$\qss
be some numbers such\dss that\dss
$x_{\dff i\dff -\dff 1}\qff <\qff y_{\dff i}\qff <\qff x_{\dff i}$\dss
for every\dss $i\qff \geq\qff 1$\nnsp.\oss
Since\vspace{3pt}
\[
\quad
p\trf(\trf y_{\dff i}\trf)
\off \in\off
G\trf(\trf F_{\dff i}\trf)
\off \subset\off
V\trf(\trf F_{\dff i}\trf)
\pff,
\]

\vspace{-9pt}
$V\trf(\trf F_{\dff i}\trf)$\sss is\dss path connected,\oss 
and\dss
$F_{\dff i}\qff \subset\qff V\trf(\trf F_{\dff i}\trf)$\nnsp,\oss
for each\qss $i\qff \geq\qff 1$\dss
we can connect\trs $p\trf(\trf y_{\dff i}\trf)$\dss by\dss a\sss path\dss
$q_{\dff i}$\dss in\sss
$V\trf(\trf F_{\dff i}\trf)$\sss with\sss some\sss point\dss 
$a_{\dff i}\qff \in\pff F_{\dff i}$\nsp.\oss
On\dss the other\dss hand,\oss
if\trs $i\qff \leq\qff n\qff -\qff 1$\nnsp,\oss
then\vspace{3pt}
\[
\quad
p\dff(\trf x_{\dff i}\trf)
\off \in\off\dff
G\trf(\trf F_{\dff i}\trf)
\qff \cap\qff  
G\trf(\trf F_{\dff i\dff +\dff 1}\trf)
\]

\vspace{-9pt}
and\dss hence\qss
$G\trf(\trf F_{\dff i}\trf)\qff \cap\qff  G\trf(\trf F_{\dff i\dff +\dff 1}\trf)
\off \neq\off
\varnothing$\nnsp.\oss
As we saw,\oss this implies\dss that\dss
$F_{\dff i}\qff \cap\qff F_{\dff i\dff +\dff 1}
\off \neq\off
\varnothing$\dss
and\dss hence\vspace{3pt}
\[
\quad
P\trf(\trf F_{\dff i}\trf)
\qff \cap\qff  
P\trf(\trf F_{\dff i\dff +\dff 1}\trf)
\off \neq\off
\varnothing
\pff.
\]

\vspace{-9pt}
Since\dss $P\trf(\trf F_{\dff i}\trf)$\dss 
and\dss $P\trf(\trf F_{\dff i\dff +\dff 1}\trf)$\dss
are path connected,\oss
we can connect\dss $a_{\dff i}$\dss with\dss $a_{\dff i\dff +\dff 1}$\dss
by\sss a\sss path\dss $r_{\dff i}$\dss in\dss
$P\trf(\trf F_{\dff i}\trf)
\qff \cup\qff  
P\trf(\trf F_{\dff i\dff +\dff 1}\trf)$\nnsp.\oss
Similarly\halfff,\oss
we can connect\dss $a$\dss with\sss $a_{\trf 1}$\dss 
by\sss a\dss path\sss $r_{\trf 0}$\dss
in\dss $P\trf(\trf F_{\trf 1}\trf)$\dss
and\dss connect\dss $a_{\dff n}$\dss with\dss $a$\dss 
by\sss a\dss path\dss $r_{\dff n}$\dss
in\trs $P\trf(\trf F_{\dff n}\trf)$\nnsp.\oss

Let\dss $y_{\trf 0}\off =\off 0$\nnsp,\pss
$y_{\trf n\dff +\dff 1}\off =\off 1$\nnsp,\oss
and\dss for\dss
$i\off =\off 0\fff,\pff 1\fff,\pff \ldots\fff,\pff n$\qss
let\qss
$f_{\dff i}\dff \colon\dff
[\trf 0\fff,\qff 1\trf]
\qff \ttoo\qff 
[\trf y_{\dff i}\fff,\qff y_{\dff i\dff +\dff 1}\trf]$\qss
be an\dss increasing\dss homeomorphism\sss and\trs
$p_{\dff i}\off =\off p\dff \circ\dff f_{\dff i}$\nsp.\oss
Then\dss for every\dss $i$\dss the path\dss 
$p_{\dff i}$\dss  
connects\dss $p\trf(\trf y_{\dff i}\trf)$\dss with\dss $p\trf(\trf y_{\dff i\dff +\dff 1}\trf)$\dss
in\dss
$G\trf(\trf F_{\dff i}\trf)\qff \cup\qff  G\trf(\trf F_{\dff i\dff +\dff 1}\trf)$\dss
and\dss 
$p$\dss is\dss homotopic\sss to\sss the product\dss 
$p_{\trf 0}\dff \cdot\dff
p_{\trf 1}\dff \cdot\dff
\ldots\dff \cdot\dff
p_{\trf n}$\nsp.\oss
For every\dss
$i\off =\off 1\fff,\pff 2\fff,\pff \ldots\fff,\pff n$\qss
the path\dss 
$q_{\dff i}^{\dff -\dff 1}\dff \cdot\qff
p_{\dff i}\dff \cdot\qff
q_{\dff i\dff +\dff 1}$\dss
connects\dss $a_{\dff i}$\dss with\dss $a_{\dff i\dff +\dff 1}$\dss
in\vspace{3pt}
\[
\quad
G\trf(\trf F_{\dff i}\trf)
\qff \cup\qff  
G\trf(\trf F_{\dff i\dff +\dff 1}\trf)
\off \subset\off
T\trf(\trf F_{\dff i}\trf)
\qff \cup\qff  
T\trf(\trf F_{\dff i\dff +\dff 1}\trf)
\pff
\]

\vspace{-9pt}
and\dss the path\dss $r_{\dff i}$\dss connects\dss 
$a_{\dff i}$\dss with\dss $a_{\dff i\dff +\dff 1}$\dss
in\vspace{3pt}
\[
\quad
P\trf(\trf F_{\dff i}\trf)
\qff \cup\qff  
P\trf(\trf F_{\dff i\dff +\dff 1}\trf)
\off \subset\off
T\trf(\trf F_{\dff i}\trf)
\qff \cup\qff  
T\trf(\trf F_{\dff i\dff +\dff 1}\trf)
\pff.
\]

\vspace{-9pt}
Since\dss
$F_{\dff i}\qff \cap\qff F_{\dff i\dff +\dff 1}
\off \neq\off
\varnothing$\nnsp,\oss
the subset\dss 
$T\trf(\trf F_{\dff i}\trf)
\qff \cup\qff  
T\trf(\trf F_{\dff i\dff +\dff 1}\trf)$\dss
is\dss simple and\dss hence\dss
$q_{\dff i}^{\dff -\dff 1}\dff \cdot\qff
p_{\dff i}\dff \cdot\qff
q_{\dff i\dff +\dff 1}$\dss
is\dss homotopic\sss to\dss $r_{\dff i}$\dss in\dss $X$\dss
relatively\dss to\sss the endpoints.\oss
Similarly\halfff,\oss the paths\dss 
$p_{\trf 0}\dff \cdot\dff q_{\trf 1}$\dss
and\dss $r_{\trf 0}$\dss
connect\dss $a$\dss with\dss $a_{\trf 1}$\dss
in\dss $G\trf(\trf F_{\dff 1}\trf)$\dss and\dss
$P\trf(\trf F_{\dff 1}\trf)$\dss respectively\halfff.\oss
Since\sss both\dss $G\trf(\trf F_{\dff 1}\trf)$\dss and\dss
$P\trf(\trf F_{\dff 1}\trf)$\dss
are contained\dss in\dss
$T\trf(\trf F_{\dff 1}\trf)$\dss
and\dss $T\trf(\trf F_{\dff 1}\trf)$\dss is\dss simple,\oss
the path\dss $p_{\trf 0}\dff \cdot\dff q_{\trf 1}$\dss
and\dss $r_{\trf 0}$\dss
are\sss homotopic in\dss $X$\dss relatively\dss to\sss the endpoints.\oss
Similarly\halfff,\pss
$q_{\dff n}^{\dff -\dff 1}\dff \cdot\dff p_{\trf n}$\dss
and\dss $r_{\trf n}$\dss
are\sss homotopic in\dss $X$\dss relatively\dss to\sss the endpoints.\oss
Clearly\halfff,\pss
$p_{\trf 0}\dff \cdot\dff
p_{\trf 1}\dff \cdot\dff
\ldots\dff \cdot\dff
p_{\trf n}$\dss
is\dss homotopic\sss 
in\dss $X$\dss relatively\dss to\sss the endpoints\sss to\vspace{3pt}
\[
\quad
p_{\trf 0}\dff \cdot\dff
q_{\trf 1}\dff \cdot\dff
q_{\trf 1}^{\dff -\dff 1}\dff \cdot\qff
p_{\trf 1}\dff \cdot\dff
q_{\trf 2}\dff \cdot\dff
q_{\trf 2}^{\dff -\dff 1}\dff \cdot\qff
p_{\trf 2}\dff \cdot\dff
q_{\trf 3}\dff \cdot\pff
\ldots\pff \cdot\dff
q_{\trf n\dff -\dff 1}^{\dff -\dff 1}\dff \cdot\qff
p_{\trf n\dff -\dff 1}\dff \cdot\dff
q_{\dff n}\dff \cdot\dff
q_{\dff n}^{\dff -\dff 1}\dff \cdot\dff p_{\trf n}
\]

\vspace{-36pt}
\[
\quad
=\off
\left(\qff
p_{\trf 0}\dff \cdot\dff
q_{\trf 1}
\qff\right)\dff \cdot\dff
\left(\qff
q_{\trf 1}^{\dff -\dff 1}\dff \cdot\qff
p_{\trf 1}\dff \cdot\dff
q_{\trf 2}
\qff\right)\dff \cdot\dff
\left(\qff
q_{\trf 2}^{\dff -\dff 1}\dff \cdot\qff
p_{\trf 2}\dff \cdot\dff
q_{\trf 3}
\qff\right)\dff \cdot\pff
\ldots\pff \cdot\dff
\left(\qff
q_{\trf n\dff -\dff 1}^{\dff -\dff 1}\dff \cdot\qff
p_{\trf n\dff -\dff 1}\dff \cdot\dff
q_{\dff n}
\qff\right)\dff \cdot\dff
\left(\qff
q_{\dff n}^{\dff -\dff 1}\dff \cdot\dff p_{\trf n}
\qff\right)
\]

\vspace{-9pt}
and\dss hence\dss is\dss homotopic\sss in\dss $X$\dss
relatively\dss to\sss the endpoints\sss to\dss
$r_{\trf 0}\dff \cdot\dff
r_{\trf 1}\dff \cdot\dff
\ldots\dff \cdot\dff
r_{\trf n}$\nsp.\oss
It\dss follows\dss that\dss $p$\dss
is\dss homotopic\sss in\dss $X$\dss
relatively\dss to\sss the endpoints\sss to\dss
$r_{\trf 0}\dff \cdot\dff
r_{\trf 1}\dff \cdot\dff
\ldots\dff \cdot\dff
r_{\trf n}$\nsp.\oss
By\dss the construction,\pss
$r_{\trf 0}\dff \cdot\dff
r_{\trf 1}\dff \cdot\dff
\ldots\dff \cdot\dff
r_{\trf n}$\qss
is\dss a\sss loop in\dss $A$\dss based\sss at\sss $a$\nnsp.\oss
Hence\sss every\dss loop\dss $p$\dss in\dss $U$\dss
based\sss at\sss $a$\dss
is\dss homotopic\sss in\dss $X$\dss
relatively\dss to\sss the endpoints\sss to
a\sss loop in\dss $A$\dss based at\dss $a$\nnsp.\oss  \eproof

\mypar{Theorem.}{pi-one-coverings-extension}
\emph{Suppose\sss that\qss
$\mathcal{F}
\off =\off
\{\qff
F_{\dff i}
\qff\}_{\dff i\dff \in\qff I}$\qss
be a closed\dss locally\dss finite covering\sss of\pss $X$\dss
by\dss path connected and\dss locally\dss path connected
subsets and\dss that\dss $a_{\dff i}\qff \in\qff F_{\dff i}$\dss
for every\dss $i\qff \in\pff I$\nnsp.\oss
Then\dss there\dss exists\dss an\dss open\dss locally\dss finite\dss covering\pss
$\mathcal{G}
\off =\off
\{\qff
G_{\dff i}
\qff\}_{\dff i\dff \in\qff I}$\qss
of\pss $X$\dss by\dss path connected subsets
such\dss that\qss
$F_{\dff i}\qff \subset\qff G_{\dff i}$\qss
for every\qss $i\qff \in\pff I$\nnsp,\oss
the nerves of\qss $\mathcal{F}$\dss and\qss $\mathcal{G}$\dss are equal,\oss
and\dss the images of\pss}\dss\vspace{3pt}
\[
\quad
\pi_{\dff 1}\dff(\trf F_{\dff i}\dff,\qff a_{\dff i}\trf)
\qff \ttoo\qff
\pi_{\dff 1}\dff(\trf X\dff,\qff a_{\dff i}\trf)
\quad\
\mbox{\emph{and}}\quad\dff\
\pi_{\dff 1}\dff(\trf G_{\dff i}\dff,\qff a_{\dff i}\trf)
\qff \ttoo\qff
\pi_{\dff 1}\dff(\trf X\dff,\qff a_{\dff i}\trf)
\]

\vspace{-9pt}
\emph{are equal\dss for every\dss $i\qff \in\pff I$\nnsp.\oss}

\proof
By\qss Theorem\qss \ref{pi-one-extension}\qss
for every\dss $i\qff \in\pff I$\dss there exists a\sss path connected\sss open set\dss $U_{\dff i}$\dss
such\dss that\dss \qss
$F_{\dff i}\qff \subset\qff U_{\dff i}$\dss
and\dss 
the images of\pss the inclusion\dss homomorphisms\dss\vspace{3pt}
\[
\quad
\pi_{\dff 1}\dff(\trf F_{\dff i}\fff,\qff a_{\dff i}\trf)
\qff \ttoo\qff
\pi_{\dff 1}\dff(\trf X\fff,\qff a_{\dff i}\trf)
\quad\
\mbox{and}\quad\
\pi_{\dff 1}\dff(\trf U_{\dff i}\fff,\qff a_{\dff i}\trf)
\qff \ttoo\qff
\pi_{\dff 1}\dff(\trf X\fff,\qff a_{\dff i}\trf)
\]

\vspace{-9pt}
are equal.\oss
Clearly\halfff,\oss if\qss $U_{\dff i}$\dss has\sss this property\halfff,\oss
then every\sss open\dss path connected\sss set\dss $G_{\dff i}$\dss
such\dss that\dss
$F_{\dff i}\qff \subset\qff G_{\dff i}\qff \subset\qff U_{\dff i}$\dss also does.\oss
Hence\sss the\sss theorem\dss follows\dss from\qss
Theorem\qss \ref{paracompact-extension}.\oss  \eproof

\mypar{Theorem.}{covering-theorem-closed}
\emph{Suppose\sss that\qss
$\mathcal{F}
\off =\off
\{\qff
F_{\dff i}
\qff\}_{\dff i\dff \in\qff I}$\qss
be a closed\dss locally\dss finite covering\sss of\pss $X$\dss
by\dss path connected and\dss locally\dss path connected
subsets.\oss
Let\trs $N$\dss be\sss the nerve of\qss the covering\qss $\mathcal{F}$\dnsp.\oss
If\qss the covering\dss $\mathcal{F}$\dss is\dss
weakly\dss boundedly\sss acyclic,\oss
then\qss
$\widehat{H}^{\dff *}\fff(\trf X \trf)
\qff \ttoo\qff
H^{\fff *}\fff(\trf X \trf)$\qss
can be factored\dss through\qss
$H^{\dff *}\fff(\trf N \trf)
\qff \ttoo\qff 
H^{\dff *}\fff(\trf X \trf)$\dnsp.\oss}

\proof
Let\dss $\mathcal{G}$\dss be\sss the covering\dss provided\dss by\qss 
Theorem\qss \ref{pi-one-coverings-extension}.\oss 
Since\dss $\mathcal{F}$\dss is\dss weakly\dss boundedly\sss acyclic,\oss
these properties imply\dss that\dss $\mathcal{G}$\dss is\dss 
an\sss open weakly\dss boundedly\sss acyclic
covering\dss with\dss the same nerve\dss $N$\dss as\dss $\mathcal{F}$\dnsp.\oss
Hence\sss the\sss theorem\dss follows\dss from\qss
Theorem\qss \ref{covering-theorem-open}.\oss  \eproof

\newpage
\mysection{Closed\qss subspaces\qss and\qss homology\qss groups}{closed-homology}

\myuppar{Introduction.}
Let\dss $X$\dss be a\dss Hausdorff\dss space 
and\sss $\mathcal{U}$\sss be a covering of\dss $X$\nnsp.\oss
The\sss goal\sss of\trs this section\dss is\dss to prove an analogue of\trs
Theorem\qss \ref{A-acyclic-open}\qss 
when\sss $\mathcal{U}$\sss is\dss  closed and\dss locally\dss finite.\oss
As\sss in\dss Section\qss \ref{extensions-coverings},\oss
this analogue\sss implies an analogue of\trs
Theorem\qss \ref{covering-theorem-open}\qss 
for such $\mathcal{U}$\nnsp.\oss
The\sss latter\dss is\dss different\dss from\trs
Theorem\qss \ref{covering-theorem-closed}\qss and complements\sss it\halfff;\oss
neither of\trs them\sss implies\sss the other\halfff.\oss
In contrast\dss with\trs Theorem\qss \ref{covering-theorem-closed},\oss
now\sss we will\sss assume\sss that\dss the closed subsets covering\sss $X$\sss
behave nicely\dss not\dss with respect\dss to\sss fundamental\dss groups,\oss
but\dss with\sss respect\dss to\sss the usual\sss singular\sss homology\sss groups.\oss

A space\sss $Z$\sss is\dss said\dss to be\qss 
\emph{homologically\dss locally\sss connected}\pss
if\trs for each\dss $n\qff \geq\qff 0$\nnsp,\oss 
each\sss $z\qff \in\qff Z$\nnsp,\oss
and each neighborhood\sss $U$\sss of\dss $z$\sss
there exists another\sss neighborhood\sss 
$V\qff \subset\qff U$\sss of\dss $z$\sss
such\dss that\dss the inclusion\sss homomorphism\sss
$H_{\fff n}\fff(\trf V,\qff \{\dff z\trf\}\trf)
\qff \ttoo\qff
H_{\fff n}\fff(\trf U,\qff \{\dff z\trf\}\trf)$\sss
is\dss equal\dss to zero.\oss
Here\sss $H_{\fff n}\fff(\trf \bullet\trf)$\sss
are\sss the usual\sss singular\sss homology\dss groups 
with\sss integer coefficients.\oss
A homologically\dss locally\sss connected space\dss is\dss also called\sss
an $H{\fff}L{\fff}C$ space,\oss
or\halfff,\oss more precisely,\oss
an\sss  $H{\fff}L{\fff}C^{\dff \infty}_{\dff \zzz}$\sss space.\oss

\myuppar{Sheaves associated\sss with singular cochains.}
For every\dss $q\qff \geq\qff 0$\dss let\dss $\mathcal{C}^{\dff q}$\dss 
be the presheaf\dss on\sss $X$\sss
assigning\dss to an\sss open subset\dss $U\qff \subset\qff X$\dss the  
vector space of\trs real-valued\sss cochains\dss $C^{\dff q}\dff(\trf U \trf)$\dss
and\dss to an\dss inclusion\dss $U\qff \subset\qff V$\dss
the restriction homomorphism\dss
$C^{\dff q}\dff(\trf V \trf)
\qff \ttoo\qff 
C^{\dff q}\dff(\trf U \trf)$\nnsp.\oss
Let\dss $\mathcal{C}^{\dff -\dff 1}$\dss be\sss the constant\dss presheaf\trs $\rrr$\nnsp.\oss
For\dss  
$q\qff \geq\qff -\qff 1$\dss
let\dss $\Gamma^{\fff q}$\dss be\sss the sheaf\dss 
associated\dss with\dss the presheaf\dss $\mathcal{C}^q$\dnsp.\oss
The maps\qss
$d_{\dff q}
\dff \colon\dff
C^{\dff q}\dff(\trf U \trf)
\qff \ttoo\qff
C^{\dff q\dff +\dff 1}\dff(\trf V \trf)$\dss
lead\dss to morphisms\dss 
$d_{\dff q}
\dff \colon\dff
\mathcal{C}^{\dff q}
\qff \ttoo\qff
\mathcal{C}^{\dff q\dff +\dff 1}$\dnsp,\oss
which,\oss in\dss turn,\oss lead\dss to morphisms\dss 
$d_{\dff q}
\dff \colon\dff
\Gamma^{\dff q}
\qff \ttoo\qff
\Gamma^{\dff q\dff +\dff 1}$\dnsp.\oss

For every\dss $Y\qff \subset\qff X$\dss
let\dss $\Gamma^{\dff q}\dff(\trf Y\trf)$\dss 
be\sss the space of\dss sections of\trs the sheaf\dss $\Gamma^{\dff q}$\dss
over\dss $Y$\nnsp.\oss
The morphisms\dss
$d_{\dff q}
\dff \colon\dff
\Gamma^{\dff q}
\qff \ttoo\qff
\Gamma^{\dff q\dff +\dff 1}$\dss
lead\dss to homomorphisms\dss
$\Gamma^{\dff q}\dff(\trf Y\trf)
\qff \ttoo\qff
\Gamma^{\dff q\dff +\dff 1}\dff(\trf Y\trf)$\dss
and\dss hence\sss to augmented\sss cochain complexes\dss
$\Gamma^{\trf \bullet}\dff(\trf Y\trf)$\nnsp.\oss
When\dss $Z\qff \subset\qff Y\qff \subset\qff X$\nnsp,\oss
the restriction of\dss sections\sss leads\sss to\sss the restriction\dss morphism\dss
$\Gamma^{\trf \bullet}\dff(\trf Y \trf)
\qff \ttoo\qff
\Gamma^{\trf \bullet}\dff(\trf Z \trf)$\nnsp.\oss
By\dss considering\sss only\sss  
subsets\dss
$Y\qff \in\qff \mathcal{U}^{\dff \cap}$\dss
we get\sss a functor\dss from
$\cat \mathcal{U}$
to augmented cochain\sss complexes,\oss
which we still\sss denote\dss by\dss $\Gamma^{\trf \bullet}$\dnsp.\oss
By\sss applying\sss constructions\sss of\qss Section\qss \ref{cohomological-leray}\qss to\qss
$A^{\bullet}\off =\off \Gamma^{\trf \bullet}$\dss
we get\dss a double complex\dss
$C^{\dff \bullet}\dff(\trf N\fff,\pff \Gamma^{\trf \bullet} \trf)$\dss
and\dss morphisms\vspace{0pt}
\[
\quad
\begin{tikzcd}[column sep=large, row sep=huge]\dis
C^{\dff \bullet}\dff(\trf N \trf) 
\arrow[r, "\dis i_{\qff \Gamma}"]
& 
T^{\dff \bullet}\dff(\trf N\fff,\pff \Gamma \trf)
&
\Gamma^{\trf \bullet}\dff(\trf X \trf)\pff,
\arrow[l, "\dis \off j_{\trf \Gamma}"']
\end{tikzcd}
\]

\vspace{-9pt}
where\dss $T^{\dff \bullet}\dff(\trf N\fff,\pff \Gamma \trf)$\dss is\dss the\sss
total\sss complex of\trs the double complex\dss
$C^{\dff \bullet}\dff(\trf N\fff,\pff \Gamma^{\trf \bullet} \trf)$\nnsp.\oss

By\sss applying\dss the same construction to open subsets of\trs $X$\dss in\dss the role
of\trs $X$\nnsp,\oss we will\dss get\dss a double complex\sss of\dss sheaves\dss
$C^{\dff \bullet}\dff(\trf N\fff,\pff \bm{\Gamma}^{\trf \bullet} \trf)$\nnsp.\oss
In\dss more details,\oss for every\sss open subset\dss $U\qff \subset\qff X$\dss 
let\vspace{4.5pt}
\[
\quad
C^{\dff p}\dff(\trf N\fff,\pff \bm{\Gamma}^{\dff q} \trf)\dff (\trf U\trf)
\off\off =\off\off
\prod\nolimits_{\qff \sigma\qff \in\qff N_{\dff p}}\qff 
\Gamma^{\dff q}\fff\left(\trf \num{\sigma}\qff \cap\qff U \trf\right)
\pff.
\]

\vspace{-9pt}
If\trs $V\qff \subset\qff U\qff \subset\qff X$\nnsp,\oss
then\dss the restriction\dss maps\dss 
$\Gamma^{\dff q}\fff\left(\trf \num{\sigma}\qff \cap\qff U \trf\right)
\qff \ttoo\qff
\Gamma^{\dff q}\fff\left(\trf \num{\sigma}\qff \cap\qff V \trf\right)$\dss
lead\dss to a map\dss\vspace{3pt}
\[
\quad
C^{\dff p}\dff(\trf N\fff,\pff \Gamma^{\dff q} \trf)\dff (\trf U\trf)
\qff \ttoo\qff
C^{\dff p}\dff(\trf N\fff,\pff \Gamma^{\dff q} \trf)\dff (\trf V\trf)
\pff.
\]

\vspace{-9pt}
These maps\sss turn\dss $C^{\dff p}\dff(\trf N\fff,\pff \bm{\Gamma}^{\dff q} \trf)$\dss
into a presheaf\halfff.\oss
Since\dss $\Gamma^{\dff q}$\dss is\dss a sheaf\halfff,\oss
$U
\off \longmapsto\off 
\Gamma^{\dff q}\fff\left(\trf \num{\sigma}\qff \cap\qff U \trf\right)$\dss
is\dss also a sheaf\halfff,\oss and\dss hence\dss
$C^{\dff p}\dff(\trf N\fff,\pff \bm{\Gamma}^{\dff q} \trf)$\nnsp,\oss
being a\sss product\sss of\dss sheaves,\oss is\dss also a sheaf\halfff.\qff\off
The maps\vspace{1.5pt}
\[
\quad
\delta_{\fff p}\dff \colon\dff
C^{\dff p}\dff(\dff N\fff,\pff \bm{\Gamma}^{\dff \bullet} \dff)
\qff \ttoo\qff
C^{\dff p\dff +\dff 1}\dff(\dff N\fff,\pff \bm{\Gamma}^{\dff \bullet} \dff)
\pff,
\]

\vspace{-10.5pt}
are defined as before.\oss

\mypar{Lemma.}{sheaf}
\emph{If\pss $\mathcal{U}$ is\dss closed\dss and\trs locally\dss finite\sss
and\qss $X$\dss is\dss paracompact\halfff,\oss
then\dss the sequence}\vspace{-1.5pt}
\begin{equation}
\label{sections-resolution}
\quad
\begin{tikzcd}[column sep=large, row sep=huge]\dis
0\arrow[r]
&
\Gamma^{\trf q}\dff(\trf X \trf) 
\arrow[r, "\dis \delta_{\trf -\dff 1}\qff"]
& 
C^{\trf 0}\dff(\trf N\fff,\pff \Gamma^{\trf q} \trf)\dff(\trf X \trf) 
\arrow[r, "\dis \delta_{\trf 0}\qff"]
& 
C^{\trf 1}\dff(\trf N\fff,\pff \Gamma^{\trf q} \trf)\dff(\trf X \trf) 
\arrow[r, "\dis \delta_{\trf 1}"] 
&   
\off \ldots
\end{tikzcd}
\end{equation}

\vspace{-9pt}
\emph{is\dss exact\qss for every\qss $q\qff \geq\qff 0$\dss
and\dss the canonical\dss morphism\qss
$\tau_{\dff \Gamma}\dff \colon\dff
\Gamma^{\trf \bullet}\dff(\trf X\trf)
\qff \ttoo\qff 
T^{\dff \bullet}\dff(\trf N\fff,\pff \Gamma \trf)$\dss
induces an\dss isomorphism\dss of\qss cohomology\dss groups.\oss}

\proof
The sheaf\dss analogue of\trs the sequence\qss (\ref{sections-resolution})\qss is\dss
the sequence\vspace{-1.5pt}
\begin{equation}
\label{sheaf-resolution}
\quad
\begin{tikzcd}[column sep=large, row sep=huge]\dis
0\arrow[r]
&
\bm{\Gamma}^{\trf q} \arrow[r, "\dis \delta_{\trf -\dff 1}\qff"]
& 
C^{\trf 0}\dff(\trf N\fff,\pff \bm{\Gamma}^{\trf q} \trf) \arrow[r, "\dis \delta_{\trf 0}\qff"]
& 
C^{\trf 1}\dff(\trf N\fff,\pff \bm{\Gamma}^{\trf q} \trf) \arrow[r, "\dis \delta_{\trf 1}"] 
&   
\off \ldots
\end{tikzcd}
\end{equation}

\vspace{-9pt}
of\dss sheaves.\oss
This\sss sequence\dss is\dss exact\halfff.\oss
In\dss fact\halfff,\oss it\dss is\dss exact\dss for every\sss sheaf\dss on\dss $X$\dss
in\dss the role of\trs $\bm{\Gamma}^{\trf q}$\dnsp.\oss
The proof\trs is\dss similar\dss to\sss the proof\dss of\qss 
Lemma\qss \ref{columns-are-exact},\oss
but\sss depends on\sss the assumption\dss that\sss $\mathcal{U}$\sss
is\dss closed and\dss locally\dss finite\qss
(it\dss works also for open coverings).\oss
See\sss the classical\dss book\qss
\cite{go}\qss by\qss R.\dss Godement\halfff,\oss
Chapter\qss II,\oss Theorem\qss 5.2.1.\oss
The sequence\qss (\ref{sections-resolution})\qss
results from\qss (\ref{sheaf-resolution})\qss by\dss taking\dss the sections over\dss $X$\nnsp.\oss
It\dss remains\sss to prove\sss that\dss this operation\dss preserves exactness.\oss

The results about\dss sections are usually\dss stated\dss for sections with supports\sss
in\sss a\sss family\dss $\Phi$\nnsp.\oss
Let\dss $\Phi$\dss be\sss the family\sss of\dss all\sss closed subsets of\trs $X$\nnsp.\oss
Then\dss the support\sss of\dss every\sss section\dss belongs\sss to\dss $\Phi$\nnsp.\oss
Since\dss $X$\dss is\dss paracompact\halfff,\oss
the family\sss of\dss all\sss closed subsets\dss is\dss a\qss
\emph{paracompactifying\qss family}\halfff.\oss 
See\qss \cite{go},\oss Section\qss II.3.2\qss
or\qss \cite{bre},\oss Section\qss I.6\qss for\sss the definitions of\dss
paracompactifying\sss families.\oss
This property\dss is\dss crucial\dss for\sss
preserving\sss exactness\sss by\dss taking sections.\oss

The sheaf\qss
$\bm{\Gamma}^{\dff q}$\dss is\dss soft\dss because\dss $X$\dss is\dss paracompact\halfff.\oss
See\qss \cite{go},\oss Chapter\qss II,\oss Example\qss 3.9.1.\oss
Since\sss the family\sss of\dss all\sss closed
subsets\dss is\dss paracompactifying,\oss this implies\dss that\dss
taking\dss the sections over\dss $X$\dss preserves exactness.\oss
See\qss \cite{go},\qss Chapter\qss II,\oss Theorem\qss 3.5.4.\oss 
Hence\sss the exactness of\qss (\ref{sheaf-resolution})\qss
implies exactness of\qss (\ref{sections-resolution}).\oss  
This proves\sss the fist\sss statement\sss of\trs the\sss lemma.\oss
The second one follows from\dss the fist\sss and\trs Theorem\qss \ref{double-complex}\qss
with\dss the rows and columns interchanged.\oss  \eproof

\myuppar{The functors of\dss global\sss sections.}
Let\sss $Z$\sss be a\sss topological\sss space.\oss
By applying\sss the above construction\sss to\sss $Z$\sss in\sss the role of\dss $X$\sss
we get\sss a sheaf\dss of\dss cochain complexes\sss 
$\Gamma^{\trf \bullet}
\off =\off
\Gamma^{\trf \bullet}_{\dff Z}$\sss on\sss $Z$\nnsp.\oss
Let\vspace{2.5pt}
\[
\quad
\gamma^{\trf \bullet}\dff(\trf Z\trf)
\off =\off\dff
\Gamma^{\trf \bullet}_{\dff Z}\trf(\trf Z\trf)
\]

\vspace{-9.5pt}
be\sss the cochain complex of\dss global\sss sections of\trs 
$\Gamma^{\trf \bullet}_{\dff Z}$\sss on\sss $Z$\nnsp.\oss
Clearly,\pss
$Z\off \longmapsto\off \gamma^{\trf \bullet}\dff(\trf Z\trf)$\sss
is\dss a\sss functor\sss from\sss topological\sss spaces\sss to cochain complexes.\oss
For\dss $Y\qff \subset\qff X$\dss the complex\dss $\gamma^{\dff \bullet}\dff(\trf Y\trf)$\sss
in\sss general\sss differs\sss from\sss $\Gamma^{\dff \bullet}\dff(\trf Y\trf)$\sss 
because\sss $\gamma^{\dff \bullet}\dff(\trf Y\trf)$\sss
is\dss defined\dss inside of\dss $Y$\nnsp,\oss
while\sss $\Gamma^{\dff \bullet}\dff(\trf Y\trf)$\sss
is\dss defined\sss in\sss terms of\dss singular\sss cochains on open subsets of\dss $X$\sss
containing\dss $Y$\dnsp.\oss
But\dss if\dss $U\qff \subset\qff X$\sss is\dss an open subset,\oss
then\sss 
$\gamma^{\dff \bullet}\dff(\trf U\trf)
\off =\off
\Gamma^{\dff \bullet}\dff(\trf U\trf)$\nnsp.\oss
In\dss particular\halfff,\pss
$\gamma^{\dff \bullet}\dff(\trf X\trf)
\off =\off
\Gamma^{\dff \bullet}\dff(\trf X\trf)$\nnsp.\oss
In\sss general,\oss the restriction of\dss singular cochains on open
sets\sss $U\qff \subset\qff X$\sss to intersections\sss $U\dff \cap\dff Y$\sss
defines a map\dss
$\Gamma^{\dff \bullet}\dff(\trf Y\trf)
\qff \ttoo\qff
\gamma^{\dff \bullet}\dff(\trf Y\trf)$\nnsp.\oss
By\sss a classical\dss theorem of\trs the sheaf\trs theory,\oss
if\trs $Z$\sss is\dss paracompact,\oss
then\dss the natural\dss homomorphism\dss
$C^{\dff \bullet}\dff(\trf Z\trf)
\qff \ttoo\qff
\gamma^{\dff \bullet}\dff(\trf Z\trf)$\dss
induces an\sss isomorphism\vspace{3pt}
\[
\quad
k\dff \colon\dff
H^{\dff n}\dff(\trf Z\trf)
\qff \ttoo\qff
H^{\dff n}\dff(\trf \gamma^{\dff \bullet}\dff(\trf Z\trf)\trf)
\pff.
\]

\vspace{-9pt}
This\dss is\dss more sophisticated\dss form of\qss
Eilenberg's\qss Theorem\qss \ref{eilenberg}.\oss
See\qss \cite{bre},\oss Section\qss I.7.\oss

\myuppar{Comparing\dss $\Gamma^{\dff \bullet}$ and\sss $\gamma^{\dff \bullet}$\dnsp.}
Suppose\sss that\sss $X$\sss is\dss paracompact.\oss
Let\sss $F$\sss be a closed subset\sss of\dss $X$\nnsp.\oss
Then\sss $F$\sss is\dss paracompact\dss
by\trs Lemma\qss \ref{subspace-paracompact}\qss
and\dss hence\sss the cohomology\sss of\dss 
$\gamma^{\dff \bullet}\dff(\trf Z\trf)$\sss are\qss 
(canonically\dss isomorphic\sss to)\qss 
the singular cohomology\sss of\dss $F$\dnsp.\oss
In\sss general,\oss there are no reasons\sss to expect\dss that\dss the map\sss
$\Gamma^{\dff \bullet}\dff(\trf F\trf)
\qff \ttoo\qff
\gamma^{\dff \bullet}\dff(\trf F\trf)$\dss
induces isomorphisms in cohomology.\oss

\mypar{Lemma.}{sheaf-subsets}
\emph{If\pss $X$ is\dss paracompact\sss and\dss 
both\dss $X$ and\dss $F$\sss are homologically\dss locally\sss connected,\oss
then\dss the morphism\dss
$\Gamma^{\dff \bullet}\dff(\trf F\trf)
\qff \ttoo\qff
\gamma^{\dff \bullet}\dff(\trf F\trf)$\dss
induces isomorphisms in cohomology.\oss}

\proof
See\qss \cite{bre},\oss Section\qss III.1,\pss
the big\sss diagram\sss on\dss p.\qss 183.\pss
In\sss more details,\oss 
let\sss $\Phi$\sss be\sss the family\sss of\dss all\sss closed
subsets of\dss $X$\nnsp.\pss
Then\dss for a sheaf\sss $\mathcal{B}$\sss on\sss $X$\trs
Bredon's\dss cohomology\sss 
$_S\dff H^{\dff \bullet}_{\dff \Phi}\dff(\qff X\dff;\qff \mathcal{B}\trf)$\sss
is\dss the cohomology\sss of\trs the complex of\dss global\sss sections of\trs
the sheaf\dss $\Gamma^{\dff \bullet}\dff \otimes\dff \mathcal{B}$\nnsp.\oss
Our\sss map\sss\vspace{3pt}\vspace{-0.4375pt}
\[
\quad
H^{\dff n}\dff(\trf \Gamma^{\dff \bullet}\dff(\trf F\trf)\trf)
\qff \ttoo\qff
H^{\dff n}\dff(\trf \gamma^{\dff \bullet}\dff(\trf F\trf)\trf)
\]

\vspace{-9pt}\vspace{-0.4375pt}
corresponds to\trs Bredon's\trs map\vspace{3pt}\vspace{-0.4375pt}
\[
\quad
f^{\dff *}\dff \colon\qff
_S\dff H^{\dff n}_{\dff \Phi}\dff(\qff X\dff;\qff \mathcal{A}_{\qff F}\trf)
\qff \ttoo\qff
_S\dff H^{\dff n}_{\dff \Phi\dff|\dff F}\dff(\qff F\dff;\qff \mathcal{A}\dff |\dff F\qff)
\]

\vspace{-9pt}\vspace{-0.4375pt}
with\sss $\mathcal{A}$ being\dss the sheaf\dss associated\sss
with\sss the constant\sss presheaf\dss $\rrr$\sss on\sss $X$\nnsp.\oss
Here\sss $\Phi\dff|\dff F$\sss is\dss the family\sss of\dss all\sss
closed subsets of\dss $F$\sss and\sss $\mathcal{A}\dff |\dff F$\sss
is\dss the restriction of\dss $\mathcal{A}$\sss to\sss $F$\dnsp.\oss
Hence\sss $\mathcal{A}\dff |\dff F$\sss is\dss the sheaf\dss associated\sss
with\sss the constant\sss presheaf\dss $\rrr$\sss on\sss $F$\dnsp.\oss
It\dss follows\dss that\vspace{3pt}
\[
\quad
_S\dff H^{\dff n}_{\dff \Phi\dff|\dff F}\dff(\qff F\dff;\qff \mathcal{A}\dff |\dff F\qff)
\off =\off
H^{\dff n}\dff(\trf \gamma^{\dff \bullet}\dff(\trf F\trf)\trf)
\pff.
\]

\vspace{-9pt}
The sheaf\dss $\mathcal{A}_{\qff F}$\sss is\dss obtained\dss by\sss
extension of\trs the sheaf\dss $\mathcal{A}\dff |\dff F$\sss
from\sss $F$\sss to\sss $X$\sss by\sss zero
and\dss hence sections of\trs
$\Gamma^{\dff \bullet}\dff \otimes\dff \mathcal{A}_{\qff F}$\sss
over\sss $X$\sss
can\sss be identified\sss with sections of\dss $\Gamma^{\dff \bullet}$\sss
over\sss $F$\dnsp.\oss 
It\dss follows\dss that\vspace{3pt}
\[
\quad
_S\dff H^{\dff n}_{\dff \Phi}\dff(\qff X\dff;\qff \mathcal{A}_{\qff F}\trf)
\off =\off
H^{\dff n}\dff(\trf \Gamma^{\dff \bullet}\dff(\trf F\trf)\trf)
\pff.
\]

\vspace{-9pt}
According\dss to\qss \cite{bre},\oss
the map\sss $f^{\dff *}$\sss is\dss an\sss isomorphism.\oss
The\sss lemma\sss follows.\oss  \eproof

\mypar{Lemma.}{comparing-total-complexes}
\emph{Suppose\sss that\dss the covering\dss $\mathcal{U}$\sss is\dss closed\sss
and\dss that\dss the space\dss $X$\sss and\sss elements of\pss
$\mathcal{U}^{\dff \cap}$\dss are\dss locally\dss homologically\sss connected.\oss
Then\dss the morphism\dss
$T^{\dff \bullet}\dff(\trf N\fff,\pff \Gamma \trf)
\qff \ttoo\qff
T^{\dff \bullet}\dff(\trf N\fff,\pff \gamma \trf)$\sss
of\qss total\sss complexes\sss induced\dss by\dss the morphisms\qss
$\Gamma^{\dff \bullet}\dff(\trf F\trf)
\qff \ttoo\qff
\gamma^{\dff \bullet}\dff(\trf F\trf)$\nnsp,\oss
where\qss $F\qff \in\qff \mathcal{U}^{\dff \cap}$\dnsp,\oss
induces an\sss isomorphism\sss in cohomology.\oss}

\proof
By\qss Lemma\qss \ref{sheaf-subsets}\qss for every\sss simplex $\sigma$ of\trs $N$\sss
the morphism\dss\vspace{3pt}
\[
\quad
\Gamma^{\dff \bullet}\dff(\trf \num{\sigma}\trf)
\qff \ttoo\qff
\gamma^{\dff \bullet}\dff(\trf \num{\sigma}\trf)
\]

\vspace{-9pt}
induces an\sss isomorphism\sss in\sss cohomology.\oss
It\dss follows\dss that\sss for every\sss $p\qff \geq\qff 0$\dss
the product\vspace{3pt}
\[
\quad
C^{\dff p}\dff(\trf N\fff,\pff \Gamma^{\dff \bullet} \trf)
\off \ttoo\off
C^{\dff p}\dff(\trf N\fff,\pff \gamma^{\dff \bullet} \trf)
\]

\vspace{-9pt}
of\trs these morphisms over\sss $\sigma\qff \in\pff N_{\dff p}$\sss
induces an\sss isomorphism\sss in\sss cohomology.\oss
It\dss remains\sss to apply\sss a well\dss known comparison\dss theorem
about\sss double complexes.\oss
See\trs Theorem\qss \ref{comparison}.\oss  \eproof

\myuppar{The homomorphism\dss
$H^{\dff *}\fff(\trf N \trf)
\qff \ttoo\qff 
H^{\dff *}\fff(\trf X \trf)$\dss
for\dss closed\dss locally\dss finite coverings.}
Suppose\sss that\dss the covering\sss $\mathcal{U}$\sss is\dss closed\sss and\dss locally\dss finite
and\sss $X$\sss is\dss paracompact.\oss
The morphisms\vspace{-1.5pt}
\[
\quad
\begin{tikzcd}[column sep=large, row sep=huge]\dis
C^{\dff \bullet}\dff(\trf N \trf) 
\arrow[r, "\dis \lambda_{\trf \Gamma}"]
& 
T^{\dff \bullet}\dff(\trf N\fff,\pff \Gamma \dff)
&
\Gamma^{\trf \bullet}\dff(\trf X \trf)\pff
\arrow[l, "\dis \off \tau_{\dff \Gamma}"']
&
C^{\dff \bullet}\dff(\trf X \trf)\pff
\arrow[l]
\end{tikzcd}
\]

\vspace{-10.5pt}
lead\dss to homomorphisms\vspace{-1.5pt}
\[
\quad
\begin{tikzcd}[column sep=large, row sep=huge]\dis
H^{\dff *}\dff(\trf N \trf) \arrow[r, "\dis \lambda_{\trf \Gamma\dff *}"]
& 
H^{\dff *}\dff(\trf N\fff,\pff \Gamma \dff)
&
H^{\dff *}\dff(\trf \Gamma^{\trf \bullet}\dff(\trf X \trf) \trf)\pff
\arrow[l, "\dis \off \tau_{\dff \Gamma\dff *}"']
&
H^{\dff *}\dff(\trf X \trf)\pff
\arrow[l, "\dis k_{\vphantom{\Gamma}}"']
\end{tikzcd}
\]

\vspace{-10.5pt}
of\dss cohomology\dss groups,\oss
where we denoted\dss by\dss
$H^{\dff *}\dff(\trf N\fff,\pff \Gamma \dff)$\dss
the cohomology\sss of\trs
$T^{\dff \bullet}\dff(\trf N\fff,\pff \Gamma \dff)$\nnsp.\oss
By\trs Lemma\qss \ref{sheaf}\qss the homomorphism\dss 
$j_{\trf \Gamma\dff *}$\dss is\dss an\dss isomorphism.\oss
Since\sss $X$\sss is\dss paracompact\sss and\dss locally\dss
homologically\sss connected,\pss
$k$\dss is\dss also an\dss isomorphism.\oss
Therefore,\oss the homomorphism\vspace{4.5pt}
\[
\quad
k^{\qff -\dff 1}
\dff \circ\qff 
\tau_{\dff \Gamma\dff *}^{\qff -\dff 1}
\pff \circ\pff 
\lambda_{\trf \Gamma\dff *}
\qff \colon\qff
H^{\dff *}\fff(\trf N \trf)
\qff \ttoo\qff 
H^{\dff *}\fff(\trf X \trf)
\]

\vspace{-7.5pt}
is\dss well\sss defined,\oss
and\dss we\sss take\sss it\dss 
as\sss the\qss \emph{canonical\dss homomorphism}\pss
$l_{\trf \mathcal{U}}\dff \colon\dff
H^{\dff *}\fff(\trf N \trf)
\qff \ttoo\qff 
H^{\dff *}\fff(\trf X \trf)$\nnsp.\oss

\myuppar{The $A^{\bullet}$\dnsp\dnsp-cohomology\sss and\dss the singular\sss cohomology.}
Suppose\sss that\sss $A^{\bullet}$\sss is\dss a\sss functor\sss of\trs generalized\sss
cochains as\sss in\dss Section\qss \ref{cohomological-leray}.\oss
Then\sss $A^{\bullet}$\dnsp\dnsp-cohomology\dss groups\sss
$\widetilde{H}^{\dff *}\dff(\trf X\trf)
\off =\off\dff
H_{\dff A}^{\dff *}\dff(\trf X\trf)$\sss
are defined.\oss
Suppose\sss that\dss the functor\sss $A^{\bullet}$\sss is\dss 
equipped\dss with a natural\dss transformation\dss 
$A^{\bullet}\qff \ttoo\qff C^{\dff \bullet}$\dnsp,\oss
leading\sss to a\sss natural\dss homomorphism\dss
$\widetilde{H}^{\dff *}\dff(\trf X\trf)
\qff \ttoo\qff
H^{\dff *}\dff(\trf X\trf)$\nnsp.\oss

\mypar{Theorem.}{A-acyclic-closed}
\emph{Suppose\sss that\dss the covering\dss $\mathcal{U}$\sss is\dss closed\sss
and\dss locally\dss finite,\oss
and\dss that\dss the space\dss $X$\sss and\sss elements of\pss
$\mathcal{U}^{\dff \cap}$\dss are\dss locally\dss homologically\sss connected.\oss
If\pss $\mathcal{U}$\dss is\dss $A^{\bullet}$\dnsp-acyclic,\oss
then\dss the\sss natural\dss homomorphism\qss
$\widetilde{H}^{\dff *}\dff(\trf X\trf)
\qff \ttoo\qff
H^{\fff *}\fff(\trf X \trf)$\dss
can\sss be\dss factored\dss through\qss
$l_{\trf \mathcal{U}}\dff \colon\dff
H^{\dff *}\fff(\trf N \trf)
\qff \ttoo\qff 
H^{\dff *}\fff(\trf X \trf)$\dnsp.\oss}

\proof
Let\dss us consider\sss the following\sss diagram of\trs total\sss complexes.\oss\vspace{6pt}
\[
\quad
\begin{tikzcd}[column sep=boom, row sep=boom]\dis
C^{\trf \bullet}\dff(\trf N \trf) 
\arrow[r]
\arrow[d, "\dis \qff ="]
& 
T^{\dff \bullet}\dff(\trf N\fff,\pff A \trf)
\arrow[d]
&
A^{\bullet}\dff(\trf X \trf)
\arrow[l]
\arrow[d]
\\
C^{\trf \bullet}\dff(\trf N \trf) 
\arrow[d, "\dis \qff ="]
\arrow[u]
\arrow[r]
& 
T^{\dff \bullet}\dff(\trf N\fff,\pff C \trf)
\arrow[d]
&
C^{\trf \bullet}\dff(\trf X \trf)
\arrow[d]
\arrow[l]
\\
C^{\trf \bullet}\dff(\trf N \trf)
\arrow[u]
\arrow[d]
\arrow[r]
& 
T^{\dff \bullet}\dff(\trf N\fff,\pff \gamma \trf)
&
\gamma^{\trf \bullet}\dff(\trf X \trf)\pff
\arrow[d]
\arrow[l]
\\
C^{\trf \bullet}\dff(\trf N \trf) 
\arrow[u, "\dis \qff ="'] 
\arrow[r]
& 
T^{\dff \bullet}\dff(\trf N\fff,\pff \Gamma \trf)
\arrow[u]
&
\Gamma^{\trf \bullet}\dff(\trf X \trf)\pff.
\arrow[l]
\arrow[u, "\dis \qff ="']
\end{tikzcd}
\]

\vspace{-3pt}
This diagram\dss leads\sss to\sss the diagram of\dss cohomology\dss
groups\vspace{3pt}
\[
\quad
\begin{tikzcd}[column sep=boomm, row sep=boom]\dis
H^{\dff *}\dff(\trf N \trf) 
\arrow[r, "\dis \lambda_{\trf A\dff *}", blue, line width=0.8pt]
\arrow[d, "\dis \qff ="]
& 
H^{\dff *}\dff(\trf N\fff,\pff A \trf)
\arrow[d]
&
H_{\trf A}^{\dff *}\dff(\trf X \trf)
\arrow[l]
\arrow[d]
\\
H^{\dff *}\dff(\trf N \trf)
\arrow[u] 
\arrow[r]
\arrow[d, "\dis \qff ="]
& 
H^{\dff *}\dff(\trf N\fff,\pff C \trf)
\arrow[d]
&
H^{\dff *}\dff(\trf X \trf)
\arrow[d, "\dis k\off \approx", blue, line width=0.8pt]
\arrow[l]
\\
H^{\dff *}\dff(\trf N \trf)
\arrow[u]
\arrow[d, "\dis \qff ="]
\arrow[r, "\dis \lambda_{\trf \gamma\dff *}"]
& 
H^{\dff *}\dff(\trf N\fff,\pff \gamma \trf)
&
H^{\dff *}\dff(\trf \gamma^{\trf \bullet}\dff(\trf X \trf) \trf)\pff.
\arrow[l, "\dis \off \tau_{\dff \gamma\dff *}"', blue, line width=0.8pt]
\arrow[d]
\\
H^{\dff *}\dff(\trf N \trf) 
\arrow[u, "\dis \qff ="']
\arrow[r, "\dis \lambda_{\trf \Gamma\dff *}"]
& 
H^{\dff *}\dff(\trf N\fff,\pff \Gamma \trf)
\arrow[u, blue, line width=0.8pt]
&
H^{\dff *}\dff(\trf \Gamma^{\trf \bullet}\dff(\trf X \trf) \trf)\pff.
\arrow[l, "\dis \off \tau_{\trf \Gamma\dff *}"', blue, line width=0.8pt]
\arrow[u, "\dis \qff ="']
\end{tikzcd}
\]

\vspace{-3pt}
Since\dss $\mathcal{U}$\dss is\dss $A^{\bullet}$\dnsp\dnsp-acyclic,\oss
the homomorphism\dss
$\lambda_{\trf A\dff *}$\dss
is\dss an\dss isomorphism\dss by\qss Lemma\qss \ref{acyclic-coverings}.\oss
By\qss Lemma\qss \ref{sheaf}\qss
the homomorphism\dss $\tau_{\dff \Gamma\dff *}$\sss is\dss
also an\sss isomorphism.\oss
As we saw,\pss $k$\sss is\dss an\sss isomorphism\dss too.\oss
Finally,\oss by\qss Lemma\qss \ref{comparing-total-complexes}\qss
the homomorphism\sss
$H^{\dff *}\dff(\trf N\fff,\pff \Gamma \trf)
\qff \ttoo\qff
H^{\dff *}\dff(\trf N\fff,\pff \gamma \trf)$\sss
is\dss an\sss isomorphism.\oss
The commutativity\sss of\trs the\sss lower\sss right\sss square\sss
implies\sss that\dss $\tau_{\dff \gamma\dff *}$\sss is\dss
also an\sss isomorphism.\oss
By\dss inverting\sss $\tau_{\dff \gamma\dff *}$\sss
and\sss $k$\sss we\sss get\dss the square\vspace{0pt}
\[
\quad
\begin{tikzcd}[column sep=boomm, row sep=boom]\dis
H^{\dff *}\dff(\trf N\fff,\pff C \trf)
\arrow[d]
&
H^{\dff *}\dff(\trf X \trf)
\arrow[l]
\\
H^{\dff *}\dff(\trf N\fff,\pff \gamma \trf)
\arrow[r]
&
H^{\dff *}\dff(\trf \gamma^{\trf \bullet}\dff(\trf X \trf) \trf)\pff,
\arrow[u]
\end{tikzcd}
\]

\vspace{-10.5pt}
commutative in\sss the sense\sss that\dss the composition of\dss its four arrows
starting at\sss $H^{\dff *}\dff(\trf X \trf)$\sss
is\dss equal\dss to\sss the identity.\oss
By\dss inverting also\dss
$\lambda_{\trf A\dff *}$\dss
and\dss $\tau_{\dff \Gamma\dff *}$\sss
we get\dss the following\sss
commutative diagram.\oss\vspace{1.5pt}
\[
\quad
\begin{tikzcd}[column sep=boomm, row sep=boom]\dis
H^{\dff *}\dff(\trf N \trf) 
\arrow[d, "\dis \qff =", red, line width=0.8pt]
& 
H^{\dff *}\dff(\trf N\fff,\pff A \trf)
\arrow[l, red, line width=0.8pt]
\arrow[d]
&
H_{\trf A}^{\dff *}\dff(\trf X \trf)
\arrow[l, red, line width=0.8pt]
\arrow[d]
\\
H^{\dff *}\dff(\trf N \trf)
\arrow[r]
\arrow[d, "\dis \qff =", red, line width=0.8pt]
& 
H^{\dff *}\dff(\trf N\fff,\pff C \trf)
\arrow[d]
&
H^{\dff *}\dff(\trf X \trf)
\arrow[l]
\\
H^{\dff *}\dff(\trf N \trf)
\arrow[d, red, line width=0.8pt, "\dis \qff ="]
\arrow[r]
& 
H^{\dff *}\dff(\trf N\fff,\pff \gamma \trf)
\arrow[r]
\arrow[d]
&
H^{\dff *}\dff(\trf \gamma^{\trf \bullet}\dff(\trf X \trf) \trf)\pff.
\arrow[u, red, line width=0.8pt]
\\
H^{\dff *}\dff(\trf N \trf) 
\arrow[r, red, line width=0.8pt]
& 
H^{\dff *}\dff(\trf N\fff,\pff \Gamma \trf)
\arrow[r, red, line width=0.8pt]
&
H^{\dff *}\dff(\trf \Gamma^{\trf \bullet}\dff(\trf X \trf) \trf)\pff.
\arrow[u, red, line width=0.8pt, "\dis \qff ="']
\end{tikzcd}
\]

\vspace{-6pt}\vspace{-0.75pt}
The commutativity\sss of\trs this diagram\sss implies\sss that\dss
$\widetilde{H}^{\dff *}\dff(\trf X\trf)
\off =\off\dff
H_{\trf A}^{\dff *}\dff(\trf X \trf)
\qff \ttoo\qff
H^{\dff *}\dff(\trf X \trf)$\sss
is\dss equal\dss to\sss the composition of\trs the red arrows.\oss
But\dss the composition of\trs the\sss last\dss four\sss red arrows\dss
is\dss nothing else\sss but\dss the canonical\dss homomorphism\dss
$l_{\trf \mathcal{U}}\dff \colon\dff
H^{\dff *}\fff(\trf N \trf)
\qff \ttoo\qff 
H^{\dff *}\fff(\trf X \trf)$\nnsp.\oss
The\sss theorem\sss follows.\oss  \eproof

\myuppar{Extensions of\dss coverings.}
Suppose\sss that\dss the covering\dss $\mathcal{U}$\sss is\dss closed\sss
and\dss locally\dss finite,\oss
and\dss that\dss the space\dss $X$\sss and\sss elements of\pss
$\mathcal{U}^{\dff \cap}$\dss are\dss locally\dss homologically\sss connected.\oss
Then\sss the covering\sss $\mathcal{U}\fff'$\sss of\dss a space\sss
$X\fff'\qff \supset\qff X$\sss constructed\sss in\trs
Theorem\qss \ref{pi-one-extensions}\qss is\dss also closed and\dss locally\dss finite.\oss
Recall\dss that\sss $X\fff'$\sss is\dss obtained\dss by\sss attaching\dss to\sss $X$\sss
some discs along\dss their\sss boundaries,\oss
and\sss elements\sss of\dss $\mathcal{U}$\sss
are obtained\dss by\sss from an element\sss of\dss $\mathcal{U}$\sss by\sss
attaching some of\trs these discs.\oss
Moreover\halfff,\oss elements of\dss $\mathcal{U}\fff'^{\trf \cap}$\sss
can\sss be obtained\dss from elements of\dss $\mathcal{U}^{\dff \cap}$\sss
in\sss the same way.\oss
Clearly,\oss if\dss a space\sss $Z$\sss is\dss
locally\dss homologically\sss connected,\oss
then\dss the result\sss of\dss attaching of\dss a collection of\dss discs\sss
to\sss $Z$\sss has\sss the same property.\oss
Therefore\sss $X\fff'$\sss and\sss elements of\pss
$\mathcal{U}\fff'^{\trf \cap}$\dss are\dss locally\dss homologically\sss connected.\oss

\mypar{Theorem.}{acyclic-closed-homology}
\emph{Suppose\sss that\dss the covering\dss $\mathcal{U}$\sss is\dss closed\sss
and\dss locally\dss finite,\oss
and\dss that\dss the space\dss $X$\sss and\sss elements of\pss
$\mathcal{U}^{\dff \cap}$\dss are\dss locally\dss homologically\sss connected.\oss
If\pss $\mathcal{U}$\dss is\dss weakly\dss boundedly\sss acyclic,\oss
then\qss 
$\widehat{H}^{\dff *}\dff(\trf X\trf)
\qff \ttoo\qff
H^{\fff *}\fff(\trf X \trf)$\dss
can\sss be\dss factored\dss through\qss
$l_{\trf \mathcal{U}}\dff \colon\dff
H^{\dff *}\fff(\trf N \trf)
\qff \ttoo\qff 
H^{\dff *}\fff(\trf X \trf)$\dnsp.\oss}

\proof
The main\sss part\sss of\trs the work\dss is\dss already\sss done.\oss
The rest\trs is\dss completely\sss similar\sss to\sss the proof\dss of\qss
Theorem\qss \ref{covering-theorem-open}.\oss
One only\sss needs\sss to refer\sss to\trs Theorem\qss \ref{A-acyclic-closed}\qss
instead of\qss Theorem\qss \ref{A-acyclic-open}\qss
and supplement\trs Corollary\qss \ref{amenable-extensions}\qss by\dss the remarks
preceding\dss the\sss theorem.\oss  \eproof

\newpage
\myappend{Double\qss complexes}{double-complexes}

\vspace{6pt}
\myuppar{Double complexes.}
A double complex\sss $K^{\dff \bullet\fff,\dff \bullet}$\sss is\dss 
a diagram of\trs the form\vspace{6.5pt}
\[
\qquad
\begin{tikzcd}[column sep=large, row sep=huge]\dis
K^{\dff 0\fff,\dff 0} \arrow[r]
\arrow[d]
& 
K^{\dff 0\fff,\dff 1} \arrow[r]
\arrow[d]
& 
K^{\dff 0\fff,\dff 2} \arrow[r] 
\arrow[d]
&   
K^{\dff 0\fff,\dff 3} \arrow[r] 
\arrow[d]
&
\off\off \ldots 
\\
K^{\dff 1\fff,\dff 0}\arrow[r]
\arrow[d]
& 
K^{\dff 1\fff,\dff 1} \arrow[r]
\arrow[d]
& 
K^{\dff 1\fff,\dff 2} \arrow[r]
\arrow[d]
&   
K^{\dff 1\fff,\dff 3} \arrow[r]
\arrow[d]
&
\off\off \ldots 
\\
K^{\dff 2\fff,\dff 0} \arrow[r]
\arrow[d]
& 
K^{\dff 2\fff,\dff 1} \arrow[r]
\arrow[d]
& 
K^{\dff 2\fff,\dff 2} \arrow[r]
\arrow[d]
&   
K^{\dff 2\fff,\dff 3} \arrow[r]
\arrow[d]
&
\off\off \ldots 
\\
K^{\dff 3\fff,\dff 0} \arrow[r]
\arrow[d]
& 
K^{\dff 3\fff,\dff 1} \arrow[r]
\arrow[d]
& 
K^{\dff 3\fff,\dff 2} \arrow[r]
\arrow[d]
&   
K^{\dff 3\fff,\dff 3} \arrow[r]
\arrow[d]
&
\off\off \ldots
\\
\ldots \vphantom{K^{\dff 3\fff,\dff 3}} &
\ldots \vphantom{K^{\dff 3\fff,\dff 3}} &
\ldots \vphantom{K^{\dff 3\fff,\dff 3}} &
\ldots \vphantom{K^{\dff 3\fff,\dff 3}} &
\off\off\off\off .
\end{tikzcd}
\]

\vspace{-6.25pt}
The horizontal\sss arrows\qss
$K^{\dff p\fff,\dff q}\qff \ttoo\qff K^{\dff p\fff,\dff q\dff +\dff 1}$\qss
and\dss the vertical\sss arrows\qss
$K^{\dff p\fff,\dff q}\qff \ttoo\qff K^{\dff p\dff +\dff 1\fff,\dff q}$\qss
are denoted\dss by\dss $d$\dss and\dss $\delta$\dss respectively\halfff,\oss
and are called\dss the\qss \emph{differentials}\qss of\trs the double complex\dss 
$K^{\dff \bullet\fff,\dff \bullet}$\dnsp.\oss
Each\dss row and each column of\trs this diagram\dss is\dss assumed\dss to
be a complex.\oss
Equivalently\halfff,\oss it\dss is\dss assumed\dss that\trs
$d\dff \circ\dff d\off =\off 0$\dss
and\qss
$\delta\dff \circ\dff \delta\off =\off 0$\nnsp.\oss

The diagram\dss is\dss assumed\dss to be commutative in one of\trs the following\dss
two senses.\oss
First\halfff,\oss one may\dss require\sss that\sss each square of\trs the diagram\dss
is\dss commutative,\oss
i.e.\qss to require\sss that\dss the\sss two differentials commute,\pss
$d\dff \circ\dff \delta\off =\off \delta\dff \circ\dff d$\nnsp.\oss
Then\dss the diagram\dss is\dss commutative in\dss the usual\sss sense.\oss
The double complexes used\dss in\qss Section\qss \ref{cohomological-leray}\qss 
are commutative in\dss this sense.\oss
Alternatively\halfff,\oss one can\dss require\sss that\dss the differentials
anti-commute,\oss i.e.\qss that\trs
$d\dff \circ\dff \delta\pff +\pff \delta\dff \circ\dff d\off =\off 0$\nnsp.\oss
The advantages of\trs this condition\dss will\dss be clear\sss in\sss a\sss moment\halfff.\oss
In order\dss to pass from one version\dss to\sss the other\dss it\dss is\dss sufficient\dss
to replace\sss differentials\dss
$\delta\off =\off \delta^{\dff p\fff,\dff q}\dff \colon\dff
K^{\dff p\fff,\dff q}\qff \ttoo\qff K^{\dff p\dff +\dff 1\fff,\dff q}$\dss
by\dss differentials\dss $(\qff -\qff 1\trf)^{\dff p}\qff \delta^{\dff p\fff,\dff q}$\dnsp.\oss

\myuppar{The\sss total\sss complex of\dss a double complex.}
Let\vspace{3pt}
\[
\quad
T^{\fff n}
\off =\off
\bigoplus_{p\qff +\qff q\qff =\qff n}\qff K^{\dff p\fff,\dff q}
\]

\vspace{-12pt}
and\dss let\dss 
$\partial\dff \colon\dff
T^{\fff n}\qff \ttoo\qff T^{\fff n\dff +\dff 1}$\dss
be\sss the map equal\dss to\dss $d\qff +\qff (\qff -\qff 1\trf)^{\dff p}\qff \delta$\dss
on\dss $K^{\dff p\fff,\dff q}$\dss if\trs the differentials commute,\oss
and\dss to\dss $d\qff +\qff \delta$\qss if\trs the differentials anti-commute.\oss
Since\dss
$d\dff \circ\dff d\off =\off 0$\dss
and\qss
$\delta\dff \circ\dff \delta\off =\off 0$\nnsp,\oss
a\sss trivial\sss computation shows\sss that\dss
in\dss both cases\qss
$\partial\dff \circ\dff \partial\off =\off 0$\nnsp.\oss
Therefore\dss $T^{\dff \bullet}$\dss together\sss with\dss $\partial$\dss
is\dss a complex.\oss
It\dss is\dss called\dss the\qss \emph{total\sss complex}\pss of\trs the double complex\dss
$K^{\dff \bullet\fff,\dff \bullet}$\dnsp.\oss

Let\dss $L^{\dff p}$\dss be the kernel of the differential\qss
$d\dff \colon\dff
K^{\fff p\fff,\dff 0}\qff \ttoo\qff K^{\fff p\fff,\dff 1}$\dnsp.\oss
Either of\trs the commutativity\sss assumptions implies\sss that\dss
$\delta$\dss maps\dss $L^{\dff p}$\dss to\dss $L^{\dff p\dff +\dff 1}$\dnsp.\oss
Therefore\qss $L^{\dff \bullet}$\qss together with the restriction of\trs the differential\qss $\delta$\qss
to\dss $L^{\dff \bullet}$\dss 
is\dss a subcomplex of the total complex\qss $T^{\dff \bullet}$\dnsp.\oss

\myapar{Theorem.}{double-complex}
\emph{If\qss each complex\qss
$\dis
\left(\qff K^{\dff p\fff,\dff \bullet},\pff d \qff\right)$\qss 
is\dss exact\halfff,\oss
then the homomorphism}\qss\vspace{3pt}\vspace{-0.125pt}
\[
\quad
H^{\fff *}\dff(\trf L^{\dff \bullet} \dff)
\qff \ttoo\qff
H^{\fff *}\dff(\trf T^{\dff \bullet} \dff)
\]

\vspace{-9pt}\vspace{-0.125pt}
\emph{induced\dss by\dss the\dss inclusion\qss
$L^{\dff \bullet}\qff \ttoo\qff T^{\dff \bullet}$\qss
is\dss an\dss isomorphism.\oss}

\proof\qss
This\dss is\dss a special case of\qss Theorem\qss 4.8.1\qss from\dss Chapter\qss I\qss of\qss
Godement's\trs book\qss \cite{go}.\oss
Its\dss standard\dss proof\trs is\dss based on\dss the properties of\dss spectral sequences 
associated\dss with\dss $K^{\dff \bullet\fff,\dff \bullet}$\dnsp.\oss

Here\dss is\dss a direct\dss proof\halfff.\oss
We may assume\sss that\dss the differentials
anti-commute.\oss
Let\dss us prove first\dss the induced\dss homomorphism\dss is\dss surjective.\oss
Let\vspace{3pt}
\[
\quad
z
\off\dff =\off
\bigoplus_{i\qff =\qff 0}^n\qff z^{\dff i}
\off\qff \in\off\qff
\bigoplus_{i\qff =\qff 0}^n\qff K^{\dff i\fff,\dff n\dff -\dff i}
\quad.
\]

\vspace{-15pt}
Then\vspace{-6pt}
\[
\quad
\partial\dff z
\off\qff \in\off\qff
\bigoplus_{k\qff =\qff 0}^n\qff K^{\dff k\fff,\dff n\dff +\dff 1\dff -\dff k}
\quad
\]

\vspace{-12pt}
and\dss the summand of\trs $\partial\dff z$\dss belonging\dss to\qss
$K^{\dff k\fff,\dff n\dff +\dff 1\dff -\dff k}$\qss is\dss equal\dss to\vspace{3pt}
\[
\quad
d\dff z^{\dff k}
\pff +\pff
\delta\dff z^{\dff k\dff -\dff 1}
\quad\
\mbox{if}\quad\
0\off <\off k\off <\off n\qff +\qff 1
\qff,
\]

\vspace{-36pt}
\[
\quad
d\dff z^{\dff 0}
\quad\
\mbox{if}\quad\
k\off =\off 0
\qff,\quad
\mbox{and}\quad\
\delta\dff z^{\dff d}
\quad\
\mbox{if}\quad\
k\off =\off n\qff +\qff 1
\qff.
\]

\vspace{-9pt}
If\trs $z$\dss
is\dss a cocycle,\oss i.e.\qss if\qss $\partial\dff z\off =\off 0$\nnsp,\oss
then\dss $d\dff z^{\dff 0}\off =\off 0$\nnsp.\oss
Since\sss the complex\dss
$\dis
(\qff K^{\dff 0\fff,\dff \bullet},\pff d \qff)$\dss 
is\dss exact\halfff,\oss
this implies\sss that\dss
$z^{\dff 0}\off =\off \partial\dff y^{\dff 0}$\dss
for some\dss
$y^{\dff 0}\qff \in\pff K^{\fff 0\fff,\dff n\dff -\dff 1}$\dnsp.\oss
Let\dss
$y\qff \in\pff T^{\dff n\dff -\dff 1}$\dss
be\sss the element\dss having\dss $y^{\dff 0}$\dss as\sss the only\dss
non-zero summand.\oss
Then\dss $z\qff -\qff \partial\dff y$\dss is\dss a cocycle representing\dss
the same cohomology\sss class as\dss $z$\dss and\dss 
the summand of\trs $z\qff -\qff \partial\dff y$\dss belonging\dss to\qss
$K^{\dff 0\fff,\dff n}$\qss is\dss equal\dss to\dss $0$\nnsp.\oss

Let\dss us\dss replace\dss $z$\dss by\dss $z\qff -\qff \partial\dff y$\dss
while keeping\dss the notation\dss $z$\nnsp.\oss
Now\dss the summand of\trs $\partial\dff z$\dss belonging\dss to\qss
$K^{\dff 1\fff,\dff n}$\qss is\dss equal\dss to\dss $d\dff z^{\dff 1}$\dss
and\dss since\dss $z$\dss is\dss a cocycle,\pss $d\dff z^{\dff 1}\off =\off 0$\nnsp.\oss
Since\dss
$\dis
(\qff K^{\dff 1\fff,\dff \bullet},\pff d \qff)$\dss 
is\dss exact\halfff,\oss
this implies\sss that\dss
$z^{\dff 1}\off =\off \partial\dff y^{\dff 1}$\dss
for some\dss
$y^{\dff 1}\qff \in\pff K^{\dff 1\fff,\dff n\dff -\dff 2}$\dnsp.\oss
Let\dss the new\dss
$y\qff \in\pff T^{\dff n\dff -\dff 2}$\dss
be\sss the element\dss having\dss $y^{\dff 1}$\dss as\sss the only\dss
non-zero summand.\oss
Then\dss $z\qff -\qff \partial\dff y$\dss represents\sss
the same cohomology\sss class as\dss $z$\dss and\dss 
the summands of\trs $z\qff -\qff \partial\dff y$\dss belonging\dss to\qss
$K^{\dff 0\fff,\dff n}$\qss and\qss
$K^{\dff 1\fff,\dff n\dff -\dff 1}$\qss are\sss equal\dss to\dss $0$\nnsp.\oss

By\sss continuing\sss in\dss this way\dss we will\dss eventually\dss reach an
element\dss $z$\dss representing\dss the same cohomology\dss class as\sss
the original\dss $z$\dss and\dss having only\sss one non-zero summand,\oss
namely\halfff,\oss the summand\dss 
$z^{\dff n}\qff \in\pff K^{\dff n\fff,\dff 0}$\dnsp.\oss
Since\sss the new\dss $z\off =\off z^{\dff n}$\dss is\dss a cocycle,\pss
$d\dff z\off =\off \delta\dff z\off =\off 0$\nnsp.\oss
Therefore\sss the new\dss $z$\dss belongs\sss to\dss $L^{\dff n}$\dss
and\dss is\dss a cocycle of\trs the complex\dss $L^{\dff \bullet}$\dnsp.\oss
Since\sss the original\dss $z$\dss was an arbitrary\sss cocycle of\trs the\sss
total\sss complex,\oss
the surjectivity\dss follows.\oss

The proof\dss of\dss injectivity\dss is\dss similar\halfff.\oss
Suppose\sss that\sss $w^{\dff n\dff +\dff 1}\qff \in\pff L^{\dff n\dff +\dff 1}$\dss
is\dss a cocycle of\trs the complex\dss $L^{\dff \bullet}$\dnsp.\oss
Let\dss $w$\dss be\sss the element\sss of\qss $T^{\dff n\dff +\dff 1}$\dss
having\dss $w^{\dff n\dff +\dff 1}$\dss as\sss the only\sss non-zero summand.\oss
Suppose\sss that\dss $w$\dss is\dss a coboundary\dss in\dss $T^{\dff \bullet}$\dnsp,\oss
i.e.\qss $w\off =\off \partial\dff z$\dss for some\dss
$z\qff \in\pff T^{\dff n}$\dnsp.\oss
Since\sss the summands of\trs $\partial\dff z$\dss in\dss 
$K^{\dff i\fff,\dff n\dff -\dff i}$\dss with\dss $i\qff >\qff 0$\dss
are equal\dss to\dss $0$\nnsp,\oss
we can apply\dss to\dss $z$\dss the same process as above.\oss
At\sss each step we replaced\dss $z$\dss by\sss an element\dss of\dss the form\dss
$z\qff -\qff \partial\dff y$\nnsp.\oss
Since\dss $\partial\dff \circ\dff \partial\off =\off 0$\nnsp,\oss
this does not\sss affect\dss the coboundary\dss $\partial\dff z$\nnsp.\oss
At\dss the\sss last\sss step we will\sss get\sss a new element\dss
$z$\dss such\dss that\dss $w\off =\off \partial\dff z$\dss 
and\dss $z$\dss has only\sss one non-zero summand,\oss
namely\halfff,\oss the summand\dss 
$z^{\dff n}\qff \in\pff K^{\dff n\fff,\dff 0}$\dnsp.\oss
Since\dss $w^{\dff n\dff +\dff 1}\qff \in\pff K^{\dff n\dff +\dff 1\fff,\dff 0}$\dss
is\dss the only\sss non-zero summand of\trs $w$\dss
and\dss $w\off =\off \partial\dff z\off =\off d\dff z\qff +\dff \delta\dff z$\nnsp,\oss
we see\sss that\dss $d\dff z^{\dff n}\off =\off 0$\dss and\dss
$w^{\dff n\dff +\dff 1}\off =\off \delta\dff z^{\dff n}$\nnsp.\oss
It\dss follows\dss that\dss $z^{\dff n}\qff \in\pff L^{\dff n}$\dss
and\dss $w^{\dff n\dff +\dff 1}$\dss is\dss a coboundary\dss in\trs $L^{\dff \bullet}$\dnsp.\oss
Since\dss $w^{\dff n\dff +\dff 1}$\dss was an arbitrary\sss cocycle of\trs
$L^{\dff \bullet}$\dss turning\sss into a coboundary\dss in\trs $T^{\dff \bullet}$\nnsp,\oss
the injectivity\dss follows.\oss  \eproof

\myuppar{Homological\dss double complexes.}
One can also consider double complex\sss $K_{\trf \bullet\fff,\dff \bullet}$\dss 
of\trs the form\vspace{9pt}
\[
\qquad
\begin{tikzcd}[column sep=large, row sep=huge]\dis
K_{\dff 0\fff,\dff 0} 
& 
K_{\dff 0\fff,\dff 1} \arrow[l]
& 
K_{\dff 0\fff,\dff 2} \arrow[l] 
&   
K_{\dff 0\fff,\dff 3} \arrow[l] 
&
\off\off \ldots \arrow[l]
\\
K_{\dff 1\fff,\dff 0}
\arrow[u]
& 
K_{\dff 1\fff,\dff 1} \arrow[l]
\arrow[u]
& 
K_{\dff 1\fff,\dff 2} \arrow[l]
\arrow[u]
&   
K_{\dff 1\fff,\dff 3} \arrow[l]
\arrow[u]
&
\off\off \ldots \arrow[l]
\\
K_{\dff 2\fff,\dff 0} 
\arrow[u]
& 
K_{\dff 2\fff,\dff 1} \arrow[l]
\arrow[u]
& 
K_{\dff 2\fff,\dff 2} \arrow[l]
\arrow[u]
&   
K_{\dff 2\fff,\dff 3} \arrow[l]
\arrow[u]
&
\off\off \ldots \arrow[l]
\\
K_{\dff 3\fff,\dff 0} 
\arrow[u]
& 
K_{\dff 3\fff,\dff 1} \arrow[l]
\arrow[u]
& 
K_{\dff 3\fff,\dff 2} \arrow[l]
\arrow[u]
&   
K_{\dff 3\fff,\dff 3} \arrow[l]
\arrow[u]
&
\off\off \ldots
\\
\ldots \vphantom{K_{\dff 3\fff,\dff 3}} \arrow[u]
&
\ldots \vphantom{K_{\dff 3\fff,\dff 3}} \arrow[u]
&
\ldots \vphantom{K_{\dff 3\fff,\dff 3}} \arrow[u]
&
\ldots \vphantom{K_{\dff 3\fff,\dff 3}} \arrow[u]
&
\off\off\off\off .
\end{tikzcd}
\]

\vspace{-6pt}
The horizontal\sss arrows\qss
$K_{\dff p\fff,\dff q}\qff \ttoo\qff K_{\dff p\fff,\dff q\dff -\dff 1}$\qss
and\dss the vertical\sss arrows\qss
$K_{\dff p\fff,\dff q}\qff \ttoo\qff K_{\dff p\dff -\dff 1\fff,\dff q}$\qss
are denoted,\oss as before,\oss by\dss $d$\dss and\dss $\delta$\dss respectively\halfff,\oss
and are called\dss the\qss \emph{differentials}\qss of\trs the double complex\dss 
$K_{\trf \bullet\fff,\dff \bullet}$\dnsp.\oss
It\dss is\dss assumed\dss that\trs
$d\dff \circ\dff d\off =\off 0$\dss
and\qss
$\delta\dff \circ\dff \delta\off =\off 0$\nnsp.\oss
Also,\oss it\dss
is\dss assumed\dss that\sss either each square\dss is\dss commutative,\oss
or\sss each square\dss is\dss anti-commutative,\oss
i.e.\qss that\sss either\dss
$d\dff \circ\dff \delta\off =\off \delta\dff \circ\dff d$\nnsp,\oss
or\dss 
$d\dff \circ\dff \delta\pff +\pff \delta\dff \circ\dff d\off =\off 0$\nnsp.\oss
The\sss total\sss complex\dss $T_{\dff \bullet}$\dss is\dss defined as before.\oss

Let\dss $L_{\dff p}$\dss be the cokernel of the differential\qss
$d\dff \colon\dff
K_{\trf p\fff,\dff 1}\qff \ttoo\qff K_{\trf p\fff,\dff 0}$\dnsp,\oss
i.e.\vspace{3pt}
\[
\quad
L_{\dff p}
\off =\off\dff
K_{\trf p\fff,\dff 0}\dff\bigl/d\dff\left(\qff K_{\trf p\fff,\dff 1} \qff\right)
\pff,
\]

\vspace{-9pt}
and\dss let\dss
$\pi\dff \colon\dff
K_{\trf p\fff,\dff 0}
\qff \ttoo\qff
L_{\dff p}$\dss
be\sss the quotient\dss map.\oss
Either of\trs the commutativity\sss assumptions\dss implies\sss that\dss
$\delta$\dss maps\dss 
$d\dff(\qff K_{\trf p\fff,\dff 1} \qff)$\dss 
to\dss 
$d\dff(\qff K_{\trf p\dff -\qff 1\fff,\dff 1} \qff)$\dss
and\dss hence induces homomorphisms
\vspace{3pt}
\[
\quad
\delta_{\trf L}\dff \colon\dff
L_{\dff p}
\qff \ttoo\qff
L_{\dff p\dff -\dff 1}
\pff.
\]

\vspace{-9pt}
Therefore\qss $L_{\dff \bullet}$\qss together with\dss 
homomorphisms\qss 
$\delta_{\trf L}$\qss
is\dss a quotient\sss complex of\trs $T_{\dff \bullet}$\dnsp.\oss

\myapar{Theorem.}{homology-double-complex}
\emph{If\qss each complex\qss
$\dis
\left(\qff K_{\trf p\fff,\dff \bullet},\pff d \qff\right)$\qss 
is\dss exact\halfff,\oss
then the homomorphism}\qss\vspace{3pt}\vspace{-0.125pt}
\[
\quad
H_{\fff *}\dff(\trf T_{\dff \bullet} \dff)
\qff \ttoo\qff
H_{\fff *}\dff(\trf L_{\dff \bullet} \dff)
\]

\vspace{-9pt}\vspace{-0.125pt}
\emph{induced\dss by\dss the\dss quotient\dss map\qss
$T_{\dff \bullet}\qff \ttoo\qff L_{\dff \bullet}$\qss
is\dss an\dss isomorphism.\oss}

\proof
We may assume\sss that\dss the differentials
commute.\oss
Let\dss us prove first\dss the induced\dss homomorphism\dss is\dss surjective.\oss
Let\vspace{3pt}
\[
\quad
z
\off\dff =\off
\bigoplus_{i\qff =\qff 0}^n\qff z_{\dff i}
\off\qff \in\off\qff
\bigoplus_{i\qff =\qff 0}^n\qff K_{\dff i\fff,\dff n\dff -\dff i}
\quad.
\]

\vspace{-12.75pt}
Then\vspace{-4.25pt}
\[
\quad
\partial\dff z
\off\qff \in\off\qff
\bigoplus_{k\qff =\qff 0}^{n\qff -\qff 1}\qff K_{\dff k\fff,\dff n\dff -\dff 1\dff -\dff k}
\quad
\]

\vspace{-12pt}
and\dss the summand of\trs $\partial\dff z$\dss belonging\dss to\qss
$K_{\dff k\fff,\dff n\dff -\dff 1\dff -\dff k}$\qss is\dss equal\dss to\vspace{3pt}
\[
\quad
d\dff z_{\dff k}
\pff +\pff
\delta\dff z^{\dff k\dff +\dff 1}
\quad\
\mbox{if}\quad\
0\off <\off k\off <\off n\qff -\qff 1
\qff,
\]

\vspace{-9pt}
for every\dss
$k\off =\off 0\fff,\pff 1\fff,\pff \ldots\fff,\pff n\qff -\qff 1$\nnsp.\oss
Suppose\sss that\dss
$w_{\dff n}\qff \in\pff K_{\dff n\fff,\qff 0}$\dss and\dss 
$\pi\trf(\dff w_{\dff n}\trf)$\dss 
is\dss a cycle of\trs $L_{\dff \bullet}$\nsp.\oss
Then\vspace{3pt}
\[
\quad
\delta\trf(\dff w_{\dff n}\trf)
\qff \in\qff
d\dff\left(\qff K_{\trf n\dff -\qff 1\fff,\dff 1} \qff\right)
\]

\vspace{-9pt}
i.e.\qss 
$\delta\trf(\dff w_{\dff n}\trf)
\off =\off 
d\dff(\trf w_{\dff n\dff -\dff 1}\trf)$\dss
for some\dss
$w_{\dff n\dff -\dff 1}\qff \in\qff K_{\trf n\dff -\qff 1\fff,\dff 1}$\nsp.\oss
Let\qss $x_{\dff n\dff -\dff 1}\off =\off \delta\trf(\trf w_{\dff n\dff -\dff 1}\trf)$\nnsp.\oss
Then\vspace{3pt}
\[
\quad
d\trf(\trf x_{\dff n\dff -\dff 1}\trf)
\off =\off
d\dff \circ\dff \delta\qff(\trf w_{\dff n\dff -\dff 1}\trf)
\]

\vspace{-36pt}
\[
\quad
\phantom{d\trf(\trf x_{\dff n\dff -\dff 1}\trf)
\off }
=\off
\delta\dff \circ\dff d\qff(\trf w_{\dff n\dff -\dff 1}\trf)
\off =\off
\delta\dff \circ\dff \delta\qff(\trf w_{\dff n}\trf)
\off =\off
0
\pff.
\]

\vspace{-9pt}
Since\trs
$(\qff K_{\trf n\dff -\dff 1\fff,\dff \bullet},\pff d \qff)$\dss
is\dss exact\halfff,\oss
it\dss follows\dss that\dss
$x_{\dff n\dff -\dff 1}
\off =\off
d\trf(\trf w_{\dff n\dff -\dff 2}\trf)$\dss
for some\dss
$w_{\dff n\dff -\dff 2}
\qff \in\pff
K_{\trf n\dff -\qff 2\fff,\dff 2}$\nsp.\oss
By\sss continuing\dss  to argue in\dss this way\dss
we will\dss get\sss elements\dss
$w_{\dff n\dff -\dff 1}\dff,\pff w_{\dff n\dff -\dff 2}\dff,\pff \ldots\dff,\pff w_{\dff 0}$\dss
such\dss that\dss
$w_{\dff n\dff -\dff i}
\qff \in\pff
K_{\trf n\dff -\qff i\fff,\dff i}$\dss
and\dss
$\delta\trf(\dff w_{\dff n\dff -\dff i}\trf)
\off =\off 
d\dff(\trf w_{\dff n\dff -\dff i\dff -\dff 1}\trf)$\dss
for every\dss $i$\nnsp.\oss
Let\vspace{3pt}
\[
\quad
w
\off\dff =\off\dff
\bigoplus_{i\qff =\qff 0}^{n}\qff w_{\dff i}
\off\qff \in\off\qff
\bigoplus_{i\qff =\qff 0}^n\qff K_{\dff i\fff,\dff n\dff -\dff i}
\quad.
\]

\vspace{-9pt}
When\dss the differentials commute,\oss
the differential\trs
$\partial\dff \colon\dff
T_{\fff n}\qff \ttoo\qff T_{\fff n\dff -\dff 1}$\dss
is\dss equal\dss to\dss $d\qff +\qff (\qff -\qff 1\trf)^{\dff p}\qff \delta$\dss
on\dss $K_{\dff p\fff,\dff q}$\nsp.\oss
It\dss follows\dss that\trs 
$\partial\dff w\off =\off 0$\nnsp,\oss
i.e.\qss $w$\dss is\dss a cycle of\trs the\sss total\sss complex.\oss
Clearly\halfff,\oss the quotient\dss map\dss
$T_{\dff \bullet}\qff \ttoo\qff L_{\dff \bullet}$\dss
maps\dss $w$\dss to\dss $\pi\trf(\dff w_{\dff n}\trf)$\nnsp.\oss
The surjectivity\dss follows.\oss

The proof\dss of\dss injectivity\dss is\dss similar\halfff.\oss
Suppose\sss that\sss\vspace{3pt}
\[
\quad
z
\off\dff =\off\dff
\bigoplus_{i\qff =\qff 0}^{n\dff -\dff 1}\qff z_{\dff i}
\off\qff \in\off\qff
\bigoplus_{i\qff =\qff 0}^n\qff K_{\dff i\fff,\dff n\dff -\dff i\dff -\dff 1}
\quad.
\]

\vspace{-9pt}
is\dss a cycle of\trs the complex\dss $T_{\dff \bullet}$\dss
such\dss that\dss $\pi\trf(\dff z_{\dff n\dff -\dff 1}\trf)$\dss
is\dss a\sss boundary\dss in\dss $L_{\dff \bullet}$\nsp.\oss
Then\dss there exists an element\dss
$w_{\dff n}\qff \in\pff K_{\dff n\fff,\qff 0}$\dss
such\dss that\trs
$\pi\trf(\dff \delta\trf(\dff w_{\dff n}\trf)\trf)
\off =\off 
\pi\trf(\dff z_{\dff n\dff -\dff 1}\trf)$\dss
and\dss hence\sss there exits an element\dss
$w_{\dff n\dff -\dff 1}\qff \in\qff K_{\trf n\dff -\qff 1\fff,\dff 1}$\dss
such\dss that\trs
$z_{\dff n\dff -\dff 1}
\off =\off  
\delta\trf(\dff w_{\dff n}\trf)
\qff -\qff
d\dff(\trf w_{\dff n\dff -\dff 1}\trf)$\nnsp.\oss
Let\dss $y\qff \in\pff T_{\dff n}$\dss
be\sss the element\dss having only\dss two non-zero summands,\oss
namely\halfff,\pss $w_{\dff n}$\dss and\dss $w_{\dff n\dff -\dff 1}$\nsp.\oss
Then\dss $z\qff +\qff \partial\trf y$\dss and\dss
$z\qff -\qff \partial\trf y$\dss are cycles and\dss for an appropriate choice of\dss sign\dss 
the summand of\trs $z\qff \pm\qff \partial\dff y$\dss belonging\dss to\qss
$K^{\dff n\dff -\dff 1\fff,\dff 0}$\qss is\dss equal\dss to\dss $0$\nnsp.\oss
Therefore,\oss it\dss is\dss sufficient\dss to prove\sss that\dss $z$\dss is\dss
a boundary\dss in\dss $T_{\dff \bullet}$\dss if\trs 
$z_{\dff n\dff -\dff 1}\off =\off 0$\nnsp.\oss

If\trs $z_{\dff n\dff -\dff 1}\off =\off 0$\nnsp,\oss
then\dss $\partial\dff z\off =\off 0$\dss implies\sss that\dss
$d\trf(\dff z_{\dff n\dff -\dff 2}\trf)\off =\off 0$\nnsp.\oss
By\dss the assumption,\oss this implies\dss that\dss
$z_{\dff n\dff -\dff 2}
\off =\off
d\trf(\dff w_{\dff n\dff -\dff 2}\trf)$\dss
for some\dss
$w_{\dff n\dff -\dff 2}\qff \in\pff K_{\dff n\dff -\dff 2\fff,\dff 2}$\nsp.\oss
Let\dss us replace\dss $z$\dss by\trs
$z\qff \pm\qff \partial\trf (\dff w_{\dff n\dff -\dff 2}\trf)$\nnsp.\oss
The new\dss $z$\dss is\dss still\sss a\sss boundary\halfff,\oss
and with an appropriate choice of\dss sign we have\dss\vspace{3pt}
\[
\quad
z_{\dff n\dff -\dff 1}
\off =\off 
z_{\dff n\dff -\dff 2}
\off =\off 
0
\pff.
\]

\vspace{-9pt}
By\sss continuing\dss to argue in\dss this way\dss we will\dss eventually\dss
get\dss $z\off =\off 0$\nnsp.\oss
Since at\sss each step we subtracted\dss from\dss $z$\dss a boundary\halfff,\oss
it\dss follows\sss that\dss
the original\sss element\dss $z$\dss is\dss also a\sss boundary\halfff.\oss
The injectivity\dss follows.\oss
This completes\sss the proof.\oss  \eproof

\myuppar{Morphisms of\dss double complexes.}
Let\sss $K^{\dff \bullet\fff,\dff \bullet}$\sss and\dss
$L^{\dff \bullet\fff,\dff \bullet}$\sss be\sss two double complexes.\oss
A\qss \emph{morphism}\qss
$f^{\dff \bullet\fff,\dff \bullet}\dff \colon\dff
K^{\dff \bullet\fff,\dff \bullet}
\qff \ttoo\qff 
L^{\dff \bullet\fff,\dff \bullet}$\sss
is\dss a family\sss of\trs homomorphisms\sss
$f^{\dff p\fff,\dff q}\dff \colon\dff
K^{\dff p\fff,\dff q}
\qff \ttoo\qff 
L^{\dff p\fff,\dff q}$\sss
commuting\sss with\sss the differentials\sss $d\fff,\pff \delta$\sss
in an obvious sense.\oss
If\trs $T^{\dff \bullet}_{\dff K}$\sss and\sss $T^{\dff \bullet}_{\dff L}$\sss
are\sss the\sss total\sss complexes of\trs $K^{\dff \bullet\fff,\dff \bullet}$\sss and\dss
$L^{\dff \bullet\fff,\dff \bullet}$\sss respectively,\oss
the\sss $f^{\dff \bullet\fff,\dff \bullet}$\sss induces a morphism\dss
$T f^{\dff \bullet}\dff \colon\dff
T^{\dff \bullet}_{\dff K}
\pff \ttoo\qff 
T^{\dff \bullet}_{\dff L}$\nsp.\oss
Also,\pss $f^{\dff \bullet\fff,\dff \bullet}$ induces for every\sss $p\qff \geq\qff 0$\sss
a morphism of\dss complexes\sss
$f^{\dff p\fff,\dff \bullet}\dff \colon\dff
K^{\dff p\fff,\dff \bullet}
\qff \ttoo\qff 
L^{\dff p\fff,\dff \bullet}$\sss
and\sss for every\sss $q\qff \geq\qff 0$\sss
a morphism of\dss complexes\sss
$f^{\dff \bullet\fff,\dff q}\dff \colon\dff
K^{\dff \bullet\fff,\dff q}
\qff \ttoo\qff 
L^{\dff \bullet\fff,\dff q}$\dnsp.\oss

\myapar{Theorem.}{comparison}
\emph{Let\dss
$f^{\dff \bullet\fff,\dff \bullet}\dff \colon\dff
K^{\dff \bullet\fff,\dff \bullet}
\qff \ttoo\qff 
L^{\dff \bullet\fff,\dff \bullet}$\sss
be a morphism of\dss double complexes.\oss
If\qss for every\sss $p\qff \geq\qff 0$\sss
the morphism of\trs complexes\dss
$f^{\dff p\fff,\dff \bullet}\dff \colon\dff
K^{\dff p\fff,\dff \bullet}
\qff \ttoo\qff 
L^{\dff p\fff,\dff \bullet}$\sss
induces an\dss isomorphism\sss in\sss cohomology,\oss
then\dss
$T f^{\dff \bullet}\dff \colon\dff
T^{\dff \bullet}_{\dff K}
\pff \ttoo\qff 
T^{\dff \bullet}_{\dff L}$\qss
also induces an\dss isomorphism\sss in\sss cohomology.\oss}

\proof
The\dss standard\dss proof\trs is\dss based on\dss comparing\sss spectral sequences 
associated\dss with\dss $K^{\dff \bullet\fff,\dff \bullet}$\dnsp,\oss
$L^{\dff \bullet\fff,\dff \bullet}$\dnsp.\oss
Here\dss is\dss a direct\sss elementary\dss proof\halfff.\oss
We will\sss assume\sss that\dss the differentials\sss $d\fff,\pff \delta$\sss 
anti-commute and\dss the differentials
of\trs total\sss complexes are defined\sss as\sss $d\qff +\qff \delta$\nnsp.\oss
Let\dss us\sss prove\sss first\dss that\dss 
the induced\dss homomorphism\dss is\dss surjective.\oss
Let\vspace{-3pt}
\[
\quad
x
\off\dff =\off
\bigoplus_{i\qff =\qff 0}^n\qff x^{\dff i}
\off\qff \in\off\qff
\bigoplus_{i\qff =\qff 0}^n\qff L^{\dff i\fff,\dff n\dff -\dff i}
\quad
\]

\vspace{-15pt}
be a cocycle.\oss
Then\dss
$d\dff x^{\trf 0}\off =\off 0$\nnsp,\pss
$\delta\dff x^{\dff n}\off =\off 0$\nnsp,\oss
and\vspace{3pt}
\[
\quad
d\dff x^{\dff i}
\pff +\pff
\delta\dff x^{\dff i\dff -\dff 1}
\off =\off
0
\quad\
\mbox{for}\quad\
0\off <\off i\off <\off n
\qff.
\]

\vspace{-9pt}
In order\sss to prove\sss the surjectivity\sss we need\dss to find a cocycle\vspace{-3pt}
\[
\quad
y
\off\dff =\off
\bigoplus_{i\qff =\qff 0}^n\qff y^{\dff i}
\off\qff \in\off\qff
\bigoplus_{i\qff =\qff 0}^n\qff K^{\dff i\fff,\dff n\dff -\dff i}
\quad
\]

\vspace{-15pt}
such\dss that\dss 
$f\dff(\trf y\trf)
\off =\off
x\qff +\qff \partial\dff z$\trs for some\vspace{-3pt}
\[
\quad
z
\off\dff =\off
\bigoplus_{i\qff =\qff 0}^n\qff z^{\dff i}
\off\qff \in\off\qff
\bigoplus_{i\qff =\qff 0}^{n\qff -\qff 1}\qff L^{\dff i\fff,\dff n\dff -\dff 1\dff -\dff i}
\quad,
\]

\vspace{-15pt}
or\halfff,\oss equivalently,\pss
$x^{\trf 0}
\off =\off
f\dff(\trf y^{\trf 0}\trf)\qff +\qff d\dff z^{\trf 0}$\nsp,\pss
$x^{\dff n}
\off =\off
f\dff(\trf y^{\dff n}\trf)\qff +\qff \delta\dff z^{\dff n\dff -\dff 1}$\nsp,\oss
and\vspace{3pt}
\[
\quad
x^{\dff i}
\off =\off
f\dff(\dff y^{\dff i}\trf)
\qff +\qff
d\dff z^{\dff i}
\pff +\pff
\delta\dff z^{\trf i\dff -\dff 1}
\quad\
\mbox{for}\quad\
0\off <\off i\off <\off n
\qff.
\]

\vspace{-9pt}
Since
{\nsp}$f^{\trf 0\fff,\dff \bullet}\dff \colon\dff
K^{\trf 0\fff,\dff \bullet}
\qff \ttoo\qff 
L^{\trf 0\fff,\dff \bullet}$\sss
induces an\sss isomorphism\sss in cohomology,\oss
there exists\sss 
$y^{\trf 0}\qff \in\pff K^{\trf 0\fff,\dff n}$\sss
and\dss
$z^{\trf 0}\qff \in\pff K^{\trf 0\fff,\dff n\dff -\dff 1}$\sss
such\dss that\dss
$f\dff(\trf y^{\trf 0}\trf)
\off =\off
x^{\trf 0}\qff +\qff d\dff z^{\trf 0}$\nsp.\oss
Suppose\sss that\sss elements\dss 
$y^{\trf 0},\pff \ldots\fff,\off y^{\dff i\dff -\dff 1}$\sss
and\dss
$z^{\trf 0},\pff \ldots\fff,\off z^{\trf i\dff -\dff 1}$\sss
with\sss the required properties are already\sss constructed.\oss
Then\vspace{3pt}
\[
\quad
d\trf \delta\dff y^{\dff i\dff -\dff 1}
\off =\off 
-\qff \delta\trf d\dff y^{\dff i\dff -\dff 1}
\off =\off
\delta\trf d\dff d\dff y^{\dff i\dff -\dff 2}
\off =\off
0
\]

\vspace{-9pt}
and\dss hence\sss
$\delta\dff y^{\dff i\dff -\dff 1}$\sss 
is\dss a cocycle of\trs the complex\dss
$K^{\dff i\fff,\dff \bullet}$\dnsp.\oss 
At\dss the same\sss time\vspace{3pt}
\[
\quad
x^{\dff i\dff -\dff 1}
\off =\off
f\dff(\dff y^{\dff i\dff -\dff 1}\trf)
\qff +\qff
d\dff z^{\dff i\dff -\dff 1}
\pff +\pff
\delta\dff z^{\trf i\dff -\dff 2}
\]

\vspace{-12pt}
and\dss hence\vspace{0pt}
\[
\quad
d\dff x^{\dff i}
\off =\off
{}-\qff \delta\dff x^{\dff i\dff -\dff 1}
\off =\off
{}-\qff \delta\dff f\dff(\dff y^{\dff i\dff -\dff 1}\trf)
\qff -\qff
\delta\trf d\dff z^{\dff i\dff -\dff 1}
\off =\off
{}-\qff f\dff(\dff \delta\dff y^{\dff i\dff -\dff 1}\trf)
\qff +\qff
d\trf \delta\trf z^{\dff i\dff -\dff 1}
\pff.
\]

\vspace{-9pt}
It\dss follows\sss that\dss
$f\dff(\dff \delta\dff y^{\dff i\dff -\dff 1}\trf)$\sss
is\dss a coboundary\dss in\dss
$L^{\dff i\fff,\dff \bullet}$\dnsp.\oss
Since
{\nsp}$f^{\dff i\fff,\dff \bullet}\dff \colon\dff
K^{\dff i\fff,\dff \bullet}
\qff \ttoo\qff 
L^{\dff i\fff,\dff \bullet}$\sss
induces an\sss isomorphism\sss in cohomology,\oss
this\sss implies\sss that\sss
$\delta\dff y^{\dff i\dff -\dff 1}$\sss 
is\dss a coboundary\dss in\dss
$K^{\dff i\fff,\dff \bullet}$\dnsp.\oss
Equivalently,\pss
$\delta\trf y^{\dff i\dff -\dff 1}
\off =\off
d\dff u^{\dff i}$\dss
for some\qss
$u^{\dff i}\qff \in\pff K^{\dff i\fff,\pff n\dff -\dff i}$\dnsp.\oss
It\dss follows\dss that\vspace{3pt}\vspace{-0.125pt}
\[
\quad
f\dff(\dff \delta\dff y^{\dff i\dff -\dff 1}\trf)
\off =\off
f\dff(\dff d\dff u^{\dff i}\trf)
\off =\off
d\dff f\dff(\dff u^{\dff i}\trf)
\]

\vspace{-12pt}
and\dss hence\vspace{-0.5pt}
\[
\quad
d\trf\left(\qff
x^{\dff i}
\qff -\qff
\delta\trf z^{\dff i\dff -\dff 1}
\qff +\qff
f\dff(\dff u^{\dff i}\trf)
\qff\right)
\off =\off
d\dff x^{\dff i}
\qff -\qff
d\trf \delta\trf z^{\dff i\dff -\dff 1}
\qff +\qff
f\dff(\dff \delta\dff y^{\dff i\dff -\dff 1}\trf)
\off =\off 
0
\pff.
\]

\vspace{-12.5pt}
Since
{\nsp}$f^{\dff i\fff,\dff \bullet}\dff \colon\dff
K^{\dff i\fff,\dff \bullet}
\qff \ttoo\qff 
L^{\dff i\fff,\dff \bullet}$\sss
induces an\sss isomorphism\sss in cohomology,\oss
this\sss implies\sss that\sss\vspace{3pt}
\[
\quad
x^{\dff i}
\qff -\qff
\delta\trf z^{\dff i\dff -\dff 1}
\qff +\qff
f\dff(\dff u^{\dff i}\trf)
\off =\off
f\dff(\dff v^{\dff i}\trf)
\qff +\qff
d\dff z^{\trf i}
\]

\vspace{-9pt}
for some\dss
$z^{\trf i}\qff \in\pff L^{\dff i\fff,\pff n\dff -\dff i}$\dnsp.\oss
Hence\vspace{2.5pt}
\[
\quad
x^{\dff i}
\off =\off
f\dff(\dff v^{\dff i}
\qff -\qff
u^{\dff i}\trf)
\qff +\qff
d\dff z^{\trf i}
\qff +\qff
\delta\trf z^{\dff i\dff -\dff 1}
\]

\vspace{-9.5pt}
and we can\sss set\dss
$y^{\dff i}
\off =\off
v^{\dff i}
\qff -\qff
u^{\dff i}$\dnsp.\oss
By continuing in\sss this way\sss we will\sss find\sss the sequences 
$y$ and\sss $z$\sss with\dss the required\dss properties.\oss

In order\sss to prove\sss the injectivity,\oss
suppose\sss that\sss $y$\sss is\dss an $n$\dnsp-cocycle in\sss the\sss total\sss
complex of\trs $K^{\dff \bullet\fff,\dff \bullet}$\sss such\dss that\dss 
$f\dff(\trf y\trf)\off =\off \partial\dff z$\nnsp.\oss
Then\dss
$f\dff(\trf y^{\trf 0}\trf)\off =\off d\dff z^{\trf 0}$\nsp,\pss
$f\dff(\trf y^{\dff n}\trf)\off =\off \delta\dff z^{\dff n\dff -\dff 1}$\nsp,\oss
and\vspace{3pt}
\[
\quad
f\dff(\dff y^{\dff i}\trf)
\off =\off
d\dff z^{\dff i}
\pff +\pff
\delta\dff z^{\trf i\dff -\dff 1}
\quad\
\mbox{for}\quad\
0\off <\off i\off <\off n
\qff.
\]

\vspace{-9pt}
We will\sss construct\sss elements\sss
$u^{\dff i}\qff \in\pff K^{\dff i\fff,\dff n\dff -\dff 1\dff -\dff i}$\sss
for\sss $0\qff \leq\qff i\qff \leq\qff n\qff -\qff 1$\sss
and\sss elements
$c^{\dff i}\qff \in\pff L^{\dff i\fff,\dff n\dff -\dff 2\dff -\dff i}$\sss
for\sss $0\qff \leq\qff i\qff \leq\qff n\qff -\qff 2$\sss
such\dss that\dss
$y^{\trf 0}\off =\off d\dff u^{\trf 0}$\nsp,\pss
$f\dff(\dff u^{\trf 0}\trf)
\off =\off 
z^{\trf 0}
\qff +\qff
d\fff c^{\trf 0}$\dnsp,\pss
$y^{\dff n}\off =\off \delta\dff z^{\dff n\dff -\dff 1}$\nsp,\oss\vspace{3pt}
\[
\quad
y^{\dff i}
\off =\off
d\fff u^{\dff i}
\pff +\pff
\delta\dff u^{\trf i\dff -\dff 1}
\quad\
\mbox{for}\quad\
0\off <\off i\off <\off n
\qff,\quad\
\mbox{and}
\]

\vspace{-36pt}
\[
\quad
f\dff(\dff u^{\dff i}\trf)
\off =\off
z^{\dff i}
\qff +\qff
d\fff c^{\dff i}
\qff +\qff
\delta\dff c^{\trf i\dff -\dff 1}
\quad\
\mbox{for}\quad\
0\off <\off i\off <\off n\qff -\qff 1
\pff.
\]

\vspace{-9pt}
To begin\sss with,\oss note\sss that\sss
since
$f^{\trf 0\fff,\dff \bullet}\dff \colon\dff
K^{\trf 0\fff,\dff \bullet}
\qff \ttoo\qff 
L^{\trf 0\fff,\dff \bullet}$\sss
induces an\sss isomorphism\sss in cohomology,\oss
there exists\sss 
$v^{\trf 0}\qff \in\pff K^{\trf 0\fff,\dff n\dff -\dff 1}$\sss
such\dss that\dss
$y^{\trf 0}
\off =\off
d\fff v^{\trf 0}$\nsp.\oss
Clearly,\pss\vspace{3pt}
\[
\quad
d\dff \left(\qff f\dff (\dff v^{\trf 0}\trf)
\qff -\qff
z^{\trf 0}
\qff\right)
\off =\off
f\dff(\trf d\fff v^{\trf 0}\trf)\qff -\qff d\dff z^{\trf 0}
\off =\off
f\dff(\trf y^{\trf 0}\trf)\qff -\qff d\dff z^{\trf 0}
\off =\off
0
\pff.
\]

\vspace{-9pt}
It\dss follows\sss that\sss
$f\dff(\trf v^{\trf 0}\trf)\qff -\qff z^{\trf 0}
\off =\off
f\dff(\trf w^{\trf 0}\trf)\qff +\qff d\dff c^{\trf 0}$\dss
and\dss hence\vspace{2.5pt}
\[
\quad
f\dff(\trf v^{\trf 0}\qff -\qff w^{\trf 0} \trf)
\off =\off
z^{\trf 0}
\qff +\qff d\dff 
c^{\trf 0}
\]

\vspace{-9.5pt}
for some\sss $w^{\trf 0}\qff \in\pff K^{\trf 0\fff,\dff n\dff -\dff 1}$\sss
such\dss that\sss $d\fff w^{\trf 0}\off =\off 0$\sss
and\sss some\sss 
$c^{\trf 0}\qff \in\pff K^{\trf 0\fff,\dff n\dff -\dff 2}$\dnsp.\oss
Hence we can set\sss 
$u^{\trf 0}
\off =\off 
v^{\trf 0}\qff -\qff w^{\trf 0}$\dnsp.\oss 
Suppose\sss that\sss we already\sss constructed\dss 
$u^{\trf 0},\pff \ldots\fff,\off u^{\dff i\dff -\dff 1}$\sss
and\dss
$c^{\trf 0},\pff \ldots\fff,\off c^{\dff i\dff -\dff 1}$\sss
with\dss the required\sss properties.\oss
Then\vspace{3pt}\vspace{-0.25pt}
\[
\quad
f\dff(\trf y^{\dff i}\trf)
\off =\off
d\dff z^{\dff i}
\pff +\pff
\delta\dff z^{\trf i\dff -\dff 1}
\off =\off
d\dff z^{\dff i}
\pff +\pff
\delta\trf f\dff(\dff u^{\dff i\dff -\dff 1}\trf)
\pff -\pff
\delta\trf d\fff c^{\dff i\dff -\dff 1}
\]

\vspace{-36pt}
\[
\quad
\phantom{f\dff(\trf y^{\dff i}\trf)
\off =\off
d\dff z^{\dff i}
\pff +\pff
\delta\dff z^{\trf i\dff -\dff 1}
\off }
=\off
d\dff z^{\dff i}
\pff +\pff
f\dff(\dff \delta\trf u^{\dff i\dff -\dff 1}\trf)
\pff +\pff
d\trf \delta\trf c^{\dff i\dff -\dff 1}
\pff.
\]

\vspace{-9pt}
It\dss follows\sss that\sss\vspace{3pt}
\[
\quad
f\dff(\trf y^{\dff i}
\qff -\qff
\delta\trf u^{\dff i\dff -\dff 1}\trf)
\off =\off
d\trf(\trf z^{\dff i}
\qff +\qff
\delta\trf c^{\dff i\dff -\dff 1}
\trf)
\pff.
\]

\vspace{-9pt}
Since $f^{\dff i\fff,\dff \bullet}\dff \colon\dff
K^{\dff i\fff,\dff \bullet}
\qff \ttoo\qff 
L^{\dff i\fff,\dff \bullet}$\sss
induces an\sss isomorphism\sss in cohomology,\oss
there exists an element\sss 
$v^{\dff i}\qff \in\pff K^{\dff i\fff,\dff n\dff -\dff 1\dff -\dff i}$\sss
such\dss that\dss
$y^{\dff i}
\qff -\qff
\delta\trf u^{\dff i\dff -\dff 1}
\off =\off
d\fff v^{\dff i}$\dnsp.\oss
Clearly,\oss\vspace{3pt}
\[
\quad
f\dff(\trf y^{\dff i}\trf)
\off =\off
f\dff(\trf \delta\trf u^{\dff i\dff -\dff 1}\trf)
\pff +\pff
f\dff(\trf d\fff v^{\dff i}\trf)
\pff.
\]

\vspace{-9pt}
By comparing\sss this with\sss the previous expression\sss for\sss
$f\dff(\trf y^{\dff i}\trf)$\sss
we conclude\sss that\vspace{3pt}
\[
\quad
d\dff z^{\dff i}
\pff +\pff
d\trf \delta\trf c^{\dff i\dff -\dff 1}
\off =\off
f\dff(\trf d\fff v^{\dff i}\trf)
\off =\off
d\dff f\dff(\trf v^{\dff i}\trf)
\]

\vspace{-12pt}
and\dss hence\vspace{0pt}
\[
\quad
d\trf\left(\qff
f\dff(\trf v^{\dff i}\trf)
\qff -\qff
z^{\dff i}
\qff -\qff
\delta\trf c^{\dff i\dff -\dff 1}
\qff\right)
\off =\off
0
\pff.
\]

\vspace{-9pt}
Since $f^{\dff i\fff,\dff \bullet}\dff \colon\dff
K^{\dff i\fff,\dff \bullet}
\qff \ttoo\qff 
L^{\dff i\fff,\dff \bullet}$\sss
induces an\sss isomorphism\sss in cohomology,\oss
there exists elements\sss 
$w^{\dff i}\qff \in\pff K^{\trf i\fff,\dff n\dff -\dff 1\dff -\dff 1}$\sss
and\dss
$c^{\dff i}\qff \in\pff L^{\dff i\fff,\dff n\dff -\dff 2\dff -\dff i}$\sss
such\dss that\dss $d\fff w^{\dff i}\off =\off 0$\sss
and\vspace{3pt}
\[
\quad
f\dff(\trf v^{\dff i}\trf)
\qff -\qff
z^{\dff i}
\qff -\qff
\delta\trf c^{\dff i\dff -\dff 1}
\off =\off
f\dff(\trf w^{\dff i}\trf)
\qff +\qff
d\fff c^{\dff i}
\pff.
\]

\vspace{-12pt}
It\dss follows\dss that\vspace{0pt}
\[
\quad
f\dff(\trf v^{\dff i}\qff -\qff w^{\dff i} \trf)
\off =\off
z^{\dff i}
\qff +\qff
d\fff c^{\dff i}
\qff +\qff
\delta\trf c^{\dff i\dff -\dff 1}
\pff
\]

\vspace{-9pt}
and\dss hence we can set\dss
$u^{\dff i}\off =\off v^{\dff i}\qff -\qff w^{\dff i}$\dnsp.\oss
By arguing\sss in\sss this way\sss we can construct\sss elements\sss $u^{\dff i}$\sss
and\sss
$c^{\dff i}$\sss for\sss $i\qff \leq\qff n\qff -\qff 2$\nnsp.\oss
The\sss last\dss two steps\sss
corresponding\dss to\dss $i\off =\off n\qff -\qff 1\fff,\pff n$\sss
are similar and can\sss be done by\sss setting\dss
$z^{\dff n}\off =\off 0$\sss
and\dss $c^{\dff n\dff -\dff 1}\off =\off c^{\dff n}\off =\off 0$\sss
in\sss the above arguments.\oss  \eproof

\begin{flushright}

December\qss 8,\oss 2020
 
https\halfff:/\!/\hspace*{-0.06em}nikolaivivanov.com

E-mail\halfff:\oss nikolai.v.ivanov{\fff}@{\dff}icloud.com

\end{flushright}

\end{document}